%% file: arxiv.tex
\pdfoutput=1
\documentclass[10pt]{article}
\usepackage{arxiv-setting}

\title{Optimal and exact recovery on the general non-uniform Hypergraph Stochastic Block Model\blfootnote{Author names are listed in alphabetical order.}}

\author{Ioana Dumitriu\thanks{Department of Mathematics, University of California, San Diego, La Jolla, CA 92093. Email: \texttt{idumitriu@ucsd.edu}}
\qquad
Hai-Xiao Wang\thanks{Department of Mathematics, University of California, San Diego, La Jolla, CA 92093. Email: \texttt{h9wang@ucsd.edu}}}

\date{This version: June 2026}

\begin{document}

\maketitle

\begin{abstract}
Consider the community detection problem in random hypergraphs under the non-uniform hypergraph stochastic block model (HSBM), where each hyperedge appears independently with some given probability depending only on the labels of its vertices. We establish, for the first time in the literature, a sharp threshold for exact recovery under this non-uniform case, subject to minor constraints; in particular, we consider the model with multiple communities. One crucial point here is that by aggregating information from all the uniform layers, we may obtain exact recovery even in cases when this may appear impossible if each layer were considered alone. Besides that, we prove a wide-ranging, information-theoretic lower bound on the number of misclassified vertices \emph{for any algorithm}, depending on a \emph{generalized Chernoff-Hellinger} divergence involving model parameters. We provide two efficient algorithms which successfully achieve exact recovery when above the threshold, and attain the lowest possible mismatch ratio when the exact recovery is impossible, proved to be optimal. The theoretical analysis of our algorithms relies on the concentration and regularization of the adjacency matrix for non-uniform random hypergraphs, which could be of independent interest. We also address some open problems regarding parameter knowledge and estimation.
\end{abstract}

\tableofcontents
\input{maintext}
\section*{Acknowledgments}
I.D. and H.X.W. acknowledge support from NSF DMS-2154099. The authors would like to thank Ery Arias-Castro and Yizhe Zhu for helpful discussions and also Nicholas Cook and Jorge Garza-Vargas for helpful comments leading to a better understanding of \Cref{ass:expected_center_separation}, during \emph{Summer School on Random Matrix Theory and Its Applications, 2023} at the Ohio State University, Columbus, OH. Part of this work was done during the \emph{Fourth ZiF Summer School, 2022} at the Center for Interdisciplinary Research, Bielefeld University, Germany. The authors thank the organizers for their hospitality. The authors thank the anonymous reviewers for detailed comments and suggestions that greatly improve the presentation of this work.
\printbibliography
\newpage
\appendix
\phantomsection
\addcontentsline{toc}{section}{Appendices}
\input{appendix}
\end{document}

%% file: maintext.tex
\section{Introduction}\label{sec:intro}
The task of \emph{community detection}, or \emph{clustering}, consists in partitioning the vertices of a graph into groups that are similarly connected \cite{Abbe2018CommunityDA}. It has become one of the central problems in network analysis and machine learning \cite{Shi2000NormalizedCA, Newman2002RandomGM, Ng2002SpectralCA, Arias2014CommunityDI}.
Random graph models, which generate community structure with a specified ground truth, have been frequently employed in many clustering algorithms, for benchmarking and proving theoretical guarantees. Mathematically, let $\gG = (\gV, \gE)$ be a graph, where the vertex set $\gV=[N]$ is composed of $K$ disjoint blocks, i.e., $\gV = \cup_{k=1}^{K} \gV_k$. The proportion of each block is denoted by $\alpha_k = |\gV_k|/|\gV|$ and we define the vector $\rvalpha = (\alpha_1, \ldots, \alpha_K)$ with $\|\rvalpha\|_1 = 1$. Let $\rvy \in [K]^{N}$ denote the membership vector of the vertices, i.e., $\ervy(v) = \ervy_{v} = k$ if the vertex $v$ belongs to block $\gV_k$. Let $\widehat{\rvy}$ denote some estimation of $\rvy$. To evaluate the accuracy of $\widehat{\rvy}$, we define the \emph{mismatch} ratio, counting the proportion of incorrectly clustered nodes:
\begin{align}\label{eqn:misratio}
   \mismatch_{N} \coloneqq \mismatch_{N}(\rvy, \widehat{\rvy}) = N^{-1}\,\,\inf_{\pi \in \mathfrak{S}_{K}} \D_{\textnormal{HD}}(\pi \circ\rvy, \widehat{\rvy}),
\end{align}
where the \textit{Hamming distance} $\D_{\textnormal{HD}}(\rvy, \widehat{\rvy})$ counts the number of different entries between $\rvy$ and $\widehat{\rvy}$. Here, $\pi \circ \rvy$, defined by $(\pi \circ \rvy)(v) = \pi(\rvy(v))$ entrywisely, is the same assignment as $\rvy$ up to some permutation $\pi$, where $\mathfrak{S}_{K}$ denotes the group of all permutations on $[K]$. Note that the expected accuracy of a random guess estimator is $\|\rvalpha\|_2^2$\footnote{By Cauchy-Schwarz inequality, $\|\rvalpha\|_2^2 \geq 1/K$, which is the accuracy lower bound without knowing $\rvalpha$.}. Consequently, the estimator $\widehat{\rvy}$ is meaningful only if $\mismatch_{N} \leq 1 - \|\rvalpha\|_2^2$. Furthermore, the recovery problem can be divided into several regimes according to $\mismatch_{N}$ \cite{Abbe2018CommunityDA}.
\begin{enumerate}
	\item Exact recovery (strong consistency): $\P(\mismatch_{N} = 0) \geq 1 - o(1)$.
	\item Almost exact recovery (weak consistency): $\P(\mismatch_{N} = o(1)) \geq 1 - o(1)$.
	\item Partial recovery: $\P(\mismatch_{N} \leq 1 - \gamma ) \geq 1 - o(1)$ for $\gamma \in (\|\rvalpha\|_2^2, 1)$.
	\item Weak recovery (detection): $\P(\mismatch_{N} \leq  1 - \|\rvalpha\|_2^2 - \Omega(1)) \geq 1 - o(1)$\footnote{The difference between partial and weak is that, the algorithm does better than random guess in weak recovery, but cannot quantify how much better as in partial recovery.}.
\end{enumerate}

The \emph{Stochastic Block Model} (SBM) is a natural generalization of the Erd\H{o}s-R\'{e}nyi graph, where vertices are densely connected within each community but sparsely connected across different communities (\emph{assortative} case). SBM was first introduced in the pioneering work \cite{Holland1983StochasticBM} for sociology research. Over the past several decades, it has been extensively studied in \cite{Bui1984GraphBA, Boppana1987EigenvaluesAG,Dyer1989TheSO, Snijders1997EstimationAP, Condon1999AlgorithmsFG, McSherry2001SpectralPO, Bickel2009ANV, Coja-Oghlan2010GraphPV, Bickel2011TheMO, Rohe2011SpectralCA, Choi2012StochasticBW, Arias2014CommunityDI, Verzelen2015CommunityDI}, driven by exploring the phase transition behaviors in various connectivity regimes. In the regime where the expected degrees grow logarithmically with respect to the number of vertices, a major breakthrough was the establishment of the exact recovery thresholds for binary \cite{Abbe2016ExactRI, Mossel2016ConsistencyTF} and multi-block case \cite{Abbe2015CommunityDI, Yun2016OptimalCR, Agarwal2017MultisectionIT}, where both necessary and sufficient conditions were provided. For the case where the expected degrees are of constant order, phase transition behavior for detection was discovered in \cite{Decelle2011AsymptoticAO, Decelle2011InferenceAP, Darst2014AlgorithmIB}, and later it was connected to the \emph{Kesten-Stigum} (KS) threshold, with the necessity \cite{Mossel2015ReconstructionAE} and sufficiency \cite{Massoulie2014CommunityDT, Mossel2018ProofOT} rigorously proved for the binary block case. For the multi-block case, it was shown that a spectral clustering algorithm based on the non-backtracking operator \cite{Bordenave2018NonbacktrackingSO}, and the acyclic belief propagation method \cite{Abbe2015DetectionIT}, succeed all the way down to the KS threshold, proving a conjecture in \cite{Krzakala2013SpectralRI}. The necessity for $K = 3, 4$ was partly established in \cite{Mossel2022ExactPT}. However, the threshold is not known when the number of communities is more than $4$, although it is known to be below the KS threshold \cite{Banks2016InformationTT, Abbe2015DetectionIT, Mossel2022ExactPT}.

The literature abounds with different methods of approach, like spectral methods \cite{McSherry2001SpectralPO, Coja-Oghlan2010GraphPV, Xu2014EdgeLI, Yun2014CommunityDV, Yun2016OptimalCR, Chin2015StochasticBM, Lei2015ConsistencyOS, Vu2018ASimple, Lei2019Unified, Abbe2020EntrywiseEA}, sphere comparison \cite{Abbe2015CommunityDI, Abbe2015RecoveringCI}, as well as \emph{semidefinite programming} (SDP) \cite{Guedon2016CommunityDC, Hajek2016AchievingEC, Jalali2016ExploitingTF, Javanmard2016PhaseTI, Montanari2016SemidefinitePO, Agarwal2017MultisectionIT, Perry2017SemidefinitePF, Bandeira2018RandomLM}. 
Recently, some more general variants of SBM have been introduced, such as degree-corrected SBM \cite{Karrer2011StochasticBM, Gao2018CommunityDI, Agterberg2022JointSC, Jin2023PhaseTF}, contextual SBM \cite{Deshpande2018ContextualSB, Lu2020ContextualSB, Abbe2022LPT}, labeled SBM \cite{Heimlicher2012CommunityDI, Xu2014EdgeLI, Yun2014CommunityDV, Yun2016OptimalCR}, and multilayer SBM \cite{Chen2022GlobalAI, Agterberg2022JointSC, Ma2023CommunityDI}. Readers may refer to \cite{Abbe2018CommunityDA} for a more detailed review.

While graphs are usually used to depict pairwise relationships among data, hypergraphs, as a generalization of graphs, have gained increasing attention in recent years due to their ability to capture higher-order interactions among vertices \cite{Benson2016HigherOO, Battiston2020NetworksBP}. They are particularly useful in modeling complex systems where relationships involve more than two entities, including but not limited to biological networks \cite{Michoel2012AlignmentAI,Tian2009HypergraphLA} with multi-way interactions, citation networks \cite{Ji2016CoauthorshipAC}, image data \cite{Wen2019LearningNH}, as well as recommendation systems \cite{Bu2010MusicRB,Li2013NewsRV} with multiple users and items. Empirically, they have been shown to be better than graph models \cite{Zhou2007LearningWH}. Besides that, hypergraphs and their spectral theory found applications in data science \cite{Jain2014ProvableTF,Zhou2021SparseRT,Harris2021DeterministicTC}, combinatorics \cite{Friedman1995SecondEO,Soma2019SpectralSO,Dumitriu2021SpectraOR} and statistical physics \cite{Caceres2021SparseSYK,Sen2018OptimizationOS}. The study of random hypergraphs has led to new insights into the structure and dynamics of complex systems, as well as the development of novel algorithms for community detection and clustering.

The \emph{Hypergraph Stochastic Block Model} (HSBM), as a generalization of the SBM, was first introduced in \cite{Ghoshdastidar2014ConsistencyOS}, where each edge of the uniform hypergraph appears independently with some given probability. In recent years, many efforts have been made to study community detection problems on random hypergraphs. For exact recovery of uniform HSBMs, the phase transition was shown to occur when the expected degree for each vertex is at least $\log(N)$, and the thresholds were given in \cite{Kim2018StochasticBM, Gaudio2023CommunityDI, Zhang2023ExactRI} by generalizing techniques in \cite{Abbe2016ExactRI, Abbe2020EntrywiseEA, Abbe2015CommunityDI}. Spectral methods were considered in \cite{Ghoshdastidar2014ConsistencyOS, Ahn2016CommunityRI, Ghoshdastidar2017ConsistencyOS,Chien2018CommunityDI, Cole2020ExactRI, Gaudio2023CommunityDI, Zhang2023ExactRI}, while SDP methods were analyzed in \cite{Kim2018StochasticBM, Lee2020RobustHC, Gaudio2023CommunityDI}. 
Meanwhile, results on almost exact and partial recovery were established in \cite{Ghoshdastidar2017ConsistencyOS, Chien2019MinimaxMR, Ke2020CommunityDF, Dumitriu2025PartialRA}. For detection of uniform HSBM, it was conjectured in \cite{Angelini2015SpectralDO} that the phase transition occurs in the regime of constant expected degrees. Spectral algorithms, which output a partition better than a random guess when above the KS threshold, were provided in \cite{Pal2021ComunityDI} for the binary community case, and in \cite{Stephan2022SparseRH} for the multi-block case. On the other hand, \cite{Gu2023WeakRT, Gu2024CommunityDI} proved that when below the KS threshold, detection is impossible for $\ell$-uniform hypergraphs with binary communities when $\ell \leq 6$ and the expected degrees are sufficiently large, thus establishing the tightness of the KS threshold in those scenarios. They also proved that in certain parameter regimes KS threshold is not tight when $\ell \geq 7$. More recently, for non-uniform HSBMs, \cite{fernandez2026achieving} proved that detection is possible whenever the combined \emph{signal-to-noise} ratio across all layers exceeds the KS threshold, and they provided a polynomial-time spectral algorithm that achieves this bound.

It should be noted that most of the results concern community detection on uniform hypergraphs, which require the same number of vertices per edge. This is a constraining and somewhat impractical assumption. As will be demonstrated after \Cref{thm:agnostic_strong_consistency}, the thresholds for the non-uniform case show that using the information from all uniform layers yields strictly better results than considering each layer alone. However, the non-uniform HSBM was less explored in the literature, with notable results in \cite{Ghoshdastidar2017ConsistencyOS, Dumitriu2025PartialRA, Alaluusua2023MultilayerHC, Philip2023NonbacktrackingSC, Wang2023ITLimits}. Generally, the results here only considered sufficiency and did not establish fundamental thresholds.

Our results are multi-folds, and work for the general, non-uniform HSBM case. All these results are new and optimal for non-uniform HSBM.
\begin{enumerate}
    \item We prove a wide-ranging, \emph{Information-Theoretic} lower bound on the number of misclassified vertices \emph{for any algorithm}, depending on a \emph{Generalized Chernoff-Hellinger} distance involving model parameters (\Cref{thm:IT_lower_bounds} (2)).
    \item The necessary condition and threshold for exact recovery are derived in \Cref{thm:IT_lower_bounds} (1).
    \item We present an efficient \Cref{alg:partition_with_prior} which, above this threshold, with full information about the model, \emph{achieves} exact recovery for all but a negligible subset of the space of problems (\Cref{thm:known_strong_consistency}).
    \item We present an efficient \Cref{alg:agnostic_partition} which, also above the threshold, \emph{achieves} exact recovery when only the number of communities is known, provided that the edge probabilities are proportional (\Cref{thm:agnostic_strong_consistency}).   
    \item Below this threshold, we prove that the provided algorithms are \emph{optimal} in Theorems \ref{thm:agnostic_optimality} and \ref{thm:known_optimality}, meaning that the number of vertices they mislabeled reaches is \emph{at most} the lower bound determined by the IT limits.
\end{enumerate}

\subsection{Non-uniform Hypergraph Stochastic Block Model}
The non-uniform HSBM was first studied in \cite{Ghoshdastidar2017ConsistencyOS}, which can be treated as a superposition of several uniform HSBMs with different model parameters. It is a more realistic model to study higher-order interaction on networks. We introduce the rigorous definition of uniform HSBM first, and extend it to non-uniform hypergraphs.
\begin{definition}[Hypergraph]\label{def:hypergraph}
    A hypergraph $\gH$ is a pair $\gH = (\gV, \gE)$, where $\gE$ denotes the set of non-empty subsets of vertex set $\gV$. We say that $\gH$ is $\ell$-uniform if every $e\in\gE$ is an $\ell$-subset of $\gV$. The degree of a vertex $v \in \gV$ is the number of hyperedges in $\gE$ containing $v$.
\end{definition}

A faithful representation of an $\ell$-uniform hypergraph is to associate it to a tensor.
\begin{definition}[Adjacency tensor]\label{def:adjacency_tensor} 
    One can associate an $\ell$-uniform hypergraph $\gH_{\ell}=(\gV, \gE_{\ell})$ to an order-$\ell$ symmetric tensor $\tA^{(\ell)}$, where the tensor entry $\etA^{(\ell)}_{e}$ denotes the presence of some $\ell$-hyperedge $e$, i.e., $\etA_{e}^{(\ell)} \coloneqq \etA_{i_1,\dots, i_{\ell}}^{(\ell)} = \indi{e \in \gE_{\ell}}$ for $i_1, \ldots, i_{\ell} \in \gV$, and $\etA_{i_1,\dots, i_{\ell}}^{(\ell)} = \etA_{i_{\pi(1)},\ldots,i_{\pi(\ell)}}^{(\ell)}$ for any permutation $\pi$ on $[\ell]$.
\end{definition}

Let $N, K \in \mathbb{N}_{+}$ denote the number of vertices and communities, respectively. Let the random variable $\rY_v$ denote the membership assignment of $v$. The vertex set $\gV$ is partitioned into classes $\gV_1, \ldots, \gV_K$ probabilistically according to a \emph{probability vector} $\rvalpha = (\alpha_1, \ldots, \alpha_K)$ with $\|\rvalpha\|_1 =1$, i.e., $\P(v \in \gV_k) = \P(\rY_{v} = k) = \alpha_k$ with $\gV_k \coloneqq \{v\in \gV\mid \rY_{v} = k\}$ for each $k\in [K]$. We focus on the scenario where the sizes of different communities are comparable.
\begin{assumption}[Comparable sizes]\label{ass:comparable_sizes}
    Assume that $\alpha_1 \geq \ldots \geq \alpha_K \geq \const$ for some universal constant $\const \in (0, 1)$, i.e., no entry of $\rvalpha$ vanishes.
\end{assumption}

\Cref{ass:comparable_sizes} is adopted throughout the dissertation. We first intorduce the uniform HSBM, which is characterized by a partition of vertices into different classes, such that all vertices in a given class are interchangeable.

\begin{definition}[Uniform HSBM]\label{def:uniform_HSBM}
With $N, K, \rvalpha$ defined above, let $\rvy\in [K]^N$ denote the membership vector where each entry of $\rvy$ is sampled independently according to $\rvalpha$. Let $\tQ^{(\ell)} \in ([0, 1]^{K})^{\otimes \ell}$ be an order-$\ell$ symmetric tensor such that $\etQ_{ \ervy_1, \ldots, \ervy_{\ell} }^{(\ell)} = \etQ_{ \ervy_{\pi(1)}, \ldots, \ervy_{\pi(\ell)} }^{(\ell)}$ for any permutation $\pi$ on $[\ell]$. Each possible edge $e = \{i_1, \ldots, i_{\ell}\}$ is generated with probability $\P( \etA_{e}^{(\ell)} = 1) = \etQ_{\rvy(e)}^{(\ell)}$, where $\rvy(e) = \{\ervy_{i_1}, \ldots, \ervy_{i_{\ell}}\}$ represents the membership sequence of edge $e$. We denote this distribution on the set of $\ell$-uniform hypergraphs by
\begin{align}\label{eqn:uniform_HSBM}
 (\rvy, \gH_{\ell}) \sim  \textnormal{HSBM}_{\ell}(N,\rvalpha, \tQ^{(\ell)})\,, \quad \ell \in \sL\,.
\end{align}
\end{definition}

\begin{figure}[h]
     \centering
    \begin{subfigure}[b]{0.45\textwidth}
         \centering
         \includegraphics[width=0.60\textwidth]{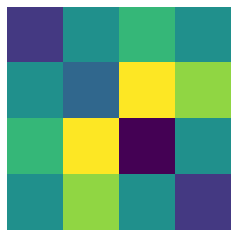}
         \caption{$\tQ^{(2)}$ is a symmetric matrix.}
     \end{subfigure}
    \begin{subfigure}[b]{0.45\textwidth}
         \centering
         \includegraphics[width=0.80\textwidth]{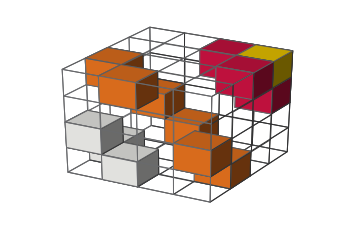}
         \caption{$\etQ^{(3)}_{124} = \etQ^{(3)}_{142} = \etQ^{(3)}_{214} = \etQ^{(3)}_{241} = \etQ^{(3)}_{412} = \etQ^{(3)}_{421}$ (orange); $\etQ^{(3)}_{344}=\etQ^{(3)}_{434} = \etQ^{(3)}_{443}$ (red).}
     \end{subfigure}
	 \centering
     \caption{Symmetry of the two probability tensors $\tQ^{(2)}$ and $\tQ^{(3)}$.}
     \label{fig:symmetry}
\end{figure}

Due to the symmetry of $\tQ^{(\ell)}$ shown in \Cref{fig:symmetry}, the probability of an edge $e \in \gE$ being present will depend only on the \emph{membership counts} of $\rvy(e)$, denoted by $\rvw =\rvw(e) = (\ervw_1, \ldots, \ervw_K)$, i.e., $\ervw_k$ nodes belong to $\gV_k$, making $\rvw$ a \emph{weak composition} of $\ell$.

\begin{definition}[Weak composition \cite{Bona2017AWalkTC}]\label{def:weak_composition}
	A sequence of non-negative integers $\rvw \coloneqq (\ervw_1, \ldots, \ervw_K)$, fulfilling $\sum_{k=1}^{K} \ervw_k = \ell$, is called a weak composition of $\ell$. Let $\WC{\ell}{K}$ denote the set of weak compositions of $\ell$ into $K$ parts, then the cardinality is $|\WC{\ell}{K}| = \binom{\ell + K - 1}{K-1}$.
\end{definition}

One can build $\tQ^{(\ell)}$ from the set of probabilities $\{ \etQ^{(\ell)}_{\rvw}\}_{\rvw \in \WC{\ell}{K}}$ for each $\ell \in \sL$, since the probability of an edge, by symmetry, depends only on its weak composition formed by membership counts. We will use $\tQ^{(\ell)}$ and $\etQ^{(\ell)}_{\rvw}$ alternately to facilitate presentation. The non-uniform HSBM can be built from uniform ones as illustrated in \cite{Ghoshdastidar2017ConsistencyOS}.
\begin{remark}\label{rem:symmetrytoWC} 
	For an $\ell$-hyperedge $e = \{i_1, \ldots, i_{\ell}\}$, let $w_k$ denote the number of vertices in $e$ with $\rY_{v} = k$, i.e., $w_k= |\{v\in e|\rY_{v} = k \}|$, then the sequence $\rvw(e) \coloneqq (\ervw_1, \ldots, w_K)$ constitutes a weak composition of $\ell$. Due to the symmetry of the probability tensor, i.e., $\etQ_{ \ervy_1, \ldots, \ervy_{\ell} }^{(\ell)} = \etQ_{ \ervy_{\pi(1)}, \ldots, \ervy_{\pi(\ell)} }^{(\ell)}$ for any $\pi \in S_{\ell}$, there exists a bijection, proved by counting, between $\rvw(e)$ and $\rvy(e)$. Recall that $\ell$-hyperedges $e_1$ and $e_2$ share the same generating probability if $\rvy(e_1) = \rvy(e_2)$, then the argument holds as well if $\rvw(e_1) = \rvw(e_2)$. Therefore, the set of $\ell$-hyperedges can be classified into different categories according to the weak compositions constituted from nodes in each edge.
\end{remark}

\begin{figure}[h]
     \centering
     \begin{subfigure}{0.3\textwidth}
         \centering
        \includegraphics[width=0.8\textwidth]{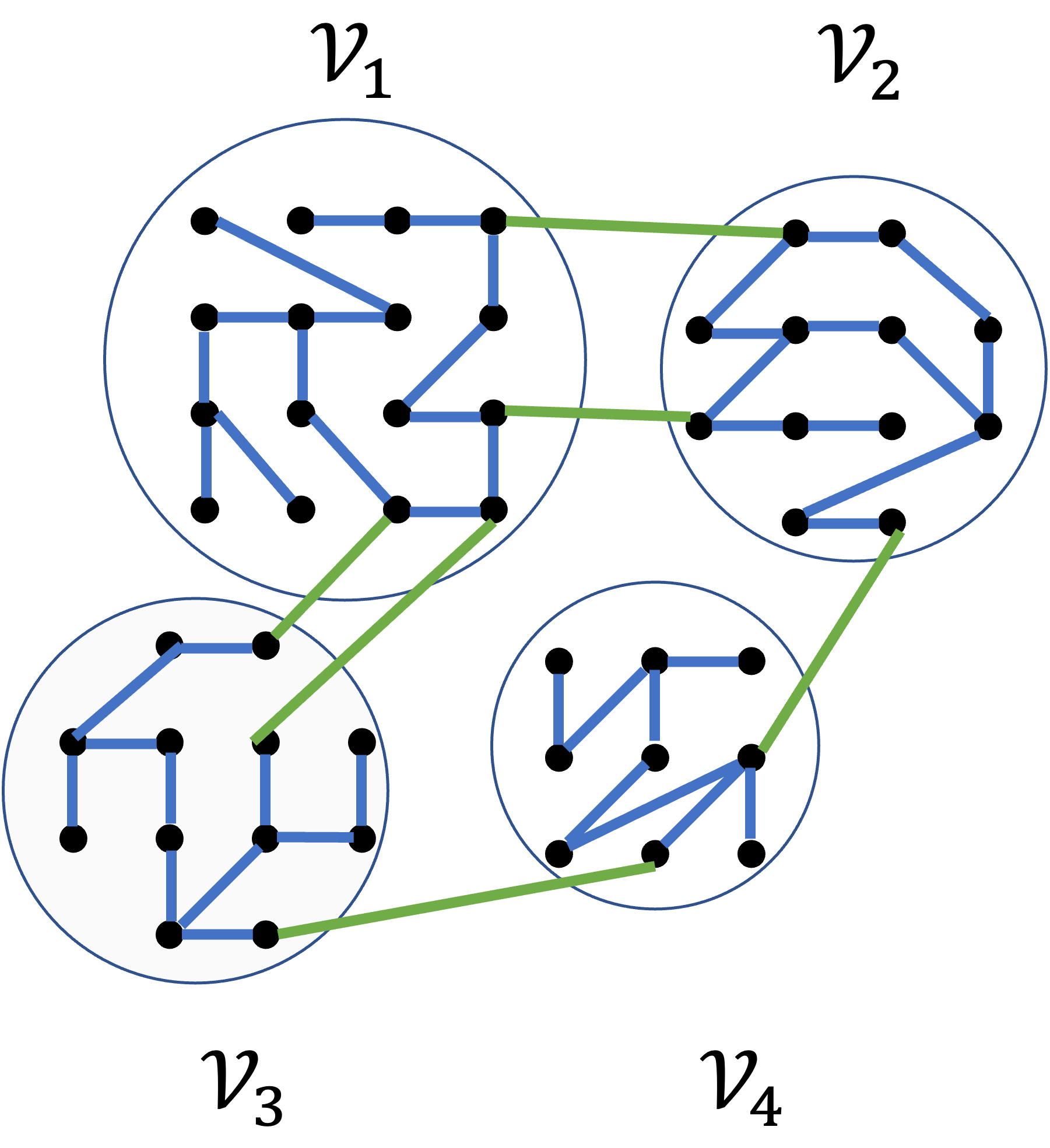}
        \subcaption{2-uniform HSBM}
     \end{subfigure}
     \begin{subfigure}{0.3\textwidth}
         \centering
        \includegraphics[width=0.8\textwidth]{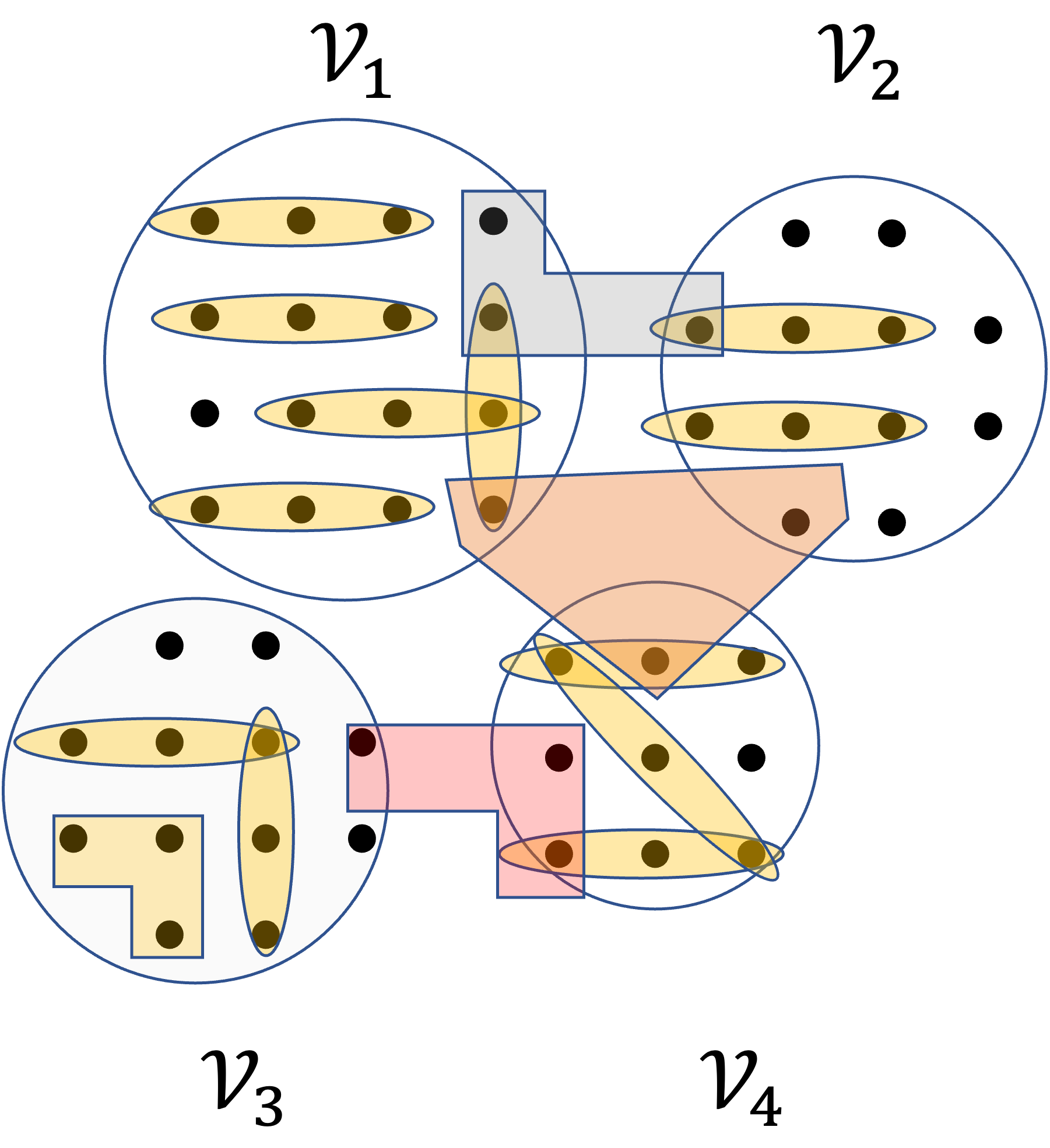}
        \subcaption{3-uniform HSBM}
     \end{subfigure}
     \begin{subfigure}{0.3\textwidth}
         \centering
        \includegraphics[width=0.8\textwidth]{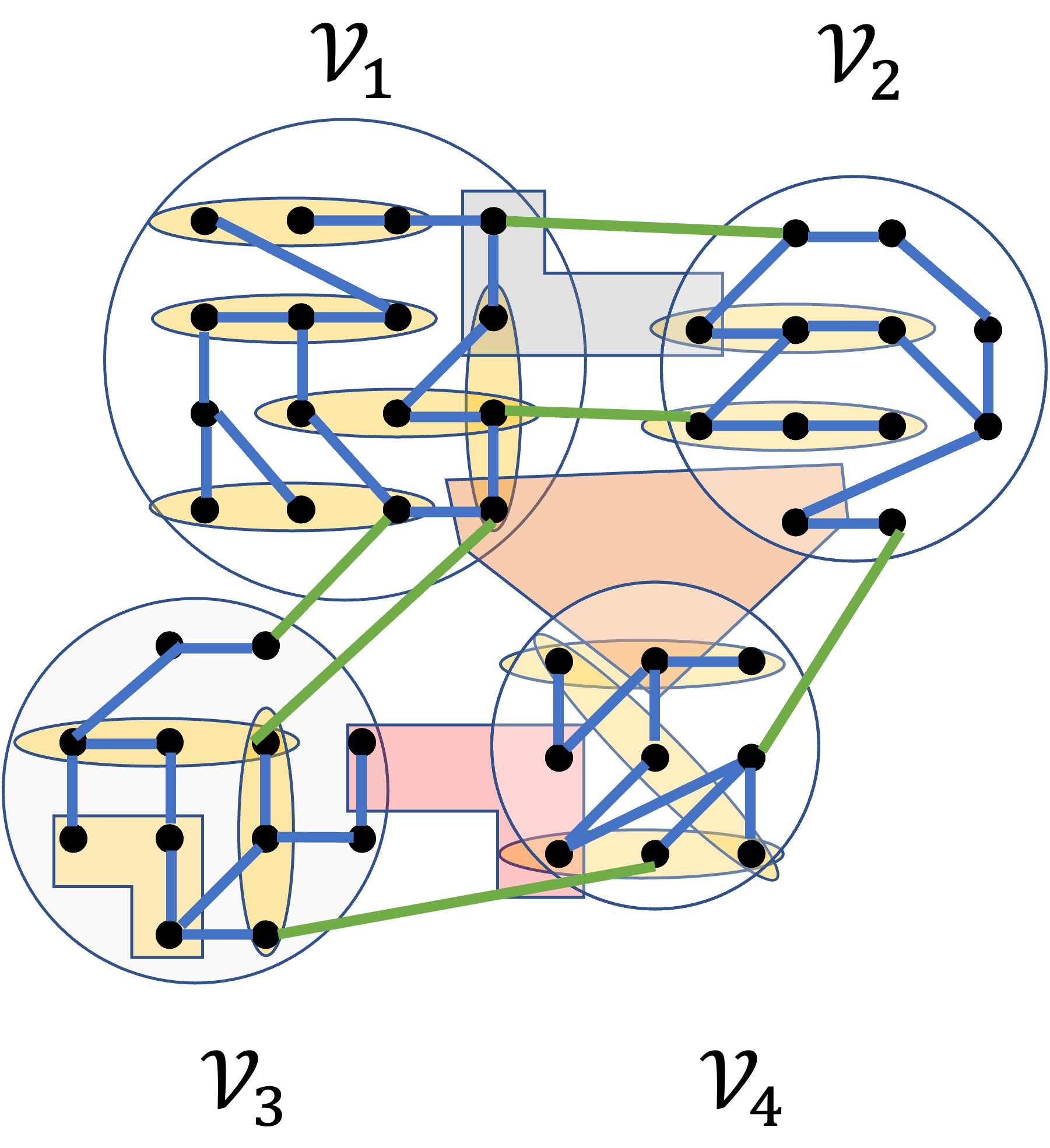}
        \subcaption{Non-uniform HSBM}
     \end{subfigure}
     \caption{Uniform \textnormal{HSBM} and non-uniform \textnormal{HSBM}.}
     \label{fig:non_uniform}
\end{figure}

\begin{definition}[Non-uniform HSBM]\label{def:non_uniform_HSBM} 
Let $\sL =\{\ell \mid \ell \geq 2, \ell \in \N\}$ be a set of integers with finite cardinality. The membership vector $\rvy$ is first sampled under $\rvalpha$, then for each $\ell \in \sL$, $\gH_{\ell}$ is independently drawn from $\textnormal{HSBM}_{\ell}(N,\rvalpha, \tQ^{(\ell)})$. The non-uniform hypergraph $\gH$ is a collection of $\ell$-uniform hypergraphs, i.e., $\gH = \cup_{\ell \in \sL} \gH_{\ell}$.
\end{definition}
Examples of $2$-uniform and $3$-uniform \textnormal{HSBM}, and an example of non-uniform \textnormal{HSBM} with $\sL = \{2, 3\}$ and $K=4$ can be seen in \Cref{fig:non_uniform} (a), \Cref{fig:non_uniform} (b) and \Cref{fig:non_uniform} (c) respectively.

Most of the computations involving tensors are NP-hard \cite{Hillar2013MostTP}. Instead, our analysis is based on the following key concept, the \emph{adjacency matrix}. 

\begin{definition}[Adjacency matrix]
For the non-uniform hypergraph $\gH$, let $\tA^{(\ell)}$ be the order-$\ell$ adjacency tensor corresponding to each underlying $\ell$-uniform hypergraph for $\ell \in\sL$. The adjacency matrix $\rmA \coloneqq [\ermA_{ij}]_{N \times N}$ of $\gH$ is defined by
\begin{align}\label{eqn:adjacency_matrix_entry}
    \ermA_{ij} \coloneqq \indi{i \neq j} \cdot \sum\limits_{\ell \in \sL}\,\,\,\,\sum\limits_{e\in \gE_{\ell},\,\, e \supset \{i, j\} }\etA^{(\ell)}_{e}\,,
\end{align}
where $\ermA_{ii} = 0$ for $i\in \gV$ since in our model each edge of size $\ell$ contains $\ell$ distinct vertices.
\end{definition}

For ease of presentation, we rewrite the notations to accommodate different regimes. Define a \emph{finite} family $\sF$ of growth rates $q:\N \mapsto \R_+$ with $\lim\limits_{N \to \infty} q(N) = \infty$, such that 
\begin{align}
\forall q_1, q_2 \in \sF,~~ \lim\limits_{N \to \infty} q_1(N)/q_2(N) \in \{0,+\infty\}~.
\end{align}
This means that $\sF$ is a family of growth rates of different magnitudes.
\begin{definition}[Edge density]
For model \ref{def:non_uniform_HSBM} , we write the probability of sampling edge $e=\{i_1, \ldots, i_{\ell}\}$ as
$\P( \etA_{e}^{(\ell)} = 1) = \etQ_{\rvy(e)}^{(\ell)} = \etQ^{(\ell)}_{\rvw}$, where $\rvy(e) \coloneqq [\ervy_{i_1}, \ldots, \ervy_{i_{\ell}}]$ is the membership vector of the vertices in $e$, and $\rvw\in \WC{\ell}{K}$ denotes the weak composition formulated by $\rvy(e)$. Furthermore, $\etQ^{(\ell)}_{\rvw}$ can be factorized as
\begin{align}
    \P( \etA_{e}^{(\ell)} = 1) = \etQ^{(\ell)}_{\rvw} = \etP_{\rvw }^{(\ell)} \cdot q_{\rvw }^{(\ell)}/\binom{N-1}{\ell - 1},\label{eqn:edge_density}
\end{align}
where $\etP_{\rvw }^{(\ell)}$ is some constant independent of $N$, and $q_{\rvw }^{(\ell)}\in \sF$. Let $q_{N}$ denote the slowest rate function among all $\{q_{\rvw }^{(\ell)}\}_{\rvw \in \WC{\ell}{K}, \ell \in \sL}$.
\end{definition}

\begin{assumption}\label{ass:prob_ratio_bound} 
There exists some absolute constant $\const_{\eqref{eqn:prob_ratio_bound} } > 1$ such that for every $N$, the probability tensors $\{\tQ^{(\ell)}\}_{\ell \in \sL}$ in model \ref{def:non_uniform_HSBM} satisfy
\begin{align} 
\const_{\eqref{eqn:prob_ratio_bound} } \coloneqq \max_{\ell \in \sL}\,\,\max_{\rvw, \rvw^{\prime} \in \WC{\ell}{K}}\, \etQ^{(\ell)}_{\rvw}/ \etQ^{(\ell)}_{\rvw^{\prime}},\label{eqn:prob_ratio_bound} 
\end{align}
which further implies $q_{\rvw }^{(\ell)} =q_{N}$ for each $\rvw \in \WC{\ell}{K}$ and $\ell \in \sL$.
\end{assumption}

\subsection{Information-Theoretic limits}
 We introduce the \emph{generalized Chernoff-Hellinger} (GCH) divergence to establish the \emph{Information-Theoretic} (IT) limits. 
\begin{definition}
For each $\rvw \in \WC{\ell - 1}{K}$, denote $k \oplus \rvw \coloneqq (\ervw_1, \ldots, w_{k-1}, w_k + 1, w_{k+1}, \ldots, w_K)$ and define the following quantity 
    \begin{align}\label{eqn:barnw}
        \overline{N}_{\rvw} \coloneqq \prod\limits_{j=1}^{K} \binom{ \lfloor \alpha_j N \rfloor }{\ervw_j}.
    \end{align}
Under \Cref{ass:prob_ratio_bound}, define the \emph{Generalized Chernoff-Hellinger} (GCH) divergence by
    \begin{align}
         &\,\D_{\rm{GCH}} = \underset{1\leq j < k \leq K}{\min} \D_{\rm{GCH}}(j,k),\label{eqn:GCH}\\
        \D_{\rm{GCH}}(j,k) \coloneqq &\, \max_{t\in [0, 1]} \sum\limits_{\ell \in \sL} \,\sum\limits_{\rvw\in \WC{\ell - 1}{K}} \frac{\overline{N}_{\rvw}}{\binom{N-1}{\ell - 1}} \big[ t\etP^{(\ell)}_{j \oplus \rvw} + (1 - t)\etP^{(\ell)}_{k \oplus \rvw} - \big(\etP^{(\ell)}_{j \oplus \rvw } \big)^{t} \cdot \big(\etP^{(\ell)}_{ k \oplus \rvw} \big)^{1 - t} \big].\notag
    \end{align}
For the more general scenario, where \Cref{ass:prob_ratio_bound} is no longer true, we define
\begin{align}
    &\,\D_{\rm{GCH}} = \underset{1\leq j < k \leq K}{\min} \D_{\rm{GCH}}(j,k) \label{eqn:GCH_not_proportional}\\
    \D_{\rm{GCH}}(j,k) \coloneqq &\, q_{N}^{-1}\max_{t\in [0, 1]} \sum\limits_{\ell \in \sL} \,\sum\limits_{\rvw\in \WC{\ell - 1}{K}} \overline{N}_{\rvw} \big[ t\etQ^{(\ell)}_{j \oplus \rvw} + (1 - t)\etQ^{(\ell)}_{k \oplus \rvw} - \big(\etQ^{(\ell)}_{j \oplus \rvw }\big)^{t} \cdot \big(\etQ^{(\ell)}_{ k \oplus \rvw} \big)^{1 - t} \big].\notag
\end{align}
Note that $\D_{\rm{GCH}}$ in \eqref{eqn:GCH} is a special case of \eqref{eqn:GCH_not_proportional} under \Cref{ass:prob_ratio_bound}. When $\const_{\eqref{eqn:prob_ratio_bound} } = 1$, leading to $\D_{\mathrm{GCH}} = 0$, HSBM is indistinguishable from a non-uniform Erd\H{o}s-R\'{e}nyi hypergraph. Therefore, the partition is undetectable. 
\end{definition}

The success of exact recovery is hindered by isolated vertices, which appear when $q_{N} =O(\log(N))$, as no algorithm could outperform random guess in the presence of isolated vertices in multiple communities. The necessary condition for exact recovery and Information-Theoretic lower bound for any algorithm are presented below.

\begin{theorem}\label{thm:IT_lower_bounds}
Let $\{\kappa_{N}\}_{N}$ be a sequence such that $\kappa_{N} \in (0, 1]$ for each $N$, and suppose
\begin{align}
   \lim\limits_{N \to \infty} \D_{\mathrm{GCH}}\cdot q_{N} /(\kappa_{N} \cdot \log(N)\,) < 1. \label{eqn:impossibleExactRecovery} 
\end{align} 
Then, given any community detection algorithm $\mathfrak{A}$, the following holds.

{\noindent}(1) With probability at least $1 - 2^{-q_{N}}$, $\mathfrak{A}$ will misclassify at least one vertex.

{\noindent}(2) When $\kappa_{N} \log(N) \to \infty$ as $N \to \infty$, it follows that $\lim\limits_{N \to \infty} (N^{\kappa_{N} } \cdot \E [\mismatch_{N}]) \geq 1$.
\end{theorem}
\begin{remark}
    Informally, (2) states that any algorithm $\mathfrak{A}$ is expected to misclassify at least $N^{1 - \kappa_{N}}$ nodes, with expectation over the sample hypergraph. 
\end{remark}

The necessary condition for exact recovery becomes $\D_{\rm{GCH}} < 1$ by taking $q_{N} = \log(N)$ and $\kappa_{N} = 1$ in \eqref{eqn:impossibleExactRecovery}. For SBM (graph)  and uniform HSBM, it was shown that exact recovery is impossible when $\D_{\rm{GCH}} < 1$ in \cite{Abbe2015CommunityDI} and \cite{Zhang2023ExactRI}, respectively, both under \Cref{ass:prob_ratio_bound}. By contrast, \Cref{thm:IT_lower_bounds} (1) is broader and its proof does not rely on \Cref{ass:prob_ratio_bound}. 

\Cref{thm:IT_lower_bounds} (2) represents a real breakthrough in the literature of HSBM. It states the lower bound on the expected mismatch ratio in regimes where exact recovery is impossible. Here are two examples: first, when $\kappa_{N} = 1$, $q_{N} = \log(N)$, $\E [\mismatch_{N}] \geq N^{-(1 + o(1))\D_{\rm{GCH}}}$; second, when $1\ll \D_{\rm{GCH}}\cdot q_{N} \ll \log(N)$ while $\kappa_{N} \to 0$ but $\kappa_{N} \log(N) \to \infty$, a lower bound on weak consistency is established. Similar results were previously shown for SBM in the binary \cite{Zhang2016MinimaxRO, Abbe2020EntrywiseEA, Abbe2022LPT} and multi-community \cite{Yun2016OptimalCR} cases. Note also that \cite{Gao2017AchievingOM, Gao2018CommunityDI} (SBM) and \cite{Chien2019MinimaxMR} (uniform HSBM) established minimax rates for mismatch ratio when certain constraints on the probability tensors $\{\tQ^{(\ell)}\}_{\ell \in \sL}$ are satisfied, which describes a subset of the entire problem space. By contrast, \Cref{thm:IT_lower_bounds} (2) is free of such constraints.

\subsection{Agnostic partition} Under \Cref{ass:prob_ratio_bound}, \Cref{alg:agnostic_partition} could achieve the exact recovery and reach the lowest expected mismatch ratio possible without knowing the probability tensors $\{\tQ^{(\ell)}\}_{\ell \in \sL}$. The number of communities $K$ is currently required as input. However, as will be discussed in \Cref{sec:communityEstimation}, one could estimate $K$ and turn \Cref{alg:agnostic_partition} into an algorithm \emph{entirely agnostic}.

There are two stages in \Cref{alg:agnostic_partition}: first, \Cref{alg:spectral_initialization} assigns all but a vanishing fraction of vertices correctly (weak consistency), then \Cref{alg:agnostic_refinement} refines the partition iteratively to reach the lowest possible mismatch ratio. We present Algorithms \ref{alg:spectral_initialization} and \ref{alg:agnostic_refinement} here, while \Cref{alg:trimming} is deferred to \Cref{sec:almostExact}.

\begin{algorithm}
\caption{\textbf{Agnostic partition}}\label{alg:agnostic_partition}
\KwData{The adjacency tensors $\{\tA^{(\ell)} \}_{\ell \in \sL}$, number of communities $K$.}


{Run \textbf{\Cref{alg:trimming} (Trimming)} on $\gH$ with $\sJ \coloneqq \sJ_1$ in \Cref{thm:regularization}  to obtain $\rmA_{\sJ_1}$}\;

{Run \textbf{\Cref{alg:spectral_initialization} (Initialization)} with input $\rmA_{\sJ_1}, K, \overline{d}$ to obtain the initial estimate $\widehat{\rvy}^{(0)}$}\;

{Run \textbf{\Cref{alg:agnostic_refinement} (Refinement)} with input $\widehat{\rvy}^{(0)}, K, \{\tA^{(\ell)} \}_{\ell \in \sL}$ and obtain the output $\widehat{\rvy}$}\;

\KwResult{$\widehat{\rvy}$}
\end{algorithm}

\subsubsection{Stage I} 
\begin{figure}[htbp]
     \centering
    \begin{subfigure}[b]{0.45\textwidth}
         \centering
         \includegraphics[width=0.9\textwidth]{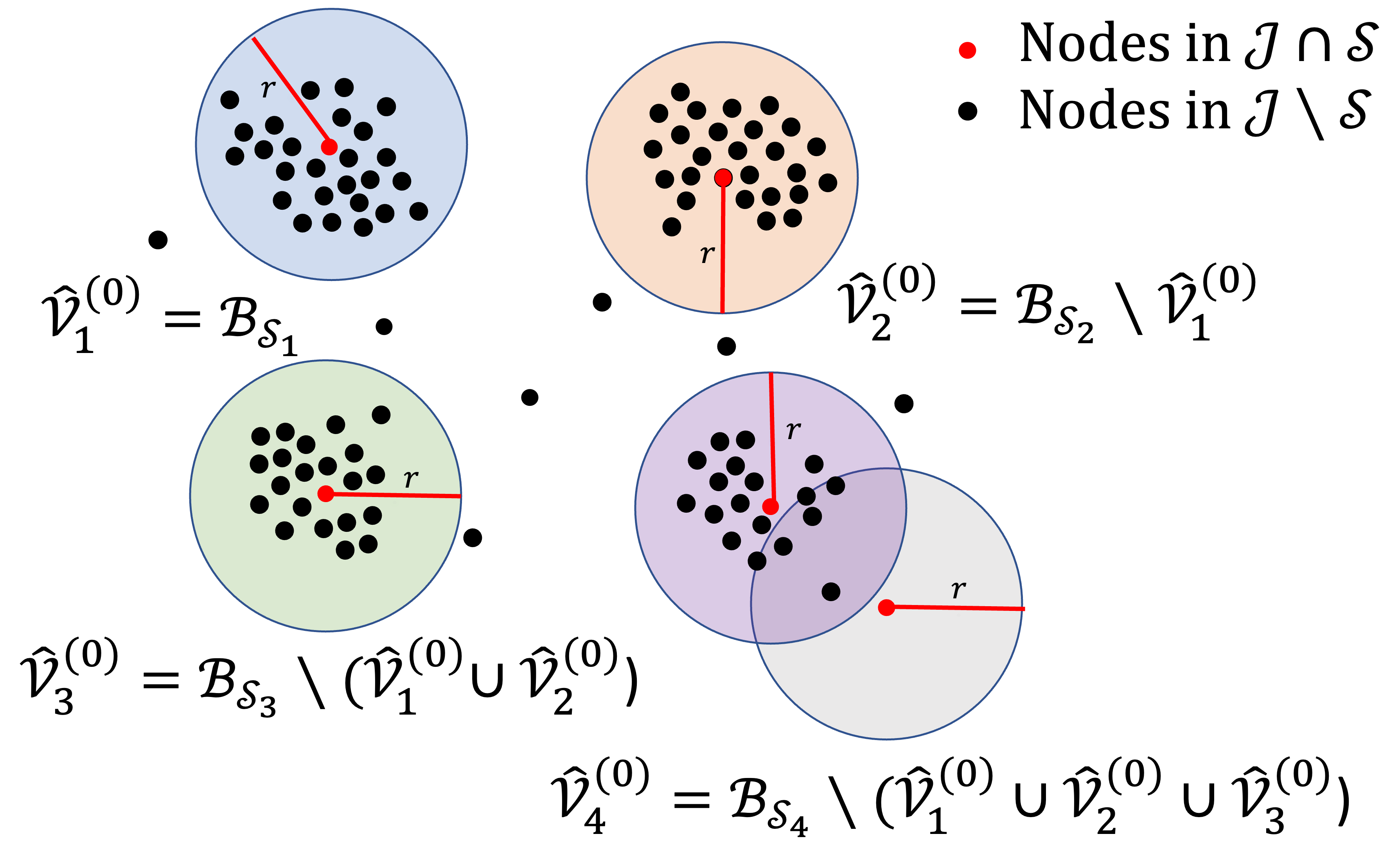}
         \caption{Line 2-7}
     \end{subfigure}
     \hspace{1cm}
    \begin{subfigure}[b]{0.38\textwidth}
         \centering
         \includegraphics[width=0.9\textwidth]{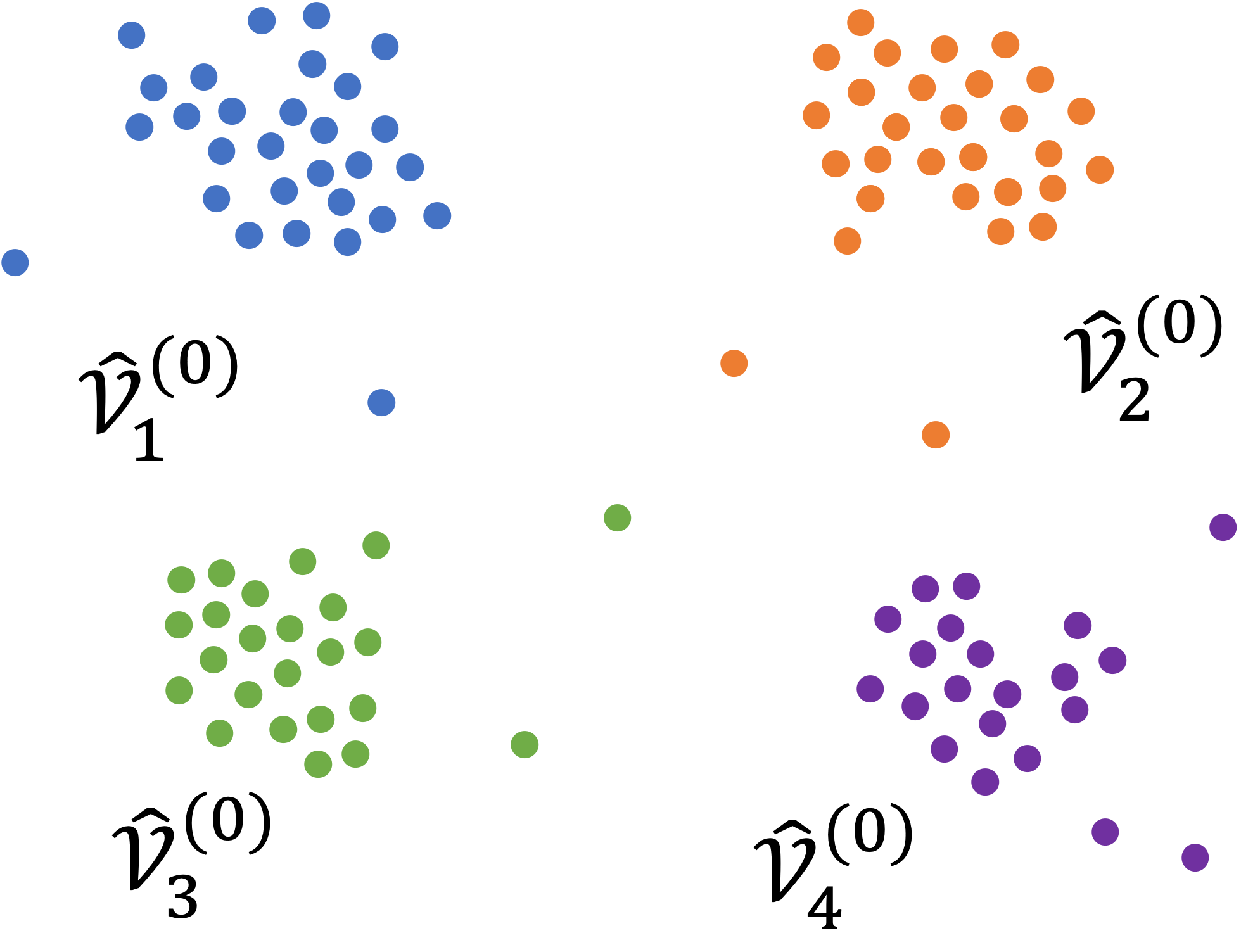}
         \caption{Line 8-12}
     \end{subfigure}
     \caption{Illustration of \Cref{alg:spectral_initialization}}\label{fig:spectral_initialization}
\end{figure}

Let $d_v$ denote the observed number of hyperedges containing $v$, and $\overline{d}$ denote the average degree of $\gH$. \Cref{fig:spectral_initialization} is an intuitive explanation of \Cref{alg:spectral_initialization}. 

\begin{algorithm}
\caption{\textbf{Spectral initialization}}\label{alg:spectral_initialization}
\KwData{$\rmA_{\sJ}$, number of communities $K$, radius $r = \overline{r}$, where $\overline{r} = [N\log(\overline{d})]^{-1}\overline{d}\,^2$ with $\overline{d} = \sum\limits_{v=1}^{N} d_v/N$.}

{Compute the rank-$K$ approximation $\rmA_{\sJ}^{(K)} = \sum\limits_{i=1}^{K}\lambda_i \rvu_i^{\sT} \rvu_i$ of $\rmA_{\sJ}$.}\ 

{Let $\setS$ be a set of $\lceil 2\log^2(N) \rceil$ nodes randomly sampled from $\sJ$ without replacement. For each $s\in \setS$, construct the ball centered at $s$ by $\ball_{r}(s) = \{w\in \sJ:  \|(\rmA_{\sJ}^{(K)})_{s:} - (\rmA_{\sJ}^{(K)})_{w:}\|_2^2 \leq r \}$.}\

{Take $\widehat{\gV}_{1}^{(0)} = \ball_{r}(s_1)$ to be the ball with most vertices, i.e., $s_{1} = \underset{s\in \setS}{\argmax}|\ball_{r}(s)|$. Break ties arbitrarily.}\

\While{$2 \leq k \leq K$ }{ \label{alg:first_loop}
    {$s_{k} = \argmax_{s\in \setS}| \ball_{r}(s)\setminus ( \bigcup_{j=1}^{k-1}\widehat{\gV}_{j}^{(0)})|$}\;
    {$\widehat{\gV}_{k}^{(0)} = \ball_{r}(s_k)\setminus ( \bigcup_{j=1}^{k-1}\widehat{\gV}_{j}^{(0)})$} \tcp*{Exclude the assigned vertices and find the remaining largest ball.}
}
\While{$v \in \sJ\setminus ( \bigcup_{k=1}^{K}\widehat{\gV}_{k}^{(0)})$ }{ \label{alg:second_loop}
    {$k = \argmin_{k\in [K]}\|(\rmA_{\sJ}^{(K)})_{s_k : } - (\rmA_{\sJ}^{(K)})_{v :}\|_2^2$}\;
    {$\widehat{\gV}_{k}^{(0)} \longleftarrow \widehat{\gV}_{k}^{(0)}\cup \{v\}$} \tcp*{Assign the remaining vertices to their nearest ball.}
}
{Randomly assign the remaining vertices $v\in \gV\setminus \sJ$ to one of the communities $\widehat{\gV}_{1}^{(0)},\ldots,\widehat{\gV}_{K}^{(0)}$}\; \label{alg:assign_vertices_outside_J}
{Obtain the initial estimate of the membership vector $\widehat{\rvy}^{(0)}$ based on $\widehat{\gV}_{1}^{(0)},\ldots,\widehat{\gV}_{K}^{(0)}$.}

\KwResult{$\widehat{\rvy}^{(0)}$}
\end{algorithm}

\begin{theorem}[Agnostic weak consistency]\label{thm:agnostic_weak_consistency}
Let $\rho_{N}$ denote the maximum expected degree among all vertices, see \eqref{eqn:rhoxi}. Suppose Assumptions \ref{ass:prob_ratio_bound}, \ref{ass:expected_center_separation} and $\rho_{N} \gtrsim \log(N)$, then with probability at least $1-2(e/2)^{-N} -2N^{-c}$ for some constant $c > 0$, \Cref{alg:spectral_initialization} achieves weak consistency by applying regularization $\sJ_1$ in \Cref{thm:regularization}. 
\end{theorem}
\begin{assumption}\label{ass:expected_center_separation} 
For each $\rvw \in \WC{\ell - 2}{K}$, let 
\begin{align}
    j \oplus i \oplus \rvw \coloneqq (\ervw_1, \ldots, \ervw_{j} + 1, \ldots, \ervw_{i} + 1, \ldots, \ervw_K )
\end{align}
denote a weak composition in $\WC{\ell}{K}$. For every distinct pair $j, k \in [K]$, there exists some $i\in [K]$ such that $\{\tQ^{(\ell)}\}_{\ell \in \sL}$ satisfies
        \begin{align}\label{eqn:expected_center_separation}
            \lim\limits_{N \to \infty} \frac{N}{\rho_{N}} \cdot \bigg| \sum\limits_{\ell \in \sL} \sum\limits_{\rvw \in \WC{\ell - 2}{K} } \overline{N}_{\rvw} \cdot (\etQ^{(\ell)}_{j \oplus i \oplus \rvw }  - \etQ^{(\ell)}_{k \oplus i \oplus \rvw} ) \bigg| \neq 0.
        \end{align}
\end{assumption}
\begin{remark}\label{rem:violate_center_separation}
Denote $\alpha^{\rvw} \coloneqq \prod_{j=1}^{K} (\alpha_j)^{\ervw_j}$ for $\rvw = (\ervw_1, \ldots, \ervw_K)$. Under \Cref{ass:prob_ratio_bound}, if there exists distinct pair $j, k\in [K]$ such that for every $i\in [K]$,
 \begin{align}
    \sum\limits_{\ell \in \sL} (\ell - 1) \sum\limits_{\rvw \in \WC{\ell - 2}{K} } \binom{\ell - 2}{\rvw}\cdot \alpha^{\rvw}  \cdot (\etP^{(\ell)}_{j \oplus i \oplus \rvw }  - \etP^{(\ell)}_{k \oplus i \oplus \rvw} ) =0, \label{eqn:expectedCenterSame}
 \end{align}
where $\binom{\ell - 2}{\rvw}\coloneqq (\ell - 2)!/(\prod_{j=1}^{K}\ervw_j!)$, satisfying $\sum_{j=1}^{K} \ervw_j = \ell - 2$, denotes the multinomial coefficient, then \Cref{ass:expected_center_separation} is violated. The set of tensors satisfying \eqref{eqn:expectedCenterSame} is a linear proper subspace of the entire parameter space. In the absence of \Cref{ass:prob_ratio_bound}, \Cref{ass:expected_center_separation} says that the expected degrees of each community, as well as the differences between communities, have to be proportional to the highest expected degree (potentially with the exception of one community). It remains true that, under \emph{HSBM}, the set of problems that cannot be solved by \Cref{alg:spectral_initialization} is \emph{asymptotically negligible} in the entire problem space.
\end{remark}

\Cref{ass:expected_center_separation} is a natural extension of the condition in \cite{Abbe2015CommunityDI} for the graph case $\sL = \{2\}$, and it also appears in \cite{Zhang2023ExactRI} under \Cref{ass:prob_ratio_bound}. However, while in the graph case this condition is truly necessary for achievability (since its absence creates two classes that are indistinguishable from each other), it is not clear whether this is the hypergraph case for \Cref{ass:expected_center_separation}. It is certainly needed for any method using the adjacency matrix to succeed, but, as shown in \cite{Zhang2023ExactRI}, due to the many parameters involved, one may construct a model with distinguishable classes via probability tensors (\Cref{fig:symmetry} (b)), but for which the expected $K\times K$ adjacency matrix has identical rows (this uniform case is yet to be solved). 

\subsubsection{Stage II}  Let $\{\widehat{\gV}^{(t)}_{k}\}_{k=1}^{K}$ denote the partitioned communities in iteration $t$ with the initial partition $\{\widehat{\gV}^{(0)}_{k}\}_{k=1}^{K}$ obtained from \Cref{alg:spectral_initialization}. For each $\ell \in \sL$ and $\rvw \in \WC{\ell}{K}$, let
\begin{align}\label{eqn:hat_Qlw}
    \widehat{\etQ}^{(\ell)}_{\rvw} \coloneqq \frac{1}{ \widehat{N}^{(0)}_{\rvw}} \sum\limits_{ \widehat{\rvy}^{(0)}(e) = \rvw} \etA^{(\ell)}_{e}\,, \quad \textnormal{  where  } \widehat{N}^{(t)}_{\rvw} = \prod\limits_{j=1}^{K} \binom{|\widehat{\gV}^{(t)}_{j}| }{\ervw_j}, \quad \forall t \geq 0,
\end{align}
be an estimation of $\etQ^{(\ell)}_{\rvw}$ in \eqref{eqn:edge_density}. Let $\widehat{\rD}^{(t)}_{v, \rvw}$ denote the number of $\ell$-hyperedges containing $v$ with the rest nodes distributed as $\rvw \in \WC{\ell - 1}{K}$ at iteration $t \geq 0$, formally,
\begin{align}\label{eqn:gEtvw}
    \widehat{\rD}^{(t)}_{v, \rvw} \coloneqq \, \big| \{e \subset [\gV]^{\ell} \mid e\ni v,\,\, \rvy(e\setminus\{v\}) = \rvw \in \WC{\ell - 1}{K}, \,\, \ervw_{j} \textnormal{ nodes in } \widehat{\gV}^{(t)}_{j}\}\, \big|.
\end{align}
Below, we provide an intuitive explanation of \Cref{alg:agnostic_refinement}. Suppose $v \in \gV_k$ and define $\widehat{\mu}^{(t)}_{v, \rvw}\coloneqq \widehat{\rD}^{(t)}_{v, \rvw}/ \widehat{N}^{(t)}_{\rvw}$, which is an estimation of $\etQ^{(\ell)}_{k\oplus \rvw}$ only using edges containing $v$. By contrast, $\widehat{\etQ}^{(\ell)}_{k \oplus \rvw}$ in \eqref{eqn:hat_Qlw} is a more accurate estimation of $\etQ^{(\ell)}_{k\oplus \rvw}$ since all edges of type $\widehat{\rvy}^{(0)}(e) = k \oplus \rvw$ are involved, enjoying sufficient concentration due to \Cref{lem:hatQmw_approx}. Let $\mathbb{H}_{\mathrm{CE}}(y, \widehat{y}) \coloneqq - y\log \widehat{y} - (1 - y) \log (1 - \widehat{y})$ denote the \emph{cross entropy} between the two probability distributions with means $y$ and $\widehat{y}$, respectively. Then \eqref{eqn:GKLEstimator} can be reformulated as
\begin{align*}
    k^{\star} =&\, \underset{1 \leq k \leq K}{\argmin} \bigg\{ \sum\limits_{\ell \in \sL} \,\, \sum\limits_{ \rvw \in \WC{\ell - 1}{K} } \widehat{N}^{(t)}_{\rvw} \cdot \Big[- \widehat{\mu}^{(t)}_{v, \rvw} \cdot \log\widehat{\etQ}^{(\ell)}_{k \oplus \rvw} - (1 - \widehat{\mu}^{(t)}_{v, \rvw}) \cdot \log( 1 - \widehat{\etQ}^{(\ell)}_{k \oplus \rvw} ) \Big] \bigg\}\\
    =&\, \underset{1 \leq k \leq K}{\argmin} \,\, \sum\limits_{\ell \in \sL} \,\, \sum\limits_{ \rvw \in \WC{\ell - 1}{K} } \widehat{N}^{(t)}_{\rvw} \cdot \mathbb{H}_{\mathrm{CE}}(\widehat{\mu}^{(t)}_{v, \rvw}, \widehat{\etQ}^{(\ell)}_{k \oplus \rvw}).
\end{align*}
Note that $\mathbb{H}_{\mathrm{CE}}(y, \widehat{y}) = 0$ if and only if $y= \widehat{y}$, and it would be relatively large when $v$ is misclassified. Thus, the correct assignment of $v$ can be found by minimizing the weighted cross-entropy terms. The theoretical guarantee is provided below.

\begin{algorithm}
\caption{\textbf{Agnostic refinement}}\label{alg:agnostic_refinement}
\KwData{$\{\tA^{(\ell)}\}_{\ell \in \sL}$, $\widehat{\rvy}^{(0)}$.}

For each $\ell \in \sL$, $\rvw \in \WC{\ell}{K}$, calculate $\widehat{\etQ}^{(\ell)}_{\rvw}$ in \eqref{eqn:hat_Qlw}. \label{alg:prob_estimate}

\While{$0 \leq t\leq \lceil\log(N) \rceil$}{
    {$\widehat{\gV}^{(t+1)}_k \leftarrow \emptyset$ for all $k\in [K]$}\;
    \While{$v\in \gV$}{
    For each $\rvw \in \WC{\ell - 1}{K}$, calculate $\widehat{\rD}^{(t)}_{v, \rvw}$ in \eqref{eqn:gEtvw}. Find \label{alg:GKL_estimator}
        \begin{align}\label{eqn:GKLEstimator}
            k^{\star}= \underset{1 \leq k \leq K}{\argmax} \bigg\{ \sum\limits_{\ell \in \sL} \,\, \sum\limits_{ \rvw \in \WC{\ell - 1}{K} } \widehat{\rD}^{(t)}_{v, \rvw} \cdot \log\widehat{\etQ}^{(\ell)}_{k \oplus \rvw} + ( \widehat{N}^{(t)}_{\rvw} - \widehat{\rD}^{(t)}_{v, \rvw}) \cdot \log( 1 - \widehat{\etQ}^{(\ell)}_{k \oplus \rvw} ) \bigg\},
        \end{align}
    {Tie broken uniformly at random}\;
    
    {$\widehat{\gV}^{(t+1)}_{k^{\star}} \leftarrow \widehat{\gV}^{(t+1)}_{k^{\star}} \cup \{v\}$.}
    }
    {$t \leftarrow t + 1$}\
}
    {Take $\widehat{\gV}_k = \widehat{\gV}^{(\lceil\log(N) \rceil)}_k$ for all $k\in [K]$ and obtain $\widehat{\rvy}$.}\
    
\KwResult{$\widehat{\rvy}$}
\end{algorithm}

\begin{theorem}[Agnostic strong consistency]\label{thm:agnostic_strong_consistency}
For some absolute $\epsilon > 0$, suppose
\begin{align}\label{eqn:exact_recovery_condition}
    \lim\limits_{N \to \infty} \D_{\mathrm{GCH}}\cdot q_{N}/\log(N) \geq 1 + \epsilon,
\end{align}
together with Assumptions \ref{ass:prob_ratio_bound} and \ref{ass:expected_center_separation}. Then with probability at least $1 -6N^{-\epsilon}$, \Cref{alg:agnostic_partition} achieves exact recovery without prior knowledge of $\{\tQ^{(\ell)}\}_{\ell \in \sL}$.
\end{theorem}

 When $q_{N} = \log(N)$, i.e. $\D_{\mathrm{GCH}} > 1$, \cite{Abbe2015RecoveringCI, Yun2016OptimalCR} (SBM) achieved exact recovery agnostically. By contrast, \cite{Zhang2023ExactRI} completed exact recovery on uniform HSBM by requiring prior knowledge of probability tensors. Therefore, our \Cref{thm:agnostic_strong_consistency} is more general than any other result in the current HSBM literature.
 
In addition, \eqref{eqn:exact_recovery_condition} indicates that the structural information of the non-uniform hypergraph is stronger than the uniform one. Exact recovery may be achievable when all uniform layers are aggregated, even if it is impossible when any of the uniform layers is used in isolation. An example with explicit formula is discussed in \cite{Wang2023ITLimits}.

\subsubsection{Optimality}
When \eqref{eqn:exact_recovery_condition} does not hold, strong consistency is impossible by \Cref{thm:IT_lower_bounds} (1). In this regime, one can only hope for weak consistency, and \Cref{thm:IT_lower_bounds} (2) provides a lower bound on the expected mismatch ratio. The following theorem shows that \Cref{alg:agnostic_partition} is optimal, as it achieves this IT lower bound without knowing $\{\tQ^{(\ell)}\}_{\ell \in \sL}$.

\begin{theorem}[Agnostic optimality]\label{thm:agnostic_optimality}
For the sequence $\kappa_{N} \in (0, 1]$ with $\kappa_{N} \log(N) \to \infty$ as $N \to \infty$ and some absolute (small) $\epsilon >0$, suppose Assumptions \ref{ass:prob_ratio_bound}, \ref{ass:expected_center_separation} hold, $\rho_{N} \gtrsim \log(N)$ and
\begin{align}
    \lim\limits_{N \to \infty} \D_{\mathrm{GCH}}\cdot q_{N} /(\kappa_{N}\cdot \log(N)\, ) \geq 1 + \epsilon.\label{eqn:optimality_condition}
\end{align}
Then with probability at least $1 -6e^{-\epsilon \cdot \kappa_{N} \log(N)}$, \Cref{alg:agnostic_partition} satisfies $\lim\limits_{N \to \infty} (N^{\kappa_{N}} \cdot \mismatch_{N}) \leq 1$.
\end{theorem}

\subsection{Partition with prior knowledge}

When the probability tensors $\{\tQ^{(\ell)}\}_{\ell \in \sL}$ are known, the regime beyond \Cref{ass:prob_ratio_bound} can be handled by \Cref{alg:partition_with_prior}, where the optimal membership assignment for each $v\in \gV$ is determined through a \emph{leave-one-out} pipeline.

\begin{algorithm}
\caption{\textbf{Partition with prior knowledge}}\label{alg:partition_with_prior}

\KwData{$\{\tA^{(\ell)}\}_{\ell \in \sL}$, $K$, probability tensors $\{\tQ^{(\ell)}\}_{\ell \in \sL}$, $\{\alpha_k\}_{k=1}^{K}$.}

\While{$v\in \gV$}{
    {Obtain $\gH_{-v}$ by removing all edges containing $v$ from $\gH$;}\

    {Run \textbf{\Cref{alg:spectral_initialization} (Initialization)} on $\gH_{-v}$ to obtain the estimation $\widehat{\rvy}^{(0)}_{-v}$ for $\rvy_{-v}$;}\label{alg:spectral_MAP}
    
    {$\widehat{\ervy}_v = \underset{k \in [K]}{\argmax} \,\, \P(\rD_v = d_v | \rY_{v} = k,\, \rmY_{-v} = \widehat{\rvy}_{-v}^{(0)} )\cdot \alpha_k $.}\label{alg:MAP_in_Alg} \tcp*{Local correction on $v$.}
}

\KwResult{$\widehat{\rvy}$}
\end{algorithm}

The essential idea of \Cref{alg:partition_with_prior} is to make full use of the \emph{Maximum A Posteriori} (MAP) estimation. We provide an intuitive explanation of the algorithm below. Let $\rvy_{-v} \in [K]^{N-1}$ denote the membership vector of the vertices other than $v$ and $\widehat{\rvy}_{-v}$ be an estimation of $\rvy_{-v}$. Let $\mathbb{H}$ and $\rmY_{-v}$ denote the laws of $\gH$ and $\rvy_{-v}$ respectively. Conditioned on events $\{\mathbb{H} = \gH \}$ and $\{ \rmY_{-v} = \widehat{\rvy}_{-v}\}$, one can prove (deferred to \Cref{sec:known_partition}) that the probability of $v$ being misclassified is minimized when applying the MAP estimator
\begin{align}\label{eqn:MAP}
    \widehat{\ervy}_{v}^{\rm{\, MAP}}\coloneqq \underset{k \in [K]}{\argmax} \,\, \P(\rY_{v} = k \mid \mathbb{H} = \gH, \,{\rmY}_{-v} = \widehat{\rvy}_{-v} )\,.
\end{align}
Similar ideas appeared in \cite{Abbe2015CommunityDI, Abbe2018CommunityDA} for graph SBM, and \cite{Chien2019MinimaxMR, Zhang2023ExactRI} for uniform HSBM. However, the posterior probability in \eqref{eqn:MAP} cannot be directly computed from the definition. Recall that $d_v$ is the observed number of hyperedges in $\gH$ containing $v$ with $\rD_v$ denoting its law. We propose \eqref{eqn:MAP_single}, an estimator equivalent to \eqref{eqn:MAP} (as proved by \Cref{lem:MAP}), which can be computed directly when knowing $\{\tQ^{(\ell)}\}_{\ell \in \sL}$ and $\{\alpha_k\}_{k=1}^{K}$. 
\begin{lemma}\label{lem:MAP}
When the estimation $\widehat{\rvy}_{-v}$ is independent of the membership assignment $\ervy_v$, the {\rm{MAP}} in \eqref{eqn:MAP} can be reformulated as follows
 \begin{align}
    \widehat{\ervy}_{v}^{\rm{\, MAP}} = &\, \underset{k \in [K]}{\argmax} \,\, \P(\rD_v = d_v \mid \rY_{v} = k,\, \rmY_{-v} = \widehat{\rvy}_{-v} )\cdot \P( \rY_{v} = k)\,.\label{eqn:MAP_single}
 \end{align}
\end{lemma}
In \Cref{alg:spectral_MAP} of \Cref{alg:partition_with_prior}, the initial estimation $\widehat{\rvy}^{(0)}_{-v}$ is obtained by running \Cref{alg:spectral_initialization} on $\gH_{-v}$. It is independent of $\ervy_{v}$ since no edges containing $v$ were involved throughout the process. Therefore by \Cref{lem:MAP}, the $\widehat{\ervy}_v$ in \Cref{alg:MAP_in_Alg} is equivalent to the best estimator MAP in \eqref{eqn:MAP_single}. 

In \cite{Abbe2015CommunityDI, Abbe2018CommunityDA}, the authors presented a two-stage strategy, where the graph $\gH$ is divided into $\gH^{(0)}$ and $\gH^{(1)}$. In the first stage, the sparser graph $\gH^{(0)}$ is used to obtain the initial estimation $\widehat{\rvy}^{(0)}$, then $\widehat{\rvy}^{(0)}_{-v}$ is obtained by removing $\widehat{\ervy}^{(0)}_v$. In the second stage, \eqref{eqn:MAP} is applied on the denser graph $\gH^{(1)}$ for local refinement. We have found their proof difficult to duplicate for the hypergraph case. In \cite{Zhang2023ExactRI}, the authors presented a similar strategy (complexity $KN^2$) and claimed (without proof) that it can be modified to work for the hypergraph. For simplicity and completion, we have chosen instead to go with a slightly slower Algorithm \ref{alg:partition_with_prior} (complexity $KN^3$), which forgoes the splitting but ensures the independence between $\ervy_v$ and $\widehat{\rvy}_{-v}$. 
\subsubsection{Stage I}
The agnostic result \Cref{thm:agnostic_weak_consistency} requires $\rho_{N} \gtrsim \log(N)$. However, with prior knowledge, \Cref{alg:spectral_initialization} can still achieve weak consistency as long as $\rho_{N} = \omega(1)$.

\begin{theorem}[Weak consistency with prior knowledge]\label{thm:known_weak_consistency}
Suppose \Cref{ass:expected_center_separation} and $\rho_{N} = \omega(1)$. If we know $\{\tQ^{(\ell)}\}_{\ell \in \sL}$ and apply $\sJ_2$ with $r = [N\log(\rho_{N})]^{-1}\rho_{N}^2$ in \Cref{thm:regularization}, \Cref{alg:spectral_initialization} achieves weak consistency with probability $1-2(e/2)^{-N} -N^{-10}$.
\end{theorem}

With knowledge of the probability tensors $\{\tQ^{(\ell)}\}_{\ell \in \sL}$, \Cref{ass:expected_center_separation} can be weakened. It is equivalent to saying that no two rows of the expected $K\times K$ adjacency matrix for the communities, when scaled by the highest degree, are asymptotically equal, as in \Cref{rem:violate_center_separation}. In fact, this can be boiled down to an "all the layers" condition, where \eqref{eqn:expected_center_separation} happens on a uniform layer by uniform layer basis. The set satisfying \eqref{eqn:expectedCenterSame} can be shrunk as follows. When \Cref{ass:expected_center_separation} fails for the non-uniform hypergraph but holds true for at least one uniform layer, we can construct a weighted adjacency matrix $\widetilde{\rmA} = \sum_{\ell \in \sL}w_{\ell} \rmA^{(\ell)}$ as an input for \Cref{alg:spectral_initialization} under some properly chosen weights $\{w_{\ell}\}_{\ell \in \sL}$ to meet \Cref{ass:expected_center_separation}, and use it to obtain an almost exact labeling in the first stage. Importantly, the second stage refinement, which gives the threshold, reverts to using tensors $\{\tA^{(\ell)}\}_{\ell \in \sL}$ and does not rely on \Cref{ass:expected_center_separation}.  Therefore at the end of stage two, provided that the problem is above the threshold, the algorithm will output an exact labeling with high probability despite the introduction of weights. Thus, the only cases not covered by \Cref{thm:known_weak_consistency} lie at the intersection of these unsolvable (so far) uniform cases.

\subsubsection{Stage II}

Compared to \Cref{alg:agnostic_refinement}, $\widehat{\ervy}_v$ in \Cref{alg:MAP_in_Alg} is able to handle the regime beyond \Cref{ass:prob_ratio_bound}. For each $\ell \in \sL$, denote $\etQ^{(\ell)}_{\max} \coloneqq \max_{\rvw\in \WC{\ell}{K}} \etQ^{(\ell)}_{\rvw}$, $\etQ^{(\ell)}_{\min} \coloneqq \min_{\rvw\in \WC{\ell}{K}} \etQ^{(\ell)}_{\rvw}$, and we write 
$\etQ^{(\ell)}_{\max} = \etP^{(\ell)}_{\max}\cdot q^{(\ell)}_{\max}/\binom{N-1}{\ell - 1}$, $\etQ^{(\ell)}_{\min} = \etP^{(\ell)}_{\min}\cdot q^{(\ell)}_{\min}/\binom{N-1}{\ell - 1}$. The theoretical guarantee of \Cref{alg:partition_with_prior} is provided below.

\begin{theorem}[Strong consistency with prior knowledge]\label{thm:known_strong_consistency}
Suppose \Cref{ass:expected_center_separation},
\begin{align}
    \log(q^{(\ell)}_{\max}/q^{(\ell)}_{\min}) /\log(q^{(\ell)}_{\max}) \ll 1,\label{eqn:Strong_consistency_prob_condition}
\end{align}
and \eqref{eqn:exact_recovery_condition} for some absolute $\epsilon > 0$. When the tensors $\{\tQ^{(\ell)}\}_{\ell \in \sL}$ are known, \Cref{alg:partition_with_prior} achieves exact recovery with probability at least $1 - N^{-\epsilon} - N^{-10}$.
\end{theorem}
An example satisfying \eqref{eqn:Strong_consistency_prob_condition} but violating \Cref{ass:prob_ratio_bound} is $q^{(\ell)}_{\rvw_1} = N^{1/6}$ while $q^{(\ell)}_{\rvw} = N^{1/6}/\log(N)$ for any $\rvw \neq \rvw_1$. The growth rates could differ only up to ploy logrithmic 

When $q_{N} = \log(N)$, the condition \eqref{eqn:exact_recovery_condition} can be simplified as $\D_{\mathrm{GCH}} > 1$. However, \eqref{eqn:exact_recovery_condition} could still be true when $1\ll q_{N} \ll \log(N)$ and $\D_{\mathrm{GCH}} \gg 1$ as in \eqref{eqn:GCH_not_proportional}. Compared with \Cref{thm:agnostic_strong_consistency}, the additional situation covered by \Cref{thm:known_strong_consistency} happens when one of the communities has some sub-critical connectivity probabilities $o(\log(N)/\binom{N-1}{\ell - 1})$ either within itself or with respect to the other communities, while all the others are super-critically connected to each other, which will means that $\D_{\mathrm{GCH}} \cdot q_{N} = \Omega(\log(N))$ satisfying \eqref{eqn:exact_recovery_condition} by \Cref{lem:GCH_max_different_order}. In graph SBM \cite{Abbe2015CommunityDI} and uniform HSBM \cite{Zhang2023ExactRI}, exact recovery was achieved only under \Cref{ass:prob_ratio_bound} when knowing the probability tensors. \Cref{thm:known_strong_consistency} removes this constraint and proves that \Cref{alg:partition_with_prior} works well under a broader regime.

\subsubsection{Optimality}
Suppose that we have the access to a weakly consistent community partition, the extra knowledge on the probability tensors helps \Cref{alg:partition_with_prior} to achieve the lowest possible expected mismatch ratio if $\rho_{N} = \omega(1)$. Compared with \Cref{thm:agnostic_optimality}, \Cref{thm:known_optimality} covers a wider scenario including $\log^{-1}(N) \ll \kappa_{N} \ll 1$ with a better convergence rate.

\begin{theorem}[Optimality with prior knowledge]\label{thm:known_optimality}
Suppose \Cref{ass:expected_center_separation}, $\rho_{N} = \omega(1)$, \eqref{eqn:Strong_consistency_prob_condition} and \eqref{eqn:optimality_condition}, where the sequence $\kappa_{N} \in (0, 1]$ satisfies $\kappa_{N} \log(N) \to \infty$ as $N \to \infty$. Then with probability at least $1 - N^{-10}$, \Cref{alg:partition_with_prior} satisfies $\lim\limits_{N \to \infty} (N^{\kappa_{N}} \cdot \E[\mismatch_{N}]) \leq 1$.
\end{theorem}

\subsection{Related literature}
\subsubsection{Graphs}
For the binary symmetric $\mathrm{SBM}(N, 2, a\frac{\log(N)}{N}, b\frac{\log(N)}{N})$, the threshold is $\D_{\rm{H}}(a, b) \coloneqq (\sqrt{a} - \sqrt{b})^2/2$, which can be obtained from \eqref{eqn:Dv} by taking $\sL =\{2\}$ and $K = 2$. Exact recovery can be achieved efficiently if $\D_{\rm{H}}(a, b) \geq 1$, but is impossible when $\D_{\rm{H}}(a, b)<1$ \cite{Abbe2016ExactRI, Mossel2016ConsistencyTF}. For multi-block $\mathrm{SBM}(N, \rvalpha, \rmP \frac{\log(N)}{N})$, exact recovery can be achieved efficiently when $\D_{\rm{CH}}(\rvalpha, \rmP) \geq 1$, but is impossible when $\D_{\rm{CH}}(\rvalpha, \rmP) < 1$ \cite{Abbe2015CommunityDI, Abbe2015RecoveringCI, Yun2016OptimalCR}, where 
\begin{align*}
   \D_{\rm{CH}}(\rvalpha, \rmP) =&\, \min_{i, j\in [K], j\neq k} \D_{\rm{CH}}(i, j),\\
   \D_{\rm{CH}}(i, j) \coloneqq &\, \max_{t\in [0, 1]}\sum\limits_{k=1}^{K} \alpha_{k} \cdot \big[ t\ermP_{ik} + (1 - t)\ermP_{jk} - (\ermP_{ik})^{t} (\ermP_{jk} )^{1 - t}\big]
\end{align*}
denotes the \emph{Chernoff-Hellinger} (CH) divergence, which can be deduced from \eqref{eqn:GCH} by taking $\sL = \{2\}$. \Cref{tab:graphSBM} reviews the related literature where exact recovery was achieved above the threshold.
\begin{table}[h]
\centering
\caption{Relevant literature for exact recovery in graph SBM} \label{tab:graphSBM}
\begin{tabular}{|P{1.2cm}|P{1.0cm}|P{1.0cm}|P{2.2cm}|P{1.8cm}|P{1.6cm}|P{2.5cm}|P{0.8cm}|}
 \hline
Graph & $K=2$ & $K\geq 3$ & First stage & Second stage & Exact recovery & Time complexity & No prior\\
\hline
\multirow{2}{*}{\cite{Abbe2016ExactRI}}& \cmark & \xmark & Spectral & Majority & \cmark & $2N^2 + N$  & \cmark \\
\cline{2-8}
& \cmark & \xmark & MLE & \xmark & \cmark & NP-hard & \xmark \\
\hline
\cite{Hajek2016AchievingEC} & \cmark & \xmark & SDP & \xmark & \cmark & $N^{3.5}$\cite{Jiang2020FasterIP}  & \cmark \\
\hline
\cite{Abbe2015CommunityDI} & \cmark & \cmark & Sphere &  Degree & \cmark & $N^{1+1/\log(c)}$ & \xmark \\
\cline{1-3} \cline{6-8}
\cite{Abbe2015RecoveringCI} & \cmark & \cmark & Comparison & Profiling & \cmark & $N^{1+1/\log(c)}$ &  \cmark \\
\hline
\cite{Mossel2016ConsistencyTF} & \cmark & \xmark & Spectral & Majority & \cmark & $2N^2 + N$ & \cmark \\
\hline
\cite{Yun2016OptimalCR} & \cmark & \cmark & Spectral & KL & \cmark & $KN^2 + KN\log(N)$ & \cmark \\
\hline
\cite{Gao2017AchievingOM} & \cmark & \cmark & Spectral & MLE & \cmark & $KN^2 + KN^3$ & \cmark \\
\hline
\cite{Abbe2020EntrywiseEA} & \cmark & \xmark & Spectral & \xmark & \cmark & $2N^2$ &  \cmark \\
\hline
\end{tabular}
\end{table}
For simplicity, the time complexity of singular value decomposition in our case is $O(N^3)$. When referring to the complexity of SDP, a factor of $\log^{O(1)}(N/\epsilon)$ is hidden, with $\epsilon$ being the accuracy of the output. If the algorithm has two stages, the entry (e.g., $N^3 + N$) represents the sum of the complexities for the first ($O(N^3)$) and second ($O(N)$) stages. If $K$ is needed as input but $\{\tQ^{(\ell)}\}_{\ell \in \sL}$ are not, the entry in column ``No prior'' is marked as \cmark. 
The second stage algorithms in \cite{Abbe2016ExactRI, Mossel2016ConsistencyTF} are based on the "Majority vote" principle. In \cite{Yun2016OptimalCR}, the estimator in second stage is aiming at minimizing \emph{Kullback-Leibler} (KL) divergence between estimated and true distribution.

\subsubsection{Hypergraphs}

Our results are most general to date, covering non-uniform hypergraphs with multiple communities. We now illustrate how they fit into existing literature. So far, all thresholds obtained in existing literature apply only to uniform hypergraphs. Our results establish, for the first time, the threshold for the non-uniform case, in the sense of showing both achievability above and impossibility below. See \Cref{tab:HSBM} for related literature.

For the binary symmetric uniform model $\mathrm{HSBM}(N, 2, a_{\ell}\frac{\log(N)}{\binom{N-1}{\ell - 1}}, b_{\ell}\frac{\log(N)}{\binom{N-1}{\ell - 1}})$, the threshold for exact recovery is given by $\D_{\rm{GH}}^{(\ell)} \coloneqq 2^{-(\ell - 1)}(\sqrt{a_{\ell}} - \sqrt{b_{\ell}})^2$ \cite{Kim2018StochasticBM}, which is a special case of \eqref{eqn:Dv} by taking $\sL = \{\ell\}$, $K = 2$ and symmetries among edges crossing two communities. Exact recovery is impossible if $\D_{\rm{GH}}^{(\ell)} <1$ \cite{Kim2018StochasticBM}, but can be achieved efficiently if $\D_{\rm{GH}}>1$ \cite{Kim2018StochasticBM, Gaudio2023CommunityDI}. For the multi-block case, the model formulation of \cite{Zhang2023ExactRI} is a bit different from our model \ref{def:uniform_HSBM}, since multisets are allowed as edges there. Our proofs work for their case as well if some steps are skipped. After proper scaling, their threshold for uniform $\mathrm{HSBM}(N, \rvalpha, \tP^{(\ell)}\log(N)/\binom{N-1}{\ell - 1} )$ can be written as $\D^{(\ell)}_{\rm{GCH}} = \underset{j, k\in [K], j\neq k}{\min} \D^{(\ell)}_{\rm{GCH}}(j,k)$, with
\begin{align*}
    \D^{(\ell)}_{\rm{GCH}}(j,k) &\, \coloneqq \max_{t\in [0, 1]} \,\sum\limits_{\rvw\in \WC{\ell - 1}{K}} \frac{\overline{N}_{\rvw}}{\binom{N-1}{\ell - 1}} \big[ t\etP^{(\ell)}_{j \oplus \rvw} + (1 - t)\etP^{(\ell)}_{k \oplus \rvw} - \big(\etP^{(\ell)}_{j \oplus \rvw } \big)^{t} \cdot \big(\etP^{(\ell)}_{ k \oplus \rvw} \big)^{1 - t} \big],
\end{align*}
which agrees with \eqref{eqn:GCH} by taking $\sL = \{\ell\}$. Exact recovery is impossible when $\D_{\rm{GCH}} < 1$, but can be achieved with the full knowledge of tensor $\tP^{(\ell)}$ and number of communities $K$ when $\D_{\rm{GCH}}> 1$ \cite{Zhang2023ExactRI}.

\begin{table}[h]
\centering
\caption{Relevant literature for exact recovery in HSBM} \label{tab:HSBM} 
\begin{tabular}{|P{1.6cm}|P{0.8cm}|P{1.0cm}|P{0.8cm}|P{0.8cm}|P{1.4cm}|P{1.2cm}|P{1.4cm}|P{2.0cm}|P{1.0cm}|}
  \hline
Hyper-graph & Uni-form & Non-uni-form & $K=2$ & $K\geq 3$ & First stage & Second stage & Exact recovery & Time complexity & No prior \\
 \hline
 \cite{Ghoshdastidar2017ConsistencyOS} & \cmark & \cmark & \cmark & \cmark & Spectral & \xmark & \xmark & $KN^2$ & \cmark \\
\hline
\cite{Kim2018StochasticBM}\tablefootnote{\cite{Gaudio2023CommunityDI} proved that SDP achieves exact recovery when $\D^{(\ell)}_{\rm{SDP}} >1$ where $\D^{(\ell)}_{\rm{SDP}} > \D^{(\ell)}_{\rm{GH}}$, leaving the area between $\D^{(\ell)}_{\rm{GH}}$ and $\D^{(\ell)}_{\rm{SDP}}$ unexplored.} & \cmark & \xmark & \cmark & \xmark &  SDP & \xmark & \cmark & $N^{3.5}$ & \cmark \\
  \hline
  \cite{Chien2019MinimaxMR} & \cmark & \xmark & \cmark & \cmark &  Spectral & MLE & \cmark & $KN^2 + KN^3$  & \cmark \\
  \hline
   \multirow{2}{*}{\cite{Gaudio2023CommunityDI}}  & \cmark & \xmark & \cmark & \xmark & Spectral & \xmark & \cmark & $2N^2$ & \cmark \\
    \cline{2-10}
     & \cmark & \xmark & \cmark & \xmark & SDP & \xmark & \cmark & $N^{3.5}$ & \cmark \\
    \hline
  \cite{Zhang2023ExactRI}  & \cmark & \xmark & \cmark & \cmark & Spectral & MAP & \cmark & $KN^2 + N$ & \xmark  \\
  \hline
  \cite{Wang2023ITLimits}  & \cmark & \cmark & \cmark & \xmark & Spectral & \xmark & \cmark & $2N^2$ & \cmark  \\ 
  \hline
  \hline
  \bf{Alg~\ref{alg:agnostic_partition}} & \cmark & \cmark & \cmark & \cmark & Spectral & KL & \cmark & $KN^2 + KN\log(N)$ & \cmark  \\
\hline
   \bf{Alg~\ref{alg:partition_with_prior}} & \cmark & \cmark & \cmark & \cmark & Spectral & MAP  & \cmark & $KN^3$ & \xmark \\
\hline
\end{tabular}
\end{table}

\subsection{Contributions}
The contribution of this paper lies on the following several aspects.

\subsubsection{Information-Theoretic limits}
\Cref{thm:IT_lower_bounds} (1) is the first result in the literature establishing the necessary condition for exact recovery under general non-uniform HSBM. Previously, the phase transition phenomenon for exact recovery was only characterized for graphs \cite{Abbe2016ExactRI, Abbe2015CommunityDI, Yun2016OptimalCR, Abbe2020EntrywiseEA, Abbe2022LPT} and uniform hypergraphs \cite{Kim2018StochasticBM, Zhang2023ExactRI}. Different from \cite{Abbe2015CommunityDI, Zhang2023ExactRI}, our proof of necessity does not rely on the additional \Cref{ass:prob_ratio_bound}.  Furthermore, for the first time in the literature of hypergraphs, \Cref{thm:IT_lower_bounds} (2) establishes the lower bound on the expected mismatch ratio for non-uniform HSBM. Previously, such results were obtained only in graph SBM for binary \cite{Abbe2020EntrywiseEA, Abbe2022LPT} and multiple \cite{Yun2016OptimalCR} communities. Minimax rates for the mismatch ratio were established for \cite{Gao2017AchievingOM, Gao2018CommunityDI} (graph SBM) and \cite{Chien2019MinimaxMR} (uniform HSBM), under specific constraints on the probability tensors $\{\tQ^{(\ell)}\}_{\ell \in \sL}$. By contrast, \Cref{thm:IT_lower_bounds} (2) is general.

\subsubsection{Weak consistency} \Cref{thm:known_weak_consistency} is the first result showing that weak consistency can be achieved under general non-uniform HSBM as long as $\rho_{N} = \omega(1)$. Previously, \cite{Dumitriu2025PartialRA} achieved weak consistency when $1 \ll \rho_{N} \ll \log(N)$, under a special case of model \ref{def:non_uniform_HSBM} . 
In \cite{Zhen2021CommunityDI}, weak consistency was achieved when $\rho_{N} \gtrsim \omega(\log N)$, a regime covered by our results. In \cite{Ke2020CommunityDF}, the authors studied the degree-corrected uniform HSBM by using a tensor power iteration method, and achieved weak consistency when the average degree is $\omega(\log^2(N))$, again covered by  our results. A way to generalize their algorithm to non-uniform hypergraphs was discussed, but the theoretical analysis remains open. 

The setting in \cite{Ghoshdastidar2017ConsistencyOS} is most similar to our model \ref{def:non_uniform_HSBM} , but it only achieved weak consistency when the minimum expected degree is of at least the order $\Omega(\log^2(N))$. Their algorithm cannot be applied to the regime $\rho_{N} \lesssim \log(N)$ due to the lack of concentration for the normalized Laplacian because of the existence of isolated vertices. We were able to overcome this issue by focusing on the adjacency matrix instead. Similar to \cite{Feige2005SpectralTA, Le2017ConcentrationAR}, we regularized the adjacency matrix by zeroing out rows and columns corresponding to vertices with large degree (\Cref{alg:trimming}) and proved a concentration result for the regularized matrix down to the bounded expected degree regime (\Cref{thm:regularization}).

\subsubsection{Strong consistency} Theorems \ref{thm:agnostic_strong_consistency} and \ref{thm:known_strong_consistency} are the first results to provide efficient algorithms for the general non-uniform case. Previously, strong consistency was only achieved under uniform HSBM for binary \cite{Kim2018StochasticBM, Gaudio2023CommunityDI}, multi-community \cite{Chien2019MinimaxMR, Zhang2023ExactRI} case, and non-uniform binary case \cite{Wang2023ITLimits}, under \Cref{ass:prob_ratio_bound}. The method in \cite{Han2022ExactCI} was not specialized to HSBM, and it does not contain any HSBM thresholds. The refinement stage is crucial for multi-block case, which brings in the sharp threshold \cite{Abbe2015CommunityDI, Yun2016OptimalCR, Chien2019MinimaxMR, Zhang2023ExactRI}. In \Cref{alg:partition_with_prior} and \cite{Zhang2023ExactRI}, the second stage is based on \emph{Maximum A Posteriori} (MAP) estimation, where the MAP achieves the minimum error and full knowledge of $\{\tQ^{(\ell)}\}_{\ell \in \sL}$ required. In \cite{Chien2019MinimaxMR}, the agnostic refinement is based on \emph{local Maximum Likelihood Estimation} (MLE), which works well for uniform hypergraphs. It remains to be seen whether this method can be generalized for non-uniform hypergraphs. Additionally, the time complexity of \cite{Chien2019MinimaxMR} is $O(N^4)$, since the spectral initialization is run $N$ times during the refinement stage. So far, \Cref{alg:agnostic_partition} is the only one that achieves strong consistency agnostically with the lowest time complexity.

\subsubsection{Optimality} Although, as noted before, many works provided algorithms achieving weak consistency, the optimality of the methods was less explored. For the first time in the literature, we are able to provide Algorithms \ref{alg:agnostic_partition} and \ref{alg:partition_with_prior} which are provably optimal in achieving weak consistency below the threshold for exact recovery, as in Theorems \ref{thm:agnostic_optimality} and \ref{thm:known_optimality}. In particular, it provides the upper bound on the expected number of misclassified nodes for non-uniform HSBM, covering the most general scenario so far. Previously, spectral methods were proved to be optimal for binary-community clustering \cite{Zhang2016MinimaxRO, Abbe2020EntrywiseEA, Abbe2022LPT} with the optimal mismatch ratio being $N^{-(\sqrt{a} - \sqrt{b})^2/2}$. The optimal mismatch ratio was also explored under multi-community SBM\cite{Gao2017AchievingOM, Gao2018CommunityDI, Zhang2023FundamentalLO}, labeled SBM \cite{Yun2016OptimalCR} and uniform HSBM \cite{Chien2019MinimaxMR}.

\subsubsection{Comparison with existing tensor methods} Many tensor-based methods contain the assumption that the $K$th smallest eigenvalue of a $\ell$-uniform tensor (perhaps resulting from Tucker decomposition) is bounded away from $0$ (in fact, the lower bounds tend to be more constraining); e.g., Assumption $1$ in \cite{Agterberg2022EstimatingHO} and (4.7) in \cite{Ke2020CommunityDF}. In general, \Cref{ass:expected_center_separation} is strictly weaker than saying that the $K$th smallest eigenvalue is nonzero: in Section ~\ref{sec:communityEstimation}, we provide an example satisfying the \Cref{ass:expected_center_separation}, but for which the least eigenvalue for the $K \times K$ ``reduced" expectation matrix is $0$, thus violating the conditions identified in \cite{Agterberg2022EstimatingHO} for tensor methods. In {\cite[Section 4.1]{Han2022ExactCI}}, the authors involves the parameter $\Delta_i^2$, which represents the minimum square $2$-norm differences between different rows of various matricizations of the expectation tensor. Their exact recovery condition, $\Delta_{min}^2/\sigma^2 \geq C ( p^{-d/2} \vee p^{-(d-1)} \log p)$, is more restrictive than ours, since it imposes a lower bound, whereas for us the lower bound is $0$. Those tensor based methods cannot be generalized to the non-uniform case directly, since an efficient method to incorporate tensors of different orders has not yet been devised. It remains to be seen if some other (non-spectral) method can be used to eliminate the remaining edge cases violating \Cref{ass:expected_center_separation}.

\subsection{Notation and preliminaries}
Let $\gH = (\gH, \gE)$ denote the hypergraph, where $\gV$ and $\gE$ denote the set of vertices and (hyper-)edges respectively. For vertex $v\in \gV$, let $d_v$ denote the observed number of hyperedges containing $v$, with $\rD_v$ denoting its law. Note that each $\ell$-hyperedge $e \ni v$ can be classified into different categories according to the distribution of the other $\ell - 1$ nodes among $\gV_{1}, \ldots, \gV_{K}$, then $\rD_v$ can be written as a summation of independent random variables due to edge-wise independence,
\begin{align}\label{eqn:Dv}
    \rD_v \coloneqq \sum\limits_{\ell \in \sL}\sum\limits_{\rvw \in \WC{\ell - 1}{K}} \rD^{(\ell)}_{v, \rvw},
\end{align}
where $\rD^{(\ell)}_{v, \rvw}$ denotes the law of $d^{(\ell)}_{v, \rvw}$, counting the number of $\ell$-hyperedges containing $v$ with the other $\ell - 1$ nodes distributed as $\rvw$ among $\gV_{1}, \ldots, \gV_{K}$. Assume that $v$ belongs to the community $\gV_k$, that is, $\rY_{v} = k$, then the distribution of membership among the blocks in $e$ is denoted by $k \oplus \rvw \coloneqq (\ervw_1, \ldots, \ervw_{k-1}, \ervw_k + 1, \ervw_{k+1}, \ldots, \ervw_K)$. According to model \ref{def:non_uniform_HSBM} , $\rD^{(\ell)}_{v, \rvw} \vert_{\rY_{v} = k} \sim {\rm{Bin}}(N_{\rvw}, \etQ^{(\ell)}_{k\oplus \rvw})$, where for each $\rvw = (\ervw_1, \ldots, \ervw_K)\in \WC{\ell - 1}{K}$, and the capacity of such $\ell$-hyperedges is denoted by
\begin{align}\label{eqn:nbw}
    N_{\rvw} \coloneqq \prod\limits_{k=1}^{K} \binom{ |\gV_k| }{\ervw_k}
\end{align}
\Cref{lem:sizeDeviation} indicates that $\overline{N}_{\rvw}$ in \eqref{eqn:barnw} is a good approximation of $N_{\rvw}$.
\begin{lemma}\label{lem:sizeDeviation}
     With probability at least $1 - N^{-\log(N)}$, $||\gV_{j}| - \alpha_{j} N| \leq \sqrt{N}\log(N)$ for each $j \in [K]$. Recall $N_{\rvw}$ from \eqref{eqn:nbw} and $\overline{N}_{\rvw}$ from \eqref{eqn:barnw}. Then for each for $\rvw \in \WC{\ell - 1}{K}$, it follows
    \begin{align}
        (1 - \log(N)/\sqrt{N} )^{\ell - 1} \leq N_{\rvw}/ \overline{N}_{\rvw} \leq (1 + \log(N)/\sqrt{N} )^{\ell - 1}.
    \end{align}
\end{lemma}

\subsection{Organization} 
This paper is organized as follows. \Cref{sec:IT_lower_bounds} contains the proof of \Cref{thm:IT_lower_bounds}. The proofs of Theorems \ref{thm:agnostic_weak_consistency}, \ref{thm:known_weak_consistency}, Theorems \ref{thm:agnostic_weak_consistency}, \ref{thm:agnostic_optimality} and Theorems \ref{thm:known_strong_consistency}, \ref{thm:known_optimality} are in Sections \ref{sec:almostExact}, \ref{sec:AgnosticRefinement} and \ref{sec:known_partition} respectively. We will discuss estimating the number of communities in \Cref{sec:communityEstimation}. Some open problems will be addressed in \Cref{sec:conclusion_general}.

\section{Information-Theoretic limits}\label{sec:IT_lower_bounds}
This section is devoted to the proof of \Cref{thm:IT_lower_bounds} with the proofs of Lemmas deferred to Appendix~\ref{app:ITLowerBounds}. 

\subsection{Impossibility of exact recovery}\label{subsec:impossibleExact}
The proof of \Cref{thm:IT_lower_bounds} (1) is presented in this section. Without loss of generality, let $j^{\star}, k^{\star} \in [K]$ be the indices where $\D_{\mathrm{GCH}}\cdot q_{N}$ in \eqref{eqn:GCH_not_proportional} attains its minimum, i.e., $(j^{\star}, k^{\star}) = \arg\min_{j, k\in [K]} \D_{\mathrm{GCH}}(j, k)\cdot q_{N}$. The main idea is to confirm the existence of many pairwise disconnected \emph{ambiguous} vertices (deferred to \Cref{def:ambiguousVertex}). With high probability, those disconnected vertices maintain the same degree profile, but half from $\gV_{j^{\star}}$ and half from $\gV_{k^{\star}}$. Consequently, no algorithm performs better than random guess, which leads to the failure of exact recovery. Similar ideas appeared in \cite{Abbe2015CommunityDI}. By contrast, the central idea of \cite{Wang2023ITLimits} is to explicitly construct some configuration where the optimal estimator Maximum A Posteriori (MAP) fails. Note that $\D_{\rm{GCH}}$ in \eqref{eqn:GCH} can be reduced to $\D_{\rm{GH}}$ \cite{Wang2023ITLimits} when the general model is restricted to the binary symmetric model. Those two distinct proofs can be viewed as a mutual verification for the correctness of the thresholds.

\subsubsection{Regime of interest} 
Before elaborating the proof sketch above, we establish the proof of impossibility for some easier cases quickly, which would simplify our discussion for the more difficult case later.

First, it suffices to consider the case $\kappa_{N} = 1$, i.e., $\D_{\mathrm{GCH}}\cdot q_{N} < \log(N)$ when under \eqref{eqn:impossibleExactRecovery}. 

Second, we restrict our attention to a sub-hypergraph described below. Denote
\begin{align}
    \const^{(\ell)}_{j^{\star}, k^{\star}, \rvw} = \max\{ \etQ^{(\ell)}_{j^{\star} \oplus \rvw}/\etQ^{(\ell)}_{k^{\star} \oplus \rvw},\,\, \etQ^{(\ell)}_{k^{\star} \oplus \rvw}/\etQ^{(\ell)}_{j^{\star} \oplus \rvw} \}
\end{align}
for each $\rvw \in \WC{\ell - 1}{K}$. When $\const^{(\ell)}_{j^{\star}, k^{\star}, \rvw} = 1$, it follows that $\etQ^{(\ell)}_{j^{\star} \oplus \rvw}=\etQ^{(\ell)}_{k^{\star} \oplus \rvw}$. The edges with the membership counts $j^{\star} \oplus \rvw$ or $k^{\star} \oplus \rvw$ have zero contribution to $\D_{\rm{GCH}}$ in \eqref{eqn:GCH_not_proportional}, and they can not be used to distinguish vertices in $\gV_{j^{\star}}$ from vertices in $\gV_{k^{\star}}$. Consequently, we can exclude them from consideration. For each $\ell \in \sL$, define the following set
\begin{align}
    \sW^{(\ell)}_{j^{\star}, k^{\star}} = \{\rvw \in \WC{\ell - 1}{K}| \const^{(\ell)}_{j^{\star}, k^{\star}, \rvw} > 1\}.\label{eqn:Wmjkstar}
\end{align}
Throughout the proof, we assume $1\ll \D_{\mathrm{GCH}}\cdot q_{N} < \log(N)$, and we focus on the type of edges with membership counts $j^{\star} \oplus \rvw$ or $k^{\star} \oplus \rvw$, where $\rvw \in \sW^{(\ell)}_{j^{\star}, k^{\star}}$ for each $\ell \in \sL$.

\subsubsection{Proof of Necessity}
Let $t^{\star} \in [0, 1]$ be the point where the maximum of $\D_{\rm{GCH}}(j^{\star}, k^{\star})\cdot q_{N}$ is attained, then $t^{\star}$ should be the critical point for the function we maximize, i.e.,
\begin{align}\label{eqn:optimalT}
    \sum\limits_{\ell \in \sL} \,\sum\limits_{\rvw\in \sW^{(\ell)}_{j^{\star}, k^{\star}}} \overline{N}_{\rvw} \big[ \etQ^{(\ell)}_{j^{\star}\oplus \rvw} - \etQ^{(\ell)}_{k^{\star} \oplus \rvw} - (\etQ^{(\ell)}_{j^{\star}\oplus \rvw})^{t^{\star}} (\etQ^{(\ell)}_{k^{\star} \oplus \rvw} )^{1 - t^{\star}} \log \frac{\etQ^{(\ell)}_{j^{\star}\oplus \rvw}}{\etQ^{(\ell)}_{k^{\star} \oplus \rvw}} \Big] = 0.
\end{align}
The $t^{\star}$ is unique since the function to the left of the equal sign in \eqref{eqn:optimalT} is strictly decreasing. For each $\rvw \in \sW^{(\ell)}_{j^{\star}, k^{\star}}$, define \textit{degree profile} as
\begin{align}\label{eqn:ambiguous}
    d_{j^{\star}, k^{\star}}^{\star}(\rvw) \coloneqq \Big \lfloor \big(\etQ^{(\ell)}_{j^{\star}\oplus \rvw} \big)^{t^{\star}} \cdot \big(\etQ^{(\ell)}_{k^{\star} \oplus \rvw} \big)^{1 - t^{\star}} \cdot \overline{N}_{\rvw} \Big\rfloor\,.
\end{align}
Intuitively, $\gV_{j^{\star}}$ and $\gV_{k^{\star}}$ are most similar among all communities. The degree profile $d^{\star}_{j^{\star}, k^{\star}}(\rvw)$ interpolates between the expected degree profiles of the vertices in $\gV_{j^{\star}}$ and $\gV_{k^{\star}}$, and, if $\gV_{j^{\star}}$ and $\gV_{k^{\star}}$ are similar enough, vertices with this degree profile are very likely to exist in both communities. If they do exist in both, they are indistinguishable: there is no way to be sure which community they belong to. We introduce the concept of \emph{ambiguity} below.
\begin{definition}[Ambiguity]\label{def:ambiguousVertex}
For $\setS \subset \gV$, let $\rD^{(\ell)}_{v, \rvw}|_{\setS}$ be the number of edges containing $v$, with the remaining $\ell - 1$ vertices all in $\gV\setminus \setS$. A vertex $v\in \setS$ is said to be $(j^{\star}, k^{\star}, \setS)$-\emph{ambiguous} between $\gV_{j^{\star}}\setminus \setS$ and $\gV_{k^{\star}}\setminus \setS$ if $\rD^{(\ell)}_{v, \rvw}|_{\setS} = d_{j^{\star}, k^{\star}}^{\star}(\rvw)$ for each $\rvw \in \sW^{(\ell)}_{j^{\star}, k^{\star}}$, $\ell \in \sL$.
\end{definition}
For distinct vertices $v, u \in \setS$, the events $\{v \textnormal{ is } (j^{\star}, k^{\star}, \setS)\textnormal{-ambiguous}\}$ and $\{u\textnormal{ is } (j^{\star}, k^{\star}, \setS)$ $\textnormal{-ambiguous}\}$ are independent, since the random variables involved are distinct. \Cref{lem:ambiguousProb} characterizes the probability of $v$ being ambiguous.

\begin{lemma}\label{lem:ambiguousProb}
Assume that vertex $v\in \gV_{k^{\star}}$, i.e., $\rY_{v} = k^{\star}$. Assume \eqref{eqn:impossibleExactRecovery} with $\kappa_{N} = 1$ and let $j^{\star}, k^{\star} \in [K]$ be the pair achieving $\D_{\mathrm{GCH}}$ in \eqref{eqn:GCH_not_proportional}, then 
\begin{align}
    \prod\limits_{\ell \in \sL} \prod\limits_{\rvw \in \sW^{(\ell)}_{j^{\star}, k^{\star}}}\P(\rD^{(\ell)}_{v, \rvw} =  d_{j^{\star}, k^{\star}}^{\star}(\rvw)) = \exp(-(1 + o(1))\cdot\D_{\rm{GCH}}\cdot q_{N}).
\end{align}
\end{lemma}
 \begin{proof}[Proof of \Cref{thm:IT_lower_bounds} (1)]
    Let $\setS \subset \gV$ be a set of vertices with size $\lfloor N q_{N}^{-2} \log^{-1}(N) \rfloor$, where each vertex is chosen randomly from $\gV$. By \Cref{lem:sizeDeviation}, with probability at least $1 - O(N^{-3})$, $|\setS\cap \gV_{j}| = \alpha_{j} N q_{N}^{-2} \log^{-1}(N)$ with deviation smaller than $O(\sqrt{N})$ for each $j\in [K]$.
Let $v^{\star} \in \setS \cap \gV_{k^{\star}}$ be a $(j^{\star}, k^{\star}, \setS)$-ambiguous vertex. Since $|\setS|/|\gV| = o(1)$, the extra condition of having the other $\ell - 1$ vertices all in $\gV \setminus \setS$ will erase only a negligible number of edges. Assume \eqref{eqn:impossibleExactRecovery} for sufficiently large $N$ with $\kappa_{N} = 1$, then $\D_{\rm{GCH}}\cdot q_{N} /\log(N) < 1 - 2\epsilon < 1$ for some absolute constant $\epsilon >0$. By following the proof of \Cref{lem:ambiguousProb}, it yields that
\begin{align}
    &\,\P(v^{\star} \textnormal{ is } (j^{\star}, k^{\star}, \setS)\textnormal{-ambiguous} ) = \prod\limits_{\ell \in \sL} \prod\limits_{\rvw \in \sW^{(\ell)}_{j^{\star}, k^{\star}}} \P(\rD^{(\ell)}_{v^{\star}, \rvw}|_{\setS} =  d_{j^{\star}, k^{\star}}^{\star}(\rvw)) \\
    = &\,\exp\big(-(1 + o(1))\cdot\D_{\rm{GCH}}\cdot q_{N}\big) \geq N^{-1 + 2 \epsilon + o(1)} \geq N^{-1 + \epsilon}.
\end{align}
By independence between ambiguous nodes and the fact $1\ll q_{N} \lesssim \log(N)$, we have
\begin{align}
   \E[\#\,\, (j^{\star}, k^{\star}, \setS)\textnormal{-ambiguous nodes in } \setS\cap \gV_{k^{\star}}] \gtrsim &\,\alpha_{k^{\star}} \lfloor N q_{N}^{-2} \log^{-1}(N) \rfloor \cdot N^{-1 + \epsilon} \notag\\
   \asymp &\, N^{\epsilon}q_{N}^{-2} \log^{-1}(N) \gg q_{N}. \label{eqn:lower_bound_{N}um_ambiguity}
\end{align}
Then by Bernstein's inequality, with probability at least $1 - e^{-N^{\epsilon}\log^{-3}(N)}$, the number of $(j^{\star}, k^{\star}, \setS)$-ambiguous nodes will concentrate around its expectation with deviation less than $N^{\epsilon}\log^{-3}(N)$. Therefore, there exists at least $q_{N}$ many $(j^{\star}, k^{\star}, \setS)$-ambiguous vertices.  

We will prove that $v^{\star}$ is isolated in $\setS$ with high probability in the considered sub-hypergraph. Recall that $\etQ^{(\ell)}_{k^{\star}\oplus\rvw}$ denotes the probability of sampling a $\ell$-hyperedge containing $v^{\star} \in \gV_{k^{\star}}$ with the other $\ell - 1$ vertices distributed as $\rvw$ among $\{\gV_{k}\}_{k=1}^{K}$. If one more constraint is added, i.e., containing $v^{\star}$ and at least one other vertex from $\setS$, then the number of such hyperedges is denoted by $\overline{N}_{\rvw}|_{\setS} = (1 + o(1)) q_{N}^{-2} \log^{-1}(N)\cdot \overline{N}_{\rvw}$, where $\overline{N}_{\rvw}$ is defined in \eqref{eqn:barnw}. The probability that $v^{\star}$ is isolated in $\setS$ in the considered sub-hypergraph is
\begin{align}\label{eqn:vstar_disconnected}
     &\,\P(v^{\star} \textnormal{ isolated in } \setS) \\
    =\prod\limits_{\ell \in \sL} &\, \prod\limits_{\rvw \in \sW^{(\ell)}_{j^{\star}, k^{\star}}} (1 - \etQ^{(\ell)}_{k^{\star}\oplus\rvw})^{\overline{N}_{\rvw}|_{\setS}} = \exp\bigg( - (1 + o(1))\sum\limits_{\ell \in \sL} \sum\limits_{\rvw \in \sW^{(\ell)}_{j^{\star}, k^{\star}}} \frac{\etQ^{(\ell)}_{k^{\star}\oplus\rvw}\cdot \overline{N}_{\rvw}}{q_{N}^2\cdot \log(N)}\bigg), \notag
\end{align}
where we applied $e^{-x} = 1 - x$ when $x = o(1)$ since $\etQ^{(\ell)}_{k^{\star}\oplus\rvw} = o(1)$. Note that we only focus on the weak compositions in $\sW^{(\ell)}_{j^{\star}, k^{\star}}$ due to the fact that $\rvw \notin \sW^{(\ell)}_{j^{\star}, k^{\star}}$ contributes $0$ to $\D_{\rm{GCH}}$. It follows by lemma~\ref{lem:GCH_max_different_order} that $\sum_{\ell \in \sL} \sum_{\rvw \in \sW^{(\ell)}_{j^{\star}, k^{\star}}}\etQ^{(\ell)}_{k^{\star}\oplus\rvw}\cdot \overline{N}_{\rvw} \lesssim \D_{\rm{GCH}} \cdot q_{N}$. Consequently, there exists some constant $\const > 0$ such that
\begin{align}
    \eqref{eqn:vstar_disconnected} = \exp\big( - \const\cdot\D_{\rm{GCH}}\cdot q_{N}\cdot [q_{N}^2\cdot \log(N)]^{-1} \big) > e^{-\const/q_{N}^2} = 1 - \const/q_{N}^2,
\end{align}
where the second inequality holds since $\D_{\rm{GCH}}\cdot q_{N} /\log(N) < 1$ by \eqref{eqn:impossibleExactRecovery}.

From \eqref{eqn:lower_bound_{N}um_ambiguity}, we can pick at least $q_{N}$ pairs of $(j^{\star}, k^{\star}, \setS)$-ambiguous vertices, where in each pair, one from $\setS \cap \gV_{j^{\star}}$ and the other one from $\setS \cap \gV_{j^{\star}}$. The probability that these $q_{N}$ pairs of vertices are all isolated from $\setS$ in the considered sub-hypergraph is at least $(1 -\const q_{N}^{-2})^{2q_{N}} \geq 1 - 4\const/q_{N}$, which tends to $1$ since $q_{N} = \omega(1)$. Consequently, those $q_{N}$ pairs of vertices are indistinguishable. 
It follows that no algorithm can perform better than random guess for \emph{each} of these vertices, and by independence, no algorithm can have probability of success better than $[ \binom{2q_{N}}{q_{N}}]^{-1} \asymp 2^{-q_{N}}$, which goes to $0$. Therefore, every algorithm will misclassify at least one vertex  with probability $1 - 2^{-q_{N}}$.
\end{proof}

\subsection{Information-Theoretic lower bound} To establish the proof of \Cref{thm:IT_lower_bounds} (2), we first introduce the following \emph{Generalized Kullback-Leibler} (GKL) divergence.
\begin{definition}
For each $\ell \in \sL$, let $\vmu^{(\ell)} = (\evmu^{(\ell)}_{\rvw})_{\rvw \in \WC{\ell - 1}{K}} \in [0, 1]^{|\WC{\ell - 1}{K}|}$ denote some sequence of numbers. Define the \emph{generalized Kullback-Leibler (GKL)} divergence
    \begin{align}
         &\,\D_{\mathrm{GKL}} =  \underset{1\leq j < k \leq K}{\min} \D_{\mathrm{GKL}} (j, k), \label{eqn:GeneralizedKL}\\
        \D_{\mathrm{GKL}} (j, k) \coloneqq &\, \underset{\{\vmu^{(\ell)}\}_{\ell \in \sL} }{\min}\,\sum\limits_{\ell \in \sL}\,\, \sum\limits_{\rvw \in \WC{\ell - 1}{K} } \overline{N}_{\rvw}\cdot \max\big\{ \D_{\mathrm{KL}}(\evmu^{(\ell)}_{\rvw}\parallel \etQ_{j \oplus \rvw}^{(\ell)}),\, \D_{\mathrm{KL}}(\evmu^{(\ell)}_{\rvw} \parallel \etQ_{k \oplus \rvw}^{(\ell)} ) \big\},\notag
    \end{align}
    where $\overline{N}_{\rvw}$ is defined in \eqref{eqn:barnw}, and
    \begin{align}\label{eqn:KLDivergence}
        \D_{\mathrm{KL}}(\evmu^{(\ell)}_{\rvw}\parallel \etQ_{j \oplus \rvw}^{(\ell)}) \coloneqq \evmu^{(\ell)}_{\rvw} \log \bigg( \frac{\evmu^{(\ell)}_{\rvw} }{\etQ^{(\ell)}_{j \oplus \rvw}} \bigg) + (1 - \evmu^{(\ell)}_{\rvw}) \log \bigg( \frac{1 - \evmu^{(\ell)}_{\rvw} }{1 - \etQ^{(\ell)}_{j \oplus \rvw}} \bigg)
    \end{align}
    denotes the KL divergence between two Bernoulli distributions with means $\evmu^{(\ell)}_{\rvw}$ and $\etQ_{j \oplus \rvw}^{(\ell)}$.
\end{definition}

\Cref{lem:equivalenveDivergence} builds the connection between $\D_{\mathrm{GCH}}$ in \eqref{eqn:GCH} and $\D_{\mathrm{GKL}}$ in \eqref{eqn:GeneralizedKL}. \Cref{lem:qMeasureExistence} establishes the existence of some sequence $\{\rvp^{(\ell)}\}_{\ell \in \sL}$ achieving $\D_{\mathrm{GKL}}$.

\begin{lemma}\label{lem:equivalenveDivergence}
    When $\etQ^{(\ell)}_{\rvw} = o(1)$ for all $\rvw\in \WC{\ell}{K}$, $\ell \in \sL$, we have $\D_{\rm{GKL}} (j, k) = (1 + o(1)) \cdot \D_{\rm{GCH}}(j,k)\cdot q_{N}$.
\end{lemma}
\begin{lemma}\label{lem:qMeasureExistence}
    Let $(j^{\star}, k^{\star}) = \arg\min_{j<k} \D_{\mathrm{GKL}}(j, k)$. There exists some probability measure sequence $\{\rvp^{(\ell)}\}_{\ell \in \sL}$ such that
    \begin{align}
        &\, \D_{\mathrm{GKL}} =\D_{\mathrm{GKL}}(j^{\star}, k^{\star}) \label{eqn:equivalenceqm}\\
        =\sum\limits_{\ell \in \sL} &\, \sum\limits_{\rvw \in \WC{\ell - 1}{K}}\overline{N}_{\rvw}\cdot  \D_{\mathrm{KL}}(\ervp^{(\ell)}_{\rvw}\parallel \etQ_{j^{\star} \oplus \rvw}^{(\ell)}) = \sum\limits_{\ell \in \sL} \sum\limits_{\rvw \in \WC{\ell - 1}{K}}\overline{N}_{\rvw} \cdot \D_{\mathrm{KL}}(\ervp^{(\ell)}_{\rvw}\parallel \etQ_{k^{\star} \oplus \rvw}^{(\ell)} ).\notag
    \end{align}
\end{lemma}

Let $\Phi$ denote the true model \ref{def:non_uniform_HSBM} and $\Psi$ denote the perturbed stochastic model \ref{def:perturbed_model}, below. The lower bound of the expected mismatch ratio in \Cref{thm:IT_lower_bounds} (2) is established by the \emph{coupling} between $\Phi$ and $\Psi$. The generalized KL divergence $\D_{\mathrm{GKL}}$ plays a fundamental role in our analysis, similar to \cite{Lai1985AsymptoticallyEA}.

\begin{definition}[Perturbed model $\Psi$]\label{def:perturbed_model}
    For $j^{\star}, k^{\star}$ in \Cref{lem:qMeasureExistence}, let $v^{\star}$ denote node with the smallest index belonging to $\gV_{j^{\star}} \cup \gV_{k^{\star}}$. The two models $\Phi$ and $\Psi$ are coupled as follows.
\begin{enumerate}
    \item The clusters $\gV_{1}, \ldots, \gV_{K}$ are generated first, thus the same under $\Phi$ and $\Psi$.
    \item For each $\ell$-edge $e\not\ni v^{\star}$, it is sampled with probability $\etQ^{(\ell)}_{\rvy(e)}$, both under $\Phi$ and $\Psi$.
    \item For each $e\ni v^{\star}$ with the other $\ell - 1$ nodes distributed as $\rvw$, it is sampled with probability $\etQ^{(\ell)}_{\ervy_{v^{\star}}\oplus \rvw}$ under $\Phi$, while with probability $\ervp^{(\ell)}_{\rvw}$ (\Cref{lem:qMeasureExistence}) under $\Psi$ respectively.
\end{enumerate}
\end{definition}
Let $\P_{\Phi}$, $\E_{\Phi}$ (resp. $\P_{\Psi}$, $\E_{\Psi}$) denote the probability measure, expectation under $\Phi$ (resp. $\Psi$). Let $\gV_{\rm{err}}$ denote the set of misclassified nodes. Those nodes are interchangeable under $\Phi$, then $\E_{\Phi}[\mismatch_{N}] = N^{-1}\E_{\Phi}[|\gV_{\rm{err}}|]$. Let $\Lambda \coloneqq \log (d\P_{\Psi}/d\P_{\Phi})$ denote the log-likelihood ratio, and
\begin{align}
    \Lambda = \sum\limits_{\ell \in \sL} \,\sum\limits_{\substack{e\in \gE_{\ell}, e\ni v^{\star}\\ \rvy(e\setminus v^{\star}) = \rvw}} \Big[\etA^{(\ell)}_{e} \log \Big(\frac{\ervp^{(\ell)}_{\rvw}}{\etQ^{(\ell)}_{\ervy_{v^{\star}}\oplus \rvw}} \Big) + (1 - \etA^{(\ell)}_{e})\log \Big(\frac{1 - \ervp^{(\ell)}_{\rvw}}{1 - \etQ^{(\ell)}_{\ervy_{v^{\star}}\oplus \rvw}} \Big) \Big].\label{eqn:Lambda}
\end{align}
Assume $\alpha_1 \geq \ldots \geq \alpha_K$ without loss of generality. Note that $j^{\star} < k^{\star}$ from \Cref{lem:qMeasureExistence}. The connection between $\Psi$ and $\Phi$ can be described quantitatively by \Cref{lem:LambdaUpperBound}.

\begin{lemma}\label{lem:LambdaUpperBound}
Let $f(N)$ be some function of $N$, then the following holds
    \begin{align}\label{eqn:LambdaUpperBound}
    \P_{\Psi}(\Lambda \leq f(N)) \leq (e^{f(N)}\cdot \E_{\Phi}[\mismatch_{N}] + \alpha_{j^{\star}}) /(\alpha_{j^{\star}} + \alpha_{k^{\star}}).
    \end{align}
\end{lemma}
\begin{proof}[Proof of \Cref{thm:IT_lower_bounds} (2)]
First, by applying $f(N) = \log(1/\E_{\Phi}[\mismatch_{N}]) - \log(2/\alpha_{k^{\star}})$ and the fact $\alpha_{j^{\star}} + \alpha_{k^{\star}} < 1$, \Cref{lem:LambdaUpperBound} yields that
\begin{align}
    \P_{\Psi}(\Lambda \leq \log(1/\E_{\Phi}[\mismatch_{N}]) - \log(2/\alpha_{k^{\star}})) \leq = 1 - (\alpha_{j^{\star}} + \alpha_{k^{\star}})^{-1}\alpha_{k^{\star}}/2 < 1 - \alpha_{k^{\star}}/4.
\end{align}
Furthermore, by Chebyshev's inequality, we have
\begin{align}
    \P_{\Psi}(\Lambda \leq \E_{\Psi}[\Lambda] + \sqrt{4\Var_{\Psi}(\Lambda)/\alpha_{k^{\star}}})\geq 1 - \alpha_{k^{\star}}/4. 
\end{align}
Consequently, the following comparison of the two upper bounds of $\Lambda$ should always be true:
\begin{align}\label{eqn:mismatch_expectedLambda}
    \log(1/\E_{\Phi}[\mismatch_{N}]) - \log(2/\alpha_{k^{\star}}) \leq \E_{\Psi}[\Lambda] + \sqrt{4\Var_{\Psi}(\Lambda)/\alpha_{k^{\star}}}\,\,.
\end{align}
\begin{lemma}\label{lem:expectation_Lambda_Psi}
Recall $1\ll \D_{\rm{GCH}} \cdot q_{N} \lesssim \log(N)$ in \eqref{eqn:GCH_not_proportional}, then the following holds:
    \begin{align}\label{eqn:expectation_Lambda_Psi}
        \E_{\Psi}[\Lambda] + \sqrt{4\Var_{\Psi}(\Lambda)/\alpha_{k^{\star}}} = (1 + o(1)) \cdot \D_{\rm{GCH}} \cdot q_{N}.
    \end{align}
\end{lemma}
Recall \eqref{eqn:impossibleExactRecovery} where $\kappa_{N} \in (0, 1]$ and $\kappa_{N} \log(N) \to \infty$ when $N \to \infty$. According \Cref{lem:expectation_Lambda_Psi}, when $\E_{\Phi}[\mismatch_{N}] < N^{-\kappa_{N} }$, the left hand side of \eqref{eqn:mismatch_expectedLambda} can be further bounded by
\begin{align}
    \kappa_{N} \cdot \log(N) - \log(2/\alpha_{k^{\star}}) \leq (1 + o(1)) \cdot \D_{\rm{GCH}} \cdot q_{N}.
\end{align}
Since $\alpha_{k^{\star}}$ is some fixed constant, the necessary condition for $\lim_{N \to \infty}(N^{\kappa_{N}} \E_{\Phi}[\mismatch_{N}] ) < 1$ is
\begin{align}
    \lim\limits_{N \to \infty} \D_{\mathrm{GCH}}\cdot q_{N}/(\kappa_{N} \cdot \log(N)\, ) \geq 1.
\end{align}
Therefore, \Cref{thm:IT_lower_bounds} (2) follows since we proved its contrapositive.
\end{proof}
\section{Stage I: Weak consistency by spectral initialization}\label{sec:almostExact}
We will show that \Cref{alg:spectral_initialization}, inspired by \cite{Coja-Oghlan2010GraphPV, Yun2014CommunityDV, Yun2016OptimalCR, Zhang2023ExactRI}, achieves almost exact recovery, i.e., $\mismatch(\rvy, \widehat{\rvy}^{(0)}) = o(1)$. Proofs of Lemmas are deferred to Appendix~\ref{app:almostExact}.

\subsection{Spectral initialization}
Recall that $d_v$ denotes the number of hyperedges containing $v$ with $\rD_v$ in \eqref{eqn:Dv} denoting its law. Let $\overline{d} = \sum_{v=1}^{N}d_v/N$ (resp. $\widetilde{d} = \max_{v\in [N]}d_v$) denote the \emph{observed average} (resp. \emph{maximum}) degree. Let $\xi_{N}$ (resp. $\rho_{N}$) denote the \emph{average} (resp. \emph{maximum}) \emph{expected} degree, defined by
\begin{align}\label{eqn:rhoxi}
    \xi_{N} \coloneqq N^{-1} \E \Big( \sum\limits_{v\in \gV} \rD_v \Big)
    , \quad \rho_{N} \coloneqq \underset{v\in \gV}{\max}\, \E(\rD_v).
\end{align}

\Cref{lem:maxProb} states the fact that $\overline{d}$ is always a good approximation of $\xi_{N}$, while $\widetilde{d}$ preserves the order of $\rho_{N}$ when $\rho_{N} \gtrsim \log(N)$.

\begin{lemma}\label{lem:maxProb}
Deterministically, $\rho_{N} \asymp \xi_{N}$. With high probability, $\overline{d} - \xi_{N} \leq \xi_{N} /\sqrt{N}$. When $\rho_{N} \gtrsim \log(N)$, $\widetilde{d} \asymp \rho_{N}$ with probability at least $1 - O(N^{-c})$ for some constant $c >0$.
\end{lemma}

\subsection{Concentration inequalities} Our analysis relies on the concentration results Theorems \ref{thm:concentration} and \ref{thm:regularization}. In fact, the results provided below holds for \emph{inhomogeneous Erd\H{o}s-R\'{e}nyi} hypergraphs \ref{def:inhomo_ER_graph}, more general than the HSBM \ref{def:non_uniform_HSBM} , which could be of independent interest.

The \emph{inhomogeneous Erd\H{o}-R\'{e}nyi} random hypergraph is a natural generalization of standard random hypergraph models, where each hyperedge is included independently with probabilities that can vary across different vertex subsets. This flexibility enables more accurate modeling of complex networks with heterogeneous interaction patterns among vertices.

\begin{definition}[Inhomogeneous Erd\H{o}s-R\`enyi hypergraph]\label{def:inhomo_ER_graph}
    Let $\tQ^{(\ell)} \in ([0, 1]^{N})^{\otimes \ell}$ be a symmetric probability tensor, i.e., $\etQ_{i_1, \ldots, i_{\ell}} = \etQ_{i_{\pi(1)}, \ldots, i_{\pi(\ell)}}$ for any permutation $\pi$ on $[\ell]$, where $\ell \geq 2$ is some finite integer. Let $\gH_{\ell} = (\gV, \gE_{\ell})$ denote inhomogeneous $\ell$-uniform Erd\H{o}s-R\'{e}nyi hypergraph associated with the probability tensor $\tQ^{(\ell)}$. Let $\tA^{(\ell)}$ denote the adjacency tensor of $\gH_{\ell}$, where each $\ell$-hyperedge $e = \{i_1, \ldots, i_{\ell}\}\subset \gV$ appears with probability $\P(\etA^{(\ell)}_{e} = 1) = \etQ^{(\ell)}_{i_1,\dots, i_{\ell}}$. Let $\sL = \{\ell \mid \ell \geq 2, \ell \in \N \}$ be a finite set of integers, where $\LM$ denotes its maximum element. A non-uniform inhomogeneous Erd\H{o}s-R\'{e}nyi hypergraph is the union of uniform ones, i.e., $\gH = \cup_{\ell \in \sL} \gH_{\ell}$.
\end{definition}

\begin{theorem}\label{thm:concentration}
Let $\gH = \cup_{\ell \in \sL} \gH_{\ell}$ be the inhomogeneous non-uniform Erd\H{o}s-R\'{e}nyi hypergraph in \Cref{def:inhomo_ER_graph}, associated with the probability tensors $\{\tQ^{(\ell)}\}_{\ell \in \sL}$. For each $\ell \in \sL$, we rescale the tensor by $\tD^{(\ell)} \coloneqq \binom{N-1}{\ell - 1}\tQ^{(\ell)}$, and denote $\etD^{(\ell)}_{\max}\coloneqq \max_{i_1, \ldots, i_{\ell} \in \gV} \etD^{(\ell)}_{i_1,\dots, i_{\ell}}$. Denote $\etD_{\max} \coloneqq \sum\limits_{\ell \in \sL} \etD^{(\ell)}_{\max}$. For some constant $\const_{\eqref{eqn:assumption_d}}>0$, suppose that
\begin{align}
    \etD_{\max} \coloneqq \sum\limits_{\ell \in \sL} \etD^{(\ell)}_{\max} \geq \const_{\eqref{eqn:assumption_d}}\cdot \log(N)\,.\label{eqn:assumption_d}
\end{align}
Then with probability at least $1-2N^{-10}- 2e^{-N}$, the adjacency matrix $\rmA$ of $\gH$ satisfies
\begin{align}
    \|\rmA - \E \rmA \| \leq \const_{\eqref{eqn:concentrateA}}\cdot \sqrt{\etD_{\max}}\,,\label{eqn:concentrateA}
\end{align}
where the constants $\const_{\eqref{eqn:concentrateA}}\coloneqq 10\LM^2 + 2\beta$, with $\beta= \beta_0 \sqrt{\beta_1} + \LM$ and $\beta_0$, $\beta_1$\footnote{The existence of such $\beta_1$ is guaranteed since $g(x) = x\log(x) - x + 1$ takes its minimum at $x=1$ and increases when $x>1$.} 
satisfying
 \begin{align}
    &\, \beta_0 = 16+32 \LM (1+e^2)+1792(1+e^{-2})\LM^2, \notag \\
    &\, \LM^{-1}\beta_1 \log(  \LM^{-1}\beta_1) -  \LM^{-1}\beta_1 + 1 > 11/\const_{\eqref{eqn:assumption_d}}\,.
\end{align} 
\end{theorem}
Similarly matrix concentration results for the graph case has appeared in the literature several times, e.g., {\cite[Theorem 5.2]{Lei2015ConsistencyOS}}, {\cite[Theorem 5.1]{Le2017ConcentrationAR}}. The main obstacle for the generalization from graph to hypergraph \ref{def:inhomo_ER_graph} is that, the entries of the adjacency matrix $\rmA$ in \eqref{eqn:adjacency_matrix_entry} have long-range dependencies, and they are no longer zero/one entries. The proofs for \Cref{thm:concentration} and \Cref{thm:regularization} are deferred to Appendix~\ref{app:proofConcentration}.

Furthermore, \Cref{alg:spectral_initialization} could still achieve weak consistency $1\ll \rho_{N} \ll \log(N)$, which is beyond regime of exact recovery. However, \Cref{thm:concentration} is not reliable in this case, and it cannot be applied directly. Instead, the concentration could be proved for some regularized versions of $\rmA$ (\Cref{def:regularize}) following the same proof strategy, leading to \Cref{thm:regularization}. This explains why the trimming stage \Cref{alg:trimming} is necessary. 
\begin{algorithm}
\caption{\textbf{Trimming}}\label{alg:trimming}
\KwData{Hypergraph $\gH = (\gV, \gE)$, regularization strategy $\sJ$}
{Obtain the adjacency matrix $\rmA$ of the hypergraph $\gH$ by \eqref{eqn:adjacency_matrix_entry}.}\

{Construct $\rmA_{\sJ}\in \R^{N \times N}$ by zeroing out the columns $\rmA_{\cdot i}$ and rows $\rmA_{i\cdot}$, where $i\notin \sJ$}\;
\KwResult{$\rmA_{\sJ}$}
\end{algorithm}
\begin{definition}[Regularization]\label{def:regularize}
For $\rmA\in \R^{N \times N}$ and index set $\sJ$, let $\rmA_{\sJ}\in \R^{N \times N}$ be the matrix obtained from $\rmA$ by zeroing out the rows and columns not in $\sJ$. Namely,
\begin{equation}\label{eqn:restricted_adjacency_matrix}
     (\rmA_{\sJ})_{ij} = \indi{i, j\in \sJ} \cdot \ermA_{ij}\,.
\end{equation}
\end{definition}

\begin{theorem}\label{thm:regularization}
Following the conventions in \Cref{thm:concentration}, let $\rmA_{\sJ}$ be the regularized version of $\rmA$, where $\sJ \in \{ \sJ_1, \sJ_2\}$ defined as below.
\begin{enumerate}
    \item Let $\sJ_1$ be the vertex set obtained by removing $\lfloor N \exp(-\overline{d}) \rfloor$ many nodes having the largest $d_v$, where $\overline{d} = N^{-1}\sum\limits_{v\in \gV}d_v$ is the sample average degree. If $\rho_{N} \gtrsim \log(N)$, 
    there exists some constant $c>0$ such that \eqref{eqn:concentrate_regularized_A} holds with probability at least $1-2(e/2)^{-N} -2N^{-c}$.
    \item Define $\sJ_2 \coloneqq \{ v\in\gV, d_v \leq 20 \LM \etD_{\max}\}$, then the inequality \eqref{eqn:concentrate_regularized_A} holds with probability at least $1-2(e/2)^{-N} -N^{-10}$ for some constant $\const_{\eqref{eqn:concentrate_regularized_A}}$.
\end{enumerate}
\begin{align}\label{eqn:concentrate_regularized_A}
    \|\rmA_{\sJ} - \E \rmA_{\sJ} \| \leq \const_{\eqref{eqn:concentrate_regularized_A}} \,\sqrt{\etD_{\max}}\,,
\end{align}
where $\const_{\eqref{eqn:concentrate_regularized_A}}\coloneqq 10\LM^2 + 2\beta$ with $\beta= \beta_0 \sqrt{\beta_1} + \LM$, $\beta_1 = \LM$ and
\begin{align*}
    \beta_0 = 16 + 32 \LM(1+e^2)+1792(1+e^{-2})\LM^2.
\end{align*}
\end{theorem}
The construction of $\sJ_1$ is purely agnostic and does not require prior knowledge. By contrast, an upper bound on the vertex degree is pre-determined for $\sJ_2$, requiring $\{\tQ^{(\ell)}\}_{\ell \in\sL}$. Therefore, $\sJ_1$ is suitable for \Cref{alg:agnostic_partition}, while $\sJ_2$ is suitable for \Cref{alg:partition_with_prior}.

\subsection{Proofs of weak consistency} 
Without loss of generality, we assume $\alpha_1 \geq \ldots \geq \alpha_K$.
\begin{lemma}\label{lem:frobeniusNormBound}
    Under the conditions of \Cref{thm:regularization}, $\|\rmA_{\sJ}^{(K)} - \E\rmA_{\sJ}\|_{\frob}^2 \lesssim \rho_{N} + \rho_{N}^2/N$ with probability at least $1 - O(N^{-c})$ for some constant $c>0$.
\end{lemma}

\begin{lemma}\label{lem:expectedCenterSeparation}
    Under \Cref{ass:expected_center_separation} and conditions in \Cref{thm:regularization}, $\|(\E \rmA_{\sJ})_{u:} - (\E\rmA_{\sJ})_{v:}\|_2^2 \gtrsim N^{-1}\rho_{N}^2$ with probability at least $1 - O(N^{-c})$ for some constant $c>0$ when $u$ and $v$ are from different communities, otherwise $\|(\E \rmA_{\sJ})_{u:} - (\E\rmA_{\sJ})_{v:}\|_2^2 \lesssim N^{-2}\rho_{N}^2$.
\end{lemma}

\begin{lemma}\label{lem:centerFar} Recall that in \Cref{alg:spectral_initialization}, $\overline{r} = [N\log(\overline{d})]^{-1}\overline{d}\,^2$, and $\setS$ is a set of $\lceil 2\log^2(N) \rceil$ nodes randomly sampled from $\sJ$ without replacement. Define the vertex sets
\begin{align}
    \sI_k &\, \coloneqq \{v\in \gV_k \cap \sJ: \| (\rmA_{\sJ}^{(K)})_{v:} - (\E\rmA_{\sJ})_{v:}\|^2_2 \leq \overline{r}/4\}\\
    \sU_k &\, \coloneqq \{v\in \gV_k \cap \sJ: \| (\rmA_{\sJ}^{(K)})_{v:} - (\E\rmA_{\sJ})_{v:} \|^2_2 \leq 4\overline{r}\}.
\end{align}
\begin{enumerate}
    \item[(a)] For all $k\in[K]$, $|\setS\cap \sI_k| \geq 1$ with probability at least $1 - O(N^{-c})$ for some $c>0$.
    \item[(b)] For each $k\in[K]$, the $s_k$ constructed in \Cref{alg:spectral_initialization} satisfies $s_k \in \sU_k$ and $|\ball_r(s_k) |\geq \alpha_k N(1 -o(1))$.
    \item[(c)] The centers $\{s_k\}_{k =1}^{K}$ constructed in \Cref{alg:spectral_initialization} are pairwise far away from each other,
    \begin{align}
        \|(\rmA_{\sJ}^{(K)})_{s_j :} -  (\rmA_{\sJ}^{(K)})_{s_k :} \|_2^2 \gtrsim N^{-1} \rho_{N}^2, \quad \textnormal{for any pair } j\neq k,\,\, j, k\in[K].
    \end{align}
\end{enumerate}
\end{lemma}

\begin{lemma}\label{lem:misclassifiedFarCenter}
In Lines \ref{alg:first_loop} and \ref{alg:second_loop} of \Cref{alg:spectral_initialization}, if $v\in \gV_i \cap \sJ$ is misclassified to $\widehat{\gV}_{j}^{(0)}$ for some $j\neq i$, then $v$ is far away from its expected center, i.e., $\|(\rmA_{\sJ}^{(K)})_{v:} - (\E\rmA_{\sJ})_{v :} \|_2^2 \gtrsim N^{-1}\rho_{N}^2$.
\end{lemma}

\begin{proof}[Proof of \Cref{thm:agnostic_weak_consistency}]
For each $k\in[K]$, \Cref{lem:centerFar} (a) ensures that at least one vertex will be picked in each $\gV_k$, (b) ensures that the constructed ball $\ball_r(s_k)$ contains almost all the vertices in $\gV_k$. Furthermore, \Cref{lem:expectedCenterSeparation} proves that the expected centers are pairwise far away from each other, while \Cref{lem:centerFar} (c) proves that the same argument holds as well for the centers $\{s_k\}_{k=1}^{K}$ of balls $\{ \ball(s_k)\}_{k=1}^{K}$. Throughout the \Cref{alg:spectral_initialization}, vertices would be misclassified during lines \ref{alg:first_loop}, \ref{alg:second_loop}, \ref{alg:assign_vertices_outside_J}. \Cref{lem:misclassifiedFarCenter} indicates that the misclassified vertices in lines \ref{alg:first_loop}, \ref{alg:second_loop} must be far away from their expected centers. The number of vertices outside $\sJ_1$ is upper bounded by $\lfloor N \exp(-\overline{d}) \rfloor$, which is at most of constant order since $\overline{d} \gtrsim \log(N)$ according to \Cref{lem:maxProb}. Then with probability at least $1-2(e/2)^{-N} -2N^{-c}$ for some $c > 0$, the mismatch ratio is bounded by
\begin{align}
    \mismatch_{N} = &\,\mismatch(\rvy, \widehat{\rvy}^{(0)}) \leq \frac{1}{N} \Bigg( \frac{\underset{v \textnormal{ misclassified}}{\sum}\|(\rmA_{\sJ}^{(K)})_{v :} - (\E\rmA_{\sJ})_{v :}\|_2^2}{ \underset{v \textnormal{ misclassified}}{\min} \|(\rmA_{\sJ}^{(K)})_{v :} - (\E\rmA_{\sJ})_{v :}\|_2^2} + \lfloor N \exp(-\overline{d}) \rfloor \Bigg)\notag\\
    \leq &\, \|\rmA_{\sJ}^{(K)} - \E\rmA_{\sJ}\|_{\frob}^2 / \Big(N \cdot \underset{v \textnormal{ misclassified}}{\min}\| (\rmA_{\sJ}^{(K)})_{v :} - (\E\rmA_{\sJ})_{v :} \|_2^2 \Big) + \exp(-\overline{d}) \notag \\
    \lesssim &\, (\rho_{N} + \rho_{N}^2/N)/(N\cdot \rho_{N}^2/N) + N^{-1} \asymp \rho_{N}^{-1} = o(1),  \quad (\textnormal{Lemmas \ref{lem:frobeniusNormBound}, \ref{lem:misclassifiedFarCenter}}) \label{eqn:weak_consistency_agnostic}
\end{align}
where the last line holds since $\rho_{N} \gtrsim \log(N)$.
\end{proof}

\begin{proof}[Proof of \Cref{thm:known_weak_consistency}]
The regularization method $\sJ_1$ relies on two key facts.
\begin{enumerate}
    \item Only $\lfloor N \exp(-\overline{d}) \rfloor = o(N)$ many nodes are removed, since $\overline{d} \gtrsim \log(N)$ by \Cref{lem:maxProb};
    \item Degrees of nodes in $\sJ_1$ are upper bounded by $\etD_{\max}$ up to some constant when $\rho_{N} \gtrsim \log(N)$, according to Lemma~\ref{lem:boundednessJ1}.
\end{enumerate}
However, the second fact is no longer valid when $1 \ll \rho_{N} \ll \log(N)$. To get rid of that, one could refer to the other trimming process $\sJ_2$ in \Cref{thm:regularization}. Deterministically, the degree of each vertex in $\sJ_2$ is upper bounded by $20 \LM \etD_{\max}$, where the construction of $\sJ_2$ relies on the prior knowledge of probability tensors $\{\tQ^{(\ell)}\}_{\ell \in \sL}$. According to Lemma~\ref{lem:maxdeg}, the number of vertices outside $\sJ_2$ is upper bounded by $N\cdot e^{-c\rho_{N}} = o(N)$ for some constant $c>0$. Then by following the same analysis above, with probability at least $1-2(e/2)^{-N} -N^{-10}$, the mismatch ratio is bounded by
\begin{align}
    \mismatch_{N} \lesssim \frac{\rho_{N} + \rho_{N}^2/{N}}{N\rho_{N}^2/N} + N^{-1} \cdot (N/e^{c\rho_{N}}) \asymp \rho_{N}^{-1} + N^{-1} + e^{-c\rho_{N}}\asymp \rho_{N}^{-1} = o(1). \label{eqn:weak_consistency_known}
\end{align}
Therefore, \Cref{alg:spectral_initialization} achieves weak consistency by applying $\sJ_2$ as long as $\rho_{N} = \omega(1)$.
\end{proof}

\section{Stage II: Agnostic iterative refinement}\label{sec:AgnosticRefinement}
This section is devoted to the proofs of \Cref{thm:agnostic_strong_consistency} and \Cref{thm:agnostic_optimality}, while the proofs of all the Lemmas are deferred to Appendix~\ref{app:AgnosticRefinement}. We first introduce \Cref{alg:agnostic_refinement}, which achieves the optimal recovery by iteratively refining the initial estimate $\widehat{\rvy}^{(0)}$ obtained from \Cref{alg:spectral_initialization}. 

\subsection{Tensor probability reconstruction}

The difference between $\etQ^{(\ell)}_{\rvw}$ and $\widehat{\etQ}^{(\ell)}_{\rvw}$ is controlled by \Cref{lem:hatQmw_approx}. It is proved that $\widehat{\etQ}^{(\ell)}_{\rvw}$ is a good estimator of $\etQ^{(\ell)}_{\rvw}$ under \Cref{ass:prob_ratio_bound}.

\begin{lemma}\label{lem:hatQmw_approx}
    For each $\ell \in \sL$ and $\rvw \in \WC{\ell}{K}$, recall $\etQ^{(\ell)}_{\rvw}$ in \eqref{eqn:edge_density} and $\widehat{\etQ}^{(\ell)}_{\rvw}$ in \eqref{eqn:hat_Qlw}. Define $\etQ^{(\ell)}_{\max} \coloneqq \max_{\rvw \in \WC{\ell}{K}} \etQ^{(\ell)}_{\rvw}$ and denote $\etQ^{(\ell)}_{\max} = \etP_{\max}^{(\ell)}q_{\max}^{(\ell)}/\binom{N-1}{\ell - 1}$ using \eqref{eqn:edge_density} where $\etP_{\max}^{(\ell)} \asymp 1$ and $1\ll q_{\max}^{(\ell)} \ll N$. The following holds with probability at least $1 - \exp(- N\log q_{\max}^{(\ell)}/q_{\max}^{(\ell)})$, 
\begin{align}
    |\widehat{\etQ}^{(\ell)}_{\rvw} - \etQ^{(\ell)}_{\rvw} | \lesssim \log q_{\max}^{(\ell)}/N^{\ell - 1}.
\end{align} 
Furthermore, $|\widehat{\etQ}^{(\ell)}_{\rvw} - \etQ^{(\ell)}_{\rvw} |/ \etQ^{(\ell)}_{\rvw} = \log q_{N} /q_{N} = o(1)$ when under \Cref{ass:prob_ratio_bound}.
\end{lemma}

The precision of the estimation $\widehat{\etQ}^{(\ell)}_{\rvw}$ is based on the fact that the number of edges corresponding to partition $\rvw$ is $q_{N}$ many for each $\rvw \in \WC{\ell}{K}$ and $\ell \in \sL$, while the number of different types $\rvw$ is finite. As long as all $\etQ^{(\ell)}_{\rvw}$ are comparable (\Cref{ass:prob_ratio_bound}), this estimation should be correct. In the absence of \Cref{ass:prob_ratio_bound}, e.g., $q_{\max}^{(\ell)} = N^{\beta}$ for some $\beta \in (0, 1)$ while $q_{\rvw}^{(\ell)} = \log\log(N)$, $|\widehat{\etQ}^{(\ell)}_{\rvw} - \etQ^{(\ell)}_{\rvw}|/ |\etQ^{(\ell)}_{\rvw}| = o(1)$ is no longer true, hence $\widehat{\etQ}^{(\ell)}_{\rvw}$ can't be accurate. 

We present briefly a discussion of \Cref{ass:prob_ratio_bound} here. According to \Cref{lem:hatQmw_approx}, \Cref{ass:prob_ratio_bound} is crucial for $\widehat{\etQ}^{(\ell)}_{\rvw}$ being an accurate estimation of $\etQ^{(\ell)}_{\rvw}$. In fact, when a mixture of Bernoulli random variables is collected where their expectations satisfy \Cref{ass:prob_ratio_bound}, it is possible to estimate each expectation accurately, since the sample average for each Bernoulli random variable (in this scenario, all but a negligible fraction of samples are chosen from their true law) would be well concentrated around its expectation with lower order deviations. However, there is no hope to estimate the lower order expectations accurately if their expectations vary in magnitude, since the fluctuation of the sample average with higher order expectation could be significantly larger than the lower order expectations. 

\subsection{Proofs of strong consistency and optimality}
 We first show that $\gV$ can be partitioned into two disjoint subsets $\gG$ and $\gV\setminus \gG$, where each $v\in \gG$ enjoys the following three properties:
\begin{description}
    \item[(G1)] Denote $\etQ^{(\ell)}_{\max} \coloneqq \max_{\rvw \in \WC{\ell}{K}}\etQ_{\rvw}^{(\ell)}$, $\etD^{(\ell)}_{\max} \coloneqq \binom{N-1}{\ell - 1}\etQ^{(\ell)}_{\max}$ and $\etD_{\max} \coloneqq \sum_{\ell \in \sL}\etD^{(\ell)}_{\max}$. The degree of $v$ is upper bounded by
        \begin{align}\label{eqn:G1}
            \rD_{v} \leq 10 \etD_{\max}.\tag{G1}
        \end{align}
    \item[(G2)] Let $\mu^{(\ell)}_{v, \rvw}\coloneqq \rD^{(\ell)}_{v, \rvw}/ N_{\rvw}$ and let $\D_{\rm{KL}}(\mu^{(\ell)}_{v, \rvw}\parallel \etQ_{k \oplus \rvw}^{(\ell)})$ be the KL divergence in \eqref{eqn:KLDivergence}. If $v\in \gV_k \cap \gG$, then the weighted sum of $\D_{\rm{KL}}(\mu^{(\ell)}_{v, \rvw}\parallel \etQ_{k \oplus \rvw}^{(\ell)})$ is much smaller than the weighted sum of $\D_{\rm{KL}}(\mu^{(\ell)}_{v, \rvw}\parallel \etQ_{j \oplus \rvw}^{(\ell)})$ for any other $j\neq k$. Formally,
        \begin{align}\label{eqn:G2}
            \sum\limits_{\ell \in \sL} \sum\limits_{ \rvw \in \WC{\ell - 1}{K} } N_{\rvw} \cdot [ \D_{\rm{KL}}(\mu^{(\ell)}_{v, \rvw}\parallel \etQ_{j \oplus \rvw}^{(\ell)}) - \D_{\rm{KL}}(\mu^{(\ell)}_{v, \rvw}\parallel \etQ_{k \oplus \rvw}^{(\ell)}) ] \geq \etD_{\max}/ \log q_{N}.\tag{G2}
        \end{align}
    \item[(G3)] Let $\gE(v, \gV\setminus \gG)$ denote the set of hyperedges containing $v$ with all the remaining nodes in $\gV\setminus \gG$. Let $\LM$ denote the maximum element of $\sL$. Let $\delta > 0$ be some absolute constant. The number of edges containing $v$ is bounded by
        \begin{align}\label{eqn:G3}
            |\gE(v, \gV\setminus \gG)| \leq \LM (\etD_{\max})^{1 + \delta}.\tag{G3}
        \end{align}
\end{description}
We explain how this construction can be done in the proof of \Cref{lem:sizeofG}. Briefly, we show that most of vertices satisfy \eqref{eqn:G1} and \eqref{eqn:G2}, and proceed to eliminate the nodes that do not satisfy \eqref{eqn:G3}. The number of vertices we need to eliminate is small.
\begin{lemma}[Size of $\gV \setminus \gG$]\label{lem:sizeofG} For the sequence $\kappa_{N} \in (0, 1]$ with $\kappa_{N} \log(N) \to \infty$ as $N \to \infty$, assume \eqref{eqn:optimality_condition} for some absolute (small) constant $\epsilon > 0$. Then $|\gV\setminus \gG| < N^{1 - \kappa_{N}}$ with probability at least $1 - 6e^{- \epsilon \kappa_{N}\log(N)}$.
\end{lemma}

Let $\widehat{\gN}^{(t)}_{jk} \coloneqq (\widehat{\gV}^{(t)}_{j} \cap \gV_k) \cap \gG$ denote the set of vertices in $\gG$ which belongs to $\gV_k$ but misclassified to $\widehat{\gV}^{(t)}_{j}$, then the misclassified vertices in $\widehat{\gV}^{(t)}_{j} \cap \gG$ can be represented as $\widehat{\gN}^{(t)}_j \coloneqq \cup_{k\neq j} \widehat{\gN}^{(t)}_{jk} = (\widehat{\gV}^{(t)}_j \setminus \gV_j) \cap \gG$.  Let $\widehat{\gN}^{(t)} \coloneqq \cup_{j=1}^{K} \widehat{\gN}^{(t)}_j = \cup_{j=1}^{K}( \widehat{\gV}_{j}^{(t)} \setminus \gV_j ) \cap \gG$ denote the set of misclassified vertices in $\gG$ at step $t$. \Cref{lem:decayofError} shows that the size of $|\widehat{\gN}^{(t)}|$ decays over consecutive iterations.
\begin{lemma}\label{lem:decayofError}
In the outer for-loop of \Cref{alg:agnostic_refinement}, the number of misclassified nodes decays
\begin{align}
    |\widehat{\gN}^{(t+1)}|/|\widehat{\gN}^{(t)}| \leq e^{-1}\cdot \log\const_{\eqref{eqn:prob_ratio_bound} }/ \const_{\eqref{eqn:prob_ratio_bound} }\,,
\end{align}
with probability at least $1-2(e/2)^{-N} -2N^{-c}$ for some constant $c > 0$.
\end{lemma}

\begin{proof}[Proof of \Cref{thm:agnostic_optimality}] Recall that $\mismatch_{N} = \mismatch(\rvy, \widehat{\rvy}^{(0)}) \lesssim (\rho_{N})^{-1}$ according to \eqref{eqn:weak_consistency_agnostic}. Note that the agnostic result \Cref{thm:agnostic_weak_consistency} requires $\rho_{N} \gtrsim \log(N)$, then by \Cref{lem:decayofError}, the number of misclassified nodes in $\gG$ after $\lceil \log(N) \rceil$ iterations will be at most 
\begin{align}
   |\widehat{\gN}^{\lceil \log(N) \rceil}|\leq &\, \mismatch_{N} N \cdot (e^{-1} \log\const_{\eqref{eqn:prob_ratio_bound} }/ \const_{\eqref{eqn:prob_ratio_bound} } )^{\lceil \log(N) \rceil} \notag\\
   \lesssim &\, N \cdot  N^{-1 - \log\const_{\eqref{eqn:prob_ratio_bound} }}/\log(N) = o(1)\,,\notag
\end{align}
where the inequality holds since $\const_{\eqref{eqn:prob_ratio_bound} } > 1$. Consequently, $\widehat{\gN}^{\lceil \log(N) \rceil} = \emptyset$, indicating that all nodes in $\gG$ are correctly assigned. Thus, the misclassified nodes are all in the set $\gV\setminus \gG$, leading to $N \cdot \mismatch_{N} \leq |\gV\setminus \gG| \leq N^{1 - \kappa_{N}}$. Therefore, the desired argument follows with probability at least $1 -6e^{-\epsilon \cdot \kappa_{N} \log(N)} - (1 + o(1))\lceil \log(N) \rceil N^{-c}$.
\end{proof}
\begin{proof}[Proof of \Cref{thm:agnostic_strong_consistency}] 
 Denote $\kappa_{N} \coloneqq \D_{\mathrm{GCH}}\cdot q_{N}/\log(N)$. In addition, assume that the limit of the sequence $\{\kappa_{N}\}_{N}$ exists, denoted by $\kappa = \lim\limits_{N \to \infty}\kappa_{N}$. There are three scenarios.
 
\noindent (a) $\kappa \in (0, 1]$. By \Cref{thm:IT_lower_bounds} (2), $\E[\mismatch_{N}] \geq N^{-(1 + o(1)) \cdot \kappa}$. According to \Cref{thm:agnostic_optimality}, $\mismatch_{N} \leq N^{-\kappa_{N}}$ with probability at least $1 -6e^{-\epsilon \cdot \kappa_{N} \log(N)} - (1 + o(1))\lceil \log(N) \rceil N^{-c}$. Therefore, we have $\E[\mismatch_{N}] = N^{-\kappa} = o(1)$, reaching the lowest expected mismatch ratio possible.

\noindent (b) $\kappa = 0$ with $\kappa_{N} \gg \log^{-1}(N)$; we still have $\E[\mismatch_{N}] \geq N^{-(1 + o(1))\kappa_{N}} = o(1)$ by \Cref{thm:IT_lower_bounds} (2), and \Cref{alg:agnostic_refinement} achieving the lowest possible expected mismatch ratio. Otherwise when $\kappa_{N} \lesssim \log^{-1}(N)$, we have $N^{-\kappa_{N}} = O(1)$, which represents the regime of partial or weak recovery, not the focus of this paper.

\noindent (c) $\kappa > 1$, equivalently \eqref{eqn:exact_recovery_condition} holds. Following the same analysis above, $\gV \setminus \gG = \emptyset$ with probability at least $1 - 6e^{- \epsilon \log(N)}$ by \Cref{lem:sizeofG}. All nodes in $\gG = \gV$ will be classified correctly by \Cref{lem:decayofError}, achieving exact recovery. Note that when $\kappa = 1$, it is not clear whether the exact recovery can be achieved or not.

Furthermore, if the sequence $\{\kappa_{N}\}_{N}$ goes to infinity as $N \to \infty$, then clearly, \Cref{alg:agnostic_partition} accomplishes exact recovery.
\end{proof}
We note that \Cref{alg:agnostic_refinement} utilizes the adjacency tensors $\{ \tA^{(\ell)}\}_{\ell \in \sL}$. The proofs of \Cref{thm:agnostic_optimality} and \Cref{thm:agnostic_strong_consistency} do not rely on \Cref{ass:expected_center_separation}. For the set of problems violating \Cref{ass:expected_center_separation}, if a weak consistent partition is obtained, it can be improved by \Cref{alg:agnostic_refinement}, achieving the lowest expected mismatch ratio possible.



\section{Partition with prior knowledge}\label{sec:known_partition}
This section is devoted to the proofs of Theorems \ref{thm:known_strong_consistency} and \ref{thm:known_optimality}, with the proofs of Lemmas deferred to \Cref{app:KnownPartition}.

\subsection{Optimality of \textnormal{MAP}}
Let $\P_{\rm{err}}(v)$ denote the probability of $v\in \gV$ being misclassified. Let $\widehat{\ervy}_{v}$, $\widehat{\rvy}_{-v}$ be some estimations of $\ervy_v$ and $\rvy_{-v}$ respectively. Let $\P(\cdot \mid \mathbb{H} = \gH,\,\rmY_{-v} = \widehat{\rvy}_{-v})$ denote the conditional probability when observing $\gH$ and $\widehat{\rvy}_{-v}$, where $\mathbb{H}$ and $\rmY_{-v}$ denote the laws of $\gH$ and $\rvy_{-v}$, respectively. Then by definition,
\begin{align}
    \P_{\rm{err}}(v) \coloneqq &\,\P( \widehat{\ervy}_{v} \neq \ervy_{v} ) = \sum\limits_{\gH, \widehat{\rvy}_{-v}} \P(\widehat{\ervy}_{v}\neq \ervy_{v} \mid \mathbb{H} = \gH,\,\rmY_{-v} = \widehat{\rvy}_{-v}) \cdot \P(\mathbb{H} = \gH, \rmY_{-v} = \widehat{\rvy}_{-v}) \notag\\
    =&\, \sum\limits_{\gH, \widehat{\rvy}_{-v}} \big[1 - \P (\widehat{\ervy}_{v} = \ervy_{v} \mid \mathbb{H} = \gH,\,\rmY_{-v} = \widehat{\rvy}_{-v}) \big]\cdot \P(\mathbb{H} = \gH,\, \rmY_{-v} = \widehat{\rvy}_{-v})\,\,.\label{eqn:errorProbability}
\end{align}
By Bayes, the term $\P(\mathbb{H} = \gH,\, \rmY_{-v} = \widehat{\rvy}_{-v})$ is irrelevant to any specific assignment of $\ervy_v$, since
\begin{align}
   &\,\P(\mathbb{H} = \gH,\, \rmY_{-v} = \widehat{\rvy}_{-v}) = \P(\mathbb{H} = \gH \mid \rmY_{-v} = \widehat{\rvy}_{-v}) \cdot \P( \rmY_{-v} = \widehat{\rvy}_{-v})\\
   = &\, \P( \rmY_{-v} = \widehat{\rvy}_{-v}) \cdot \sum\limits_{j\in [K]} \P(\mathbb{H} = \gH \mid \rY_{v} = j,\, \rmY_{-v} = \widehat{\rvy}_{-v}) \cdot \P( \rY_{v} = j)
\end{align}
Therefore, $\P_{\rm{err}}(v)$ is minimized only if the MAP estimator $\widehat{\ervy}_{v}^{\rm{\, MAP}}$ in \eqref{eqn:MAP} is substituted into \eqref{eqn:errorProbability}. However, $\widehat{\ervy}_{v}^{\rm{\, MAP}}$ in \eqref{eqn:MAP} cannot be computed directly in practice as explained. Instead, \Cref{alg:partition_with_prior} utilized the simplified version \eqref{eqn:MAP_single} inspired by \Cref{lem:MAP}. Recall $\rD_v = \sum_{\ell, \rvw}\rD^{(\ell)}_{v, \rvw}$ in \eqref{eqn:Dv} and $d_v = \sum_{\ell, \rvw}d^{(\ell)}_{v, \rvw}$, where $\rD^{(\ell)}_{v, \rvw} \vert_{\rY_{v} = k} \sim {\rm{Bin}}(N_{\rvw}, \etQ^{(\ell)}_{k\oplus \rvw})$ and $d^{(\ell)}_{v, \rvw}$ denote its realization. From simplicity, we denote
\begin{subequations}
\begin{align}
    \widehat{\P}_{k}(\rD_v = d_v) &\, \coloneqq \prod\limits_{\ell \in \sL}\prod\limits_{\rvw\in \WC{\ell - 1}{K} }\P(\rD^{(\ell)}_{v, \rvw} = d^{(\ell)}_{v, \rvw} | \rY_{v} = k,\, \rmY_{-v} = \widehat{\rvy}_{-v}), \label{eqn:conditional_hat_Dvk} \\
    \P_{k}(\rD_v = d_v) &\, \coloneqq \prod\limits_{\ell \in \sL}\prod\limits_{\rvw\in \WC{\ell - 1}{K} }\P(\rD^{(\ell)}_{v, \rvw} = d^{(\ell)}_{v, \rvw} | \rY_{v} = k,\, \rmY_{-v} = \widehat{\rvy}_{-v}).\label{eqn:conditionalDvk}
\end{align}
\end{subequations}
According to \Cref{lem:MAP} and edgewise independence, the probabilities in \eqref{eqn:errorProbability} that are irrelevant to $v$ cancels out. By substituting $\widehat{\ervy}_{v}^{\rm{\, MAP}}$ in \eqref{eqn:MAP_single} into \eqref{eqn:errorProbability}, it follows that\footnote{The summation in \eqref{eqn:errorProbabilityDv}, \eqref{eqn:upper_lower_bound_error_prob} and 
\eqref{eqn:circledTerms} is over all possible realizations $\{d^{(\ell)}_{v, \rvw}\}^{\ell \in \sL}_{\rvw \in \WC{\ell - 1}{K}}$.}
    \begin{align}
        \P_{\rm{err}}(v) = &\,\sum \P(\widehat{\ervy}_{v}^{\rm{\, MAP}} \neq k \big| \rD_v = d_v, \,\rmY_{-v} = \widehat{\rvy}_{-v}) \cdot \P_k(\rD_v = d_v) \notag\\
        =&\, \sum \P_k(\rD_v = d_v) \cdot \P \big( \exists j \neq k, \textnormal{s.t. } \alpha_j \widehat{\P}_j(\rD_v = d_v) > \alpha_k\widehat{\P}_k (\rD_v = d_v) \big). \label{eqn:errorProbabilityDv}
    \end{align}
Furthermore,  $\P_{\rm{err}}(v)$ admits the upper and lower bounds presented in \Cref{lem:boundsErrorProb}.

\begin{lemma}\label{lem:boundsErrorProb}
Assume $v\in \gV_k$ for some $k\in [K]$. Define the following quantity
\begin{align}
    \widehat{\rR}_{k} \coloneqq \sum\limits_{j\in[K], \, j\neq k} \P \big( \exists j \neq k, \textnormal{s.t. } \alpha_j \widehat{\P}_j(\rD_v = d_v) > \alpha_k\widehat{\P}_k (\rD_v = d_v) \big).\label{eqn:ProbComparation}
\end{align}
Then the upper and lower bounds of $\P_{\rm{err}}(v)$ are given by
    \begin{align}
        (K-1)^{-1}\sum\P_k(\rD_v = d_v) \cdot \widehat{\rR}_{k} \, \leq \, \P_{\rm{err}}(v) \, \leq \, \sum \P_k(\rD_v = d_v) \cdot \widehat{\rR}_{k}.\label{eqn:upper_lower_bound_error_prob} 
    \end{align}
\end{lemma}

\subsection{Error decomposition} As discussed previously, $\P_{\rm{err}}(v)$ can be characterized from both direction, as long as $\widehat{\rR}_{k}$ is evaluated accurately. The difficulty lies on the fact that \eqref{eqn:conditional_hat_Dvk} is not directly accessible for analysis due to non-equivalence between $\widehat{\rvy}_{-v}$ and $\rvy_{-v}$, and we do not know how far it deviates from the true probability \eqref{eqn:conditionalDvk}. Fortunately, the ratio $\eqref{eqn:conditional_hat_Dvk}/\eqref{eqn:conditionalDvk}$ can be bounded from both sides in almost all cases with some allowable errors, thus \eqref{eqn:conditional_hat_Dvk} can be replaced by \eqref{eqn:conditionalDvk} when evaluating $\widehat{\rR}_{k}$ in \eqref{eqn:ProbComparation}.

We now present the proof sketch. Let $\gE_{\ell}(v)$ denote the set of $\ell$-edges containing $v$, while $\gE_{\ell}(v, \widehat{\gN}^{(0)})$ be its subset where each edge contains $v$ and at least one node from $\widehat{\gN}^{(0)}$, where $\widehat{\gN}^{(0)}$ in \Cref{lem:decayofError} denotes the set of misclassified nodes after initial estimation. Note that $\rvy(e) = \widehat{\rvy}^{(0)}(e)$ for $e\notin \gE_{\ell}(v, \widehat{\gN}^{(0)})$, then the desired ratio is
\begin{align}
    \frac{\eqref{eqn:conditional_hat_Dvk}}{\eqref{eqn:conditionalDvk}} 
    = &\, \prod\limits_{\ell \in \sL} \,\, \prod\limits_{ e \in \gE_{\ell}(v, \,\, \widehat{\gN}^{(0)})} \, \Bigg[ \frac{ \widehat{\P}_k(\etA_{e}^{(\ell)} = 1 ) }{ \P_k( \etA_{e}^{(\ell)} = 1) } \Bigg]^{\etA_{e}^{(\ell)}} \cdot \Bigg[ \frac{ \widehat{\P}_k(\etA_{e}^{(\ell)} = 0 ) }{ \P_k( \etA_{e}^{(\ell)} = 0) } \Bigg]^{1 - \etA_{e}^{(\ell)}},
\end{align}
since the probabilities for edges $e\notin \gE_{\ell}(v, \widehat{\gN}^{(0)})$ cancel out. We say that the realization $\gE_{\ell}(v)$ is \emph{good} if the cardinality of its subset $\gE_{\ell}(v, \widehat{\gN}^{(0)})$ is relatively small. We say that the edge set $\gE(v) = \cup_{\ell \in \sL}\,\gE_{\ell}(v)$ is good if $\gE_{\ell}(v)$ is good for each $\ell \in \sL$.
\begin{definition}[Good realization]\label{def:good_degree}
Recall $\etQ^{(\ell)}_{\max}$, $q_{\max}^{(\ell)}$ defined in \Cref{lem:hatQmw_approx}. We say that the realization $\gE(v)$ is \emph{good} if $\rG (v) \coloneqq \prod_{\ell \in \sL}\rG_{\ell}(v) = 1$, where for each $\ell \in \sL$
  \begin{align}
     \rG_{\ell}(v) = \indi{|\gE_{\ell}(v, \widehat{\gN}^{(0)})| \leq 4 q^{(\ell)}_{\max} / \log( q^{(\ell)}_{\max})}. \label{eqn:Dvl_good}
 \end{align}
 Conversely, $\gE(v)$ is \emph{Not good} if $\rG(v) = 0$, i.e., $\rG_{\ell}(v) = 0$ for at least one $\ell \in \sL$.
\end{definition}
According to Bayes, we have $\widehat{\rR}_{k} = \widehat{\rR}_{k | 0} \cdot\P(\rG(v) = 0) + \widehat{\rR}_{k | 1} \cdot \P(\rG(v) = 1)$, where
\begin{align}
    \widehat{\rR}_{k| \ell} \coloneqq \sum\limits_{j\in[K], \, j\neq k} \P \big[ \exists j \neq k, \textnormal{s.t. } \alpha_j \widehat{\P}_j(\rD_v = d_v) > \alpha_k\widehat{\P}_k (\rD_v = d_v) \big \vert \rG(v) = \ell\big],\label{eqn:ConditionedProbComparation}
\end{align}
for $\ell\in \{0, 1\}$. Consequently, $\P_{\rm{err}}(v)$ in \eqref{eqn:errorProbabilityDv} is upper bounded by $\P_{\rm{err}}(v) \leq \circled{0} + \circled{1}$, where
\begin{align}
    \circled{$\ell$} = \sum \P_k(\rD_v = d_v) \cdot \widehat{\rR}_{k | \ell} \cdot \P(\rG(v) = \ell). \label{eqn:circledTerms}
\end{align}

\noindent \textbf{Bound for} \circled{0} The probability of $\gE(v)$ being \emph{Not good} can be controlled by $\D_{\rm{GCH}}\cdot q_{N}$.
\begin{lemma}\label{lem:badProbUpperBound}
For $v\in \gV_k$, we have $\P_k(\rG(v) = 0) \lesssim \exp(- (4 + o(1))\cdot \D_{\rm{GCH}} \cdot q_{N})$.
\end{lemma}
Since $\widehat{\rR}_{k|0} \leq K - 1$ as defined in \eqref{eqn:ConditionedProbComparation} and $\sum_{d_v}\P_k(\rD_v = d_v) = 1$, it then follows that
\begin{align}
    \circled{0} \lesssim \exp(- (4 + o(1))\cdot \D_{\rm{GCH}} \cdot q_{N}). \label{eqn:0UpperBound}
\end{align}

\noindent \textbf{Bound for} \circled{1}. When $\rG(v) = 1$, the ratio $\eqref{eqn:conditional_hat_Dvk}/\eqref{eqn:conditionalDvk}$ can be bounded from both sides as proved by \Cref{lem:boundsofHatProbability}, which further leads to an upper bound on $\widehat{\rR}_{k | 1}$, as well as \circled{1}.

\begin{lemma}\label{lem:boundsofHatProbability}
Suppose \eqref{eqn:Strong_consistency_prob_condition}. Conditioned on $v\in \gV_k$ and $\rG(v) = 1$, there exists some $U_k \in \R^{+}$ with $\log U_k \ll \D_{\mathrm{GCH}}(j, k) \cdot q_{N}$ for any $j\neq k$, such that $U^{-1}_k \leq \widehat{\P}_k( \rD_v = d_v )/\P_k(\rD_v = d_v) \leq U_k$ with probability at least $1 - N^{-10}$.
\end{lemma}
With some arguments deferred to \Cref{app:KnownPartition}, with probability $1 - N^{-10}$, \circled{1} is upper bounded by
\begin{align}
     \circled{1} \lesssim \exp(- (1 + o(1)) \cdot \D_{\rm{GCH}}\cdot q_{N}).\label{eqn:1UpperBound}
\end{align}

\noindent \textbf{Bound for} $\P_{\rm{err}}(v)$. Combining \eqref{eqn:upper_lower_bound_error_prob} and upper bounds above for \circled{0} and \circled{1}, it follows that
 \begin{align}
     \P_{\rm{err}}(v) \leq \circled{0} + \circled{1} \lesssim (1 + o(1)) \cdot e^{- (1 + o(1)) \cdot \D_{\rm{GCH}} \cdot q_{N}}, \label{eqn:UpperBoundErrorProb}
 \end{align}
 with probability at least $1 - N^{-10}$, where the last inequality holds since $\D_{\rm{GCH}} \cdot q_{N} \to \infty$ as $N \to \infty$. The term $\D_{\rm{GCH}} \cdot q_{N}$ plays a crucial role in understanding the limit of \Cref{alg:partition_with_prior}.
\subsection{Algorithm correctness} 
The complete proofs are presented in this subsection.

\begin{proof}[Proof of \Cref{thm:known_optimality}] Suppose \eqref{eqn:optimality_condition} where the sequence $\kappa_{N} \in (0, 1]$ with $\kappa_{N} \log(N) \to \infty$, then there exists some universal (small) constant $\epsilon >0$ such that $\D_{\rm{GCH}} \cdot q_{N} \geq (1 + \epsilon)\kappa_{N}\log(N)$. According to \eqref{eqn:UpperBoundErrorProb}, $\P_{\rm{err}}(v) \lesssim (1 + o(1)) \cdot N^{-(1 + \epsilon)\kappa_{N}}$ for each $v\in \gV$, then with probability at least $1 - N^{-10}$, the expected mismatch ratio is upper bounded by
\begin{align}
\E[\mismatch_{N}] = \E[\# \textnormal{ of misclassified nodes}]/N 
     =\sum\P_{\rm{err}}(v)/N \lesssim (1 + o(1)) \cdot N^{-(1 + \epsilon)\kappa_{N}}.
\end{align}
An argument similar to the proof of \Cref{thm:agnostic_strong_consistency}, together with \Cref{thm:IT_lower_bounds} (2), shows that \Cref{alg:partition_with_prior} reaches the lowest possible expected mismatch ratio.
 \end{proof}
 
\begin{proof}[Proof of \Cref{thm:known_strong_consistency}] 
Note that $\D_{\rm{GCH}} \cdot q_{N} \geq (1 + \epsilon)\log(N)$ when taking $\kappa_{N} = 1$. Then the probability of $\widehat{\ervy}_{v}^{\rm{\, MAP}}$ in \eqref{eqn:MAP_single} failing exact recovery is upper bounded by
\begin{align}
    \P \big( \exists v\in\gV,\,\, \widehat{\ervy}_{v}^{\rm{\, MAP}} \neq \ervy_v \big) \leq \sum \P_{\rm{err}}(v) \leq N \cdot (1 + o(1)) \cdot N^{-(1 + \epsilon)\kappa_{N}} \leq (1 + o(1)) \cdot N^{-\epsilon}.
\end{align}
Therefore, \Cref{alg:partition_with_prior} achieves exact recovery with probability at least $1 - N^{-\epsilon} - N^{-10}$.
\end{proof}

\section{Estimation for number of communities}\label{sec:communityEstimation}

For the practical purpose, it is better to assume as little prior knowledge as possible. To make \Cref{alg:agnostic_partition} completely agnostic, we propose the following \Cref{alg:communityEstimation} for the estimation of the number of communities.
\begin{algorithm}
\caption{\textbf{Number of Communities}}\label{alg:communityEstimation}
\KwData{ $\{\tA^{(\ell)}\}_{\ell \in \sL}$}
{Compute the degree $d_v$ for each vertex $v\in \gV$ and find $\widetilde{d} = \max_{v\in\gV}d_v$}\;
{Construct the adjacency matrix $\rmA$ and compute the eigenvalues $\lambda_1(\rmA) \geq \ldots \geq \lambda_{N}(\rmA)$}\;
{Find the first eigenvalue $\lambda_r(\rmA)$ such that $\lambda_r \leq (\widetilde{d})^{3/4}$}\;
\KwResult{$\widehat{K} = r-1$}
\end{algorithm}

The functionality of \Cref{alg:communityEstimation} is based on \Cref{ass:eigenvalueGap} (\emph{eigenvalue separation}) and \Cref{ass:prob_ratio_bound} (\emph{proportional degree}) under the exact recovery regime \eqref{eqn:edge_density}. On the other hand, one may refer to \cite{Stephan2024Community} for estimating the number of communities under the bounded degree regime. Recall that $\rho_{N}, \xi_{N}$ in \eqref{eqn:rhoxi} represent the maximum and average expected degrees, respectively. Let $\zeta_{N} \coloneqq \min_{v\in \gV}\E (\rD_{v})$ denote the minimum expected degree.


\begin{assumption}\label{ass:eigenvalueGap}
    Define the matrix $\E \widetilde{\rmA}$ where $\E \widetilde{\ermA}_{ij} = \E \ermA_{ij}$ for $i\neq j$ and $\E \widetilde{\ermA}_{ii} = \E \ermA_{ij}$ when $i$ and $j$ are from the same block. The eigenvalues $\lambda_1(\E \widetilde{\rmA}) \geq \ldots \geq \lambda_{N}(\E \widetilde{\rmA})$ satisfy $\lambda_K(\E \widetilde{\rmA}) \gtrsim \zeta_{N}$ and $\lambda_{K+1}(\E \widetilde{\rmA}) = 0$.
\end{assumption}

\begin{lemma}\label{lem:eigenvaluesEA}
   Under Assumptions \ref{ass:prob_ratio_bound}, \ref{ass:eigenvalueGap}, $\lambda_k(\E \rmA) \gtrsim \rho_{N}$ for $k\in[K]$, and $\lambda_k(\E \rmA) \lesssim N^{-1}\rho_{N}$ for $k\in \gV\setminus[K]$.
\end{lemma}

\begin{proof}[Proof of \Cref{lem:eigenvaluesEA}]
    Define the diagonal matrix $\rmD \coloneqq \E \widetilde{\rmA} - \E\rmA$, then $\lambda_k(\rmD) = \E \widetilde{\ermA}_{kk} \lesssim N^{-1}\rho_{N}$. Note the facts that $\zeta_{N} \asymp \rho_{N}$ by \Cref{ass:prob_ratio_bound}, and $\lambda_K(\E \widetilde{\rmA}) \gtrsim \zeta_{N} \gg N^{-1}\rho_{N}$ by \Cref{ass:eigenvalueGap}, then the desired result follows by Weyl's inequality \ref{lem:weyl}, since $\lambda_{k}(\E \widetilde{\rmA}) - \|\rmD\| \leq \lambda_k(\E\rmA) \leq \lambda_{k}(\E \widetilde{\rmA}) + \|\rmD\|$ for all $k\in \gV$.
\end{proof}
\begin{proof}[Correctness of \Cref{alg:communityEstimation}]
    By Weyl's inequality \ref{lem:weyl}, for any $k\in \gV$, the eigenvalue $\lambda_k(\E \rmA)$ satisfies
    \begin{align}
        \lambda_{k}(\E \rmA) - \|\rmA - \E \rmA\| \leq  \lambda_{k}(\rmA) \leq \lambda_{k}(\E \rmA) + \|\rmA - \E \rmA\|.
    \end{align}
By \Cref{lem:maxProb} $\widetilde{d} \asymp \etD_{\max} \asymp \rho_{N}$, then for sufficiently large $N$, by \Cref{lem:eigenvaluesEA} and \Cref{thm:concentration}
\begin{align}
    \lambda_{K}(\rmA) &\, \gtrsim \rho_{N} - \const_{\eqref{eqn:concentrateA}}\cdot \sqrt{\etD_{\max}} \, \gtrsim \rho_{N} - \sqrt{\rho_{N}} \gg (\rho_{N})^{3/4} \asymp (\widetilde{d})^{3/4},\\
    \lambda_{K+1}(\rmA) &\, \leq N^{-1}\rho_{N} +  \const_{\eqref{eqn:concentrateA}}\cdot \sqrt{\etD_{\max}} \, \lesssim \sqrt{\rho_{N}} \ll (\rho_{N})^{3/4} \asymp (\widetilde{d})^{3/4}.
\end{align}
Therefore, the accuracy of $\widehat{K}$ in \Cref{alg:communityEstimation} is guaranteed by the separation between $\lambda_{K}(\rmA)$ and $\lambda_{K+1}(\rmA)$.
\end{proof}

We have seen that \Cref{alg:communityEstimation} works when all the degrees are proportional. One may be interested in exploring modifications of \Cref{alg:communityEstimation} for scenarios not covered by Assumptions \ref{ass:prob_ratio_bound}, \ref{ass:eigenvalueGap}, such as the presence of communities with outstandingly large degrees.

\subsection{Eigenvalue separation}
In particular, we want to emphasize that \Cref{ass:eigenvalueGap} is stronger than \Cref{ass:expected_center_separation}. Consider the following example, where $K = 3$, $\LM =2$ with $\alpha_1 = \alpha_2 = \alpha_3 = 1/3$, and 
\begin{align}
    \tQ^{(2)} = \frac{\log(N)}{N}\begin{bmatrix}
    4 & 3 & 7\\
    3 & 5 & 8 \\
    7 & 8 & 15
    \end{bmatrix},\quad \E\widetilde{\rmA} = \tQ^{(2)} \otimes (\ones_{N/3}\,\ones_{N/3}^{\sT}),
\end{align}
where $\otimes$ denotes the \emph{Kronecker} product. Obviously, \Cref{ass:expected_center_separation} is satisfied. However, $\lambda_3(\E\widetilde{\rmA}) = 0$ since $\tQ^{(2)}_{3:} = \tQ^{(2)}_{1:} + \tQ^{(2)}_{2:}$. As we can see, a mixture of assortative ($\tQ^{(2)}_{3:}$) and disassortative  ($\tQ^{(2)}_{1:}, \tQ^{(2)}_{2:}$) clusters may result in violating \Cref{ass:eigenvalueGap} but satisfying \Cref{ass:expected_center_separation}. We further note that this would never happen when all communities are purely assortative (diagonally dominant) or purely disassortative due to the symmetry requirement for probability tensors.

As far as we know, even in the graph case, \Cref{ass:eigenvalueGap} is crucial in most algorithms aiming to estimate for the number of communities. For example, \cite[Theorem 1]{Chen2014NetworkCF} applies a \emph{Network cross-validation} algorithm. This implicitly requires \Cref{ass:eigenvalueGap}, as it relies on Davis-Kahan \cite[Theorem 4]{Yu2015UsefulVO} to obtain an upper bound for the singular subspace perturbation, which is meaningful only if the $K$th largest singular value is at least of order $\log(N)$. 

As we just mentioned, \Cref{ass:expected_center_separation}, which \Cref{alg:agnostic_partition} relies on, is weaker than \Cref{ass:eigenvalueGap}. It would be interesting to close the gap between the those two assumptions, i.e., to find whether it is possible to give an accurate estimation the number of communities when the underlying generating model is a mixture of assortative and disassortative clusters.

\subsection{Proportional degrees pitfalls and possible fixes}
 \Cref{ass:prob_ratio_bound} is crucial for \Cref{alg:communityEstimation}. When this is dropped, and we are under the regime $\rho_{N} \gg \zeta_{N} \gtrsim \log(N)$, \Cref{alg:communityEstimation} tends to underestimate the number of communities due to the possibility that the upper bound $(\widetilde{d})^{3/4}$ is much larger than $\lambda_{K}(\rmA)$, leading to the problem that sparse clusters would be absorbed by denser clusters. One could strengthen \Cref{ass:eigenvalueGap} to satisfy $\lambda_K(\E \widetilde{\rmA}) \gtrsim \rho_{N}$ to make full use of \Cref{alg:communityEstimation}, though not so realistic in practice.

Instead, one could start with a lax overestimation of the number of communities, and then slightly modify \Cref{alg:spectral_initialization} to correct the overestimate. The intuition is that all but a vanishing fraction of the rows corresponding to vertices in $\gV_k$ would concentrate around their expected center, thus their corresponding vertices will be located in the ball $\ball_r(s_k)$, while $\widehat{\gV}^{(0)}_i  = o(N)$ for $i > K$ by \Cref{lem:centerFar}. Also, \Cref{alg:spectral_initialization} works agnostically without \Cref{ass:prob_ratio_bound} when $\rho_{N} \gtrsim \log(N)$ as proved in \Cref{sec:almostExact}. Thus, we propose \Cref{alg:communityEstimationIterative} as an informal alternative; here we pick some large enough $i$ as an initial guess and iteratively decrease $i$ until the partition has all classes proportional. However, we leave the criterion for deciding if a class is too small up to the reader, and we offer no guarantees. Note also that
\Cref{alg:communityEstimationIterative} is computationally inefficient: it has time complexity $O(N^2\log^2(N))$ since at least $\log(N)$ eigenpairs should be computed.

\begin{algorithm}[h]
\caption{\textbf{Number of Communities}}\label{alg:communityEstimationIterative}
\KwData{$\rmA_{\sJ}$, radius $r = [N\log(\rho_{N})]^{-1}\rho_{N}^2$ with $\rho_{N} = \widetilde{d} = \max_{v\in\gV}d_v$, $i = \log(N)$}
\While{$i \geq 2$ }{
{Compute the rank-$i$ approximation $\rmA_{\sJ}^{(i)} = \sum\limits_{j=1}^{i}\lambda_j \rvu_j^{\sT} \rvu_j$ of $\rmA_{\sJ}$.} 

{Construct balls of vertices centered at $s\in \setS$ by $\ball_{r}(s) = \{w\in \sJ:  \|(\rmA_{\sJ}^{(K)})_{s:} - (\rmA_{\sJ}^{(K)})_{w:}\|_2^2 \leq r \}$, where $\setS$ is randomly sampled from $\sJ$ with $|\setS| = \log^2(N)$.}

{Let $\ball_{r}(s_1)$ denote the ball with most vertices and take $\widehat{\gV}_{1}^{(0)} = \ball_{r}(s_1)$.}

\While{$2 \leq k \leq i$ }{
    {$s_{k} = \argmax_{s\in \setS}| \ball_{r}(s)\setminus ( \bigcup_{j=1}^{k-1}\widehat{\gV}_{j}^{(0)})|$}\;
    {$\widehat{\gV}_{k}^{(0)} = \ball_{r}(s_k)\setminus ( \bigcup_{j=1}^{k-1}\widehat{\gV}_{j}^{(0)})$}\;
}
\uIf{$|\widehat{\gV}^{(0)}_{i}| \ll N$ \tcp*{Informal condition, should be replaced with an explicit criterion.} }{
    $i \leftarrow i - 1$
    }
  \Else{
    {Quit the loop}\;
    }
}
\KwResult{$\widehat{K} = i$}
\end{algorithm} 


\section{Concluding remarks}\label{sec:conclusion_general}
In this paper, we characterize the phase transition of the exact recovery under non-uniform HSBM by formulating the sharp threshold and establishing necessary and sufficient conditions. Algorithms \ref{alg:agnostic_partition} and \ref{alg:partition_with_prior} achieve strong consistency within the stated assumptions. At the same time, algorithms provided are proved to be optimal, reaching the lowest possible mismatch ratio. Community detection on random hypergraphs remains a rich area for exploration. Besides finding an algorithm for correct estimation for the number of communities, we mention a few other interesting directions below.

First, \Cref{ass:expected_center_separation} is crucial for any spectral algorithm based on $\rmA$, since it is required for separation between expected centers. It will be interesting to see whether this condition can be removed, or at least weakened, such that the weak consistency could be achieved by some other non-spectral method. Perhaps in the vein of \cite{Ke2020CommunityDF, Han2022ExactCI, Agterberg2022EstimatingHO}, one could achieve weak consistency for uniform hypergraphs. Additionally, the generalization of those methods to handle the non-uniform hypergraphs when given the tensors $\{\tA^{(\ell)}\}$ with different orders remains open.

Second, it is known that exact recovery is  achievable \emph{at} the threshold in the case of graph SBM, i.e., when $\D_{\rm{H}}(a, b) = (\sqrt{a} - \sqrt{b})^2/2 = 1$ \cite{Mossel2016ConsistencyTF} and $\D_{\rm{CH}}(\rvalpha, \rmP) = 1$ \cite{Abbe2015CommunityDI}. However, it is still unclear whether exact recovery can be achieved when $\D_{\rm{GCH}} = 1$ for model \ref{def:non_uniform_HSBM}.

Last but not least, it will be interesting to see if model \ref{def:non_uniform_HSBM} can be generalized to have other degree distributions \cite{Gao2018CommunityDI, Hu2022MultiwaySC}, different edge labels \cite{Yun2016OptimalCR} and contextual feature vectors \cite{Chen2022GlobalAI, Ma2023CommunityDI}, which would be of particular interest for real-world applications.

%% file: appendix.tex

\section{The Kahn-Szemer\'{e}di approach and proofs of concentration}\label{app:proofConcentration}

\subsection{Proof outline of Theorem~\ref{thm:concentration}}\label{subsec:outline_concentration}
We follow the approach in \cite{Friedman1989OnTS, Lei2015ConsistencyOS, Cook2018SizeBC} to establish the proof of \Cref{thm:concentration}. According to Courant-Fischer-Weyl minimax principle, the spectral norm of any $N \times N$ Hermitian matrix $\rmW$ can be reformulated as
\begin{align}\label{eqn:operator_norm}
    \|\rmW\| \coloneqq\sup_{\rvx \in \mathbb{S}^{N - 1}} \|\rmW \rvx\|_2 = \sup_{\rvx \in \mathbb{S}^{N - 1}} |\rvx^{\sT}\rmW \rvx|\,,
\end{align}
where $\S^{N - 1} \coloneqq \{\rvx\in \R^{N}: \|\rvx\|_2 = 1\}$ denotes the unit sphere, and the second equality makes use of the Hermitian property of $\rmW$. As indicated by a short continuity argument (\Cref{lem:operator_norm_net}), if one wants to bound the supremum in \eqref{eqn:operator_norm}, it suffices to control $\rvx^{\sT}\rmW \rvx$ for all $\rvx$ in a suitable $\epsilon$-net $\sN_{\epsilon}$ of $\S^{N - 1}$ with cardinality bounded by $|\sN_{\epsilon}|\leq (1+ \frac{2}{\epsilon})^{N}$, while the existence of such $\sN_{\epsilon}$ is guranteed by compactness of $\S^{N - 1}$ and \cite[Corollary 4.2.13]{Vershynin2018HighDP}.

\begin{lemma}[{\cite[Lemma 4.4.1]{Vershynin2018HighDP}}]\label{lem:operator_norm_net}
Let $\rmW$ be any Hermitian $N \times N$ matrix and $\sN_{\epsilon}$ be an $\epsilon$-net on the unit sphere $\mathbb{S}^{N - 1}$ with $\epsilon\in (0,1)$, then 
$
\|\rmW\| \leq \frac{1}{1-\epsilon}\sup_{\rvx\in \sN_{\epsilon}} |\langle \rmW \rvx, \rvx\rangle |.
$
\end{lemma}
\begin{proof}[Proof of \Cref{lem:operator_norm_net}]
An $\epsilon$-net for a compact metric space $(\mathcal{X}, d)$ is a finite subset $\sN$ of $\mathcal{X}$ such that for each point $\rvx \in \mathcal{X}$, there is a point $\rvx \in \sN \subset \mathcal{X}$ with $d(\rvx, \rvx) \leq \epsilon$. Note that $\sN_{\epsilon}$ is an $\epsilon$-net on the unit sphere $\mathbb{S}^{N - 1}$, then for any $\rvx \in \mathbb{S}^{N - 1}$, there exists some $\rvx_0 \in \sN_{\epsilon}$ such that $\|\rvx - \rvx_0\|_2 \leq \epsilon$, and we have
\begin{equation*}
   \|\rmW \rvx\|_2 - \|\rmW \rvx_0\|_2 \leq \|\rmW \rvx - \rmW \rvx_0\|_2 \leq \|\rmW\|_2 \cdot \|\rvx - \rvx_0\|_2 \leq \epsilon \|\rmW\|\, ,
\end{equation*}
leading to $\|\rmW \rvx\|_2 - \epsilon \|\rmW\| \leq \|\rmW \rvx_0\|_2$. If the supremum is taken over $\rvx$ on the left hand side, then 
\begin{equation*}
   \sup_{\rvx \in \mathbb{S}^{N - 1}} ( \|\rmW \rvx\|_2 - \epsilon \|\rmW\|) = (1 - \epsilon) \|\rmW\| \leq \|\rmW \rvx_0\|_2 \leq \sup_{\rvy \in \sN_{\epsilon}} \|\rmW \rvy\|_2. 
\end{equation*}
Consequently, $ \|\rmW\| \leq \frac{1}{1 - \epsilon} \sup_{\rvx \in \sN_{\epsilon}} \|\rmW \rvx\|_2 = \frac{1}{1 - \epsilon} \sup_{\rvx\in \sN_{\epsilon}} |\langle \rmW \rvx, \rvx \rangle|$.
\end{proof}
Define $\rmW \coloneqq \rmA - \E \rmA$ and take $\sN\coloneqq \sN_{\frac{1}{2}}$ with its size bounded by $|\sN| \leq 5^{N}$, then
\begin{align}\label{eqn:epsilon_net}
    \|\rmA - \E \rmA\| \eqqcolon \|\rmW\| \leq 2\sup_{\rvx \in \sN }| \rvx^{\sT}\rmW\rvx| \leq 2\cdot 5^{N}\cdot  |\rvx^{\sT}\rmW\rvx|, \textnormal{ for some } \rvx\in \sN.
\end{align}
Therefore, upper bound of $\|\rmA - \E \rmA\|$ will be obtained once the concentration of random variable $\rvx^{\sT}\rmW\rvx$ is demonstrated for each fixed $\rvx \coloneqq [\ervx_1, \ervx_2, \cdots, \ervx_{N}]^{\sT}\in \sN$. Roughly speaking, one might seek upper bound on $|\rvx^{\sT}\rmW\rvx|$ of order $O(\sqrt{\etD_{\max}})$ which holds with high probability, then a union bound over $\sN$ establishes \Cref{thm:concentration}, as long as concentration of $\rvx^{\sT}\rmW\rvx$ is sufficient to beat the cardinality of the net $\sN$.

The problem now reduces to prove $|\rvx^{\sT}\rmW \rvx = \sum_{i}\sum_{j}\ermW_{ij}\ervx_i \ervx_j| = O(\sqrt{\etD_{\max}})$ with high probability. However, two main difficulties lie on
\begin{enumerate}
    \item the dependence between matrix elements $\ermW_{ij}$,
    \item not every pair $\ermW_{ij}\ervx_i \ervx_j$ enjoys the sufficient concentration to beat the cardinality of $\sN$.
\end{enumerate}
To get rid of those difficulties, we follow Kahn-Szemer\'{e}di's argument {\cite[Theorem 2.2]{Friedman1989OnTS}} to split $\rvx^{\sT} \rmW \rvx$ into \textit{light} and \textit{heavy} pairs according to $\max_{(i, j)\in [N]^2} |\ervx_i \ervx_j|$,
\begin{align}\label{eqn:lightheavypair}
    &\light(\rvx)\coloneqq \bigg\{ (i,j): |\ervx_i \ervx_j|\leq {\frac{\sqrt{\etD_{\max}}}{N}}\bigg\}, \quad \heavy(\rvx) \coloneqq \bigg\{ (i,j): |\ervx_i \ervx_j|> { \frac{\sqrt{\etD_{\max}} }{N} }\bigg\}\,,
\end{align}
where $\etD_{\max} \coloneqq \sum_{\ell \in \sL}\etD_{\max}^{(\ell)}$. By triangle inequality, \eqref{eqn:epsilon_net} adapts the following decomposition, where contribution from \textit{light} and \textit{heavy} couples can be bounded separately,
\begin{align}
    \|\rmA - \E \rmA\| \leq &\, 2\sup_{\rvx \in \sN }| \<\rmW\rvx, \rvx\>| = 2\sup_{\rvx \in \sN } \bigg| \sum_{(i,j)\in [N^2]} \ermW_{ij} \ervx_i \ervx_j \bigg|\notag \\
    \leq &\,2 \Bigg( \sup_{\rvx \in \sN }\bigg|\sum_{(i,j)\in \light(\rvx)} \ermW_{ij} \ervx_i \ervx_j\bigg| + \sup_{\rvx\in \sN }\bigg|\sum_{(i,j)\in \heavy(\rvx)} \ermW_{ij} \ervx_i \ervx_j \bigg| \Bigg)\label{eqn:contribution}\,.
\end{align}

As will be shown later, the contribution of \textit{light} pairs is bounded by $\alpha \sqrt{\etD_{\max}}$ with probability at least $1-2e^{-N}$ according to \eqref{eqn:lightbound}, while contribution of \textit{heavy} pairs is bounded by $\beta\sqrt{\etD_{\max}}$ with probability at least $1-2N^{-10}$ if $\etD_{\max} \geq c\log(N)$ for some constant $c>0$, as in \Cref{lem:contribution_of_heavy}. Therefore, 
\begin{align*}
    \|\rmA - \E \rmA \|\leq \const_{\eqref{eqn:concentrateA}}\sqrt{\etD_{\max}},
\end{align*}
holds with probability at least $1-2e^{-N}-2N^{-10}$, where $\const_{\eqref{eqn:concentrateA}}\coloneqq 2(\alpha + \beta)$ with $\alpha=5\LM^{2} $, $\beta= \beta_0 \sqrt{\beta_1} + \LM$,
 \begin{align*}
    &\, \beta_0 = 16+32 \LM (1+e^2)+1792(1+e^{-2})\LM^{2} ,\quad \frac{\beta_1}{\LM} \log \Big( \frac{\beta_1}{\LM} \Big) - \frac{\beta_1}{\LM} + 1 >\frac{11}{c},
\end{align*} 
which completes the proof of \Cref{thm:concentration}. The remaining of this section is organized as follows. The bounds for contributions of \textit{light} and \textit{heavy} pairs are provided in \Cref{subsec:light}, \Cref{subsec:heavy}, and the proof of \Cref{thm:regularization} will be established in \Cref{subsec:regularization}.

\subsection{Light couples}\label{subsec:light}
The contribution from \textit{light} pairs can be decomposed into summation of contributions from independent hyperedges. An application of Bernstein's inequality (\Cref{lem:Bernstein}) provides us with the concentration we need to beat the cardinality $|\sN| = 5^{N}$.

For each $\ell$-hyperedge $e\in \gE_{\ell}$, define $\etW^{(\ell)}_{e} \coloneqq \etA^{(\ell)}_e - \E\etA^{(\ell)}_e$. Then for each fixed $\rvx\in \mathbb{S}^{N - 1}$, the contribution from light pair can be reformulated as
\begin{align}
       \sum_{(i,j)\in \light(\rvx)} \ermW_{ij} \ervx_i \ervx_j =&\, \sum_{(i,j)\in \light(\rvx)} \bigg( \sum_{\ell \in \sL} \sum_{\substack{e\in \gE_{\ell}\\e \supset \{i,j\}} } \etW^{(\ell)}_{e} \bigg) \ervx_i \ervx_j = \sum_{\ell \in \sL} \sum_{e\in \gE_{\ell}} \etW^{(\ell)}_{e} \cdot \sum_{\substack{(i,j)\in \light(\rvx)\\ \{i,j\}\subset e} }\ervx_i \ervx_j\notag\\
       =&\, \sum_{\ell \in \sL} \sum_{e\in  \gE_{\ell}} \etY^{(\ell)}_{e},\label{eq:lightW}
\end{align}
where we denote
$
    \etY^{(\ell)}_{e}\coloneqq \etW^{(\ell)}_{e} \big( \sum_{\substack{(i,j)\in \light(\rvx)\\ \{i,j\}\subset e} }\ervx_i \ervx_j \big)
$.
Note that $e\supset (i, j)$ is an $\ell$-hyperedge, then the number of choices for the pair $(i, j)$ under some fixed $e$ is at most $\ell^{2}$. By the  definition of light pair in \eqref{eqn:lightheavypair}, $\etY^{(\ell)}_{e}$ is upper bounded by
$|\etY^{(\ell)}_{e}|\leq \ell^{2}\frac{\sqrt{\etD_{\max}} }{N} \leq  \LM^{2}\frac{\sqrt{\etD_{\max}} }{N}$ for all $\ell \in \sL$. Moreover, \eqref{eq:lightW} is a sum of independent, mean-zero ($\E \etY^{(\ell)}_{e} = 0$) random variables, and its second moment is bounded by
\begin{align*}
    &\,\sum_{\ell \in \sL} \sum_{e\in \gE_{\ell}} \E [(\etY^{(\ell)}_{e})^2] \coloneqq\sum_{\ell \in \sL} \sum_{e\in \gE_{\ell}} \Bigg[ \E[(\etW^{(\ell)}_{e})^2] \Bigg( \sum_{\substack{(i,j)\in \light(\rvx)\\ \{i,j\}\subset  e}}\ervx_i \ervx_j \Bigg)^2 \Bigg] \\
    \leq &\, \sum_{\ell \in \sL} \sum_{e\in \gE_{\ell}} \Bigg[ \E[ \etA^{(\ell)}_e ] \cdot \ell^{2} \Bigg(\sum_{\substack{(i,j)\in \sL(x)\\ \{i,j\}\subset e}} \ervx_i^2 \ervx_j^2 \Bigg) \Bigg] \leq \sum_{\ell \in \sL} \binom{N - 2}{\ell - 2} \cdot \frac{\etD^{(\ell)}_{\max}}{\binom{N - 1}{\ell - 1}} \cdot  \ell^{2}\sum_{(i,j)\in \sL(x)}\ervx_i^2 \ervx_j^2 \\
    \leq &\, \sum_{\ell \in \sL} \frac{\ell^{2}(\ell - 1)}{N - \ell + 2} \cdot \etD^{(\ell)}_{\max} \leq \frac{2}{N}\sum_{\ell \in \sL} \ell^{3} \cdot \etD^{(\ell)}_{\max}\leq \frac{2\LM^{3}}{N} \etD_{\max}\,,
\end{align*}
where in the last line, we use the facts $\etD_{\max} = \sum_{\ell \in \sL} \etD^{(\ell)}_{\max}$, $\rvx\in \mathbb{S}^{N - 1}$, $\rvx^{\sT} \rvx = 1$ and
\begin{align*}
    \sum_{(i,j)\in \sL(x)}\ervx_i^2 \ervx_j^2 \leq \sum_{(i,j)\in [N^2]}\ervx_i^2 \ervx_j^2 = \sum_{i=1}^{N} \ervx_i^2 \Big(\sum_{j=1}^{N} \ervx_j^2 \Big) = \sum_{i=1}^{N} \ervx_i^2 =1.
\end{align*}
Then Bernstein's inequality (\Cref{lem:Bernstein}) implies that for any $\alpha>0$,
\begin{subequations}
\begin{align*}
    &\,\P \bigg( \bigg| \sum_{(i,j)\in\light(\rvx)} \ermW_{ij} \ervx_i \ervx_j\bigg| \geq \alpha \sqrt{\etD_{\max} }\bigg) = \P \bigg( \bigg| \sum_{\ell \in \sL} \sum_{e\in  \gE_{\ell}} \etY^{(\ell)}_{e} \bigg| \geq \alpha\sqrt{\etD_{\max} } \bigg) \\
    \leq &\, 2\exp \bigg( -\frac{\frac{1}{2}\alpha^2 \etD_{\max} } {\frac{2\LM^{3}}{N} \etD_{\max} + \frac{1}{3} \LM^{2} \frac{\sqrt{\etD_{\max} }}{N} \cdot \alpha\sqrt{\etD_{\max}} } \bigg)
    \leq 2\exp \bigg( -\frac{\alpha^2 N}{4\LM^{3}+\frac{2\alpha}{3}\LM^{2} } \bigg).
\end{align*}
\end{subequations}
Consequently by a union bound,
    \begin{align}
        &\, \P \bigg(  \sup_{\rvx\in \sN } \bigg| \sum_{(i,j)\in\light(\rvx)} \ermW_{ij} \ervx_i \ervx_j \bigg|\geq \alpha\sqrt{\etD_{\max} }\Bigg) \notag \\ 
        \leq &\,|\sN| \cdot \P \bigg( \bigg| \sum_{(i,j)\in\light(\rvx)} \ermW_{ij} \ervx_i \ervx_j\bigg| \geq \alpha\sqrt{\etD_{\max}}\bigg) \notag\\
        \leq &\, 2\exp \bigg( \log(5) \cdot N -\frac{\alpha^2 N}{4\LM^{3}+\frac{2\alpha}{3}\LM^{2} }\bigg)\leq 2e^{-N}\, ,\label{eqn:lightbound}
    \end{align}
where the last inequality holds by taking $\alpha=5\LM^{2}$. 

\subsection{Heavy couples}\label{subsec:heavy}
\subsubsection{General Strategy}
The contribution of \textit{heavy} couples can be bounded separately,
\begin{align}
    \bigg|\sum_{(i,j)\in \heavy(\rvx)} \ermW_{ij} \ervx_i \ervx_j \bigg|  \leq \bigg|\sum_{(i,j)\in \heavy(\rvx)} \ermA_{ij} \ervx_i \ervx_j \bigg| + \bigg|\sum_{(i,j)\in \heavy(\rvx)} (\E \ermA)_{ij}\, \ervx_i \ervx_j \bigg|.
\end{align}
Meanwhile, for any pair $(i,j)\in [N]^2$,
\begin{align}\label{eqn:EA_upperbound}
    (\E \ermA)_{ij}\leq \sum_{\ell \in \sL} \binom{N - 2}{\ell - 2}\frac{\etD_{\max}^{(\ell)}}{\binom{N - 1}{\ell - 1}} \leq  \sum_{\ell \in \sL}\frac{(\ell - 1)}{N - 1} \cdot \etD_{\max}^{(\ell)}= \frac{\LM  \etD_{\max}}{N - 1},
\end{align}
and by definition of \textit{heavy} pair
\begin{align}\label{eqn:heavyexp}
   \bigg|\sum_{(i,j)\in \heavy(\rvx)} (\E \ermA)_{ij}\, \ervx_i \ervx_j \bigg| = &\, \bigg|\sum_{(i,j)\in \heavy(\rvx)} (\E \ermA)_{ij} \frac{ \ervx_i^2 \ervx_j^2}{\ervx_i \ervx_j}\bigg| \leq\sum_{(i,j)\in \heavy(\rvx)} \frac{\LM  \etD_{\max}}{N - 1} \frac{\ervx_i^2 \ervx_j^2}{|\ervx_i \ervx_j|} \notag\\
   \leq &\, \LM \sqrt{\etD_{\max}}\sum_{(i,j)\in \heavy(\rvx)} \ervx_i^2 \ervx_j^2\leq \LM\sqrt{\etD_{\max}},
\end{align}
where we use the fact $\sum_{(i,j)\in \heavy(\rvx)} \ervx_i^2 \ervx_j^2 \leq 1$ again. Therefore it suffices to show that, with high probability, 
\begin{align}\label{eqn:heavyA}
    \bigg| \sum_{(i,j)\in \sH(x)} \ermA_{ij} \ervx_i \ervx_j  \bigg| =  O \big(\sqrt{\etD_{\max}} \big).
\end{align}
Unfortunately, the contribution of \textit{heavy} couples does not enjoy sufficient concentration to beat the cardinality of the net $\sN$, which was the strategy we used in \Cref{subsec:light}. Instead, the key idea here is to prove the \textit{discrepancy property} (\Cref{def:Discrepancy}) holds with high probability for the associated random regular graph. Essentially, the edge counts
$
    \gE_{\rmA}(\setS, \setT) \coloneqq \sum_{u \in \setS} \sum_{v \in \setT} = \ones_{\setS}^{\sT} \rmA \ones_{\setT}
$
are not much larger than their expectation, uniformly over choices of $\setS, \setT\subset [N]$, which can be proved by using tail estimates of random variables $\gE_{\rmA}(\setS, \setT)$. Conditioning on the event that \textit{discrepancy} property holds, one can show that the contribution of \textit{heavy} couples is deterministically of order $O(\sqrt{\etD_{\max}})$ by \Cref{lem:heavybound}, as long as row and column sums of $\rmA$ are bounded by $\sqrt{\etD_{\max}}$ up to some constant.

\begin{definition}[Discrepancy property, \bf{DP}]\label{def:Discrepancy}
Let $\rmA$ be an $N \times N$ matrix with non-negative entries. For subsets $\setS,\setT\subset [N]$, define
\begin{align*}
    \gE_{\rmA}(\setS,\setT) \coloneqq \sum_{i\in \setS}\sum_{j\in \setT} \ermA_{ij}.
\end{align*}
We say that $\rmA$ has the discrepancy property with parameter $\delta>0$, $\kappa_1>1, \kappa_2\geq 0$, denoted by {\bf{DP}}$(\delta,\kappa_1,\kappa_2)$, if for all non-empty $\setS,\setT\subset [N]$, at least one of the following hold:
\begin{enumerate}
    \item $\gE_{\rmA}(\setS,\setT)\leq \kappa_1 \delta |\setS| |\setT|$;
    \item  $\gE_{\rmA}(\setS,\setT) \cdot \log \big(\frac{\gE_{\rmA}(\setS,\setT)}{\delta |\setS|\cdot |\setT|}\big)\leq \kappa_2 (|\setS|\vee |\setT|)\cdot \log \big(\frac{eN}{|\setS|\vee |\setT|}\big)$.
\end{enumerate}
\end{definition}

\begin{lemma}[{\cite[Lemma 6.6]{Cook2018SizeBC}}, {\bf{DP}} $\implies$ heavy couples are small]\label{lem:heavybound}
Let $\rmA$ be a non-negative symmetric $N \times N$ matrix. Suppose that
\begin{enumerate}
    \item all row and column sums of $\rmA$ are bounded by $d$,
    \item $\rmA$ has $\textnormal{\bf{DP}}(\delta,\kappa_1,\kappa_2)$ with $\delta=\frac{Cd}{N}$ for some $C>0,\kappa_1>1,\kappa_2\geq 0$,
\end{enumerate}
then for any $x\in \mathbb{S}^{N - 1}$, deterministically,
\begin{align*}
     \bigg| \sum_{(i,j)\in \sH(x)} \ermA_{ij} \ervx_i \ervx_j  \bigg| \leq \beta_0\sqrt{d},
\end{align*}
where
$
    \beta_0 =16+32C(1+\kappa_1)+64\kappa_2  (1+\frac{2}{\kappa_1\log\kappa_1})
$.
\end{lemma}

\subsubsection{Our model} To prove the \textit{discrepancy} property in our model, we follow the discrepancy analysis in \cite{Feige2005SpectralTA, Cook2018SizeBC}, and consider the weighted graph associated with the adjacency matrix $\rmA$. By \Cref{lem:DPforAdjacency}, the \textit{discrepancy} property follows with high probability if \textit{uniform upper tail} property (\Cref{Def:UUTP}) holds. Therefore, our tasks reduce to prove that row sums of $\rmA$ are bounded and {\bf{UUTP}} holds with high probability, which are shown in \Cref{lem:maxdeg} and \Cref{lem:UUTPforA} respectively.

\begin{definition}[Uniform upper tail property, \textbf{UUTP}]\label{Def:UUTP}
Let $\rmA$ be a random symmetric $N \times N$ matrix with non-negative entries. An $N \times N$ matrix $\rmQ$ is associated to $\rmA$ through $f_{\rmQ}(\rmA) \coloneqq \sum_{i,j=1}^{N} \ermQ_{ij} \ermA_{ij}$. Define
\begin{align*}
    \mu \coloneqq f_{\rmQ}(\E \rmA) = \sum_{i,j=1}^{N} \ermQ_{ij}\E \ermA_{ij}, \quad \widetilde{\sigma}^2 \coloneqq f_{\rmQ \circ \rmQ }(\E \rmA) = \sum_{i,j=1}^{N} \ermQ_{ij}^2 \E \ermA_{ij},
\end{align*}
where $\circ$ denotes the Hadamard (entrywise) matrix product. We say that $\rmA$ satisfies the uniform upper tail property {\bf{UUTP}}$(c_0,\gamma_0)$ with $c_0>0,\gamma_0\geq 0$, if for any $a, t>0$ and any symmetric $N \times N$ matrix $\rmQ$ satisfying $\ermQ_{ij}\in [0,a]$ for all $i, j\in [N]$, we have
\begin{align*}
    \P \Big(  f_{\rmQ}(\rmA)\geq (1 + \gamma_0)\mu +t \Big)\leq \exp \bigg( -c_0 \frac{\widetilde{\sigma}^2}{a^2} h\Big( \frac{at}{\widetilde{\sigma}^2}\Big)\bigg).
\end{align*}
where the function $h(x):=(1 + x)\log(1+x) - x$ for all $x>-1$.
\end{definition}

\begin{lemma}[{\cite[Lemma 6.4]{Cook2018SizeBC}}, {\bf{UUTP}} $\implies$ {\bf{DP}}]\label{lem:DPforAdjacency}
Let $\rmA$ be a symmetric $N \times N$ random matrix with non-negative entries. Assume that 
\begin{enumerate}
    \item $\E \ermA_{ij}\leq \delta$ for all $i,j\in [N]$ for some $\delta>0$,
    \item $\rmA$ has {\bf{UUTP}}$(c_0,\gamma_0)$ with parameter $c_0,\gamma_0 >0$.
\end{enumerate}
Then for any $\theta>0$, the \textit{discrepancy} property {\bf{DP}}$(\delta,\kappa_1,\kappa_2)$ holds for $\rmA$ with probability at least $1 - N^{-\theta}$ with 
$
    \kappa_1 = e^2(1+\gamma_0)^2,  \kappa_2 = \frac{2}{c_0}(1+\gamma_0)(\theta + 4).
$
\end{lemma}

\begin{lemma}\label{lem:maxdeg}
Denote $\etD_{\max} \coloneqq c q_{N}$ for some constant $c>0$, where $q_{N}$ is the magnitude parameter in \eqref{eqn:edge_density}. Then for some large enough constant $\beta_1 > 0$ \footnote{The existence of such $\beta_1$ is guaranteed since $g(x) = x\log(x) - x + 1$ takes its minimum at $x=1$ and increases when $x>1$.}, satisfying $\LM^{-1}\beta_1 \log(  \LM^{-1}\beta_1) -  \LM^{-1}\beta_1 + 1 > 11/c$, there exists some $\widetilde{c} > 0$ such that
\begin{align}\label{eqn:bounded_each_degree}
    \P \big( \rD_v > \beta_1 \etD_{\max}\big) \leq e^{-\widetilde{c}\rho_{N}}.
\end{align}
At the same time, the following holds with probability at least $1 - Ne^{-11q_{N}}$,
\begin{align}\label{eqn:bounded_max_degree}
    \max_{1\leq i\leq n} \big( \sum_{j=1}^{N} \ermA_{ij} \big)\leq \beta_1 \etD_{\max}.
\end{align}
In particular, when $q_{N} = \log(N)$, \eqref{eqn:bounded_max_degree} holds with probability at least $1 - N^{-10}$.
\end{lemma}
\begin{proof}[Proof of \Cref{lem:maxdeg}]
Recall $\rD_v$ in \eqref{eqn:Dv}, and for each fixed $v\in[N]$, we have
\begin{align}\label{eqn:rowsum_edge}
    \sum_{j=1}^{N} \ermA_{vj} = \sum_{j=1}^{N} \sum_{\ell \in \sL} \, \sum_{\substack{e\in \gE_{\ell}\\ e\supset\{v, j\}} } \etA^{(\ell)}_e = \sum_{\ell \in \sL}\sum_{ \substack{e\in \gE_{\ell}\\ e\ni v }} (\ell - 1)\etA^{(\ell)}_e \asymp \sum_{\ell \in \sL} \sum_{ e\ni v} \etA^{(\ell)}_e = \rD_v,
\end{align}
since each $\ell \in \sL$ is constant. From now on, we work on $\sum_{j=1}^{N} \ermA_{vj}$ for convenience. By Markov inequality (\Cref{lem:Markov}) and independence between hyperedges, \eqref{eqn:bounded_each_degree} can be proved as follows
\begin{align*}
    &\, \P \bigg( \sum_{j=1}^{N} \ermA_{vj} > \beta_1 \etD_{\max}\bigg) = \P \bigg( \sum_{\ell \in \sL} \sum_{ \substack{e\in \gE_{\ell}\\ e\ni v}} (\ell - 1)\etA^{(\ell)}_e > \beta_1 \etD_{\max}\bigg)\\
    \leq &\, \inf_{\theta >0} \bigg( \exp(-\theta \beta_1 \etD_{\max} ) \prod_{\ell \in \sL} \prod_{ \substack{e\in \gE_{\ell}\\ e\ni i}} \E \exp[\theta (\ell - 1) \etA^{(\ell)}_e] \bigg)\\
    \leq &\, \inf_{\theta >0} \bigg( \exp(-\theta \beta_1 \etD_{\max} ) \prod_{\ell \in \sL} \prod_{ \substack{e\in \gE_{\ell}\\ e\ni i}} \bigg[ \frac{\etD_{\max}^{(\ell)} }{ \binom{N - 1}{\ell - 1} } e^{\theta(\ell - 1)} + 1 - \frac{\etD_{\max}^{(\ell)} }{ \binom{N - 1}{\ell - 1} } \bigg ] \bigg)\\
    \leq &\, \inf_{\theta >0} \bigg( \exp(-\theta \beta_1 \etD_{\max} ) \prod_{\ell \in \sL} \prod_{ \substack{e\in \gE_{\ell}\\ e\ni i}} \exp \bigg[ \frac{\etD_{\max}^{(\ell)} }{ \binom{N - 1}{\ell - 1} } ( e^{\theta(\ell - 1)} - 1) \bigg] \bigg), \quad (1 + x \leq e^{x}) \\
    \leq &\, \inf_{\theta >0} \bigg( \exp(-\theta \beta_1 \etD_{\max} ) \cdot  \exp \bigg[ \sum_{\ell \in \sL} \etD_{\max}^{(\ell)} ( e^{\theta(\ell - 1)} - 1) \bigg] \bigg), \, (\textnormal{at most } \binom{N - 1}{\ell - 1} \textnormal{ edges contains } i)\\
    \leq &\, \inf_{\theta > 0} \exp[ -\etD_{\max}(\theta \beta_1 - e^{\theta \LM} + 1) ] = \exp \bigg[ -\etD_{\max} \Big( \frac{\beta_1}{\LM} \log \Big( \frac{\beta_1}{\LM} \Big) - \frac{\beta_1}{\LM} + 1 \Big)\bigg] \\
    \leq&\, \exp \bigg[ - q_{N} \cdot c \Big( \frac{\beta_1}{\LM} \log \Big( \frac{\beta_1}{\LM} \Big) - \frac{\beta_1}{\LM} + 1 \Big)\bigg] \leq e^{-11q_{N}} = e^{- \widetilde{c}\rho_{N}},
\end{align*}
where the last line holds since $\etD_{\max} \asymp \rho_{N}$ as proved in \Cref{lem:maxProb}. The proof of \eqref{eqn:bounded_max_degree} is completed by a union bound using \eqref{eqn:bounded_each_degree}
\begin{align*}
    \P \Big(\max_{1\leq i\leq n} \Big(\sum_{j=1}^{N} \ermA_{ij} \Big) > \beta_1 \etD_{\max} \Big) \leq \sum_{i=1}^{N} \, \P \Big( \sum_{j=1}^{N} \ermA_{ij} > \beta_1 \etD_{\max}\Big) \leq Ne^{-11q_{N}}.
\end{align*}
The last part follows directly by plugging $q_{N} = \log(N)$.
\end{proof}

\begin{lemma}[UUTP from Non-uniform Inhomogeneous Hypergraph]\label{lem:UUTPforA}
Let $\rmA$ be the adjacency matrix of $\gH = \cup_{\ell \in \sL}\gH_m$, then $\rmA$ satisfies {\bf{UUTP}}$(c_0,\gamma_0)$ with $c_0 = \LM^{-2}, \gamma_0 = 0$.
\end{lemma}
\begin{proof}[Proof of \Cref{lem:UUTPforA}]
Note that
\begin{align*}
    &\,f_{\rmQ}(\rmA) - \mu = \sum_{i,j=1}^{N} \ermQ_{ij}(\ermA_{ij} - \E \ermA_{ij}) = \sum_{i,j=1}^{N} \ermQ_{ij} \ermW_{ij}= \sum_{i,j=1}^{N} \ermQ_{ij} \bigg( \sum_{\ell \in \sL}\,\, \sum_{\substack{e\in \gE_{\ell}\\ e \supset \{i,j\}} } \etW^{(\ell)}_{e}  \bigg) \\
    =&\,\sum_{\ell \in \sL} \,\, \sum_{e\in \gE_{\ell}} \etW^{(\ell)}_{e} \bigg( \sum_{\{i,j\}\subset e} \ermQ_{ij}\bigg) =\sum_{\ell \in \sL} \,\,\sum_{e\in \gE_{\ell}}\etZ^{(\ell)}_{e}\,,
\end{align*}
where $\etZ^{(\ell)}_{e} = \etW^{(\ell)}_{e} \big( \sum_{\{i,j\} \subset e} \ermQ_{ij}\big)$ are independent centered random variables upper bounded by
$
    |\etZ^{(\ell)}_{e}| \leq \sum_{ \{i,j\} \subset e} \ermQ_{ij}\leq \LM^{2}a
$
for each $ \ell \in \sL$ since $\ermQ_{ij} \in [0, a]$. Moreover, the variance can be written as
\begin{align*}
    &\, \sum_{\ell \in \sL} \sum_{e\in \gE_{\ell}} \E(\etZ^{(\ell)}_{e})^2 \coloneqq \sum_{\ell \in \sL} \,\, \sum_{e\in \gE_{\ell}} \E(\etW^{(\ell)}_{e})^2 \bigg( \sum_{\{i,j\}\subset e} \ermQ_{ij}\bigg)^2\\
    \leq &\,\sum_{\ell \in \sL} \,\,\sum_{e\in \gE_{\ell}}\E [\etA^{(\ell)}_{e}] \cdot \ell^{2} \sum_{ \{i,j \}\subset e } \ermQ_{ij}^2
    \leq \LM^{2} \sum_{i,j=1}^{N} \ermQ_{ij}^2 \E \ermA_{ij} = \LM^{2} \widetilde{\sigma}^2.
\end{align*}
where the last inequality holds since by definition $\E \ermA_{ij} = \sum_{\ell \in \sL} \,\,\sum_{\substack{e\in \gE_{\ell}\\ \{i,j\}\subset e}}\E [\etA^{(\ell)}_{e}]$. Then by Bennett's inequality (\Cref{lem:Bennett}), we obtain
\begin{align*}
    \P (f_{\rmQ}(\rmA)- \mu\geq t) \leq \exp\Big( -\frac{\widetilde{\sigma}^2}{\LM^{2} a^2} h \Big( \frac{at}{\widetilde{\sigma}^2}\Big) \Big),
\end{align*}
where the last line holds since $xh(\frac{1}{x}) = (1+x)\log(1 + \frac{1}{x}) - 1$ decreases when $x > 0$.
\end{proof}

\subsubsection{Contribution of Equation~(\ref{eqn:heavyA})}
\begin{lemma}\label{lem:contribution_of_heavy}
If $\etD_{\max} \geq c\log(N)$ for some constant $c>0$, there exists some constant $\beta$ such that with probability at least $1-2N^{-10}$,
 \begin{align}\label{eqn:betasqrtd}
  \Bigg|\sum_{(i,j)\in \mathcal H(x)}  (\ermA_{ij} - \E \ermA_{ij}) \ervx_i \ervx_j\Bigg|\leq \beta \sqrt{\etD_{\max}}.   
 \end{align}
\end{lemma}
\begin{proof}[Proof of \Cref{lem:contribution_of_heavy}]
Note that $\E \ermA_{ij}$ is bounded (\eqref{eqn:EA_upperbound}) and $\rmA$ satisfies $\textnormal{\bf{UUTP}}(\LM^{-2}, 0)$ (\Cref{lem:UUTPforA}), then by \Cref{lem:DPforAdjacency}, for any $\theta>0$, the \textit{discrepancy} property $\textnormal{\bf{DP}}(\delta, \kappa_1,\kappa_2)$ holds for $\rmA$ with probability at least $1-N^{-\theta}$, where
\begin{align*}
    \delta = \frac{\LM}{N} \etD_{\max}, \quad \kappa_1 = e^2(1+\gamma_0)^2 = e^2,  \quad \kappa_2 = \frac{2}{c_0}(1+\gamma_0)(\theta + 4) = 2\LM^{2}(\theta + 4).
\end{align*}
Let $\mathcal{E}_1$ be the event that all row sums of $\rmA$ are bounded by $\beta_1 \etD_{\max}$, and $\mathcal{E}_2$ be the event that $\textnormal{\textbf{DP}}(\delta, \kappa_1,\kappa_2)$ holds for $\rmA$, then $\P (\mathcal{E}_1 \cap \mathcal{E}_2)\geq 1- 2N^{-10}$ where we take $\theta = 10$ for convenience. On the event $\mathcal{E}_1 \cap \mathcal{E}_2$, by \Cref{lem:heavybound}, \eqref{eqn:betasqrtd} holds with $\beta= \beta_0 \sqrt{\beta_1} + \LM$, where $\beta_0, \beta_1$ satisfying
 \begin{align*}
    &\, \beta_0= 16 + 32\LM (1 + \kappa_1) + 64\kappa_2(1 + \frac{2}{\kappa_1 \log(\kappa_1)}) =16+32\LM (1+e^2)+1792(1+e^{-2})\LM^{2},\\
    &\, \frac{\beta_1}{\LM} \log \Big( \frac{\beta_1}{\LM} \Big) - \frac{\beta_1}{\LM} + 1 > \frac{11}{c}.
\end{align*} 

\end{proof}

\subsection{Proof of Theorem~\ref{thm:regularization}} \label{subsec:regularization}
Following the strategy in \Cref{subsec:outline_concentration}, for any $\sJ\subset [N]$,
\begin{align}\label{eqn:contribution_regularized}
    &\,\|(\rmA - \E \rmA)_{\sJ}\| \\
   \leq  2 \Bigg( \sup_{\rvx \in \sN }\bigg| &\, \sum_{(i,j)\in \light(\rvx)} {\big[ (\rmA - \E \rmA)_{\sJ}\big] }_{ij} \ervx_i \ervx_j\bigg| + \sup_{\rvx\in \sN }\bigg|\sum_{(i,j)\in \heavy(\rvx)} {\big[ (\rmA - \E \rmA)_{\sJ}\big] }_{ij} \ervx_i \ervx_j \bigg| \Bigg)\,,\notag
\end{align}
we want to bound the contribution of \textit{light} and \textit{heavy} pairs separately.

\subsubsection{Light couples}
By \eqref{eqn:lightbound}, with probability at least $1-2e^{-N}$,
\begin{align}\label{eq:J_light}
    \sup_{\rvx\in \sN }  \bigg| \sum_{(i,j)\in \light(\rvx)}  {\big[ (\rmA - \E \rmA)_{\sJ}\big] }_{ij}\, \ervx_i \ervx_j \bigg|\leq 5\LM^{2}\sqrt{\etD_{\max}}\,.
\end{align}
A union bound proves that \eqref{eq:J_light} holds for any given subset $\sJ\subset [N]$ with probability at least $1-2(e/2)^{-N}$, since the number of possible choices for $\sJ$ is at most $2^{N}$, including the options $\sJ \coloneqq\sJ_1, \sJ_2$ of our interests.

\subsubsection{Heavy couples}
Our goal is to prove
\begin{align}\label{eq:betasqrtd_regularized}
    \sup_{\rvx\in \sN }  \bigg| \sum_{(i,j)\in \heavy(\rvx)}  {\big[ (\rmA - \E \rmA)_{\sJ}\big] }_{ij}\, \ervx_i \ervx_j \bigg|\leq \beta\sqrt{\etD_{\max}}\,,
\end{align}
Similar to \eqref{eqn:heavyexp}, deterministically for $\sJ\coloneqq \sJ_1, \sJ_2$. 
\begin{align}\label{eq:reg2}
     \bigg|\sum_{(i,j)\in \heavy(\rvx)} {\big[ (\E \rmA)_{\sJ}\big]}_{ij} \bigg|\leq \LM\sqrt{\etD_{\max}}\,.
\end{align}
 To control the contribution of \textit{heavy} pairs, it suffices to show that \textit{discrepancy} holds for $\rmA_{\sJ}$ and all row sums of $\rmA_{\sJ}$ are bounded by $\beta_1 \etD_{\max}$,  according to \Cref{lem:heavybound}.
 
 To verify the \textit{discrepancy} property, it suffices to show that {\bf{UUTP}}$(\LM^{-2}, 0)$ holds for $\rmA_{\sJ}$, which is true since $\rmA_{\sJ}$ is obtained from $\rmA$ by restricting rows and columns to $\sJ$. 
 Then by \Cref{lem:DPforAdjacency}, for any $\theta>0$, the \textit{discrepancy} property $\textnormal{\bf{DP}}(\delta, \kappa_1,\kappa_2)$ holds for $\rmA_{\sJ}$ with probability at least $1-N^{-\theta}$, where
\begin{align*}
    \delta = \frac{\LM}{N} \etD_{\max}, \quad \kappa_1 = e^2(1+\gamma_0)^2 = e^2,  \quad \kappa_2 = \frac{2}{c_0}(1+\gamma_0)(\theta + 4) = 2\LM^{2}(\theta + 4).
\end{align*}
The remaining is to verify the boundedness of row sums of $\rmA_{\sJ}$.
\begin{itemize}
    \item For $\sJ_1$, \Cref{lem:boundednessJ1} proves that row sums of each $v\in \sJ_1$ is bounded by $\beta_1 \etD_{\max}$ with probability at least $1 - N^{-c}$ for some constant $c>0$ under the regime $\rho_{N} \gtrsim \log(N)$.
    \item Deterministically, it is true for $v\in \sJ_2 \coloneqq \{ i\in[N], \sum_{j=1}^{N} \ermA_{ij} \leq \beta_1 \etD_{\max}\}$ by taking $\beta_1 = \LM$.
\end{itemize}
Let $\mathcal{E}_1$ be the event that all row sums of $\rmA_{\sJ}$ are bounded by $\beta_1 \etD_{\max}$, while $\mathcal{E}_2$ be the event that $\textnormal{\textbf{DP}}(\delta, \kappa_1,\kappa_2)$ holds for $\rmA_{\sJ}$, then
\begin{itemize}
    \item For $\sJ_1$, $\P (\mathcal{E}_1 \cap \mathcal{E}_2)\geq 1- 2N^{-c}$, where we take $\theta = c$ for convenience.
    \item For $\sJ_2$, $\P (\mathcal{E}_1 \cap \mathcal{E}_2)\geq 1- N^{-10}$, where we take $\theta = 10$ for convenience.
\end{itemize}
On the event $\mathcal{E}_1 \cap \mathcal{E}_2$, \eqref{eq:betasqrtd_regularized} holds with $\beta= \beta_0 \sqrt{\beta_1} + \LM$, $\beta_1 = \LM$,
 \begin{align*}
    \beta_0= 16 + 32\LM (1 + \kappa_1) + 64\kappa_2(1 + \frac{2}{\kappa_1 \log(\kappa_1)}) =16+32\LM (1+e^2)+1792(1+e^{-2})\LM^{2}.
\end{align*}
\begin{lemma}\label{lem:boundednessJ1}
If $\rho_{N} \gtrsim \log(N)$, then there exists some constant $c>0$ such that row sums of each $v\in \sJ_1$ is bounded by $\beta_1 \etD_{\max}$ with probability at least $1 - N^{-c}$.
\end{lemma}
\begin{proof}[Proof of \Cref{lem:boundednessJ1}]
Recall $\rD_v$ in \eqref{eqn:Dv}, which can be further written as
\begin{align}
    \rD_v = \sum_{\ell \in \sL} \sum_{ e\ni v} \etA^{(\ell)}_e \asymp \sum_{\ell \in \sL} \sum_{ e\ni v} (\ell - 1)\etA^{(\ell)}_e = \sum_{j=1}^{N} \ermA_{vj},
\end{align}
since $\ell \geq 2$ is some constant for each $\ell \in \sL$. From now on, we work on $\sum_{j=1}^{N} \ermA_{vj}$ for convenience. In the following, we calculate the probability that the maximum vertex degree in $\sJ_1$ is greater than $\beta_1 \etD_{\max}$. By the definition of $\sJ_1$,
    \begin{align*}
        &\, \P \Big( \max_{v\in \sJ_1} \Big( \sum_{j=1}^{N} \ermA_{vj} \Big) > \beta_1 \etD_{\max} \Big) = \P \Big( \# \Big\{ v: \sum_{j=1}^{N} \ermA_{vj} > \beta_1 \etD_{\max} \Big\} >\lfloor N \exp(-\overline{d}) \rfloor \Big)\\
        \overset{\textnormal{(Markov)}}{\leq} &\, \frac{ 1 }{\lfloor N \exp(-\overline{d}) \rfloor + 1} \E \Big( \# \Big\{ v: \sum_{j=1}^{N} \ermA_{vj} > \beta_1 \etD_{\max}\Big\} \Big) = \frac{\sum_{v\in[N]} \E \mathbbm{1}\Big\{ v: \sum_{j=1}^{N} \ermA_{vj} > \beta_1 \etD_{\max} \Big\} }{\lfloor N \exp(-\overline{d}) \rfloor + 1}\\
            \leq &\, \frac{N}{\lfloor N \exp(- \overline{d}) \rfloor + 1} \, \max_{v\in [N]}\, \P \Big( \sum_{j=1}^{N} \ermA_{vj} > \beta_1 \etD_{\max} \Big) \leq \frac{1}{\exp(- \overline{d}) } \, \max_{v\in [N]}\,  \P \Big( \sum_{j=1}^{N} \ermA_{vj} > \beta_1 \etD_{\max} \Big).
    \end{align*}
Note that $\etD_{\max} \asymp \rho_{N}$ as proved in \Cref{lem:maxProb}, then $\rho_{N} \gtrsim \log(N)$ is equivalent to say $\etD_{\max} \geq c\log(N)$ for some constant $c>0$.
    \begin{itemize}
        \item For the numerator, by going through the same procedure in \Cref{lem:maxdeg} and a union bound, one can prove that
        \begin{align*}
            &\, \max_{v\in [N]}\,  \P \Big( \sum_{j=1}^{N} \ermA_{vj} > \beta_1 \etD_{\max} \Big) \leq \sum_{v\in [N]}\,  \P \Big( \sum_{j=1}^{N} \ermA_{vj} > \beta_1 \etD_{\max} \Big) \\
            \leq &\, \exp \bigg[ -\log(N) \cdot c \Big( \frac{\beta_1}{\LM} \log \Big( \frac{\beta_1}{\LM} \Big) - \frac{\beta_1}{\LM} + 1 \Big) + \log(N)\bigg]. 
        \end{align*}   

    \item For the denominator, $\overline{d}\leq \xi_{N} + \xi_{N}/\sqrt{N}$ with probability tending to $1$ according to \Cref{lem:maxProb}.
    
\end{itemize}
Consequently, when $\rho_{N} \gtrsim \log(N)$, by taking $\beta_1$ large enough satisfying $\frac{\beta_1}{\LM} \log \Big( \frac{\beta_1}{\LM} \Big) - \frac{\beta_1}{\LM} + 1 - \frac{1}{c} - \frac{\xi_{N}}{c\log(N)}\geq 1$, we then have
    \begin{align*}
        &\,\P \Big( \max_{v\in \sJ_1} \Big( \sum_{j=1}^{N} \ermA_{vj} \Big) > \beta_1 \etD_{\max} \Big) \leq \exp(\overline{d}) \cdot \,\max_{v\in [N]} \, \P \Big( \sum_{j=1}^{N} \ermA_{vj} > \beta_1 \etD_{\max} \Big) \\
        \leq&\, \exp \bigg[ -\log(N) \cdot c \Big( \frac{\beta_1}{\LM} \log \Big( \frac{\beta_1}{\LM} \Big) - \frac{\beta_1}{\LM} + 1 - \frac{1}{c} - \frac{\xi_{N}}{c\log(N)}\Big)\bigg] \leq N^{-c},
    \end{align*}
    which completes the proof. The existence of such $\beta_1$ is guaranteed since the function $f(x) = x\log(x) - x$ takes its minimum at $x = 1$ with minimum $f(1) = -1$ and $c^{-1} = O(1)$, $\xi_{N}/[c\log(N)] \lesssim \xi_{N}/\rho_{N} = O(1)$. The requirement $\rho_{N} \gtrsim \log(N)$ is crucial since one needs to beat the union bound as in the proof of \eqref{eqn:bounded_max_degree}.
\end{proof}

\subsubsection{Putting together}
Following the discussions above, \eqref{eqn:contribution_regularized} adapts the following upper bound
\begin{align*}
    \|(\rmA - \E \rmA)_{\sJ}\| \leq \const_{\eqref{thm:regularization}} \sqrt{\etD_{\max}}\,,
\end{align*}
where the constant $\const_{\eqref{thm:regularization}} \coloneqq 2(\alpha + \beta)$, with $\alpha = 5\LM^{2}$, $\beta = \beta_0 \sqrt{\beta_1} + \LM$, $\beta_0 =16+32\LM (1+e^2)+1792(1+e^{-2})\LM^{2}$, and $\beta_1 = \LM$.

\section{Proofs in Section \ref{sec:IT_lower_bounds}}\label{app:ITLowerBounds}
\begin{proof}[Proof of \Cref{lem:sizeDeviation}]
    This can be established by Hoeffding \Cref{lem:Hoeffding}.
\end{proof}

\begin{proof}[Proof of \Cref{lem:ambiguousProb}]
Note that $\rY_{v} = k^{\star}$. Without loss of generality, we assume $\etQ^{(\ell)}_{k^{\star} \oplus \rvw} \lesssim \etQ^{(\ell)}_{j^{\star} \oplus \rvw}$, then $\etQ^{(\ell)}_{k^{\star} \oplus \rvw} \lesssim \D_{\rm{GCH}} \cdot q_{N}$ according to \Cref{lem:GCH_max_different_order}. Also, $N_{\rvw} \asymp N^{\ell - 1}$ and $\etQ^{(\ell)}_{k^{\star} \oplus \rvw} \asymp q^{(\ell)}_{k^{\star} \oplus \rvw}/N^{\ell - 1} = o(1)$, then $N_{\rvw}\log( 1 - \etQ^{(\ell)}_{k^{\star} \oplus \rvw} ) = - (1 + o(1))\cdot N_{\rvw}\etQ^{(\ell)}_{k^{\star} \oplus \rvw}$. At the same time, $1 \ll d_{j^{\star}, k^{\star}}^{\star}(\rvw) \ll N$ in \eqref{eqn:ambiguous}, then by Stirling's approximation \Cref{lem:stirling}, for large enough $N$,
    \begin{align*}
        &\,\P[ \rD^{(\ell)}_{v, \rvw} = d_{j^{\star}, k^{\star}}^{\star}(\rvw) ] = \binom{N_{\rvw}}{d_{j^{\star}, k^{\star}}^{\star}(\rvw)} ( \etQ^{(\ell)}_{k^{\star} \oplus \rvw} )^{d_{j^{\star}, k^{\star}}^{\star}(\rvw)} (1 - \etQ^{(\ell)}_{k^{\star} \oplus \rvw} )^{N_{\rvw} - d_{j^{\star}, k^{\star}}^{\star}(\rvw)}\\
        = &\, \exp\Big\{ \log \binom{N_{\rvw}}{d_{j^{\star}, k^{\star}}^{\star}(\rvw)} + d_{j^{\star}, k^{\star}}^{\star}(\rvw) \cdot \log ( \etQ^{(\ell)}_{k^{\star} \oplus \rvw}) + \big[ N_{\rvw} - d_{j^{\star}, k^{\star}}^{\star}(\rvw) \big] \cdot \log ( 1 - \etQ^{(\ell)}_{k^{\star} \oplus \rvw}) \Big\}\\
        = &\, \exp\Big\{ d_{j^{\star}, k^{\star}}^{\star}(\rvw) \log\Big( \frac{N_{\rvw}}{d_{j^{\star}, k^{\star}}^{\star}(\rvw)} \Big) + d_{j^{\star}, k^{\star}}^{\star}(\rvw) + d_{j^{\star}, k^{\star}}^{\star}(\rvw) \log \Big( \frac{\etQ^{(\ell)}_{k^{\star} \oplus \rvw}}{1 - \etQ^{(\ell)}_{k^{\star} \oplus \rvw} } \Big) \\
        &\, \quad \quad \quad \quad \quad \quad \quad \quad \quad \quad \quad \quad \quad \quad \quad \quad \quad + N_{\rvw} \cdot \log( 1 - \etQ^{(\ell)}_{k^{\star} \oplus \rvw} ) + o(q^{(\ell)}_{k^{\star} \oplus \rvw})  \Big\}\\
        = &\, \exp\Big\{ d_{j^{\star}, k^{\star}}^{\star}(\rvw) \cdot \log \Big( \frac{ N_{\rvw} \etQ^{(\ell)}_{k^{\star} \oplus \rvw}}{d_{j^{\star}, k^{\star}}^{\star}(\rvw)\cdot (1 - \etQ^{(\ell)}_{k^{\star} \oplus \rvw})} \Big) - N_{\rvw} \cdot \etQ^{(\ell)}_{k^{\star} \oplus \rvw} + d_{j^{\star}, k^{\star}}^{\star}(\rvw)  + o(q^{(\ell)}_{k^{\star} \oplus \rvw})\Big\}\\
        = &\, \exp\Big\{ -\overline{N}_{\rvw} \Big[ t^{\star}\big(\etQ^{(\ell)}_{j^{\star} \oplus \rvw} \big)^{t^{\star}} \big(\etQ^{(\ell)}_{k^{\star} \oplus \rvw} \big)^{1 - t^{\star}} \log \frac{ \etQ^{(\ell)}_{j^{\star} \oplus \rvw}}{\etQ^{(\ell)}_{k^{\star} \oplus \rvw}} + \etQ^{(\ell)}_{k^{\star} \oplus \rvw} - \big(\etQ^{(\ell)}_{j^{\star} \oplus \rvw} \big)^{t^{\star}} \big(\etQ^{(\ell)}_{k^{\star} \oplus \rvw} \big)^{1 - t^{\star}} + o(\etQ^{(\ell)}_{k^{\star} \oplus \rvw})\Big]\Big\}
        \end{align*}
    where in the last inequality, we plug in $d_{j^{\star}, k^{\star}}^{\star}(\rvw)$ in \eqref{eqn:ambiguous}, and replace $N_{\rvw}$ with $\overline{N}_{\rvw}$ according to \Cref{lem:sizeDeviation}. Therefore, by hyperedge independence, we have
    \begin{subequations}
        \begin{align*}
            &\, \prod_{\ell \in \sL} \prod_{ \rvw \in \sW^{(\ell)}_{j^{\star}, k^{\star}} } \P[ \rD^{(\ell)}_{v, \rvw} = d_{j^{\star}, k^{\star}}^{\star}(\rvw) ]\\
            =&\, \exp\Big\{ -\sum_{\ell \in \sL} \sum_{ \rvw \in \sW^{(\ell)}_{j^{\star}, k^{\star}} } \overline{N}_{\rvw} \Big[ t^{\star}\big(\etQ^{(\ell)}_{j^{\star} \oplus \rvw} \big)^{t^{\star}} \big(\etQ^{(\ell)}_{k^{\star} \oplus \rvw} \big)^{1 - t^{\star}} \log\bigg( \frac{ \etQ^{(\ell)}_{j^{\star} \oplus \rvw}}{\etQ^{(\ell)}_{k^{\star} \oplus \rvw}} \bigg) \\
            &\, \quad \quad \quad \quad \quad \quad \quad \quad \quad \quad \quad \quad \quad \quad \quad  + \etQ^{(\ell)}_{k^{\star} \oplus \rvw} - \big(\etQ^{(\ell)}_{j^{\star} \oplus \rvw} \big)^{t^{\star}} \big(\etQ^{(\ell)}_{k^{\star} \oplus \rvw} \big)^{1 - t^{\star}} + o(\etQ^{(\ell)}_{k^{\star} \oplus \rvw}) \Big] \Big\}\\
            =&\, \exp\Big\{ -\sum_{\ell \in \sL} \sum_{ \rvw \in \sW^{(\ell)}_{j^{\star}, k^{\star}} } \overline{N}_{\rvw}\cdot \Big[ t^{\star}\etQ^{(\ell)}_{j^{\star} \oplus \rvw} + (1 - t^{\star})\etQ^{(\ell)}_{k^{\star} \oplus \rvw} - \big(\etQ^{(\ell)}_{j^{\star} \oplus \rvw} \big)^{t^{\star}} \big(\etQ^{(\ell)}_{k^{\star} \oplus \rvw} \big)^{1 - t^{\star}} + o(\etQ^{(\ell)}_{k^{\star} \oplus \rvw}) \Big] \Big\} \\
            =&\, \exp(-(1 + o(1))\cdot\D_{\rm{GCH}}(j, k)\cdot q_{N}),
        \end{align*}
    \end{subequations}
    where the second to last equality holds according to the first order condition \eqref{eqn:optimalT}.
\end{proof}

\begin{proof}[Proof of \Cref{lem:equivalenveDivergence}]
    We rewrite $\D_{\textnormal{GKL}} (j, k)$ as
    \begin{align}\label{eqn:objectiveDivergence}
        &\underset{\{\rvy^{(\ell)}\}_{\ell \in \sL} }{\min} \sum_{\ell \in \sL} \, \sum_{\rvw \in \WC{\ell - 1}{K}} \overline{N}_{\rvw}\cdot \Big\{ \ervy^{(\ell)}_{\rvw} \log \frac{\ervy^{(\ell)}_{\rvw} }{\etQ^{(\ell)}_{k \oplus \rvw}} + (1 - \ervy^{(\ell)}_{\rvw}) \log \frac{1 - \ervy^{(\ell)}_{\rvw} }{1 - \etQ^{(\ell)}_{k \oplus \rvw}} \Big\} \\
        &\quad \quad  \textnormal{s.t.} \sum_{\ell \in \sL} \, \sum_{\rvw \in \WC{\ell - 1}{K} } \overline{N}_{\rvw} \cdot \big\{ \D_{\textnormal{KL}}(\ervy^{(\ell)}_{\rvw}\parallel \etQ_{j \oplus \rvw}^{(\ell)}) - \D_{\textnormal{KL}}(\ervy^{(\ell)}_{\rvw}\parallel \etQ_{k \oplus \rvw}^{(\ell)}) \big\} \leq 0\notag.
    \end{align}
By \cite{Boyd2004ConvexO}, the associated Lagrangian is defined as
\begin{align*}
    L(\{\rvy^{(\ell)}\}_{\ell \in \sL}\,, \lambda) \coloneqq &\, \sum_{\ell \in \sL} \, \sum_{\rvw \in \WC{\ell - 1}{K}} \overline{N}_{\rvw}\cdot \Big\{ \ervy^{(\ell)}_{\rvw} \log \frac{\ervy^{(\ell)}_{\rvw} }{\etQ^{(\ell)}_{k \oplus \rvw}} + (1 - \ervy^{(\ell)}_{\rvw}) \log \frac{1 - \ervy^{(\ell)}_{\rvw} }{1 - \etQ^{(\ell)}_{k \oplus \rvw}} \Big\} \\
    &\,+ \lambda \sum_{\ell \in \sL} \, \sum_{\rvw \in \WC{\ell - 1}{K} } \overline{N}_{\rvw} \cdot \Big\{ \ervy^{(\ell)}_{\rvw} \log \frac{\etQ^{(\ell)}_{k \oplus \rvw}}{\etQ^{(\ell)}_{j \oplus \rvw}} + (1 - \ervy^{(\ell)}_{\rvw}) \log \frac{1 - \etQ^{(\ell)}_{k \oplus \rvw}}{1 - \etQ^{(\ell)}_{j \oplus \rvw}} \Big\}.
\end{align*}
The optimal choice of $\{\rvy^{(\ell)}\}_{\ell \in \sL}$, which minimizes $L(\{\rvy^{(\ell)}\}_{\ell \in \sL}\,, \lambda)$, satisfies the first order condition
\begin{align*}
    \frac{\partial L(\{\rvy^{(\ell)}\}_{\ell \in \sL}, \lambda) }{\partial \ervy^{(\ell)}_{\rvw}} = \overline{N}_{\rvw} \cdot \Big[ \Big( \log \frac{\ervy^{(\ell)}_{\rvw} }{\etQ^{(\ell)}_{k \oplus \rvw} } - \log \frac{1 - \ervy^{(\ell)}_{\rvw}}{1 - \etQ^{(\ell)}_{k \oplus \rvw} } \Big) + \lambda \Big( \log \frac{\etQ^{(\ell)}_{k \oplus \rvw}}{\etQ^{(\ell)}_{j \oplus \rvw}} - \log \frac{1 - \etQ^{(\ell)}_{k \oplus \rvw}}{1 - \etQ^{(\ell)}_{j \oplus \rvw}} \Big) \Big] =0.
\end{align*}
Note that $\etQ^{(\ell)}_{k \oplus \rvw} = o(1)$ in \eqref{eqn:edge_density}, then the solution to \eqref{eqn:objectiveDivergence} has to satisfy $\ervy^{(\ell)}_{\rvw} = o(1)$ for each $\rvw \in \WC{\ell - 1}{K}$, otherwise optimal value of \eqref{eqn:objectiveDivergence} goes to infinity since $\ervy^{(\ell)}_{\rvw} \log \frac{\ervy^{(\ell)}_{\rvw} }{\etQ^{(\ell)}_{k \oplus \rvw}}$ would diverge. Thus $\log \frac{1 - \ervy^{(\ell)}_{\rvw}}{1 - \etQ^{(\ell)}_{k \oplus \rvw} }$ and $\log \frac{1 - \etQ^{(\ell)}_{k \oplus \rvw}}{1 - \etQ^{(\ell)}_{j \oplus \rvw}}$ converge to $0$ as $N$ grows. Then \eqref{eqn:objectiveDivergence} is minimized at
\begin{align*}
    \ervy^{(\ell)}_{\rvw} = (\etQ^{(\ell)}_{j \oplus \rvw} )^{\lambda}\cdot (\etQ_{k \oplus \rvw}^{(\ell)} )^{1 - \lambda}\cdot (1 + o(1))\,, \, \textnormal{for each }\, \rvw \in \WC{\ell - 1}{K}, \,\, \ell \in \sL.
\end{align*}
Consider the Lagrange dual function
\begin{align*}
    &\,g(\lambda) \coloneqq \underset{\{\rvy^{(\ell)}\}_{\ell \in \sL}}{\min}  L(\{\rvy^{(\ell)}\}_{\ell \in \sL},\,\, \lambda)\\
    =&\, \underset{\{\rvy^{(\ell)}\}_{\ell \in \sL}}{\min} \sum_{\ell \in \sL} \, \sum_{\rvw \in \WC{\ell - 1}{K} } \overline{N}_{\rvw} \cdot \Big\{ \ervy^{(\ell)}_{\rvw} \log \frac{(\etQ^{(\ell)}_{k \oplus \rvw})^{\lambda - 1 } \ervy^{(\ell)}_{\rvw} }{(\etQ^{(\ell)}_{j \oplus \rvw})^{\lambda} } + (1 - \ervy^{(\ell)}_{\rvw}) \log \frac{ (1 - \etQ^{(\ell)}_{k \oplus \rvw})^{\lambda - 1} (1 - \ervy^{(\ell)}_{\rvw}) }{ (1 - \etQ^{(\ell)}_{j \oplus \rvw} )^{\lambda} } \Big\}\\
    =&\, \underset{\{\rvy^{(\ell)}\}_{\ell \in \sL}}{\min} \sum_{\ell \in \sL} \, \sum_{\rvw \in \WC{\ell - 1}{K} } \overline{N}_{\rvw} \cdot (1 - \ervy^{(\ell)}_{\rvw}) \cdot \Big[ (\lambda - 1) \log(1 - \etQ^{(\ell)}_{k \oplus \rvw}) + \log(1 - \ervy^{(\ell)}_{\rvw}) - \lambda \log(1 - \etQ^{(\ell)}_{j \oplus \rvw} ) \Big]\\
    =&\, \sum_{\ell \in \sL} \, \sum_{\rvw \in \WC{\ell - 1}{K} } \overline{N}_{\rvw} \cdot \Big[ \lambda  \etQ^{(\ell)}_{j \oplus \rvw} + (1 - \lambda) \etQ^{(\ell)}_{k \oplus \rvw} - (\etQ^{(\ell)}_{j \oplus \rvw} )^{\lambda}\cdot (\etQ^{(\ell)}_{k \oplus \rvw} )^{1 - \lambda} \Big],
\end{align*}
where in the last two equalities we use the facts $\ervy^{(\ell)}_{\rvw} = o(1)$, $\etQ^{(\ell)}_{j \oplus \rvw} = o(1)$ and $\lim_{x\to 0} \log(1 + x) = x$. By taking maximum of $g(\lambda)$, we have
\begin{align*}
    &\,\D_{\textnormal{GKL}} (j, k) = \eqref{eqn:objectiveDivergence} = \underset{\lambda\in[0, 1]}{\max} \,\,g(\lambda)\\
    =&\, \underset{\lambda\in[0, 1]}{\max} \sum_{\ell \in \sL} \, \sum_{\rvw \in \WC{\ell - 1}{K} } \, \overline{N}_{\rvw} \cdot \Big[ \lambda  \etQ^{(\ell)}_{j \oplus \rvw} + (1 - \lambda) \etQ^{(\ell)}_{k \oplus \rvw} - (\etQ^{(\ell)}_{j \oplus \rvw} )^{\lambda}\cdot (\etQ_{k \oplus \rvw}^{(\ell)} )^{1 - \lambda} \Big]\\
    \overset{\eqref{eqn:edge_density}}{=} &\, q_{N} \cdot \underset{\lambda\in[0, 1]}{\max} \sum_{\ell \in \sL} \, \sum_{\rvw \in \WC{\ell - 1}{K} } \, \frac{\overline{N}_{\rvw}}{\binom{N - 1}{\ell - 1}} \cdot \Big[ \lambda  \etP^{(\ell)}_{j \oplus \rvw} + (1 - \lambda) \etP^{(\ell)}_{k \oplus \rvw} - (\etP_{j \oplus \rvw} )^{\lambda}\cdot (\etP_{k \oplus \rvw}^{(\ell)} )^{1 - \lambda} \Big]\\
    =&\, (1 + o(1))\cdot \D_{\rm{GCH}}(j,k) \cdot q_{N}\,.
\end{align*}
\end{proof}

\begin{proof}[Proof of \Cref{lem:qMeasureExistence}]
    Proof by contradiction. Suppose that there is no such $\{\rvp^{(\ell)}\}_{\ell \in \sL}$ satisfying \eqref{eqn:equivalenceqm}. Then for each $\{\rvp^{(\ell)}\}_{\ell \in \sL}$, without loss of generality, we assume
    \begin{align*}
        \D_{\mathrm{GKL}}(j^{\star}, k^{\star}) =&\, \sum_{\ell \in \sL} \sum_{\rvw \in \WC{\ell - 1}{K}} \overline{N}_{\rvw} \cdot  \D_{\mathrm{KL}}(\ervp^{(\ell)}_{\rvw}\parallel \etQ_{j^{\star} \oplus \rvw}^{(\ell)}) \\
        >&\, \sum_{\ell \in \sL} \sum_{\rvw \in \WC{\ell - 1}{K}}\overline{N}_{\rvw}\cdot \D_{\mathrm{KL}}(\ervp^{(\ell)}_{\rvw}\parallel \etQ_{k^{\star} \oplus \rvw}^{(\ell)} ).
    \end{align*}
    We pick a sequence $\{\rvp^{(\ell)}\}_{\ell \in \sL}$, where $\rvp^{(\ell)} = (\ervp^{(\ell)}_{\rvw})_{\rvw \in \WC{\ell - 1}{K}}$ for each $\ell$. There exists some $\rvw_0 \in \WC{\ell - 1}{K}$ such that $\D_{\mathrm{KL}}(\ervp^{(\ell)}_{\rvw_0}\parallel \etQ_{j^{\star} \oplus \rvw_0}^{(\ell)}) > \D_{\mathrm{KL}}(\ervp^{(\ell)}_{\rvw_0}\parallel \etQ_{k^{\star} \oplus \rvw_0}^{(\ell)}) \geq 0$, due to the non-negativity of KL divergence, with their difference $f(\ervp^{(\ell)}_{\rvw_0}) >0$ defined as
    \begin{align*}
        f(\ervp) \coloneqq \D_{\mathrm{KL}}(\ervp \parallel \etQ_{j^{\star} \oplus \rvw_0}^{(\ell)}) - \D_{\mathrm{KL}}(\ervp \parallel \etQ_{k^{\star} \oplus \rvw_0}^{(\ell)}) = \ervp \cdot \log \frac{\etQ_{k^{\star} \oplus \rvw_0}^{(\ell)}}{\etQ_{j^{\star} \oplus \rvw_0}^{(\ell)}} + (1 - \ervp) \cdot \log \frac{1 - \etQ_{k^{\star} \oplus \rvw_0}^{(\ell)}}{1 - \etQ_{j^{\star} \oplus \rvw_0}^{(\ell)}}.
    \end{align*}
    By continuity of $f(\ervp)$ and $\D_{\mathrm{KL}}(\ervp \parallel \etQ_{j^{\star} \oplus \rvw_0}^{(\ell)})$, one is able to find some $\ervt^{(\ell)}_{\rvw_0}$ such that $\D_{\mathrm{KL}}(\ervt^{(\ell)}_{\rvw_0} \parallel \etQ_{j^{\star} \oplus \rvw_0}^{(\ell)}) < \D_{\mathrm{KL}}(\ervp^{(\ell)}_{\rvw_0} \parallel \etQ_{j^{\star} \oplus \rvw_0}^{(\ell)})$ but $f(\ervt^{(\ell)}_{\rvw_0})$ remains positive. As a result, for this fixed $\ell \in \sL$, there exists some sequence $( \ervt^{(\ell)}_{\rvw} )_{\rvw \in \WC{\ell - 1}{K}}$ with $\ervt^{(\ell)}_{\rvw} = \ervp^{(\ell)}_{\rvw}$ for all $\rvw \neq \rvw_0$, and $\ervt^{(\ell)}_{\rvw_0}$ satisfying the property as before. Consequently, one could construct a new sequence $\{\rvt^{(\ell)}\}_{\ell \in \sL}$ satisfying
    \begin{align*}
        \D_{\mathrm{GKL}}(j^{\star}, k^{\star}) &\, = \sum_{\ell \in \sL} \sum_{\rvw \in \WC{\ell - 1}{K}} \overline{N}_{\rvw} \cdot  \D_{\mathrm{KL}}(\ervp^{(\ell)}_{\rvw}\parallel \etQ_{j^{\star} \oplus \rvw}^{(\ell)})\\
        &\, > \sum_{\ell \in \sL} \sum_{\rvw \in \WC{\ell - 1}{K}} \overline{N}_{\rvw} \cdot \D_{\mathrm{KL}}(\ervt^{(\ell)}_{\rvw}\parallel \etQ_{j^{\star} \oplus \rvw}^{(\ell)}) \\
        &\, > \sum_{\ell \in \sL} \sum_{\rvw \in \WC{\ell - 1}{K}} \overline{N}_{\rvw} \cdot \D_{\mathrm{KL}}(\ervt^{(\ell)}_{\rvw}\parallel \etQ_{k^{\star} \oplus \rvw}^{(\ell)}), 
    \end{align*}
    which contradicts to the minimum property in definition of $\D_{\mathrm{GKL}}(j^{\star}, k^{\star})$ in \eqref{eqn:GeneralizedKL}.
\end{proof}

\begin{proof}[Proof of \Cref{lem:LambdaUpperBound}]
Note that the nodes are interchangeable under $\Phi$, and
\begin{align}\label{eqn:Lambda_bounded_fn}
    \P_{\Psi}(\Lambda \leq f(N)) = \P_{\Psi}(\Lambda \leq f(N),\,\, v^{\star} \in \gV_{\rm{err}}) + \P_{\Psi}(\Lambda \leq f(N),\,\, v^{\star} \notin \gV_{\rm{err}}).
\end{align}
The first term on the right hand side of \eqref{eqn:Lambda_bounded_fn} can be upper bounded by the following \emph{change of measure} argument,
\begin{align*}
    &\,\P_{\Psi}(\Lambda \leq f(N),\,\, v^{\star} \in \gV_{\rm{err}}) = \int_{\Lambda \leq f(N), v^{\star} \in \gV_{\rm{err}} } d\P_{\Psi} = \int_{\Lambda \leq f(N), v^{\star} \in \gV_{\rm{err}}} \exp(\Lambda) d\P_{\Phi}\\
    \leq &\, \exp(f(N)) \cdot \P_{\Phi} (\Lambda \leq f(N), v^{\star} \in \gV_{\rm{err}} ) \leq \exp(f(N)) \cdot \P_{\Phi} ( v^{\star} \in \gV_{\rm{err}} ) \leq \frac{e^{f(N)} \cdot \E_{\Phi}[\mismatch_{N}]}{\alpha_{j^{\star}} + \alpha_{k^{\star}}},
\end{align*}
where the last inequality holds since we can not distinguish $v^{\star}$ and any other node $v\in \gV_{\ervy_{v^{\star}}}$, indeed 
\begin{align*}
    \P_{\Phi}(v^{\star} \in \gV_{\rm{err}}) =&\, \P_{\Phi}(v\in \gV_{\rm{err}}|v\in \gV_{j^{\star}} \cup \gV_{k^{\star}}) \\
    \leq &\, \frac{\P_{\Phi}(v\in \gV_{\rm{err}})}{\P_{\Phi}(v\in \gV_{j^{\star}} \cup \gV_{k^{\star}})} = \frac{\E_{\Phi}[|\gV_{\rm{err}}|]}{(\alpha_{j^{\star}} + \alpha_{k^{\star}})n} = \frac{\E_{\Phi}[\mismatch_{N}]}{\alpha_{j^{\star}} + \alpha_{k^{\star}}}.
\end{align*}

Now, we focus on the second term on the right hand side of \eqref{eqn:Lambda_bounded_fn}. Under the perturbed model $\Psi$ in \ref{def:perturbed_model}, the observed edges do not depend on whether $v^{\star}$ belongs to $\gV_{j^{\star}}$ or $\gV_{k^{\star}}$, since the generation of edge $e \ni v^{\star}$ only relies on the distribution of $\rvw \in \WC{\ell - 1}{K}$ for the rest $\ell - 1$ nodes. To determine the community membership of $v^{\star}$, one can only use the information from the edges containing $v^{\star}$. Therefore, it doesn't matter whether $v^{\star}$ truly belongs to $\gV_{j^{\star}}$ or $\gV_{k^{\star}}$ when assigning $v^{\star}$ to $\widehat{\gV}_{j^{\star}}$. By symmetry, the two conditional probabilities should be the same, and the following holds
\begin{align}\label{eqn:equaling_error_prob}
    \P_{\Psi}(v^{\star}\in \widehat{\gV}_{j^{\star}}| v^{\star}\in \gV_{j^{\star}}) =&\, \P_{\Psi}(v^{\star}\in \widehat{\gV}_{j^{\star}}| v^{\star}\in \gV_{k^{\star}}), \\
    \P_{\Psi}(v^{\star}\in \widehat{\gV}_{k^{\star}}| v^{\star}\in \gV_{j^{\star}}) =&\, \P_{\Psi}(v^{\star}\in \widehat{\gV}_{k^{\star}}| v^{\star}\in \gV_{k^{\star}}).\notag
\end{align}
Then the second term on the right hand side of \eqref{eqn:Lambda_bounded_fn} can be bounded as follows
\begin{align*}
    &\, \P_{\Psi}(\Lambda \leq f(N),\,\, v^{\star} \notin \gV_{\rm{err}} ) \\
    \leq &\,  \P_{\Psi}(v^{\star} \notin \gV_{\rm{err}}) \\
    = &\, \frac{\alpha_{j^{\star}}}{\alpha_{j^{\star}} + \alpha_{k^{\star}}}\P_{\Psi}(v^{\star}\in \widehat{\gV}_{j}| v^{\star} \in \gV_{j^{\star}}) + \frac{\alpha_{k^{\star}}}{\alpha_{j^{\star}} + \alpha_{k^{\star}}} \P_{\Psi}(v^{\star}\in \widehat{\gV}_{k}| v^{\star} \in \gV_{k^{\star}}) \\
    \overset{\eqref{eqn:equaling_error_prob}}{=} &\, \frac{\alpha_{j^{\star}}}{\alpha_{j^{\star}} + \alpha_{k^{\star}}}\P_{\Psi}(v^{\star}\in \widehat{\gV}_{j}| v^{\star} \in \gV_{j^{\star}}) + \frac{\alpha_{k^{\star}}}{\alpha_{j^{\star}} + \alpha_{k^{\star}}} \P_{\Psi}(v^{\star}\in \widehat{\gV}_{k}| v^{\star} \in \gV_{j^{\star}}) \\
    \leq &\, \frac{\alpha_{j^{\star}}}{\alpha_{j^{\star}} + \alpha_{k^{\star}}} \Big( \P_{\Psi}(v^{\star}\in \widehat{\gV}_{j}| v^{\star} \in \gV_{j^{\star}}) + \P_{\Psi}(v^{\star}\in \widehat{\gV}_{k}| v^{\star} \in \gV_{j^{\star}}) \Big) \leq \frac{\alpha_{j^{\star}}}{\alpha_{j^{\star}} + \alpha_{k^{\star}}},
\end{align*}
where the last line holds since $j^{\star} < k^{\star}$ and $\alpha_1 \geq \ldots \geq \alpha_K$ implies $\alpha_{j^{\star}} \geq \alpha_{k^{\star}}$, together with the fact $\P_{\Psi}(v^{\star}\in \widehat{\gV}_{j^{\star}}| v^{\star}\in \gV_{j^{\star}}) + \P_{\Psi}(v^{\star}\in \widehat{\gV}_{k^{\star}}| v^{\star}\in \gV_{j^{\star}}) \leq 1$. Therefore the desired result follows.
\end{proof}

\begin{proof}[Proof of \Cref{lem:expectation_Lambda_Psi}]
Since $\Lambda$ in \eqref{eqn:Lambda} is a summation of Bernoulli random variables, a standard concentration argument shows that $\sqrt{\Var_{\Psi}(\Lambda)}/\E_{\Psi}(\Lambda) = o(1)$, thus $\E_{\Psi}(\Lambda)$ dominates. We present the details by evaluating $\E_{\Psi}[\Lambda]$ and $\Var_{\Psi}(\Lambda)$. First, the randomness of $\Lambda$ comes from edges $\etA^{(\ell)}_e$, then by linearity of expectation,
\begin{align*}
    \E_{\Psi} [\Lambda] = &\,\sum_{\ell \in \sL} \,\, \sum_{\substack{e\in \gE_{\ell},\,\, e\ni v^{\star}\\ \rvy(e\setminus v^{\star}) = \rvw \in \WC{\ell - 1}{K}}} \E \bigg( \etA^{(\ell)}_e \log\Big( \frac{\ervp^{(\ell)}_{\rvw}}{\etQ^{(\ell)}_{\ervy_{v^{\star}}\oplus \rvw}} \Big) + (1 - \etA^{(\ell)}_e) \log\Big( \frac{1 - \ervp^{(\ell)}_{\rvw}}{1 - \etQ^{(\ell)}_{\ervy_{v^{\star}}\oplus \rvw}} \Big) \bigg)\\
    = &\,\sum_{\ell \in \sL} \,\, \sum_{\substack{e\in \gE_{\ell},\,\, e\ni v^{\star}\\ \rvy(e\setminus v^{\star}) = \rvw \in \WC{\ell - 1}{K}}} \bigg( \ervp^{(\ell)}_{\rvw} \log\Big( \frac{\ervp^{(\ell)}_{\rvw}}{\etQ^{(\ell)}_{\ervy_{v^{\star}}\oplus \rvw}} \Big) + (1 - \ervp^{(\ell)}_{\rvw}) \log\Big( \frac{1 - \ervp^{(\ell)}_{\rvw}}{1 - \etQ^{(\ell)}_{\ervy_{v^{\star}}\oplus \rvw}} \Big) \bigg)\\
    =&\, \sum_{\ell \in \sL} \sum_{\ell \in \WC{\ell - 1}{K}} N_{\rvw} \cdot \D_{\mathrm{KL}}(\ervp^{(\ell)}_{\rvw}\parallel \etQ_{\ervy_{v^{\star}} \oplus \rvw}^{(\ell)} ) = (1 + o(1)) \cdot \D_{\mathrm{GKL}} = (1 + o(1)) \cdot \D_{\rm{GCH}} \cdot q_{N},
\end{align*}
where the last two equalities hold by Lemmas \ref{lem:sizeDeviation}, \ref{lem:equivalenveDivergence}, \ref{lem:qMeasureExistence}. On the other hand, using the facts $\Var(c \etA^{(\ell)}_{e}) = c^2\Var( \etA^{(\ell)}_{e})$ and $\Var(\etA^{(\ell)}_{e} + c) = \Var(\etA^{(\ell)}_{e})$ for any deterministic $c \in \R$, we have
\begin{align*}
    \Var_{\Psi}(\Lambda) = \sum_{\ell \in \sL} \sum_{\ell \in \WC{\ell - 1}{K}} N_{\rvw} \ervp^{(\ell)}_{\rvw}(1 - \ervp^{(\ell)}_{\rvw}) \cdot \bigg[\log\bigg( \frac{\ervp^{(\ell)}_{\rvw} \cdot (1 - \etQ^{(\ell)}_{\ervy_{v^{\star}}\oplus \rvw}) }{ \etQ^{(\ell)}_{\ervy_{v^{\star}}\oplus \rvw} \cdot (1 - \ervp^{(\ell)}_{\rvw}) } \bigg) \bigg]^2.
\end{align*}
Note that $\etQ_{\rvw}^{(\ell)} = o(1)$ for each $\rvw \in \WC{\ell}{K}$, $\ell \in \sL$, then by the argument above, the so constructed $\{\rvp^{(\ell)}\}_{\ell \in \sL}$ in \Cref{lem:qMeasureExistence} satisfies $\ervp_{\rvw}^{(\ell)} = o(1)$ and $N_{\rvw}\ervp^{(\ell)}_{\rvw} \gg 1$. 
Then the desired result follows since
\begin{align*}
    \sqrt{\Var_{\Psi}(\Lambda)}/\E_{\Psi} [\Lambda] \lesssim \frac{1}{\sqrt{\D_{\rm{GCH}} \cdot q_{N}}} = o(1).
\end{align*}
\end{proof}

\section{Proofs in Section \ref{sec:almostExact}}\label{app:almostExact}
\begin{proof}[Proof of \Cref{lem:maxProb}]
    For the first part, $\xi_{N} \asymp \rho_{N}$ is trivially true when \Cref{ass:prob_ratio_bound} holds. On the other hand, due to the fact $\alpha_1, \ldots, \alpha_K = \Omega(1)$, the conclusion $\xi_{N} \asymp \rho_{N}$ still holds since the dense part inside the summation $\xi_{N} = N^{-1}\sum_{v=1}^{N}\E[\rD_v]$ distinguishes itself.
    
    For the second part, by Bernstein \ref{lem:Bernstein}, we have
    \begin{align*}
        &\,\P(\overline{d} - \xi_{N} \geq \xi_{N}/\sqrt{N}) = \P\Big(\sum_{v\in \gV}\sum_{\ell \in \sL} \sum_{\rvw \in\WC{\ell - 1}{K} } (\rD_{v, \rvw}^{(\ell)} - \E \rD_{v, \rvw}^{(\ell)} ) \geq \sqrt{N} \xi_{N} \Big)\\
        \leq &\, 2\exp\Big( -\frac{N\cdot \xi_{N}^2/2}{ M\binom{K + \LM - 1}{\ell - 1}\xi_{N} + \sqrt{N}\xi_{N}/3}\Big).
    \end{align*}
    Then $\overline{d} - \xi_{N} \leq \xi_{N}/\sqrt{N}$ with probability at least $1 - e^{-N^{1/2}}$ as long as $\xi_{N} = \omega(1)$.
    
    For the last part, we assume $q_{N} \geq \log(N)$ in \eqref{eqn:edge_density} from now on. We are going to present some large deviation results. Let $\widetilde{d} = \max_{v\in\gV} d_v$ denote realized maximum degree. Define the random variable $\widetilde{\rD}$ by $\widetilde{\rD} = \max_{v\in\gV} \rD_v$, then $\E \widetilde{\rD} = \E \max_{v\in\gV} \rD_v$. Define $\rW_v \coloneqq \rD_v - \E \rD_v$, then following \eqref{eqn:Dv}, $\rW_v = \sum_{\ell \in \sL} \sum_{\rvw \in\WC{\ell - 1}{K} } \rW_{v, \rvw}^{(\ell)}$ where $\rW_{v, \rvw}^{(\ell)}$ are mean zero random variables, and define
\begin{align}
    \widetilde{\rW} \coloneqq \max_{v\in\gV} \rW_v \leq \sum_{\ell \in \sL} \sum_{\rvw \in\WC{\ell - 1}{K} } \max_{v\in\gV} \rW_{v, \rvw}^{(\ell)}.
\end{align}
Note that $\rD^{(\ell)}_{v, \rvw} \sim \mathrm{Bin}(N_{\rvw}, \etQ^{(\ell)}_{k \oplus \rvw})$ when $\rY_{v} = k$, then $\sigma^2 = \E(\rW_{v, \rvw}^{(\ell)})^2 = N_{\rvw}\etQ^{(\ell)}_{k \oplus \rvw}(1 - \etQ^{(\ell)}_{k \oplus \rvw}) \leq \E \rD^{(\ell)}_{v, \rvw}  < \E\widetilde{\rD}$. By Bernstein \ref{lem:Bernstein}, for any $\epsilon >0$, we have
\begin{align}
    \P[ \rW_{v, \rvw}^{(\ell)} \geq (1 - \epsilon) \E\widetilde{\rD}] \leq \exp \Big( - \frac{3(1 - \epsilon)^2 (\E\widetilde{\rD})^2}{6\sigma^2 + 2\sqrt{(1 - \epsilon)\E\widetilde{\rD}} } \Big) \leq \exp \Big( - \frac{1}{2}(1 - \epsilon)^2 \E\widetilde{\rD} \Big)
\end{align}
Note that $\E\widetilde{\rD} \geq c q_{N}$ for some $c>0$ under regime \eqref{eqn:edge_density}. By the union bounds
\begin{align*}
    &\, \P(\widetilde{\rW} \geq (1 - \epsilon) \E\widetilde{\rD} ) \leq \sum_{\ell \in \sL} \sum_{\rvw \in\WC{\ell - 1}{K} } \P( \max_{v\in\gV} \rW_{v, \rvw}^{(\ell)} \geq (1 - \epsilon) \E\widetilde{\rD} ) \\
    \leq &\, \sum_{\ell \in \sL} \sum_{\rvw \in\WC{\ell - 1}{K} }\sum_{v=1}^{N} \P[ \rW_{v, \rvw}^{(\ell)} \geq (1 - \epsilon) \E\widetilde{\rD} ] \leq (\ell - 1) \cdot \binom{\ell+K-1}{K-1} \exp \Big( \log(N) - \frac{1}{2}(1 - \epsilon)^2 \E\widetilde{\rD} \Big)\\
    \lesssim &\, \exp \big\{ - \log(N) \cdot [ (1 - \epsilon)^2 c^2/2 -1 ] \big\} \overset{ N \to \infty}{\longrightarrow } 0, \quad \textnormal{ when } c > \sqrt{2}/(1 - \epsilon),
\end{align*}
where the last inequality holds when we choose $q_{N} = \log(N)$ and $c > \sqrt{2}/(1 - \epsilon)$. Note that $\D_{\rm{GCH}} > 1$ implies $c > 2^{\ell - 1}$ for $\ell \in\sL$, then it just a matter of choosing proper $\epsilon$. Similarly, $\widetilde{\rW} \geq -(1 - \epsilon) \E\widetilde{\rD}$ with probability tending to $1$.
Note that $|\widetilde{\rD} - \E\widetilde{\rD}| \leq |\widetilde{\rW}| \leq (1 - \epsilon) \E\widetilde{\rD}$ with high probability, then
\begin{align}
   \epsilon \E\widetilde{\rD} \leq \widetilde{\rD} \leq 2 \E\widetilde{\rD}.
\end{align}

We now refer to the magnitude of $\E\widetilde{\rD}$. By Jensen's inequality, for any $t>0$, we have
\begin{align}
    \exp(t\E\widetilde{\rD}) \leq \E \exp(t\widetilde{\rD}) = \E \exp(t\max_{v\in\gV} \rD_v) \leq \prod_{\ell \in \sL} \prod_{\rvw \in \WC{\ell - 1}{K}} \E \exp(t\max_{v\in\gV} \rD^{(\ell)}_{v, \rvw}),
\end{align}
since $\widetilde{\rD} \leq \sum_{\ell \in \sL} \sum_{\rvw \in \WC{\ell - 1}{K}} \max_{v\in\gV} \rD^{(\ell)}_{v, \rvw}$. By the union bound,
\begin{align}
    &\, \E \exp(t\max_{v\in\gV} \rD^{(\ell)}_{v, \rvw}) \leq \sum_{v=1}^{N} \E\exp(t\rD^{(\ell)}_{v, \rvw}) = \sum_{v=1}^{N} [ \etQ^{(\ell)}_{\rvy(v)\oplus \rvw} e^{t} + (1 - \etQ^{(\ell)}_{\rvy(v)\oplus \rvw})]^{N_{\rvw}}\\
    \leq &\, N [ 1 + \etQ^{(\ell)}_{\max} (e^{t} -1)]^{N_{\rvw}} \leq \exp[ \log(N) + N_{\rvw}\etQ^{(\ell)}_{\max} (e^{t} -1)],
\end{align}
since by definition $\etQ^{(\ell)}_{\max} = \max_{\rvw \in \WC{\ell}{K}} \etQ^{(\ell)}_{\rvw}$ and $1 + x \leq e^{x}$. We adapt the notation in \eqref{eqn:edge_density} and \eqref{thm:concentration}, then $\etD_{\max} = \sum_{\ell \in \sL} N_{\rvw}\etQ^{(\ell)}_{\max}$. Note that $N_{\rvw}\etQ^{(\ell)}_{\max} \leq \etP^{(\ell)}_{\max} q_{N}$, By taking $t = 1$, we have
\begin{align}
    \E\widetilde{\rD} \leq &\, \sum_{\ell \in \sL} \sum_{\rvw \in \WC{\ell - 1}{K}}  (\log(N) + (e-1)N_{\rvw}\etQ^{(\ell)}_{\max}) \\
    \leq &\, \sum_{\ell \in \sL} |\WC{\ell - 1}{K}|\cdot  (\log(N) + (e-1)N_{\rvw}\etQ^{(\ell)}_{\max})\lesssim \etD_{\max}\,,
\end{align}
where $\sum_{\ell \in \sL}N_{\rvw}\etQ^{(\ell)}_{\max} \gtrsim \log(N)$ since $q_{N} = \log(N)$ and $|\WC{\ell - 1}{K}| =\binom{\ell+K-1}{K-1} = O(1)$. At the same time, $\E\widetilde{\rD} \geq \max_{\ell \in \sL} N_{\rvw}\etQ^{(\ell)}_{\max} \asymp \etD_{\max}$. Therefore $\widetilde{d} \asymp \etD_{\max}$ with high probability. The desired result follows if $\etD_{\max} \asymp \rho_{N}$, which is obviously true when \Cref{ass:prob_ratio_bound} holds. When the assumption is not true, where some specific configuration $\rvw \in \WC{\ell}{K}$ dominates for some $\ell \in \sL$, then the vertex with maximum expected degree should certainly contain this type of edge. At the same time, $\etD_{\max}$ should contain the dominant type $\rvw$ as well. Therefore $\etD_{\max} \asymp \rho_{N} \asymp \widetilde{d}$ with high probability.
\end{proof}

\begin{proof}[Proof of \Cref{lem:frobeniusNormBound}]
Define the matrix $\E \widetilde{\rmA}$ where $\E \widetilde{\ermA}_{ij} = \E \ermA_{ij}$ for $i\neq j$ and $\E \widetilde{\ermA}_{ii} = \E \ermA_{ij}$ for some $j\in \gV$ in the same community as $i$, i.e., $\ervy_i = \ervy_j$, then $\mathrm{rank}(\E \widetilde{\rmA}) = K$. By triangle inequality and $(a + b)^2 \leq 2a^2 + 2b^2$,
\begin{align}
    \|\rmA_{\sJ}^{(K)} - \E\rmA_{\sJ}\|_{\frob}^2 \leq 2\|\rmA_{\sJ}^{(K)} - \E \widetilde{\rmA}_{\sJ}\|_{\frob}^2
 + 2\|\E \rmA_{\sJ} - \E \widetilde{\rmA}_{\sJ}\|_{\frob}^2
\end{align}
The matrices $\E \rmA_{\sJ}$ and $\E \widetilde{\rmA}_{\sJ}$ only differ in diagonal elements while $\E \widetilde{\ermA}_{ii} \lesssim N^{-1}\etD_{\max}$, then $\|\E \rmA_{\sJ} - \E \widetilde{\rmA}_{\sJ}\|_{\frob}^2 \lesssim N^{-1}\etD_{\max}^2$. Meanwhile, the rank of $\rmA_{\sJ}^{(K)} - \E \widetilde{\rmA}_{\sJ}$ is at most $2K$, then by rank inequality
\begin{align}
    \|\rmA_{\sJ}^{(K)} - \E \widetilde{\rmA}_{\sJ}\|_{\frob}^2 \leq 2K \|\rmA_{\sJ}^{(K)} - \E\widetilde{\rmA}_{\sJ}\|_{2}^2 \leq 4K\big( \|\rmA_{\sJ}^{(K)} - \rmA_{\sJ} \|_2^2 + \| \rmA_{\sJ} -\E\widetilde{\rmA}_{\sJ}\|_{2}^2 \big).
\end{align}
By Eckart–Young–Mirsky \Cref{lem:EYM}, $\|\rmA_{\sJ} - \rmA_{\sJ}^{(K)} \|_2 = \lambda_{K+1} \leq \|\rmA_{\sJ} - \E \widetilde{\rmA}_{\sJ} \|_2$. Also, $\rmA_{\sJ} -\E\widetilde{\rmA}_{\sJ} = \rmA_{\sJ} - \E\rmA_{\sJ} - \textnormal{diag}(\E \ermA_{ii})$. Then desired the argument holds since
\begin{align*}
        &\, \|\rmA_{\sJ}^{(K)} - \E\rmA_{\sJ}\|_{\frob}^2 \lesssim 8K \big( \|\rmA_{\sJ}^{(K)} - \rmA_{\sJ} \|_2^2 + \| \rmA_{\sJ} -\E\widetilde{\rmA}_{\sJ}\|_{2}^2 \big) + 2\etD_{\max}^2/N \\
        \leq &\, 16 K \cdot \| \rmA_{\sJ} -\E\widetilde{\rmA}_{\sJ}\|_{2}^2 + 2\etD_{\max}^2/N \lesssim 32K \cdot \| \rmA_{\sJ} -\E \rmA_{\sJ}\|_{2}^2 + (32K + 2)\etD_{\max}^2/N\\
        \overset{\textnormal{\Cref{thm:regularization}} }{\leq} &\, 32K \const_{\eqref{eqn:concentrate_regularized_A}}^2 \etD_{\max} + (32K + 2)\, \etD_{\max}^2/N \asymp \rho_{N} + \rho_{N}^2/N,
    \end{align*}
    where the last step holds since $\rho_{N} \asymp \etD_{\max}$, as proved in \Cref{lem:maxProb}.
\end{proof}

\begin{proof}[Proof of \Cref{lem:expectedCenterSeparation}]
    First, we discuss the case both $u, v\in \sJ$, and this is a high probability event. By construction, the number of vertices not in $\sJ_1$ is upper bounded by $\lfloor N\exp(-\overline{d})\rfloor \lesssim 1$ since $\rho_{N} \gtrsim \log(N)$ and $\overline{d} \asymp \log(N)$ with high probability by \Cref{lem:maxProb}. Thus vertex $l\in \sJ_1$ with probability at least $1 - N^{-1}$. Also, by taking $\beta_1 = 20 \LM$ in \eqref{eqn:bounded_each_degree}, the number of vertices not in $\sJ_2$ is at most $N e^{-\widetilde{c}\rho_{N}}$ for some constant $\widetilde{c} > 0$, thus $l \in \sJ_2$ with probability at least $1 - e^{-\widetilde{c}\rho_{N}}$. By definition, 
    \begin{align}
        \|(\E \rmA_{\sJ})_{u:} - (\E\rmA_{\sJ})_{v:}\|_2^2 = \sum_{l\in \gV} [(\E \rmA_{\sJ})_{ul} - (\E\rmA_{\sJ})_{vl}]^2.
    \end{align}
    For vertex $l$ with $\ervy_l = t$ for some $t\in[K]$, entry $(\E \rmA_{\sJ})_{ul}$ denotes the expected number of hyperedges containing $u$ and $l$. Hyperedges can be classified into different categories according to the distribution of the remaining $\ell - 2$ nodes among $K$ communities. Assuming that $u\in \gV_j$ and $v\in \gV_k$ with $j\neq k$, the membership distribution of $e$ can be represented by $j \oplus l \oplus \rvw$ for each $\rvw \in \WC{\ell - 2}{K}$ with $N_{\rvw} \asymp N^{-(\ell - 2)}$. For each $l\in \sJ$, the difference between $(\E \rmA_{\sJ})_{ul}$ and $(\E\rmA_{\sJ})_{vl}$ is 
    \begin{align}
        (\E \rmA_{\sJ})_{ul} - (\E\rmA_{\sJ})_{vl} = \sum_{\ell \in \sL} \sum_{\rvw \in \WC{\ell - 2}{K} } N_{\rvw} \cdot (\etQ^{(\ell)}_{j \oplus l \oplus \rvw }  - \etQ^{(\ell)}_{k \oplus l \oplus \rvw} ).
    \end{align}
    At the same time, $N_{\rvw}$ can be replaced by $\overline{N}_{\rvw}$ with high probability by \Cref{lem:sizeDeviation}. Then according to \Cref{ass:expected_center_separation}, it follows that with probability at least $1 - N^{-c}$,
    \begin{align}
        \|(\E \rmA_{\sJ})_{u:} - (\E\rmA_{\sJ})_{v:}\|_2^2 \gtrsim N \big[ (\E \rmA_{\sJ})_{ul} - (\E\rmA_{\sJ})_{vl} \big]^2 \gtrsim N^{-1} \rho_{N}^2.
    \end{align}
    On the other hand, when $u$ and $v$ are from the same cluster, i.e., $u\in \gV_j$ and $v\in \gV_k$ with $j= k$, since $(\E \rmA_{\sJ})_{uu} = (\E \rmA_{\sJ})_{vv} = 0$ due to the fact that self loops are not allowed in model \ref{def:non_uniform_HSBM}, it then follows
    \begin{align}
        \|(\E \rmA_{\sJ})_{u:} - (\E\rmA_{\sJ})_{v:}\|_2^2 =[ (\E \rmA_{\sJ})_{uv}]^2 + [\E \rmA_{\sJ})_{vu}]^2 \lesssim N^{-2}\rho^2_{N}\,.
    \end{align}

    We now turn to the discussion when $u \in \sJ$ but $v\notin \sJ$, and the conclusion holds similarly. As for the case $u \notin \sJ$, $v\notin \sJ$, it follows directly from the definition of $\E \rmA_{\sJ}$ that $\|(\E \rmA_{\sJ})_{u:} - (\E\rmA_{\sJ})_{v:}\|_2^2 = 0$.
\end{proof}

\begin{proof}[Proof of \Cref{lem:centerFar}]
The algorithm inputs $\overline{r} = \lfloor N \log(\overline{d})\rfloor^{-1} \overline{d}^2$. Note that $\overline{d} \asymp \rho_{N}$ with high probability at least $1 - e^{-N^{1/2}}$ by \Cref{lem:maxProb}, then $\overline{r} \asymp r \coloneqq \lfloor N \log(\rho_{N})\rfloor^{-1} \rho_{N}^2$. For now on, we will work on $r$ instead of $\overline{r}$.

Define $\sO_k \coloneqq \{v\in \gV_k \cap \sJ: \| (\rmA_{\sJ}^{(K)})_{v:} - (\E\rmA_{\sJ})_{v:} \|^2_2 > 4r\}$, while $\sI_k$, $\sU_k$, $\sO_k$ and $\ball_{r}(v)$ enjoy the following properties.
\begin{enumerate}[label=(\roman*)]
    \item\label{property:IinBall} For each $k\in [K]$, $\sI_k \subset  \ball_{r}(v)\subset \sU_k$ for any $v\in \sI_k \cap \setS$, since $(\E\rmA_{\sJ})_{v:}=(\E\rmA_{\sJ})_{w:}$ for all $w \in \sI_k$,
    \begin{align}
        \|(\rmA_{\sJ}^{(K)})_{v:} - (\rmA_{\sJ}^{(K)})_{w:}\|^2_2 \leq 2 \|(\rmA_{\sJ}^{(K)})_{v:} - (\E\rmA_{\sJ})_{v:}\|^2_2 + 2 \| (\E\rmA_{\sJ})_{w:} - (\rmA_{\sJ}^{(K)})_{w:}\|^2_2 \leq r.\notag
    \end{align}
    
    \item\label{property:BallinU} For each $k\in [K]$, $\ball_{r}(v)\subset \sU_k$ for any $v\in \sI_k \cap \setS$.
    
    \item\label{property:AllinI} $|\sI_k| \geq \alpha_k N(1 - o(1))$ for each $k \in [K]$, which follows from the fact that almost all nodes in $\gV_k \cap \sJ$ are $r/4$-close to its expected center, since the number of vertices outside $\sI_k$ is bounded by
    \begin{align}
        \frac{|(\gV_k \cap \sJ)\setminus \sI_k|}{N} \leq &\, \frac{1}{N} \cdot \frac{ \underset{v\in (\gV_k \cap \sJ) \setminus \sI_k }{\sum} \|(\rmA_{\sJ}^{(K)})_{v:} - (\E\rmA_{\sJ})_{v:} \|_{2}^2 }{ \underset{v\in (\gV_k \cap \sJ) \setminus \sI_k }{\min}\|(\rmA_{\sJ}^{(K)})_{v:} - (\E\rmA_{\sJ})_{v:} \|^2_2 }\notag
        \leq \frac{1}{N} \frac{\|\rmA_{\sJ}^{(K)} - \E\rmA_{\sJ}\|_{\frob}^2 }{r/4} \\
        \lesssim &\, \frac{1}{N}\frac{\rho_{N} + \rho_{N}^2/N}{[N\log(\rho_{N})]^{-1}\rho_{N}^2} = \frac{\log(\rho_{N})}{\rho_{N}} + \frac{\log(\rho_{N})}{N} = o(1),\notag
    \end{align}
    where the last two equalities hold since $\|\rmA_{\sJ}^{(K)} - \E\rmA_{\sJ}\|_{\frob}^2 \lesssim \rho_{N} + \rho_{N}^2/N$ (\Cref{lem:frobeniusNormBound}), $r = [N\log(\rho_{N})]^{-1}\rho_{N}^2$ and $1 \ll \rho_{N} \ll N$. Thus $|\sI_k| = |\gV_k \cap \sJ | - |(\gV_k \cap \sJ) \setminus \sI_k| \geq \alpha_k N(1 - o(1))$.
    
    \item\label{property:SizeofU} Following from \ref{property:AllinI}, for all $k\in [K]$,
    \begin{align}
        |\sU_k| \leq N - |\cup_{j\neq k} \sU_j| \leq N - \sum_{j \neq k}|\sI_j| \leq \alpha_k N \Big[ 1 + O\Big(\frac{\log(\rho_{N})}{\rho_{N}}\Big) \Big] = \alpha_k N(1 + o(1)).\notag
    \end{align}
    
    \item\label{property:LeastinO} $|\cup_{k=1}^{K}\sO_k| \leq o(N)$, since most nodes are in $\cup_{k=1}^{K}\sI_k$, and
    \begin{align}
        |\cup_{k=1}^{K}\sO_k| \leq \frac{\sum_{v\in \sJ} \|( \rmA_{\sJ}^{(K)})_{v:} - (\E\rmA_{\sJ})_{v:} \|_{2}^2 }{ \underset{v\in (\cup_{k=1}^{K}\sO_k) }{\min}\|(\rmA_{\sJ}^{(K)})_{v:} - (\E\rmA_{\sJ})_{v:} \|^2_2 }\leq \frac{\|\rmA_{\sJ}^{(K)} - \E\rmA_{\sJ}\|_{\frob}^2}{4r} = o(N)\notag
    \end{align}
    
    \item\label{property:BallSizeinO} $|\ball_r(v)| = o(N)$ for all $v\in \sO_k \cap \setS$. First, $\ball_{r}(v) \cap \sU_k = \emptyset$ because the distance between $v\in \sO_k$ and $w\in \sI_k$ is larger than $r$, due to the facts $(\E\rmA_{\sJ})_{v:} = (\E\rmA_{\sJ})_{w:}$, $(x - y)^2 \geq \frac{1}{2}(x - z)^2 - (y-z)^2$, and
    \begin{align}
       \|(\rmA_{\sJ}^{(K)})_{v:} - (\rmA_{\sJ}^{(K)})_{w:}\|^2_2 \geq \frac{1}{2}\|(\rmA_{\sJ}^{(K)})_{v:} - (\E\rmA_{\sJ})_{v:}\|^2_2 - \|(\rmA_{\sJ}^{(K)})_{w:} - (\E\rmA_{\sJ})_{w:}\|^2_2 > r.\notag
    \end{align}
    Then as a consequence of \ref{property:LeastinO},  $|\ball_r(v)| \leq |\cup_{k=1}^{K}\sO_k| = o(N)$.
    
    \item\label{property:DisjointBalls} $\ball_{r}(u)\cap \ball_{r}(v) = \emptyset$ for any $u\in \sU_j$ and $v\in \sU_k$ with $j\neq k$. This follows the fact that $|\ball_{r}(v)\cap \sU_j | = 0$ when $|\ball_{r}(v)\cap \sU_k |\geq 1$, since the distance between expected centers $\| (\E \rmA_{\sJ})_{v:} - (\E\rmA_{\sJ})_{u:}\|^2_2 \gtrsim N^{-1}\rho_{N}^2$ (\Cref{lem:expectedCenterSeparation}) is much larger than the radius $4r \asymp [N\log(\rho_{N})]^{-1}\rho_{N}^2$ where $\rho_{N} = \omega(1)$.
\end{enumerate}
We now prove those three arguments using the properties above.

{\noindent} (a). For each $k\in [K]$, the probability that a randomly selected node $s\in \gV$ does not belong to $\sI_k$ is 
    \begin{align}
    1 - \frac{|\sI_k|}{N} \overset{\ref{property:AllinI}}{\leq} 1 - \alpha_k(1 -o(1)) = \underset{{j\neq k}}{\sum}\alpha_j + o(1) < 1,
    \end{align}
    where $1 > \alpha_1 \geq \cdots \geq \alpha_K > 0$ are some constants by \Cref{def:uniform_HSBM}, then the probability that there exists at least one node $s\in \setS \cap \sI_k$ is
    \begin{align}
        1 - \Big( 1 - \frac{|\sI_k|}{N}\Big)^{\log(N)} \geq 1 - \Big(\sum_{j\neq k}\alpha_j + o(1) \Big)^{\log(N)} = 1 - N^{-\log( \frac{1}{\sum_{j\neq k}\alpha_j + o(1)})},
    \end{align}
    A simple union bound completes the argument. 

{\noindent} (b). We prove part (b) by induction.
\begin{enumerate}
    \item[(1)] $k=1$. By part $(a)$, there exists a node $v_1$ in $\sI_1 \cap \setS$, satisfying $\ball_r(v_1) \supset \sI_1$ due to \ref{property:IinBall}. By algorithm procedure, $|\ball_r(s_1)|$ is at least $|\ball_r(v_1)|$, then $|\ball_r(s_1)| \geq |\ball_r(v_1)| \geq |\sI_1| \geq \alpha_1 N(1 - o(1))$ thanks to \ref{property:AllinI}. To prove $s_1 \in \sU_1$, one should verify the failure of other possibilities, namely, for $v\notin \sU_1$, the cardinality of $\ball_r(v)$ is too small to make $v$ selected as a center.
    \begin{itemize}
        \item For any $v\in \sO_k \cap \setS$ for each $k\in [K]$, $|\ball_r(v)| = o(N) < |\ball_r(s_1)|$ by property \ref{property:BallSizeinO}.
        \item For any $v\in \sU_k \cap \setS$ with $k \neq 1$, $|\ball_r(v)| \overset{\ref{property:BallinU}}{\leq} |\sU_k| \overset{\ref{property:SizeofU}}{\leq} \alpha_k N(1 + o(1)) < \alpha_1 N$ since $\alpha_1 > \alpha_k$.
    \end{itemize}
    
    \item[(2)] $2 \leq k \leq K$. Suppose $s_j \in \sU_j$ for all $j\in [k-1]$. Similarly by \ref{property:IinBall}, there exists a node $v_k\in \sU_k \cap \setS$ satisfying $\ball_r(v_k) \supset \sU_k$. By induction hypothesis $|\widehat{\gV}_{j}^{(0)}| = |\ball_{r}(v_j)\setminus (\cup_{l=1}^{j-1}\widehat{\gV}_{l}^{(0)}) | \geq \alpha_j N(1 - o(1))$ for $j\in [k-1]$. Then by \ref{property:DisjointBalls}, $\ball_r(v_k)\setminus (\bigcup_{l=1}^{k-1}\widehat{\gV}_{l}^{(0)}) = \ball_r(v_k)$, and $|\ball_r(s_k)| \geq |\ball_r(v_k)| \geq |\sU_k| \geq \alpha_k N(1 - o(1))$ thanks to \ref{property:AllinI}. To prove $s_k\in \sU_k$, one should verify the failure of other possibilities, namely, for $v\notin \sU_k$, the cardinality of so obtained set $\ball_r(v)\setminus (\cup_{l=1}^{k-1}\widehat{\gV}_{l}^{(0)})$ is too small to make $v$ selected as a center.
\begin{itemize}
    \item For any $v\in \sO_k \cap \setS$ for all $k\in [K]$, $|\ball_r(v)| = o(N) < |\ball_r(s_k)|$ by property \ref{property:BallSizeinO}.
    
    \item For any $v\in \sU_j \cap \setS$ where $j < k$, note that $|\ball_r(s_j)| \geq \alpha_j N(1 - o(1))$ by induction, then
    \begin{align}
        |\ball_r(v)\setminus (\cup_{l=1}^{k-1}\widehat{\gV}_{l}^{(0)}) | &\, \overset{ \ref{property:DisjointBalls}}{\leq} |\ball_r(v)\setminus \widehat{\gV}_{j}^{(0)}| \overset{\ref{property:BallinU}}{\leq} |\sU_j| - |\ball_{r}(s_j)| \overset{\ref{property:SizeofU}}{=} o(N) \\
        &\,\,\, < |\ball_r(s_k)| = \alpha_k N (1 + o(1)).
    \end{align}
    
    \item For any $v\in \sU_j \cap \setS$ where $j > k$,
    \begin{align}
        |\ball_r(v)\setminus (\cup_{l=1}^{k-1}\widehat{\gV}_{l}^{(0)}) | &\, \overset{\ref{property:DisjointBalls}}{\leq} |\ball_r(v)| \overset{\ref{property:BallinU}}{\leq} |\sU_j|  \overset{\ref{property:SizeofU}}{=} \alpha_j N(1 + o(1)) \\
        &\,\,\,< \alpha_k N(1 + o(1)) =|\ball_r(s_k)|,
    \end{align}
    where the last inequality holds since $\alpha_k > \alpha_j$.
\end{itemize}
\end{enumerate}

{\noindent}(c). Consequently, $\{s_k\}_{k=1}^{K}$ are pairwise far from each other, i.e., for any pair $j\neq k$
\begin{subequations}
\begin{align*}
    &\, \|(\rmA_{\sJ}^{(K)})_{s_j :} -  (\rmA_{\sJ}^{(K)})_{s_k :} \|_2^2
    \geq \frac{1}{2}\|(\rmA_{\sJ}^{(K)})_{s_j :} -  (\E\rmA_{\sJ})_{s_k :} \|_2^2 - \|(\E\rmA_{\sJ})_{s_k :} -  (\rmA_{\sJ}^{(K)})_{s_k :} \|_2^2\\
    \geq&\, \frac{1}{2}\left(\frac{1}{2} \|(\E\rmA_{\sJ})_{s_j : } -  (\E\rmA_{\sJ})_{s_k :} \|_2^2 - \|(\rmA_{\sJ}^{(K)})_{s_j :} -  (\E\rmA_{\sJ})_{s_j :} \|_2^2\right) - \|(\E\rmA_{\sJ})_{s_k :} -  (\rmA_{\sJ}^{(K)})_{s_k :} \|_2^2 \gtrsim \frac{\rho_{N}^2}{N},
\end{align*}
\end{subequations}
where the last equality holds by since \Cref{lem:expectedCenterSeparation} $\| (\E \rmA_{\sJ})_{s_j :} - (\E\rmA_{\sJ})_{s_k :}\|^2_2 \gtrsim N^{-1} \rho_{N}^2 \gg  [N\log(\rho_{N})]^{-1} \rho_{N}^2 \asymp 4r \geq \|(\E\rmA_{\sJ})_{s_j :} -  (\rmA_{\sJ})_{s_j :} \|_2^2$ with $\rho_{N} = \omega(1)$ by the property that $s_j \in \sU_j$ for each $j\in [K]$.
\end{proof}

\begin{proof}[Proof of \Cref{lem:misclassifiedFarCenter}]
Note that $(\E\rmA_{\sJ})_{v:} = (\E\rmA_{\sJ})_{s_i :}$ when $v\in \gV_i \cap \sJ$ where $s_i$ is the center of $\widehat{\gV}_{i}^{(0)}$ obtained in \Cref{alg:first_loop}. There are two scenarios where $v\in \gV_i\cap \sJ$ would be misclassified into $\widehat{\gV}_j^{(0)}$ ($j \neq i$).
\begin{enumerate}
    \item In \Cref{alg:first_loop} of \Cref{alg:spectral_initialization}, if $v$ is $r$-close to center of $\widehat{\gV}_j^{(0)}$, i.e., $\|(\rmA_{\sJ}^{(K)})_{v:} - (\rmA_{\sJ}^{(K)})_{s_j :}\|_2^2 \leq r$, then
\begin{subequations}
    \begin{align*}
        &\, \|(\rmA_{\sJ}^{(K)})_{v:} - (\E\rmA_{\sJ})_{s_i :} \|_2^2 \geq \frac{1}{2} \|(\E\rmA_{\sJ})_{s_i :} - (\E\rmA_{\sJ})_{s_j :}\|_2^2 - \|(\rmA_{\sJ}^{(K)})_{v:} - (\E\rmA_{\sJ})_{s_j :}\|_2^2\\
        \geq &\, \frac{1}{2} \|(\E\rmA_{\sJ})_{s_i :} - (\E\rmA_{\sJ})_{s_j :}\|_2^2 - 2\|(\rmA_{\sJ}^{(K)})_{v:} - (\rmA_{\sJ}^{(K)})_{s_j :}\|_2^2 - 2\|(\rmA_{\sJ}^{(K)})_{s_j :} - (\E\rmA_{\sJ})_{s_j :}\|_2^2 \gtrsim \frac{\rho_{N}^2}{N},
    \end{align*}
\end{subequations}
where $\|(\E\rmA_{\sJ})_{s_i :} - (\E\rmA_{\sJ})_{s_j :}\|_2^2 \gtrsim N^{-1}\rho_{N}^2$ by \Cref{lem:expectedCenterSeparation}, and $\|(\rmA_{\sJ}^{(K)})_{s_j :} -  (\E\rmA_{\sJ})_{s_j :} \|_2^2 \leq 4r \asymp [N\log(\rho_{N})]^{-1}\rho_{N}^2 \ll \rho_{N}^2/N$ since $\rho_{N} = \omega(1)$ and $s_j \in \sU_j$ as proved in \Cref{lem:centerFar}.

    \item In \Cref{alg:second_loop} of \Cref{alg:spectral_initialization}, when $v$ is closer to the center of $\widehat{\gV}_j^{(0)}$ than $\widehat{\gV}_i^{(0)}$, i.e., $\|(\rmA_{\sJ}^{(K)})_{s_j :} - (\rmA_{\sJ}^{(K)})_{v:}\|_2^2 \leq \|(\rmA_{\sJ}^{(K)})_{s_i :} - (\rmA_{\sJ}^{(K)})_{v:}\|_2^2$. One can verify that $v$ is far from its expected center, namely,
\begin{align*}
    \|(\rmA_{\sJ}^{(K)})_{v:} - (\E\rmA_{\sJ})_{s_i :}\|_2^2 \geq \frac{1}{2}\|(\rmA_{\sJ}^{(K)})_{s_i :} - (\rmA_{\sJ}^{(K)})_{v:}\|_2^2 - \| (\E\rmA_{\sJ})_{s_i :} - (\rmA_{\sJ}^{(K)})_{s_i :}\|_2^2 \gtrsim \frac{\rho_{N}^2}{N},
\end{align*}
where $\| (\E\rmA_{\sJ})_{s_i :} - (\rmA_{\sJ}^{(K)})_{s_i :}\|_2^2  \leq 4r \asymp [N\log(\rho_{N})]^{-1}\rho_{N}^2$ and
\begin{align*}
    \|(\rmA_{\sJ}^{(K)})_{s_i :} - (\rmA_{\sJ}^{(K)})_{v:}\|_2^2 \geq&\, \frac{1}{2} \Big( \|(\rmA_{\sJ}^{(K)})_{s_i :} - (\rmA_{\sJ}^{(K)})_{v:}\|_2^2 + \|(\rmA_{\sJ}^{(K)})_{s_j :} - (\rmA_{\sJ}^{(K)})_{v:}\|_2^2 \Big) \\
    \geq&\, \frac{1}{4} \|(\rmA_{\sJ}^{(K)})_{s_i :} - (\rmA_{\sJ}^{(K)})_{s_j :}\|_2^2 \overset{\textnormal{\Cref{lem:centerFar} (c)}}{\gtrsim} \frac{\rho_{N}^2}{N},
\end{align*}
where the last line holds since $(x - y)^2 \leq 2(x - z)^2 + 2(y - z)^2$.
\end{enumerate}
\end{proof}

\section{Proofs in Section \ref{sec:AgnosticRefinement}}\label{app:AgnosticRefinement}

\begin{proof}[Proof of \Cref{lem:hatQmw_approx}]
At step $0$, let $\widehat{\rD}^{(0)}_{\rvw}$ denote the number of $\ell$-hyperedges with 
\[
\widehat{\rvy}^{(0)}(e) = \rvw = (\ervw_1, \ldots, \ervw_K),
\]
meaning that $\ervw_k$ nodes in $\widehat{\gV}_{k}^{(0)}$ for each $k\in [K]$. Recall $\widehat{\etQ}^{(\ell)}_{\rvw} = \widehat{\rD}^{(0)}_{\rvw} /\widehat{N}^{(0)}_{\rvw}$ in \eqref{eqn:hat_Qlw}. Consequently, for any $\rvw \in \WC{\ell}{K}$, we have
\begin{align}
    |\widehat{\etQ}^{(\ell)}_{\rvw} - \etQ^{(\ell)}_{\rvw} | = \frac{\big| \widehat{\rD}^{(0)}_{\rvw} - \etQ^{(\ell)}_{\rvw} \cdot \widehat{N}^{(0)}_{\rvw} \big| }{\widehat{N}^{(0)}_{\rvw}} \leq \frac{1}{\widehat{N}^{(0)}_{\rvw}} \Big( \big| \widehat{\rD}^{(0)}_{\rvw} - \E \widehat{\rD}^{(0)}_{\rvw} \big| + \big|\E \widehat{\rD}^{(0)}_{\rvw} - \etQ^{(\ell)}_{\rvw} \cdot \widehat{N}^{(0)}_{\rvw} \big| \Big).
\end{align}
We will bound the two terms in the numerator separately. For the first term $\big| \widehat{\rD}^{(0)}_{\rvw} - \E \widehat{\rD}^{(0)}_{\rvw} \big|$, denote $\etQ_{\min}^{(\ell)} \coloneqq \min_{\rvw \in \WC{\ell}{K}} \etQ_{\rvw}^{(\ell)}$, $\etQ_{\max}^{(\ell)} \coloneqq \max_{\rvw \in \WC{\ell}{K}} \etQ_{\rvw}^{(\ell)}$ for each $\ell \in \sL$. Using \eqref{eqn:edge_density}, we write $\etQ_{\min}^{(\ell)} = \etP_{\min}^{(\ell)} q_{\min}^{(\ell)}/\binom{N - 1}{\ell - 1}$, $\etQ_{\max}^{(\ell)} = \etP_{\max}^{(\ell)} q_{\max}^{(\ell)}/\binom{N - 1}{\ell - 1}$ where $\etP_{\min}^{(\ell)} \asymp \etP_{\max}^{(\ell)} \asymp 1$ and $q_{\min}^{(\ell)} \lesssim q_{\max}^{(\ell)}$. Since $\widehat{N}^{(0)}_{\rvw} \asymp N^{\ell}$, we then have
\begin{align}
     N q_{\min}^{(\ell)} \asymp \widehat{N}^{(0)}_{\rvw} \cdot \etQ_{\max}^{(\ell)} \leq  \E \widehat{\rD}^{(0)}_{\rvw} \leq \widehat{N}^{(0)}_{\rvw} \cdot \etQ_{\max}^{(\ell)} \asymp N q_{\max}^{(\ell)}
\end{align}
Let $\widehat{\gV}$ denote the collection of partitions $\{\widehat{\gV}_k\}_{k=1}^{K}$ such that 
\begin{align}
    \Big|\bigcup_{k=1}^{K} (\gV_k \setminus \widehat{\gV}_k) \Big| \leq \mismatch_{N}\cdot N, \quad \forall\, \{\widehat{\gV}_k\}_{k=1}^{K} \in \widehat{\gV},
\end{align}
where $\mismatch_{N}$ is the mismatch ratio defined in \eqref{eqn:misratio} and $\mismatch_{N}\lesssim (\rho_{N})^{-1}$ with probability at least $1 - N^{-c}$ for some constant $c> 0$ by \eqref{eqn:weak_consistency_agnostic}. For any $\{\widehat{\gV}^{(0)}_k\}_{k=1}^{K} \in \widehat{\gV}$, thanks to edgewise independence, Bernstein (\Cref{lem:Bernstein}) shows that
\begin{align*}
    \P \big( \big| \widehat{\rD}^{(0)}_{\rvw} - \E \widehat{\rD}^{(0)}_{\rvw} \big| \geq N\log q_{\max}^{(\ell)} \big) \leq &\, 2 \exp \bigg( - \frac{ \big(N\log q_{\max}^{(\ell)}\big)^2/2}{Nq_{\max}^{(\ell)} + N/3\cdot\log q_{\max}^{(\ell)}} \bigg) \\
    =&\, 2\exp\big(-N(\log q_{\max}^{(\ell)})^2/(2q_{\max}^{(\ell)}) \big),
\end{align*}
Meanwhile, the upper bound on the cardinality of $\widehat{\gV}$ can be computed as follows: first $\mismatch_{N} N$ many vertices (wrongly classified) are selected, then each of them can be assigned to $K$ different communities. Thus by \Cref{lem:stirling}, for some universal constant $\const > 0$, the following holds
\begin{align}
    |\widehat{\gV}|\leq \binom{N}{\mismatch_{N} N} K^{\mismatch_{N} N} \leq \Big( \frac{eN K}{\mismatch_{N} N} \Big)^{\mismatch_{N} N} = \exp\Big( \mismatch_{N} N\log\Big(\frac{eK}{\mismatch_{N}} \Big)\Big) \lesssim \exp \Big(\const \frac{N\log\rho_{N}}{\rho_{N}} \Big)
\end{align}
Note that $\rho_{N} \asymp \sum_{\ell \in \sL} q^{(\ell)}_{\max}$ and the function $h(x) = x/ \log(x)$ increases when $x \geq e$, then a simple union bound shows that the following holds with probability at least $1 - N^{-c}$
\begin{align}
    &\, \P \Big( \exists \{\widehat{\gV}^{(0)}_k\}_{k=1}^{K} \in \widehat{\gV}\,\, \textnormal{s.t.}\,  \big| \widehat{\rD}^{(0)}_{\rvw} - \E \widehat{\rD}^{(0)}_{\rvw} \big| \geq \cdot N\log q_{\max}^{(\ell)} \Big) \\
    \leq &\, |\widehat{\gV}| \binom{K}{2} \cdot  2\exp \Big( - N \frac{[\log q_{\max}^{(\ell)}]^2}{2q_{\max}^{(\ell)}} \Big) \lesssim \exp \Big( - N\frac{\log q_{\max}^{(\ell)}}{q_{\max}^{(\ell)}} \Big).
\end{align}

Consider the second term $\E \widehat{\rD}^{(0)}_{\rvw} = \sum_{\widehat{\rvy}^{(0)}(e) = \rvw} \etQ^{(\ell)}_{\rvy(e)}$. For hyperedge $e \subset \gV \setminus \widehat{\gN}^{(0)}$, i.e., $\widehat{\rvy}^{(0)}(e) = \rvw = \rvy(e)$, the contribution of such $e$ inside the $|\E \widehat{\rD}^{(0)}_{\rvw} - \etQ^{(\ell)}_{\rvw} \cdot \widehat{N}^{(0)}_{\rvw}|$ cancels out. Thus it suffices to consider contributions from $e$ which contains at least one node from $\widehat{\gN}^{(0)}$, and the number of such hyperedges is at most $\mismatch_{N} N \cdot\binom{N - 1}{\ell - 1}$. According to \eqref{eqn:weak_consistency_agnostic}, $\mismatch_{N} \lesssim (\rho_{N})^{-1}$ where $\rho_{N} \asymp \sum_{\ell \in \sL}q_{\max}^{(\ell)}$, then the following holds with probability at least $1 - N^{-c}$ for some $c > 0$
    \begin{align}
        |\E  \widehat{\rD}^{(0)}_{\rvw} - \etQ^{(\ell)}_{\rvw} \cdot \widehat{N}^{(0)}_{\rvw} | \leq \etQ_{\max}^{(\ell)} \cdot \mismatch_{N} N \binom{N - 1}{\ell - 1} \lesssim q_{\max}^{(\ell)} N/\rho_{N} \lesssim N.
    \end{align}
Therefore, with probability at least $1 - \exp(- N\log q_{\max}^{(\ell)} /q_{\max}^{(\ell)} )$,
\begin{align}
    |\widehat{\etQ}^{(\ell)}_{\rvw} - \etQ^{(\ell)}_{\rvw} | \lesssim \frac{ N\log q_{\max}^{(\ell)}  + N}{\widehat{N}^{(0)}_{\rvw}} \lesssim \frac{\log q_{\max}^{(\ell)}}{N^{\ell - 1}}.
\end{align}
The second part follows easily since $\etQ^{(\ell)}_{\rvw} \asymp q_{N} N^{-m+1}$ for each $\ell \in \sL$ when under \Cref{ass:prob_ratio_bound}.
\end{proof}

\begin{lemma}\label{lem:notG1G2}
Recall that \eqref{eqn:optimality_condition} holds for some absolute constant $\epsilon >0$, where $\kappa_{N} \log(N) \to \infty$ as $N \to \infty$. Then the number of nodes which don't satisfy either \eqref{eqn:G1} or \eqref{eqn:G2} is at most $N^{1 - \kappa_{N}}/3$ with probability at least $1 - 6e^{- \epsilon\kappa_{N} \log(N)}$.
\end{lemma}
\begin{proof}[Proof of \Cref{lem:notG1G2}]
\textbf{Probability of satisfying} \eqref{eqn:G1}.
For each $\rvw\in \WC{\ell - 1}{K}$, $\rD^{(\ell)}_{v, \rvw} \big\vert_{\rY_{v} = k} \sim \textnormal{Bin}(N_{\rvw}, \etQ^{(\ell)}_{k \oplus \rvw})$ where $N_{\rvw}$ is defined in \eqref{eqn:nbw}. Denote $\etQ^{(\ell)}_{\max} = \max_{\rvw \in \WC{\ell}{K}}\etQ^{(\ell)}_{\rvw}$, then 
\begin{align}
    \sum_{\ell \in \sL} \sum_{\rvw \in\WC{\ell - 1}{K}} N_{\rvw} \etQ^{(\ell)}_{k \oplus \rvw} \lesssim \sum_{\ell \in \sL} \sum_{\rvw \in\WC{\ell - 1}{K}} N_{\rvw}\etQ^{(\ell)}_{\max} \lesssim \sum_{\ell \in \sL} \sum_{\rvw \in\WC{\ell - 1}{K}} \binom{N - 1}{\ell - 1}\cdot \etQ^{(\ell)}_{\max} = \etD_{\max}.\notag
\end{align}
Thanks to hyperedges independence and Markov inequality, for each $v\in \gG$, we have
\begin{align*}
    &\, \P(\rD_v > 10\etD_{\max}) = \P \Big( \sum_{\ell \in \sL}\sum_{\rvw \in\WC{\ell - 1}{K}} \rD^{(\ell)}_{v, \rvw} > 10 \etD_{\max} \Big) \\\
    \leq &\, \inf_{\theta >0} \Big[\exp(-10 \theta \etD_{\max}) \prod_{\ell \in \sL} \prod_{\rvw \in\WC{\ell - 1}{K}} \E \exp(\theta \cdot \rD^{(\ell)}_{v, \rvw}) \Big]\\
    \leq &\, \inf_{\theta >0} \exp(-10 \theta \etD_{\max}) \prod_{\ell \in \sL} \prod_{\rvw \in\WC{\ell - 1}{K}} \exp\big[ (1 - \etQ^{(\ell)}_{k \oplus \rvw}) + \etQ^{(\ell)}_{k \oplus \rvw}e^{\theta} \big]^{N_{\rvw}}, \quad (1 + x \leq e^{x})\\
    \leq &\, \inf_{\theta >0} \exp( -10 \theta \etD_{\max}) \exp\Big( (e^{\theta} - 1) \sum_{\ell \in \sL} \sum_{\rvw \in\WC{\ell - 1}{K}} N_{\rvw} \etQ^{(\ell)}_{k \oplus \rvw}  \Big)\\
    \leq &\, \inf_{\theta >0} \exp \big( (e^{\theta} -1 -10 \theta) \cdot \etD_{\max} \big) \leq e^{-10 \etD_{\max}},
\end{align*}
where the last line holds by taking $\theta = 2$, and $\P(\rD_v \leq 10 \etD_{\max}) \leq 1 - e^{-10 \etD_{\max}} \to 1$ since $\etD_{\max} = \omega(1)$.


\textbf{Probability of satisfying} \eqref{eqn:G2}
Note that $\D_{\textnormal{GKL}}(j, k) = (1 + o(1))\D_{\textnormal{GCH}}(j, k)\cdot q_{N}$ by \Cref{lem:equivalenveDivergence}. We claim that it suffices to prove the following holds with probability at least $1 - 2 e^{-\frac{1}{2}\log(N/s)}$,
\begin{align}
    \sum_{\ell \in \sL} \sum_{\rvw \in \WC{\ell - 1}{K} } N_{\rvw}\cdot \D_{\mathrm{KL}}(\mu^{(\ell)}_{v, \rvw}\parallel \etQ_{k \oplus \rvw}^{(\ell)}) \leq \Big(1 - \frac{\log(N)}{\sqrt{N}} \Big) \cdot \D_{\textnormal{GKL}}(j, k) -\frac{\etD_{\max}}{\log q_{N}},\label{eqn:KLUpperBound}
\end{align}
where the proof is deferred to \Cref{lem:ProbG1andG2}. Assuming \eqref{eqn:KLUpperBound}, subsequent logic is indicated as follows.
\begin{enumerate}
    \item Remind that $(1 - \log(N)/\sqrt{N} )^{\ell - 1} \leq N_{\rvw}/ \overline{N}_{\rvw}$ by \Cref{lem:sizeDeviation}, then by \eqref{eqn:KLUpperBound},
\begin{align}
    &\, \sum_{\ell \in \sL} \sum_{\rvw \in \WC{\ell - 1}{K} } \overline{N}_{\rvw}\cdot \D_{\textnormal{KL}}(\mu^{(\ell)}_{v, \rvw}\parallel \etQ_{k \oplus \rvw}^{(\ell)}) \notag \\
    \leq &\, \Big(1 - \frac{\log(N)}{\sqrt{N}} \Big)^{-1} \sum_{\ell \in \sL} \sum_{\rvw \in \WC{\ell - 1}{K} } N_{\rvw}\cdot \D_{\textnormal{KL}}(\mu^{(\ell)}_{v, \rvw}\parallel \etQ_{k \oplus \rvw}^{(\ell)}) < \D_{\textnormal{GKL}}(j, k).\notag
\end{align}

    \item Then by \eqref{eqn:GeneralizedKL}, $\sum_{\ell \in \sL} \sum_{\rvw \in \WC{\ell - 1}{K} } \overline{N}_{\rvw}\cdot \D_{\textnormal{KL}}(\mu^{(\ell)}_{v, \rvw}\parallel \etQ_{j \oplus \rvw}^{(\ell)}) \geq \D_{\textnormal{GKL}}(j, k)$, since part (1) and
    \begin{align}
       \sum_{\ell \in \sL} \sum_{\rvw \in \WC{\ell - 1}{K} } \overline{N}_{\rvw}\cdot \max\big\{ \D_{\textnormal{KL}}(\mu^{(\ell)}_{v, \rvw} \parallel \etQ_{j \oplus \rvw}^{(\ell)}),\, \D_{\textnormal{KL}}(\mu^{(\ell)}_{v, \rvw}\parallel \etQ_{k \oplus \rvw}^{(\ell)}) \big\} \geq \D_{\textnormal{GKL}}(j, k). \notag
    \end{align}

    \item Then we use part (2) and $(1 - \log(N)/\sqrt{N} )^{\ell - 1} \leq N_{\rvw}/ \overline{N}_{\rvw}$ again,
    \begin{align}
        &\, \sum_{\ell \in \sL} \sum_{\rvw \in \WC{\ell - 1}{K} } N_{\rvw}\cdot \D_{\textnormal{KL}}(\mu^{(\ell)}_{v, \rvw}\parallel \etQ_{j \oplus \rvw}^{(\ell)}) \notag \\
        \geq&\, \Big(1 - \frac{\log(N)}{\sqrt{N}} \Big)\cdot \sum_{\ell \in \sL} \sum_{\rvw \in \WC{\ell - 1}{K} } \overline{N}_{\rvw}\cdot \D_{\textnormal{KL}}(\mu^{(\ell)}_{v, \rvw}\parallel \etQ_{j \oplus \rvw}^{(\ell)}) \geq \Big(1 - \frac{\log(N)}{\sqrt{N}} \Big) \cdot \D_{\textnormal{GKL}}(j, k). \notag
    \end{align}

    \item Therefore by part (3) and \eqref{eqn:KLUpperBound}, $\sum_{\ell \in \sL} \, \sum_{ \rvw \in \WC{\ell - 1}{K} } N_{\rvw} \cdot \big[ \D_{\rm{KL}}(\mu^{(\ell)}_{v, \rvw}\parallel \etQ_{j \oplus \rvw}^{(\ell)}) - \D_{\rm{KL}}(\mu^{(\ell)}_{v, \rvw}\parallel \etQ_{k \oplus \rvw}^{(\ell)}) \big] \geq \etD_{\max}/\log q_{N}$, which proves \eqref{eqn:G2}.
\end{enumerate}
Let $s = N^{1 - \kappa_{N}}$. According to \Cref{lem:ProbG1andG2}, the probability that there are more than $s/3$ nodes in $\gV$ not satisfying either \eqref{eqn:G1} or \eqref{eqn:G2} is upper bounded by Markov inequality as follows,
\begin{align}
    &\, \frac{\E[\# \textnormal{ of nodes that do not satisfy either \eqref{eqn:G1} or \eqref{eqn:G2}}] }{s/3} \leq \frac{N\cdot 2e^{-\D_{\rm{GCH}} \cdot q_{N}}}{s/3} \notag\\
    \leq &\, 6\exp(- \D_{\rm{GCH}} \cdot q_{N} + \log(N/s)) \leq 6\exp(- \epsilon \cdot \kappa_{N}\cdot \log(N)) \to 0,\notag
\end{align}
where the last inequality holds since \eqref{eqn:optimality_condition} and $\kappa_{N}\cdot \log(N) \to \infty$.
\end{proof}

\begin{lemma}\label{lem:ProbG1andG2}
    Under \Cref{ass:prob_ratio_bound}, with probability at least $ 1 - 2e^{-\D_{\rm{GCH}} \cdot q_{N}}$, vertex $v$ satisfies \eqref{eqn:G1} and \eqref{eqn:G2} simultaneously.
\end{lemma}
\begin{proof}[Proof of \Cref{lem:ProbG1andG2}]
    Let $\sX$ denote the set of degree profiles satisfying \eqref{eqn:G1}, i.e.,
    \begin{align}
        \sX = \Big\{ \{\rvx^{(\ell)}\}_{\ell \in \sL} \Big|&\,\,\, \rvx^{(\ell)} \in \N^{|\WC{\ell - 1}{K}|}, \textnormal{ s.t. } x_{\rvw}^{(\ell)} = \omega(1) \textnormal{ for each } m, \rvw, \textnormal{ and} \sum_{\ell \in \sL} \sum_{\rvw \in \WC{\ell - 1}{K}} x_{\rvw}^{(\ell)} \leq 10 \etD_{\max}, \Big\}.\notag
    \end{align}
Since the cardinality $|\WC{\ell - 1}{K}| = \binom{\ell+K-1}{K-1}$, for simplicity, we denote
\begin{align}
    T = \sum_{\ell \in \sL} \binom{\ell+K-1}{K-1} = O(1).
\end{align}
Let $\{\rD^{(\ell)}_{v, \rvw} = x^{(\ell)}_{\rvw} \}^{\ell}_{\rvw}$ denote the event $\rD^{(\ell)}_{v, \rvw} = x^{(\ell)}_{\rvw}$ for each $\ell \in \sL$, $\rvw \in \WC{\ell - 1}{K}$, then
\begin{align}
    \P( \{\rD^{(\ell)}_{v, \rvw} = x^{(\ell)}_{\rvw} \}^{\ell}_{\rvw} ) = \prod_{\ell \in \sL} \prod_{\rvw \in \WC{\ell - 1}{K} } \binom{N_{\rvw}}{x^{(\ell)}_{\rvw}} \cdot[ \etQ^{(\ell)}_{k \oplus \rvw}]^{x^{(\ell)}_{\rvw}} \cdot (1 - \etQ^{(\ell)}_{k \oplus \rvw})^{ N_{\rvw} - x^{(\ell)}_{\rvw}}.
\end{align}
Note that $\mu^{(\ell)}_{v, \rvw}\coloneqq \rD^{(\ell)}_{v, \rvw}/ N_{\rvw}$ where $N_{\rvw} \coloneqq \prod_{k=1}^{K}\binom{|\gV_{k}|}{\ervw_{k}} \asymp N^{\ell - 1}$. For any $\const > 0$, by Markov inequality, we have
\begin{align}
    &\, \P\bigg( \Big\{ \sum_{\ell \in \sL} \sum_{\rvw \in \WC{\ell - 1}{K} } N_{\rvw}\cdot \D_{\textnormal{KL}}(\mu^{(\ell)}_{v, \rvw}\parallel \etQ_{k \oplus \rvw}^{(\ell)}) > \const \Big\} \cap \{\rD_v \leq 10 \etD_{\max} \} \bigg)\notag\\
    =&\,  \sum_{\{\rvx^{(\ell)}\}\in \sX } \P( \{\rD^{(\ell)}_{v, \rvw} = x^{(\ell)}_{\rvw} \}^{\ell}_{\rvw} ) \cdot \P\bigg( \sum_{\ell \in \sL} \sum_{\rvw \in \WC{\ell - 1}{K} } N_{\rvw}\cdot \D_{\textnormal{KL}}( \mu^{(\ell)}_{v, \rvw} \parallel \etQ_{k \oplus \rvw}^{(\ell)}) > \const \Big \vert \{\rD^{(\ell)}_{v, \rvw} = x^{(\ell)}_{\rvw} \}^{\ell}_{\rvw} \bigg)\notag\\
    \leq &\, \sum_{\{\rvx^{(\ell)}\}\in \sX } \P( \{\rD^{(\ell)}_{v, \rvw} = x^{(\ell)}_{\rvw} \}^{\ell}_{\rvw} ) \cdot e^{-\const} \cdot \E\Big( \exp \Big [\sum_{\ell \in \sL} \sum_{\rvw \in \WC{\ell - 1}{K} }N_{\rvw}\cdot \D_{\textnormal{KL}}( \mu^{(\ell)}_{v, \rvw} \parallel \etQ_{k \oplus \rvw}^{(\ell)}) \Big]\Big \vert \{\rD^{(\ell)}_{v, \rvw} = x^{(\ell)}_{\rvw} \}^{\ell}_{\rvw} \Big) \notag\\
    \leq &\,  e^{-\const} \sum_{\{\rvx^{(\ell)}\}\in \sX } \prod_{\ell \in \sL} \prod_{\rvw \in \WC{\ell - 1}{K} } \P(\rD^{(\ell)}_{v, \rvw} = x^{(\ell)}_{\rvw}  ) \cdot \exp\Big( x^{(\ell)}_{\rvw}  \log \frac{x^{(\ell)}_{\rvw}/N_{\rvw} }{\etQ^{(\ell)}_{k \oplus \rvw}} + (N_{\rvw} - x^{(\ell)}_{\rvw}) \log \frac{1 - x^{(\ell)}_{\rvw}/N_{\rvw} }{1 - \etQ^{(\ell)}_{k \oplus \rvw}}\Big)\notag\\
    \lesssim &\,  e^{-\const} \exp\Big( \frac{T}{2}\log(\etD_{\max})\Big),\label{eqn:upperBoundProbG1G2}
    \end{align}
    where the proof of \eqref{eqn:upperBoundProbG1G2} is deferred later. By \Cref{ass:prob_ratio_bound}, $\const_{\eqref{eqn:prob_ratio_bound} } = O(1)$ and $ \etD_{\max} \asymp \const_{\eqref{eqn:prob_ratio_bound} } \D_{\rm{GCH}}(j,k)q_{N}$, then $\D_{\rm{GKL}}(j, k) = \D_{\rm{GCH}}(j,k) \cdot q_{N} \gg \frac{\etD_{\max}}{\log q_{N}}$ by \Cref{lem:equivalenveDivergence}. We take
    \begin{align}
        \const =&\, \Big(1 - \frac{\log(N)}{\sqrt{N}} \Big) \cdot \D_{\rm{GKL}}(j, k) - \frac{\etD_{\max}}{\log q_{N}} \overset{\textnormal{\Cref{lem:equivalenveDivergence}}}{=} (1 - o(1))\cdot \D_{\rm{GCH}}(j,k) q_{N} \asymp \frac{\etD_{\max}}{ \const_{\eqref{eqn:prob_ratio_bound} }} \gg \log(\etD_{\max}),\notag
    \end{align}
    leading to $e^{-\const}\exp\Big( \frac{T}{2}\log(\etD_{\max})\Big) = e^{- (1 - o(1))\cdot \const}$. Then we are able to establish \eqref{eqn:KLUpperBound},
    \begin{align*}
        &\, \P\bigg( \Big\{ \sum_{\ell \in \sL} \sum_{\rvw \in \WC{\ell - 1}{K} } N_{\rvw}\cdot \D_{\textnormal{KL}}(\mu^{(\ell)}_{v, \rvw}\parallel \etQ_{k \oplus \rvw}^{(\ell)}) \leq  \Big(1 - \frac{\log(N)}{\sqrt{N}} \Big) \D_{\textnormal{GKL}}(j, k) - \frac{\etD_{\max}}{\log q_{N}} \Big\} \bigcap \{\rD_v \leq 10 \etD_{\max} \} \bigg)\\
        \geq &\, \P(\rD_v \leq 10 \etD_{\max}) \cdot \bigg( 1 - e^{-\const}\exp\Big( \frac{T}{2}\log(\etD_{\max})\Big) \bigg) \\
        \geq &\, \P(\rD_v \leq 10 \etD_{\max}) - e^{-\const}\exp\Big( \frac{T}{2}\log(\etD_{\max})\Big) \\
        \geq &\, 1 - e^{-10 \etD_{\max}} - e^{- \const} =  1 - e^{-10 \etD_{\max}} - e^{- \D_{\rm{GCH}} \cdot q_{N}} \geq 1 - 2e^{-\D_{\rm{GCH}} \cdot q_{N}},
    \end{align*}
    where the last line holds since $\D_{\rm{GCH}} = \min_{1\leq j, k \leq K} \D_{\rm{GCH}} (j, k)$, $\const \gg T \log(\etD_{\max})$ and $\etD_{\max} \geq \D_{\rm{GCH}} \cdot q_{N}$ as shown in \Cref{lem:GCH_max_different_order}.
\end{proof}

\begin{proof}[Proof of inequality \eqref{eqn:upperBoundProbG1G2}]
    Note that $\rD^{(\ell)}_{v, \rvw}\big|_{\rY_{v} = k} \sim \mathrm{Bin}(N_{\rvw}, \etQ_{k \oplus \rvw}^{(\ell)})$, then by \Cref{lem:stirling},
\begin{align*}
          \P(\rD^{(\ell)}_{v, \rvw} = x^{(\ell)}_{\rvw}) = &\, \binom{N_{\rvw}}{x^{(\ell)}_{\rvw}} \exp\Big( x^{(\ell)}_{\rvw} \log \etQ^{(\ell)}_{\rvw} + (N_{\rvw} - x^{(\ell)}_{\rvw})\log(1 - \etQ^{(\ell)}_{\rvw}) \Big)\\
         \log \binom{N_{\rvw}}{x^{(\ell)}_{\rvw}} = &\, -\frac{1}{2}\log(2\pi x_{\rvw}^{(\ell)}) - \frac{1}{2}  \log \Big(1 - \frac{x_{\rvw}^{(\ell)}}{N_{\rvw}} \Big) + N_{\rvw} \log N_{\rvw} - x^{(\ell)}_{\rvw} \log x^{(\ell)}_{\rvw} \\
         &\, \quad - (N_{\rvw} - x^{(\ell)}_{\rvw}) \log (N_{\rvw}) - (N_{\rvw} - x^{(\ell)}_{\rvw}) \log \Big(1 - \frac{x_{\rvw}^{(\ell)}}{N_{\rvw}} \Big) + o(1)
    \end{align*}
Thus for each $\ell \in \sL$ and $\rvw \in \WC{\ell - 1}{K}$, we have
\begin{align*}
    &\,\P(\rD^{(\ell)}_{v, \rvw} = x^{(\ell)}_{\rvw}  ) \cdot \exp\Big( x^{(\ell)}_{\rvw}  \log \frac{x^{(\ell)}_{\rvw}/N_{\rvw} }{\etQ^{(\ell)}_{k \oplus \rvw}} + (N_{\rvw} - x^{(\ell)}_{\rvw}) \log \frac{1 - x^{(\ell)}_{\rvw}/N_{\rvw} }{1 - \etQ^{(\ell)}_{k \oplus \rvw}}\Big)\\
    = &\, \exp \Big( -\frac{1}{2}\log(2\pi x^{(\ell)}_{\rvw}) - \frac{1}{2}\log(1 - \frac{x^{(\ell)}_{\rvw}}{N_{\rvw}}) + o(1)\Big) = \exp \Big( -\frac{1}{2}\log(2\pi x^{(\ell)}_{\rvw}) + o(1)\Big),
\end{align*}
where the last equality holds since $1 \ll x^{(\ell)}_{\rvw} \lesssim \etD_{\max}^{(\ell)} \ll N^{\ell - 1}$.  Also, for sufficiently large $N$, the cardinality of $\sX$ is bounded by
\begin{align}
    |\sX| \leq \binom{10\etD_{\max}}{T} \leq \frac{(10\etD_{\max})^{T}}{T!} \leq \exp \Big( T\log(\etD_{\max}) \Big),
\end{align}
where $T = O(1)$, then
\begin{align*}
    &\,\sum_{\{\rvx^{(\ell)}\}\in \sX } \prod_{\ell \in \sL} \prod_{\rvw \in \WC{\ell - 1}{K} } \P(\rD^{(\ell)}_{v, \rvw} = x^{(\ell)}_{\rvw}  ) \cdot \exp\Big( x^{(\ell)}_{\rvw}  \log \frac{x^{(\ell)}_{\rvw}/N_{\rvw} }{\etQ^{(\ell)}_{k \oplus \rvw}} + (N_{\rvw} - x^{(\ell)}_{\rvw}) \log \frac{1 - x^{(\ell)}_{\rvw}/N_{\rvw} }{1 - \etQ^{(\ell)}_{k \oplus \rvw}}\Big)\\
    \leq &\, |\sX| \cdot \exp \Big( -\frac{T}{2}\log(2\pi x^{(\ell)}_{\rvw}) + o(1)\Big) \leq \exp\Big( \frac{T}{2}\log(\etD_{\max})\Big).
\end{align*}

\end{proof}

\begin{proof}[Proof of \Cref{lem:sizeofG}]
Let $\gE_{\ell}(v, \setS)$ denote the set of $\ell$-hyperedges containing $v$ with the rest $\ell - 1$ nodes in $\setS\subset \gV$, and $\gE_{\ell}(\setS, \setS) \coloneqq \cup_{v\in \setS} \gE_{\ell}(v, \setS)$ denote the set of $\ell$-hyperedges with all vertices in $\setS\subset \gV$, where the repeated edges are dropped, and $\gE(\setS, \setS)\coloneqq \cup_{\ell \in \sL}\gE_{\ell}(\setS, \setS)$ similarly. To conclude the proof, we build the sequence of sets $\{\setS^{(t)} \subset \gV\}_{0\leq t \leq t^{\star}}$ as follows.
\begin{itemize}
    \item $\setS^{(0)} = \setS$, where $\setS$ denotes the set of nodes which do not satisfy \eqref{eqn:G1} or \eqref{eqn:G2}. By \Cref{lem:notG1G2},  $|\setS|\leq N^{1 - \kappa_{N}}/3$ with probability at least $1 - 6e^{- \epsilon\kappa_{N}\log(N)}$.
    \item For $t\geq 1$, $\setS^{(t)} = \setS^{(t-1)} \cup \{v_t\}$ if there exists $v_t\in \gV$ s.t. $|\gE (v_t, \setS^{(t-1)})| > \LM (\etD_{\max})^{1 + \delta}$.
    \item The sequence ends if no such vertex can be found in $\gV \setminus \setS^{(t)}$, and we obtain $\setS^{(t^{\star})}$.
\end{itemize}

Denote $s = N^{1 - \kappa_{N}}$. Assume that $|\setS^{(t^{\star})}| \geq s$. Consider the step $\tilde{t} = s - |\setS^{(0)}|$, then $\tilde{t} \geq 2s/3$. By construction, each of the nodes added so far contributes more than $\LM(\etD_{\max})^{1 + \delta}$ in $|\gE(\setS^{(\tilde{t})}, \setS^{(\tilde{t})})|$. At that point, $|\setS^{(\widetilde{t})}| = s$. The number of edges 
$|\gE(\setS^{(\tilde{t})}, \setS^{(\tilde{t})})| \geq 2s/3 \cdot \LM(\etD_{\max})^{1 + \delta} \geq 4s/3 (\etD_{\max})^{1 + \delta} > s(\etD_{\max})^{1 + \delta}$ since $\LM \geq 2$. On the other hand, with high probability, there cannot be any such set $|\setS^{(\tilde{t})}|$ by \Cref{lem:fewOutsideEdges}. Therefore, $|\gV\setminus \gG| < s$.
\end{proof}

\begin{lemma}\label{lem:fewOutsideEdges}
    Let $\delta > 0$ be an absolute constant. Denote $s = N^{1 - \kappa_{N}}$ where $\kappa_{N} \in (0, 1]$ satisfies \eqref{eqn:optimality_condition}. Recall $\etD_{\max}$ in \eqref{eqn:G1}. Then with probability at least $1 - e^{-s(\etD_{\max})^{1 + \delta}/4}$, there is no subset $\setS\subset \gV$ with $|\setS| = s$ such that $|\gE(\setS, \setS)| \geq s (\etD_{\max})^{1 + \delta}$ when $N$ is large enough.
\end{lemma}
\begin{proof}[Proof of \Cref{lem:fewOutsideEdges}]
    For any given subset $\setS\subset \gV$ with $|\setS| = s$, by Markov's inequality,
    \begin{align*}
    &\,\P( |\gE(\setS, \setS)| \geq s(\etD_{\max})^{1 + \delta} ) \\
    \leq&\, \inf_{\beta >0} e^{-\beta s(\etD_{\max})^{1 + \delta}}\E[\exp(\beta \cdot |\gE(\setS, \setS)|) ] = \inf_{\beta >0} e^{-\beta s(\etD_{\max})^{1 + \delta}} \prod_{\ell \in \sL} \E[ \exp(\beta \cdot |\gE_{\ell}(\setS, \setS)|) ] \\
    \leq &\, \inf_{\beta >0} e^{-\beta s(\etD_{\max})^{1 + \delta}} \prod_{\ell \in \sL} [ \etQ^{(\ell)}_{\max}\cdot e^{\beta} + (1 - \etQ^{(\ell)}_{\max}) ]^{\binom{s}{\ell} },\quad ( \textnormal{edge independence and } 1 + x \leq e^{x})\\
   \leq &\, \inf_{\beta >0} \exp\Big[ \sum_{\ell \in \sL} \binom{s}{\ell}\etQ^{(\ell)}_{\max}(e^{\beta} -1) - \beta s (\etD_{\max})^{1 + \delta} \Big]\\
   = &\, \exp \Big\{ - \Big[ \sum_{\ell \in \sL} \binom{s}{\ell}\etQ^{(\ell)}_{\max} - s(\etD_{\max})^{1 + \delta} \Big] - s(\etD_{\max})^{1 + \delta} \cdot \log \Big( \frac{s(\etD_{\max})^{1 + \delta}}{\sum_{\ell \in \sL} \binom{s}{\ell}\etQ^{(\ell)}_{\max}} \Big)\Big\}\\
   \leq &\, \exp \Big[ - \frac{s}{2}(\etD_{\max})^{1 + \delta}\Big],
\end{align*}
where the second to last inequality holds since the function $h(\beta) = a(e^{\beta} - 1) - b\beta$ achieves its infimum at $\beta = \log(\frac{b}{a})$, and the last inequality holds since $o(1) = \binom{s}{\ell}\etQ^{(\ell)}_{\max} \ll s(\etD_{\max})^{1 + \delta}$ and $\binom{N}{\ell}\etQ^{(\ell)}_{\max} \asymp \etD_{\max}$. The number of subsets $\setS\subset \gV$ with size $s$ is $\binom{N}{s} \leq (\frac{eN}{s})^{s}$. Then by Markov inequality, the probability of existing a subset $\setS\subset \gV$ with $|S| = s$ such that $|\gE(\setS, \setS)| > s (\etD_{\max})^{1 + \delta}$ is upper bounded by
\begin{align*}
    &\, \P\big( \big|\{ \setS: |\gE(\setS, \setS)| > s (\etD_{\max})^{1 + \delta} \textnormal{ and } |\setS| = s\} \big| \geq 1 \big) \\
    \leq &\, \E \big[ \,\, \big|\{ \setS: |\gE(\setS, \setS)| > s (\etD_{\max})^{1 + \delta} \textnormal{ and } |\setS| = s\}\big| \,\, \big],
\end{align*}
where the term on the right hand side can be further bounded by
\begin{align*}
    &\,\E \big[ \,\, \big|\{ \setS: |\gE(\setS, \setS)| > s (\etD_{\max})^{1 + \delta} \textnormal{ and } |\setS| = s\}\big| \,\, \big] \\
    =&\, \binom{N}{s} \cdot \P( \{ |\gE(\setS, \setS)| \geq s (\etD_{\max})^{1 + \delta} \textnormal{ and } |S| = s\} ) \\
    \leq &\,  \Big(\frac{eN}{s} \Big)^{s}\cdot \exp\Big( - \frac{s}{2} (\etD_{\max})^{1 + \delta} \Big) \quad \textnormal{ (plug in $s = N^{1 - \kappa_{N}}$)}\\
    = &\, \exp\Big( -\frac{s}{2}(\etD_{\max})^{1 + \delta} + s\log\frac{eN}{s} \Big) = \exp\Big( -\frac{N^{1 - \kappa_{N}}}{2} (\etD_{\max})^{1 + \delta} \big[1 - \frac{2 \kappa_{N} \log(N) + 1}{(\etD_{\max})^{1 + \delta}} \big] \Big) \\
    \lesssim &\, \exp\Big( -\frac{N^{1 - \kappa_{N}}}{4}(\etD_{\max})^{1 + \delta}\Big),
\end{align*}
where the last inequality holds since $\etD_{\max} \gtrsim \D_{\rm{GCH}} \cdot q_{N} \geq (1 + \epsilon) \kappa_{N} \log(N)$ by \eqref{eqn:optimality_condition}, and
\begin{align*}
    \frac{2\kappa_{N} \log(N) + 1}{(\etD_{\max})^{1 + \delta}} \lesssim \frac{2\kappa_{N}\log(N) + 1}{(\D_{\rm{GCH}} \cdot q_{N})^{1 + \delta}} \lesssim \frac{1}{[\kappa_{N}\log(N)]^{\delta}} = o(1) < \frac{1}{2},
\end{align*}
since $\kappa_{N} \log(N) \to \infty$ as $N \to \infty$.
\end{proof}

\begin{proof}[Proof of \Cref{lem:decayofError}]
Let $\{\widehat{\gV}^{(t)}_{k} \}_{k=1}^{K}$ denote the partition at time $t$. Remind that $\widehat{\gN}^{(t)} \coloneqq \cup_{j=1}^{K} \widehat{\gN}^{(t)}_j = \cup_{j=1}^{K}( \widehat{\gV}_{j}^{(t)} \setminus \gV_j ) \cap \gG$, where $\widehat{\gN}^{(t)}_j \coloneqq \cup_{k\neq j} \widehat{\gN}^{(t)}_{jk} = (\widehat{\gV}^{(t)}_j \setminus \gV_j) \cap \gG$, and $\widehat{\gN}^{(t)}_{jk} \coloneqq (\widehat{\gV}^{(t)}_{j} \cap \gV_k) \cap \gG$, representing the set of vertices in $\gG$ which belong to $\gV_k$ but are misclassified in $\widehat{\gV}^{(t)}_{j}$.

For node $v\in \widehat{\gN}^{(t+1)}_{jk}$ ($v\in \gV_k$ but misclassified in $\widehat{\gV}^{(t+1)}_j$), by \Cref{alg:agnostic_refinement}, it has to satisfy
\begin{align}
    0 \leq \sum_{\ell \in \sL} \sum_{ \rvw \in \WC{\ell - 1}{K} } \Big( \widehat{\rD}^{(t)}_{v, \rvw} \cdot \log \frac{\widehat{\etQ}^{(\ell)}_{j \oplus \rvw}}{\widehat{\etQ}^{(\ell)}_{k \oplus \rvw}}  + [ \widehat{N}^{(t)}_{\rvw} - \widehat{\rD}^{(t)}_{v, \rvw}] \cdot \log \frac{( 1 - \widehat{\etQ}^{(\ell)}_{j \oplus \rvw} )}{( 1 - \widehat{\etQ}^{(\ell)}_{k \oplus \rvw} )} \Big)
\end{align}
where $\widehat{\rD}^{(t)}_{v, \rvw}$ denotes the number $\ell$-hyperedges containing $v$ with the rest $\ell - 1$ nodes distributed as $\rvw$ at step $t$, i.e., $\ervw_j$ nodes from $\widehat{\gV}^{(t)}_{j}$ for each $j\in [K]$. Let $\widehat{N}^{(t)}_{\rvw}$ denote the capacity of such hyperedges. Then with the proofs of \eqref{eqn:step1}, \eqref{eqn:step2}, \eqref{eqn:step3} \eqref{eqn:step4} deferred to next several pages, we have, 
\begin{subequations}
    \begin{align}
        0 \leq&\, \sum_{1 \leq j\neq k \leq K} \sum_{v\in \widehat{\gN}^{(t+1)}_{jk}} \sum_{\ell \in \sL} \, \sum_{ \rvw \in \WC{\ell - 1}{K} } \Big( \widehat{\rD}^{(t)}_{v, \rvw} \cdot \log \frac{\widehat{\etQ}^{(\ell)}_{j \oplus \rvw}}{\widehat{\etQ}^{(\ell)}_{k \oplus \rvw}}  + ( \widehat{N}^{(t)}_{\rvw} - \widehat{\rD}^{(t)}_{v, \rvw}) \cdot \log \frac{( 1 - \widehat{\etQ}^{(\ell)}_{j \oplus \rvw} )}{( 1 - \widehat{\etQ}^{(\ell)}_{k \oplus \rvw} )} \Big) \notag \\
        =&\, \sum_{1 \leq j\neq k \leq K} \sum_{v\in \widehat{\gN}^{(t+1)}_{jk}} \sum_{\ell \in \sL} \, \sum_{ \rvw \in \WC{\ell - 1}{K} } \bigg\{ \widehat{\rD}^{(t)}_{v, \rvw} \cdot \Big( \log \frac{\etQ^{(\ell)}_{j \oplus \rvw} }{\etQ^{(\ell)}_{k \oplus \rvw}} + \log \frac{\widehat{\etQ}^{(\ell)}_{j \oplus \rvw}}{\etQ^{(\ell)}_{j \oplus \rvw} } - \log \frac{\widehat{\etQ}^{(\ell)}_{k \oplus \rvw}}{\etQ^{(\ell)}_{k \oplus \rvw}} \Big)  \notag \\
        & \quad \quad + (\widehat{N}^{(t)}_{\rvw} - \widehat{\rD}^{(t)}_{v, \rvw}) \cdot \Big[ \log \frac{( 1 - \etQ^{(\ell)}_{j \oplus \rvw} )}{( 1 - \etQ^{(\ell)}_{k \oplus \rvw} )} + \log \frac{( 1 - \widehat{\etQ}^{(\ell)}_{j \oplus \rvw} )}{( 1 - \etQ^{(\ell)}_{j \oplus \rvw} )}  - \log \frac{( 1 - \widehat{\etQ}^{(\ell)}_{k \oplus \rvw} )}{( 1 - \etQ^{(\ell)}_{k \oplus \rvw} )} \Big] \bigg\} \notag \\        
        \leq &\, \sum_{1 \leq j\neq k \leq K} \sum_{v\in \widehat{\gN}^{(t+1)}_{jk}} \sum_{\ell \in \sL} \, \sum_{ \rvw \in \WC{\ell - 1}{K} } \Big\{ \widehat{\rD}^{(t)}_{v, \rvw} \cdot \log \frac{\etQ^{(\ell)}_{j \oplus \rvw} }{\etQ^{(\ell)}_{k \oplus \rvw}} + ( \widehat{N}^{(t)}_{\rvw} - \widehat{\rD}^{(t)}_{v, \rvw}) \cdot \log \frac{( 1 - \etQ^{(\ell)}_{j \oplus \rvw} )}{( 1 - \etQ^{(\ell)}_{k \oplus \rvw} )} \Big\}  \notag \\
        &\, \quad \quad + |\widehat{\gN}^{(t+1)}| \cdot \const_{\eqref{eqn:step1}} \cdot \log q_{N} \label{eqn:step1}\\
        \leq &\, \sum_{1 \leq j\neq k \leq K} \sum_{v\in \widehat{\gN}^{(t+1)}_{jk}} \sum_{\ell \in \sL} \, \sum_{ \rvw \in \WC{\ell - 1}{K} } \Big\{ \rD^{(\ell)}_{v, \rvw} \cdot \log \frac{\etQ^{(\ell)}_{j \oplus \rvw} }{\etQ^{(\ell)}_{k \oplus \rvw}} + (N_{\rvw} - \rD^{(\ell)}_{v, \rvw}) \cdot \log \frac{( 1 - \etQ^{(\ell)}_{j \oplus \rvw} )}{( 1 - \etQ^{(\ell)}_{k \oplus \rvw} )} \Big\}  \notag \\
        &\, \quad \quad + 2\log \const_{\eqref{eqn:prob_ratio_bound} } \sum_{v\in \widehat{\gN}^{(t+1)}} \sum_{\ell \in \sL} \sum_{ \rvw \in \WC{\ell - 1}{K} } |\gE_{\ell}(v, \widehat{\gN}^{(t)}, \rvw)| \notag \\
        &\, \quad \quad + |\widehat{\gN}^{(t+1)}| \cdot (  \const_{\eqref{eqn:step1}} \cdot \log q_{N} + \const_{\eqref{eqn:step2}} ) \label{eqn:step2} \\
        \leq &\, -|\widehat{\gN}^{(t+1)}|\cdot \frac{\etD_{\max}}{\log q_{N}} + 2\log \const_{\eqref{eqn:prob_ratio_bound} } \sum_{v\in \widehat{\gN}^{(t+1)}} \sum_{\ell \in \sL} \sum_{ \rvw \in \WC{\ell - 1}{K} } |\gE_{\ell}(v, \widehat{\gN}^{(t)}, \rvw)| \notag\\
        &\, \quad \quad + |\widehat{\gN}^{(t+1)}| \cdot (  \const_{\eqref{eqn:step1}} \cdot \log q_{N} + \const_{\eqref{eqn:step2}} ) \label{eqn:step3}\\ 
        \leq &\, - |\widehat{\gN}^{(t+1)}| \frac{\etD_{\max}}{\log q_{N}} + \const_{\eqref{eqn:step4}} \cdot \sqrt{ |\widehat{\gN}^{(t)}|\cdot |\widehat{\gN}^{(t+1)}| \cdot \etD_{\max}} \notag\\
        &\,+ |\widehat{\gN}^{(t+1)}| \cdot \bigg( \const_{\eqref{eqn:step1}}\cdot \log q_{N} + \const_{\eqref{eqn:step2}} + 2 \log \const_{\eqref{eqn:prob_ratio_bound} } \cdot \sum_{\ell \in \sL} \binom{\ell+K-1}{\ell} \bigg), \label{eqn:step4}
    \end{align}
\end{subequations}
where $|\gE_{\ell}(v, \widehat{\gN}^{(t)}, \rvw)|$ in \eqref{eqn:step2} denotes the number of $\ell$-hyperedges $e$ containing $v$ where the rest $\ell - 1$ nodes distributed as $\rvw$ among $\{\widehat{\gV}_{k}^{(t)}\}_{k=1}^{K}$ with at least one of those $\ell - 1$ nodes in $\widehat{\gN}^{(t)}$, and $\const_{\eqref{eqn:step4}} \coloneqq 2 \const_{\eqref{eqn:concentrateA}} \log \const_{\eqref{eqn:prob_ratio_bound} }$. Note that $\const_{\eqref{eqn:step1}}, \const_{\eqref{eqn:step2}}, \const_{\eqref{eqn:prob_ratio_bound} }, \LM, K = O(1)$, then
\begin{align}
    \frac{\etD_{\max}}{\log q_{N}} \gg \const_{\eqref{eqn:step1}} \cdot \log q_{N} + \const_{\eqref{eqn:step2}} + 2 \log \const_{\eqref{eqn:prob_ratio_bound} } \cdot \sum_{\ell \in \sL} \binom{\ell+K-1}{\ell}.
\end{align}
At the same time, $\etD_{\max} \asymp q_{N}$ by the assumption of \Cref{thm:concentration}. Then for sufficiently large $N$, we have
\begin{align}
    \frac{|\widehat{\gN}^{(t+1)}|}{|\widehat{\gN}^{(t)}|} \leq 2 \const_{\eqref{eqn:concentrateA}} \log \const_{\eqref{eqn:prob_ratio_bound} } \cdot \frac{[\log q_{N}]^2}{\etD_{\max}} \leq \frac{2 \const_{\eqref{eqn:concentrateA}} \log\const_{\eqref{eqn:prob_ratio_bound} }}{\const_{\eqref{eqn:prob_ratio_bound} }} \frac{[\log q_{N}]^2}{q_{N}} \leq \frac{\log\const_{\eqref{eqn:prob_ratio_bound} }}{\const_{\eqref{eqn:prob_ratio_bound} }} \frac{1}{e}\,.\notag
\end{align}
\end{proof}

\begin{proof}[Proof of \eqref{eqn:step1}]
Note that for any $\widehat{q}, q>0$ with $\widehat{q} < 2q$, the following holds
\begin{align}
    \Big| \log \frac{\widehat{q}}{q} \Big| \leq \log\Big( \frac{q}{q - |\widehat{q} - q|} \Big) \leq \frac{|\widehat{q} -q|}{q - |\widehat{q} - q|},
\end{align}
where the first inequality follows by discussions on cases $2q>\widehat{q} > q$ and $ \widehat{q} <q$, and the second comes from $\log(1 + x) \leq x$. Thus one may focus on bounding $|\widehat{\etQ}^{(\ell)}_{k \oplus \rvw} - \etQ^{(\ell)}_{k \oplus \rvw}|$. According to \Cref{lem:hatQmw_approx}, under \Cref{ass:prob_ratio_bound}, the following holds with probability at least $1 - \exp(-N\log q_{N}/q_{N})$ for each $\ell \in \sL$,
\begin{align*}
    \Big|\log \frac{\widehat{\etQ}^{(\ell)}_{k \oplus \rvw} }{\etQ^{(\ell)}_{k \oplus \rvw}} \Big| \leq &\, \frac{|\widehat{\etQ}^{(\ell)}_{k \oplus \rvw} - \etQ^{(\ell)}_{k \oplus \rvw} |}{ \etQ^{(\ell)}_{k \oplus \rvw} - |\widehat{\etQ}^{(\ell)}_{k \oplus \rvw} - \etQ^{(\ell)}_{k \oplus \rvw} | } \lesssim \frac{\log q_{N}}{q_{N} },\\
    \Big|\log \frac{( 1 - \widehat{\etQ}^{(\ell)}_{k \oplus \rvw} )}{( 1 - \etQ^{(\ell)}_{k \oplus \rvw} )} \Big| \leq &\, \frac{|\widehat{\etQ}^{(\ell)}_{k \oplus \rvw} - \etQ^{(\ell)}_{k \oplus \rvw} |}{ (1 - \etQ^{(\ell)}_{k \oplus \rvw}) - |\widehat{\etQ}^{(\ell)}_{k \oplus \rvw} - \etQ^{(\ell)}_{k \oplus \rvw} | } \lesssim \frac{\log q_{N}}{N^{\ell - 1}}.
\end{align*}
Meanwhile, $\etD_{\max} = \sum_{\ell \in \sL} \etD_{\max}^{(\ell)} \asymp q_{N}$ under \Cref{ass:prob_ratio_bound}, and $\widehat{\rD}^{(t)}_{v, \rvw} \lesssim \etD_{\max} \ll N^{\ell - 1}$ by {\eqref{eqn:G1}} with $\widehat{N}^{(t)}_{\rvw} \coloneqq \prod_{k=1}^{K} \binom{|\widehat{\gV}^{(t)}_{k}| }{\ervw_{k}} \asymp N^{\ell - 1}$. Then $(\widehat{N}^{(t)}_{\rvw} - \widehat{\rD}^{(t)}_{v, \rvw}) \asymp N^{\ell - 1}$, and
\begin{align*}
    &\, \sum_{\ell \in \sL} \, \sum_{ \rvw \in \WC{\ell - 1}{K} } \widehat{\rD}^{(t)}_{v, \rvw} \cdot \Big| \log \frac{\widehat{\etQ}^{(\ell)}_{j \oplus \rvw}}{\etQ^{(\ell)}_{j \oplus \rvw} } - \log \frac{\widehat{\etQ}^{(\ell)}_{k \oplus \rvw}}{\etQ^{(\ell)}_{k \oplus \rvw}} \Big| \lesssim \log q_{N}, \\
    &\, \sum_{\ell \in \sL} \, \sum_{ \rvw \in \WC{\ell - 1}{K} } ( \widehat{N}^{(t)}_{\rvw} - \widehat{\rD}^{(t)}_{v, \rvw}) \cdot \Big| \log \frac{( 1 - \widehat{\etQ}^{(\ell)}_{j \oplus \rvw} )}{( 1 - \etQ^{(\ell)}_{j \oplus \rvw} )}  - \log \frac{( 1 - \widehat{\etQ}^{(\ell)}_{k \oplus \rvw} )}{( 1 - \etQ^{(\ell)}_{k \oplus \rvw} )} \Big| \lesssim \log q_{N}.
\end{align*}
Consequently, with probability at least $1 - \exp(- N\log q_{N}/q_{N})$,
\begin{align*}
    &\,\sum_{j\neq k} \sum_{v\in \widehat{\gN}^{(t+1)}_{jk}} \sum_{\ell \in \sL} \, \sum_{ \rvw \in \WC{\ell - 1}{K} } \bigg\{ \widehat{\rD}^{(t)}_{v, \rvw} \cdot \Big| \log \frac{\widehat{\etQ}^{(\ell)}_{j \oplus \rvw}}{\etQ^{(\ell)}_{j \oplus \rvw} } - \log \frac{\widehat{\etQ}^{(\ell)}_{k \oplus \rvw}}{\etQ^{(\ell)}_{k \oplus \rvw}} \Big| \\
    &\, \quad \quad \quad \quad \quad \quad \quad \quad \quad \quad \quad \quad + (\widehat{N}^{(t)}_{\rvw} - \widehat{\rD}^{(t)}_{v, \rvw}) \cdot \Big| \log \frac{( 1 - \widehat{\etQ}^{(\ell)}_{j \oplus \rvw} )}{( 1 - \etQ^{(\ell)}_{j \oplus \rvw} )}  - \log \frac{( 1 - \widehat{\etQ}^{(\ell)}_{k \oplus \rvw} )}{( 1 - \etQ^{(\ell)}_{k \oplus \rvw} )} \Big| \bigg\}\\
    \leq &\,\const_{\eqref{eqn:step1}} |\widehat{\gN}^{(t+1)}|\cdot \log q_{N},
\end{align*}
where $\const_{\eqref{eqn:step1} } >0$ is some constant.
\end{proof}

\begin{proof}[Proof of \eqref{eqn:step2}]
We want to replace $\widehat{\rD}^{(t)}_{v, \rvw}$ with $\rD^{(\ell)}_{v, \rvw}$ in this step, where $\rD^{(\ell)}_{v, \rvw}$ (resp. $\widehat{\rD}^{(t)}_{v, \rvw}$) denotes the number of $\ell$-hyperedges containing $v$ with the rest $\ell - 1$ nodes distributed as $\rvw$, i.e., $w_l$ nodes in $\gV_l$ (resp. $\widehat{\gV}^{(t)}_l$) for each $l\in[K]$. Note that $\widehat{\rD}^{(t)}_{v, \rvw} \leq \rD^{(\ell)}_{v, \rvw} + |\widehat{\rD}^{(t)}_{v, \rvw} - \rD^{(\ell)}_{v, \rvw}|$ and $ ( \widehat{N}^{(t)}_{\rvw} - \widehat{\rD}^{(t)}_{v, \rvw}) \leq ( N_{\rvw} - \rD^{(\ell)}_{v, \rvw}) + |( \widehat{N}^{(t)}_{\rvw} - \widehat{\rD}^{(t)}_{v, \rvw}) -( N_{\rvw} - \rD^{(\ell)}_{v, \rvw}) |$, then it suffices to control the contribution from deviation terms.
\begin{itemize}
    \item For hyperedge $e \ni v$ with $\rvy(e) = \widehat{\rvy}^{(t)}(e)$, the contribution of such $e$ inside $|\widehat{\rD}^{(t)}_{v, \rvw} - \rD^{(\ell)}_{v, \rvw}|$ cancels out. Thus it suffices to consider $e\in \gE_{\ell}(v, \widehat{\gN}^{(t)}, \rvw)$, where $\gE_{\ell}(v, \widehat{\gN}^{(t)}, \rvw)$ denotes the set of $\ell$-hyperedges $e$ containing $v$ where the rest $\ell - 1$ nodes distributed as $\rvw$ among $\{\widehat{\gV}_{k}^{(t)}\}_{k=1}^{K}$ with at least one of those $\ell - 1$ nodes in $\widehat{\gN}^{(t)}$. Then under the regime \eqref{eqn:edge_density} and \Cref{ass:prob_ratio_bound},
\begin{align*}
    &\,\sum_{1 \leq j\neq k \leq K} \sum_{v\in \widehat{\gN}^{(t+1)}_{jk}} \sum_{\ell \in \sL} \, \sum_{ \rvw \in \WC{\ell - 1}{K} } |\widehat{\rD}^{(t)}_{v, \rvw} - \rD^{(\ell)}_{v, \rvw}| \cdot \log \frac{\etQ^{(\ell)}_{j \oplus \rvw} }{\etQ^{(\ell)}_{k \oplus \rvw}} \\
    \leq &\, \log \const_{\eqref{eqn:prob_ratio_bound} } \cdot \sum_{v\in \widehat{\gN}^{(t+1)}} \sum_{\ell \in \sL}\, \sum_{ \rvw \in \WC{\ell - 1}{K} } |\gE_{\ell}(v, \widehat{\gN}^{(t)}, \rvw)|.
\end{align*}

\item Meanwhile, $(1 - \mismatch_{N} )^{\ell - 1} \lesssim \Big( \frac{|\widehat{\gV}_{k}^{(t)} | }{|\gV_{k}|} \Big)^{\ell - 1} \asymp \frac{\widehat{N}^{(t)}_{\rvw}}{N_{\rvw}} \lesssim (1 + \mismatch_{N} )^{\ell - 1}$, where $\mismatch_{N} \lesssim (\rho_{N})^{-1}$ with probability at least $1 - O(N^{-c})$ for some $c > 0$ according to \eqref{eqn:weak_consistency_agnostic}, then we have
\begin{align*}
  |\widehat{N}^{(t)}_{\rvw} - N_{\rvw}| = N_{\rvw}\cdot \Big| 1 - \frac{\widehat{N}^{(t)}_{\rvw}}{N_{\rvw}}\Big| \leq [1 - (1 - \mismatch_{N})^{\ell - 1}]N_{\rvw} \leq (\ell - 1)\mismatch_{N} N_{\rvw} \lesssim N^{\ell - 1}(\rho_{N})^{-1}.
\end{align*}
Also, $\log \frac{1 - \etQ^{(\ell)}_{j \oplus \rvw} }{1 - \etQ^{(\ell)}_{k \oplus \rvw}} = \log(1 + \frac{\etQ^{(\ell)}_{k \oplus \rvw} - \etQ^{(\ell)}_{j \oplus \rvw} }{1 - \etQ^{(\ell)}_{k \oplus \rvw}}) = \frac{\etQ^{(\ell)}_{k \oplus \rvw} - \etQ^{(\ell)}_{j \oplus \rvw} }{1 - \etQ^{(\ell)}_{k \oplus \rvw}} (1 + o(1)) \asymp N^{-(\ell - 1)}\rho_{N}$, then there exists some constant $\const_{\eqref{eqn:step2}} > 0$ such that
\begin{align*}
    &\,\sum_{1 \leq j\neq k \leq K} \sum_{v\in \widehat{\gN}^{(t+1)}_{jk}} \sum_{\ell \in \sL} \, \sum_{ \rvw \in \WC{\ell - 1}{K} }  \big|( \widehat{N}^{(t)}_{\rvw} - \widehat{\rD}^{(t)}_{v, \rvw}) -( N_{\rvw} - \rD^{(\ell)}_{v, \rvw}) \big| \cdot \log \frac{1 - \etQ^{(\ell)}_{j \oplus \rvw} }{1 - \etQ^{(\ell)}_{k \oplus \rvw}}\\
   \leq &\, \sum_{v\in \widehat{\gN}^{(t+1)}} \sum_{\ell \in \sL} \, \sum_{ \rvw \in \WC{\ell - 1}{K} } |\widehat{N}^{(t)}_{\rvw} - N_{\rvw}| \cdot \frac{\etQ^{(\ell)}_{k \oplus \rvw} - \etQ^{(\ell)}_{j \oplus \rvw} }{1 - \etQ^{(\ell)}_{k \oplus \rvw}} + \sum_{v\in \widehat{\gN}^{(t+1)}} \sum_{\ell \in \sL} \, \sum_{ \rvw \in \WC{\ell - 1}{K} } | \widehat{\rD}^{(t)}_{v, \rvw} - \rD^{(\ell)}_{v, \rvw})| \cdot \frac{\etQ^{(\ell)}_{k \oplus \rvw} - \etQ^{(\ell)}_{j \oplus \rvw} }{1 - \etQ^{(\ell)}_{k \oplus \rvw}}\\
   \leq &\, \const_{\eqref{eqn:step2}} \cdot|\widehat{\gN}^{(t+1)}| + \log \const_{\eqref{eqn:prob_ratio_bound} } \cdot \sum_{v\in \widehat{\gN}^{(t+1)}} \sum_{\ell \in \sL}\, \sum_{ \rvw \in \WC{\ell - 1}{K} } |\gE_{\ell}(v, \widehat{\gN}^{(t)}, \rvw)|.
\end{align*}
\end{itemize}
\end{proof}

\begin{proof}[Proof of \eqref{eqn:step3}]
For $\mu^{(\ell)}_{v, \rvw}\coloneqq \rD^{(\ell)}_{v, \rvw}/ N_{\rvw}$, $N_{\rvw}\asymp N^{\ell - 1}$ in \eqref{eqn:nbw} for $\rvw \in \WC{\ell - 1}{K}$, then
\begin{align}
    &\, \rD^{(\ell)}_{v, \rvw} \cdot \log \frac{\etQ^{(\ell)}_{j \oplus \rvw} }{\etQ^{(\ell)}_{k \oplus \rvw}} + (N_{\rvw} - \rD^{(\ell)}_{v, \rvw}) \cdot \log \frac{( 1 - \etQ^{(\ell)}_{j \oplus \rvw} )}{( 1 - \etQ^{(\ell)}_{k \oplus \rvw} )} = N_{\rvw} \cdot \big[ \D_{\rm{KL}}(\mu^{(\ell)}_{v, \rvw}\parallel \etQ_{k \oplus \rvw}^{(\ell)}) - \D_{\rm{KL}}(\mu^{(\ell)}_{v, \rvw}\parallel \etQ_{j \oplus \rvw}^{(\ell)}) \big],\notag
\end{align}
where $\D_{\rm{KL}}(\mu^{(\ell)}_{v, \rvw}\parallel \etQ_{j \oplus \rvw}^{(\ell)})$ denotes the KL divergence defined in \eqref{eqn:KLDivergence}. It suffices to show
\begin{align}
    \sum_{\ell \in \sL} \, \sum_{ \rvw \in \WC{\ell - 1}{K} } N_{\rvw} \cdot [ \D_{\rm{KL}}(\mu^{(\ell)}_{v, \rvw}\parallel \etQ_{k \oplus \rvw}^{(\ell)}) - \D_{\rm{KL}}(\mu^{(\ell)}_{v, \rvw}\parallel \etQ_{j \oplus \rvw}^{(\ell)}) ] \leq - \frac{\etD_{\max}}{\log q_{N}}
\end{align}
for each $v\in \widehat{\gN}^{(t)}$, which holds true according to \eqref{eqn:G2} since $\widehat{\gN}^{(t)} \subset \gG$.
\end{proof}

\begin{proof}[Proof of \eqref{eqn:step4}]
Note that $\rho_{N} \asymp \etD_{\max} = \sum_{\ell \in \sL} \etD_{\max}^{(\ell)} \asymp q_{N}$ under \Cref{ass:prob_ratio_bound}, where $\etD_{\max}^{(\ell)} = \binom{N - 1}{\ell - 1} \etQ_{\max}^{(\ell)}$, then
\begin{align}
    \E |\gE_{\ell}(v, \widehat{\gN}^{(t)}, \rvw)| \leq \etQ^{(\ell)}_{\max} \cdot \mismatch_{N} N \cdot \binom{N - 2}{\ell - 2} \leq \mismatch_{N} q_{N} \lesssim 1,
\end{align}
where the last inequality holds since and $\mismatch_{N} \lesssim (\rho_{N})^{-1}$ by \eqref{eqn:weak_consistency_agnostic}. Then by a union bound, the following holds
\begin{align}
    \underset{v\in\gV}{\max}{ \sum_{\ell \in \sL} \sum_{ \rvw \in \WC{\ell - 1}{K} } \E |\gE_{\ell}(v, \widehat{\gN}^{(t)}, \rvw)|} \leq \sum_{\ell \in \sL} \binom{\ell+K-1}{\ell} 
\end{align}
Therefore by \Cref{thm:concentration},
    \begin{align*}
        &\,\sum_{v\in \widehat{\gN}^{(t+1)}} \sum_{\ell \in \sL} \,\sum_{ \rvw \in \WC{\ell - 1}{K} }\, |\gE_{\ell}(v, \widehat{\gN}^{(t)}, \rvw)| \\
        =&\, \sum_{v\in \widehat{\gN}^{(t+1)}} \sum_{\ell \in \sL} \,\sum_{ \rvw \in \WC{\ell - 1}{K} }\,  \big[ |\gE_{\ell}(v, \widehat{\gN}^{(t)}, \rvw)| - \E |\gE_{\ell}(v, \widehat{\gN}^{(t)}, \rvw)| \big] + \sum_{v\in \widehat{\gN}^{(t+1)}} \sum_{\ell \in \sL} \,\sum_{ \rvw \in \WC{\ell - 1}{K} }\,  \E |\gE_{\ell}(v, \widehat{\gN}^{(t)}, \rvw)|\\
        \leq &\, \ones^{\sT}_{\widehat{\gN}^{(t+1)}} \cdot (\rmA - \E \rmA) \cdot \ones_{\widehat{\gN}^{(t)}} + |\widehat{\gN}^{(t+1)}| \cdot \sum_{\ell \in \sL} \binom{\ell+K-1}{\ell} \\
        \leq &\, \|\ones_{\widehat{\gN}^{(t+1)}}\|_2 \cdot \| \rmA - \E \rmA\|_2 \cdot \|\ones_{\widehat{\gN}^{(t)}} \|_2 + |\widehat{\gN}^{(t+1)}|\cdot \sum_{\ell \in \sL} \binom{\ell + K - 1}{\ell} \\
        \leq &\, \const_{\eqref{eqn:concentrateA}} \sqrt{ |\widehat{\gN}^{(t)}|\cdot |\widehat{\gN}^{(t+1)}| \cdot \etD_{\max} \,\,} \,\, + |\widehat{\gN}^{(t+1)}|\cdot \sum_{\ell \in \sL} \binom{\ell+K-1}{\ell}.
    \end{align*}
Define $\const_{\eqref{eqn:step4}} \coloneqq 2 \const_{\eqref{eqn:concentrateA}} \log \const_{\eqref{eqn:prob_ratio_bound} }$, then by \Cref{thm:concentration}, the following holds
    \begin{align}
        &\, 2\log \const_{\eqref{eqn:prob_ratio_bound} } \sum_{v\in \widehat{\gN}^{(t+1)}} \sum_{\ell \in \sL} \sum_{ \rvw } |\gE_{\ell}(v, \widehat{\gN}^{(t)}, \rvw)|  + |\widehat{\gN}^{(t+1)}|\cdot (  \const_{\eqref{eqn:step1}}  \log q_{N} + \const_{\eqref{eqn:step2}} ) \notag \\
        \leq &\, \const_{\eqref{eqn:step4}} \cdot \sqrt{ |\widehat{\gN}^{(t)}|\cdot |\widehat{\gN}^{(t+1)}| \cdot \etD_{\max}} + |\widehat{\gN}^{(t+1)}| \cdot \Bigg(  \const_{\eqref{eqn:step1}}  \log q_{N} + \const_{\eqref{eqn:step2}} + 2 \log \const_{\eqref{eqn:prob_ratio_bound} } \cdot \sum_{\ell \in \sL} \binom{\ell+K-1}{\ell} \Bigg)\,\,,\notag
    \end{align}
    with probability at least $1-2N^{-10}- 2e^{-N}$.
\end{proof}

\section{Proofs in Section \ref{sec:known_partition}}\label{app:KnownPartition}

\begin{proof}[Proof of \Cref{lem:MAP}]
Recall that $\widehat{\rvy}_{-v}$ is some estimation of $\rvy_{-v}$, while  $\mathbb{H}$ and $\rmY_{-v}$ denote the laws of $\gH$ and $\rvy_{-v}$, respectively. By definition, the posterior probability is
\begin{align}
    \P (\rY_{v} = k \mid \mathbb{H} = \gH, \,\rmY_{-v} = \widehat{\rvy}_{-v}) = \frac{\P (\rY_{v} = k, \, \mathbb{H} = \gH, \, \rmY_{-v} = \widehat{\rvy}_{-v})}{\P (\mathbb{H} = \gH, \,\rmY_{-v} = \widehat{\rvy}_{-v})}.
\end{align}

We first evaluate the numerator. Let $d_v$ be the observed number of edges in $\gH$ containing $v$ with $\rD_v$ denoting its law. Let $\gH_{-v}$ be the hypergraph obtained by deleting edges containing $v$ from $\gH$, with $\mathbb{H}_{-v}$ denoting its law. Let $\{ \rY_{v} = k, \rmY_{-v} =  \widehat{\rvy}_{-v}\}$ denote the event that the membership vector is $\rY_{v} = k$ with the others being $\widehat{\rvy}_{-v}$. Thanks to the edgewise independence between $\gH_{-v}$ and $\rD_v$, the following holds
\begin{align*}
    &\,\P (\mathbb{H} = \gH, \, \rY_{v} = k,\,\rmY_{-v} = \widehat{\rvy}_{-v} ) = \P(\mathbb{H} = \gH \mid \rY_{v} = k,\, \rmY_{-v} = \widehat{\rvy}_{-v} )  \cdot \P (\rY_{v} = k,\, \rmY_{-v} = \widehat{\rvy}_{-v})\\
    =&\, \P(\rD_v = d_{v} \mid \rY_{v} = k,\, \rmY_{-v} = \widehat{\rvy}_{-v}) \cdot \P(\mathbb{H}_{-v} = \gH_{-v} \mid \rY_{v} = k,\, \rmY_{-v} = \widehat{\rvy}_{-v})  \cdot \P (\rY_{v} = k,\,\rmY_{-v} =  \widehat{\rvy}_{-v})\\
    =&\, \P(\rD_v = d_{v} \mid \rY_{v} = k,\, \rmY_{-v} = \widehat{\rvy}_{-v}) \cdot \P(\mathbb{H}_{-v} = \gH_{-v},\, \rY_{v} = k,\, \rmY_{-v} = \widehat{\rvy}_{-v} )\\
    =&\, \P(\rD_v = d_{v} \mid \rY_{v} = k,\, \rmY_{-v} = \widehat{\rvy}_{-v} ) \cdot \P ( \mathbb{H}_{-v} = \gH_{-v}, \,\rmY_{-v} = \widehat{\rvy}_{-v} \mid \rY_{v} = k) \cdot \P (\rY_{v} = k)\\
    =&\, \P(\rD_v = d_{v} \mid \rY_{v} = k,\, \rmY_{-v} = \widehat{\rvy}_{-v} ) \cdot \P ( \mathbb{H}_{-v} = \gH_{-v}, \,\rmY_{-v} = \widehat{\rvy}_{-v})\cdot \P (\rY_{v} = k),
\end{align*}
where the last equality holds due to two types of independence: first, $\gH_{-v}$ is independent of the assignment of $\{\rY_{v} = k\}$ according to the generating process of model \ref{def:non_uniform_HSBM}; second, the estimation $\widehat{\rvy}_{-v}$ is independent of the assignment of $\{\rY_{v} = k\}$ as well due to assumption. Then by edgewise independence between $\rD_v$ and $\gH_{-v}$, and the derivation above, we have
\begin{align}
    &\,\P (\rY_{v} = k \mid \mathbb{H} = \gH, \,\rmY_{-v} = \widehat{\rvy}_{-v} ) = \frac{\P (\rY_{v} = k,\, \mathbb{H} = \gH, \,\rmY_{-v} = \widehat{\rvy}_{-v}) }{ \P(\mathbb{H} = \gH, \,\rmY_{-v} = \widehat{\rvy}_{-v}) }\notag\\
    =&\, \frac{\P(\rD_v = d_{v} \mid \rY_{v} = k,\, \rmY_{-v} = \widehat{\rvy}_{-v} ) \cdot \P ( \mathbb{H}_{-v} = \gH_{-v}, \,\rmY_{-v} = \widehat{\rvy}_{-v})\cdot \P (\rY_{v} = k)}{ \P(\rD_v = d_v, \,\rmY_{-v} = \widehat{\rvy}_{-v} )\cdot \P(\mathbb{H}_{-v} = \gH_{-v}, \,\rmY_{-v} = \widehat{\rvy}_{-v} )} \notag\\
    =&\, \P(\rD_v = d_{v} \mid \rY_{v} = k,\, \rmY_{-v} = \widehat{\rvy}_{-v}) \cdot \P (\rY_{v} = k) \cdot [\P(\rD_v = d_v, \,\rmY_{-v} = \widehat{\rvy}_{-v})]^{-1}.\notag
\end{align}
Note that the denominator $\P(\rD_v = d_v, \,\rmY_{-v} = \widehat{\rvy}_{-v})$ can be written as follows
\begin{align}
    &\,\P(\rD_v = d_v, \,\rmY_{-v} = \widehat{\rvy}_{-v}) = \P(\rD_v = d_v \mid \rmY_{-v} = \widehat{\rvy}_{-v})\cdot \P(\rmY_{-v} = \widehat{\rvy}_{-v})\notag\\
    = &\, \P(\rmY_{-v} = \widehat{\rvy}_{-v}) \cdot \sum_{j\in [K]} \P(\rD_v = d_v\mid \rY_{v} = j, \,\, \rmY_{-v} = \widehat{\rvy}_{-v}),\notag
\end{align}
which is irrelevant to any specific assignment of $\{\rY_{v} = k\}$. Then it can be factored from the MAP, and the proof of the desired argument follows directly
\begin{align}
    \widehat{\ervy}_{v}^{\rm{\, MAP}} =&\, \underset{k \in [K]}{\argmax} \,\, \P( \rY_{v} = k\mid\mathbb{H} = \gH,\, \rmY_{-v} = \widehat{\rvy}_{-v}) \notag\\
    =&\, \underset{k \in [K]}{\argmax} \,\, \P(\rD_v = d_{v} | \rY_{v} = k,\, \rmY_{-v} = \widehat{\rvy}_{-v}) \cdot \P (\rY_{v} = k) \cdot [\P(\rD_v = d_v, \,\rmY_{-v} = \widehat{\rvy}_{-v})]^{-1} \notag \\
    =&\, \underset{k \in [K]}{\argmax} \,\, \P(\rD_v = d_v | \rY_{v} = k,\, \rmY_{-v} = \widehat{\rvy}_{-v} )\cdot \P(\rY_{v} = k ).\notag
\end{align}
Another way to verify the derivation above is to directly employ the independence between $\{\rY_{v} = k\}$ and $\{\rmY_{-v} = \widehat{\rvy}_{-v}\}$, i.e., $\P( \rY_{v} = k,\, \rmY_{-v} = \widehat{\rvy}_{-v}) = \P( \rY_{v} = k) \cdot \P( \rmY_{-v} = \widehat{\rvy}_{-v})$, which then gives us
\begin{align*}
    &\,\P(\rD_v = d_{v} | \rY_{v} = k,\, \rmY_{-v} = \widehat{\rvy}_{-v}) \cdot \P (\rY_{v} = k) \cdot [\P(\rD_v = d_v, \,\rmY_{-v} = \widehat{\rvy}_{-v} )]^{-1}\\
    =&\, \frac{\P( \rY_{v} = k,\,\rD_v = d_v ,\, \rmY_{-v} = \widehat{\rvy}_{-v} )}{\P( \rD_v = d_v ,\, \rmY_{-v} = \widehat{\rvy}_{-v} )} \cdot \frac{\P (\rY_{v} = k)}{\P( \rY_{v} = k,\, \rmY_{-v} = \widehat{\rvy}_{-v})}\\
    = &\,\P( \rY_{v} = k|\rD_v = d_v ,\, \rmY_{-v} = \widehat{\rvy}_{-v} ) \cdot[ \P( \rmY_{-v} = \widehat{\rvy}_{-v} )]^{-1},
\end{align*}
where both $\P( \rmY_{-v} = \widehat{\rvy}_{-v})$ and $\P(\rD_v = d_v, \,\rmY_{-v} = \widehat{\rvy}_{-v})$ on two sides of the equation can be factored out from MAP since they are not relevant to any specific assignment $\{\rY_{v} = k\}$. Then the \rm{MAP} becomes
\begin{align*}
    \widehat{\ervy}_{v}^{\rm{\, MAP}}=&\, \underset{k \in [K]}{\argmax} \,\, \P( \rY_{v} = k|\rD_v = d_v ,\, \rmY_{-v} = \widehat{\rvy}_{-v} )\\
    =&\, \underset{k \in [K]}{\argmax} \,\, \P(\rD_v = d_v | \rY_{v} = k,\, \rmY_{-v} = \widehat{\rvy}_{-v})\cdot \P(\rY_{v} = k ).
\end{align*}
\end{proof}

\begin{proof}[Proof of \Cref{lem:boundsErrorProb}]
For ease of presentation, define the event $\widehat{\gE}_{jk}$ to be
\begin{align}
    \widehat{\gE}_{jk} \coloneqq \mathbbm{1}\Big\{  \,\widehat{\P}_j(\rD_v = d_v)\cdot \alpha_j > \widehat{\P}_k(\rD_v = d_v)\cdot \alpha_k \Big\} 
\end{align}
Equivalently, assuming $\rY_{v} = k$, the probability of $v$ being misclassified is
\begin{align}
    \P_{\rm{err}}(v) = \sum_{d_v} \P_k(\rD_v = d_v) \cdot \P \big( \exists j \neq k, \textnormal{ s.t. } \{ \widehat{\gE}_{jk} = 1\} \big).
\end{align}
Then the upper bound in \eqref{eqn:upper_lower_bound_error_prob} in follows by a simple union bound, i.e.,
\begin{align}
    \P \big( \exists j \neq k, \textnormal{ s.t. } \{ \widehat{\gE}_{jk} = 1\} \big) \leq \sum_{j=1, j\neq k}^{K}\P( \widehat{\gE}_{jk} = 1) = \sum_{j=1, j\neq k}^{K} \Big(1 - \P( \widehat{\gE}_{jk} = 0) \Big).
\end{align}
On the other hand, $\P( \widehat{\gE}_{jk} = 0) \geq \prod_{j=1, j\neq k}^{K}\P( \widehat{\gE}_{jk} = 0)$ for any $j\in [K]$ with $j\neq k$, then 
\begin{align}
   &\, \sum_{j=1, j\neq k}^{K}  \P( \widehat{\gE}_{jk} = 0) \geq (K-1) \prod_{j=1, j\neq k}^{K}\P( \widehat{\gE}_{jk} = 0) \notag \\
   \iff &\, \sum_{j=1, j\neq k}^{K} \Big(1 - \P( \widehat{\gE}_{jk} = 0) \Big) \leq (K-1) \bigg( 1 - \prod_{j=1, j\neq k}^{K}\P( \widehat{\gE}_{jk} = 0) \bigg)\notag
\end{align}
Then the lower bound in\eqref{eqn:upper_lower_bound_error_prob} follows since
    \begin{align}
       \P \big( \exists j \neq k, \textnormal{ s.t. } \{ \widehat{\gE}_{jk} = 1\} \big) = 1 - \prod_{j=1, j\neq k}^{K} \P( \widehat{\gE}_{jk} = 0) \geq \frac{1}{K - 1}\sum_{j=1, j\neq k}^{K} \Big(1 - \P( \widehat{\gE}_{jk} = 0) \Big).\notag
    \end{align}
\end{proof}

\begin{proof}[Proof of \Cref{lem:boundsofHatProbability}]
According to \eqref{eqn:weak_consistency_known}, the mismatch ratio after the initial stage satisfies $\mismatch_{N} \lesssim \rho_{N}^{-1} = o(1)$ with probability at least $1 - N^{-10}$. Recall that $\widehat{\gN}^{(0)}$ in \Cref{lem:decayofError} denotes the set of nodes misclassified by \Cref{alg:spectral_initialization}, then $|\widehat{\gN}^{(0)}| \lesssim \mismatch_{N} N$. Let $\gE_{\ell}(v, \widehat{\gN}^{(0)})$ denote the set of $\ell$-hyperedges $e$ containing $v$ and at least one node in $e$ belongs to $\widehat{\gN}^{(0)}$, then the capacity of such $\ell$-hyperedge is at most $N_v^{(\ell)} \coloneqq \mismatch_{N} N \binom{N - 2}{\ell - 2} \lesssim N^{\ell - 1}/\rho_{N}$. Without loss of generality, we assume $v \in \gV_k$. Note that $\WC{\ell - 1}{K}$ is a finite set, then there exist $\rvw_1, \rvw_2, \rvw_3, \rvw_4 \in \WC{\ell - 1}{K}$ such that
\begin{align}
     \const^{(\ell)}_{k} \coloneqq \max_{\rvw, \rvw^{\prime} \in \WC{\ell - 1}{K}} \frac{\etQ^{(\ell)}_{k \oplus \rvw }}{\etQ^{(\ell)}_{k \oplus \rvw^{\prime}}} = \frac{\etQ^{(\ell)}_{k \oplus \rvw_1}}{\etQ^{(\ell)}_{k \oplus \rvw_2}},\quad \quad \frac{1 - \etQ^{(\ell)}_{k \oplus \rvw_3 }}{1 - \etQ^{(\ell)}_{k \oplus \rvw_4}} = \max_{\rvw, \rvw^{\prime} \in \WC{\ell - 1}{K}} \frac{1 - \etQ^{(\ell)}_{k \oplus \rvw }}{1 - \etQ^{(\ell)}_{k \oplus \rvw^{\prime}}},
\end{align}
where we can derive the facts $\etQ^{(\ell)}_{k \oplus \rvw_1} \geq \etQ^{(\ell)}_{k \oplus \rvw_2}$ and $\etQ^{(\ell)}_{k \oplus \rvw_4} \geq \etQ^{(\ell)}_{k \oplus \rvw_3}$. Then for edge $e \in \gE_{\ell}(v, \widehat{\gN}^{(0)})$, the ratio $\widehat{\P}_k(\etA_{e}^{(\ell)} = 1)/\P_k( \etA_{e}^{(\ell)} = 1)$ can be bounded as below
\begin{align}
    \frac{\etQ^{(\ell)}_{k \oplus \rvw_2 }}{\etQ^{(\ell)}_{k \oplus \rvw_1}} \leq&\, \frac{\widehat{\P}_k(\etA_{e}^{(\ell)} = 1) }{\P_k( \etA_{e}^{(\ell)} = 1)} \leq \frac{\etQ^{(\ell)}_{k \oplus \rvw_1 }}{\etQ^{(\ell)}_{k \oplus \rvw_2}}\,, \quad \frac{1 - \etQ^{(\ell)}_{k \oplus \rvw_4 }}{1 - \etQ^{(\ell)}_{k \oplus \rvw_3}} \leq \frac{\widehat{\P}_k(\etA_{e}^{(\ell)} = 0) }{\P_k( \etA_{e}^{(\ell)} = 0)} \leq \frac{1 - \etQ^{(\ell)}_{k \oplus \rvw_3 }}{1 - \etQ^{(\ell)}_{k \oplus \rvw_4}} \,.
\end{align}
By \Cref{def:good_degree}, $\rG(v) = 1$ if $\rG_{\ell}(v) = 1$ for each $\ell \in \sL$, i.e., $|\gE_{\ell}(v, \widehat{\gN}^{(0)})| \leq c_{N}^{(\ell)}$ where we denote $c_{N}^{(\ell)} = 4q^{(\ell)}_{\max}/\log(q^{(\ell)}_{\max})$ for simplicity. Consequently, the ratio $\widehat{\P}_k(\rD_v = d_v)/\P_k(\rD_v = d_v)$ can be bounded from both sides by $U_k$ and $1/U_k$ respectively, where
\begin{subequations}
    \begin{align}
        \frac{\widehat{\P}_k(\rD_v = d_v)}{\P_k(\rD_v = d_v)} \leq &\, \prod_{\ell \in \sL} \Bigg( \frac{\etQ^{(\ell)}_{k \oplus \rvw_1 }}{\etQ^{(\ell)}_{k \oplus \rvw_2}} \Bigg)^{c_{N}^{(\ell)}} \Bigg( \frac{1 - \etQ^{(\ell)}_{k \oplus \rvw_3 }}{1 - \etQ^{(\ell)}_{k \oplus \rvw_4}} \Bigg)^{N^{(\ell)}_{v} - c_{N}^{(\ell)}} \eqqcolon U_k\,,\\
        \frac{\widehat{\P}_k(\rD_v = d_v)}{\P_k(\rD_v = d_v)} \geq &\, \prod_{\ell \in \sL} \Bigg( \frac{\etQ^{(\ell)}_{k \oplus \rvw_2 }}{\etQ^{(\ell)}_{k \oplus \rvw_1}} \Bigg)^{c_{N}^{(\ell)}} \Bigg( \frac{1 - \etQ^{(\ell)}_{k \oplus \rvw_4 }}{1 - \etQ^{(\ell)}_{k \oplus \rvw_3}} \Bigg)^{N^{(\ell)}_{v} - c_{N}^{(\ell)}} \eqqcolon 1/U_k\,.
    \end{align}
\end{subequations}
Furthermore, we rewrite the terms above in the following way
    \begin{align}
       &\,\bigg( \frac{\etQ^{(\ell)}_{k \oplus \rvw_1 }}{\etQ^{(\ell)}_{k \oplus \rvw_2}} \bigg)^{c_{N}^{(\ell)}} = \exp\bigg(\log\Big( \frac{\etQ^{(\ell)}_{k \oplus \rvw_1 }}{\etQ^{(\ell)}_{k \oplus \rvw_2}}\Big)\cdot c_{N}^{(\ell)} \bigg) \,\\
        &\, \bigg( \frac{1 - \etQ^{(\ell)}_{k \oplus \rvw_3 }}{1 - \etQ^{(\ell)}_{k \oplus \rvw_4}} \bigg)^{N^{(\ell)}_{v} - c_{N}^{(\ell)}} = \exp\bigg(\log\Big(1 + \frac{\etQ^{(\ell)}_{k \oplus \rvw_4} - \etQ^{(\ell)}_{k \oplus \rvw_3 }}{1 - \etQ^{(\ell)}_{k \oplus \rvw_4}}\Big)\cdot (N^{(\ell)}_{v} - c_{N}^{(\ell)}) \bigg).
    \end{align}
Moreover, by using the facts that $\log(1 + x) = x$ for $x = o(1)$, we have
    \begin{align*}
        \log U_k = &\, \sum_{\ell \in \sL} \bigg( \log\Big(1 + \frac{\etQ^{(\ell)}_{k \oplus \rvw_4} - \etQ^{(\ell)}_{k \oplus \rvw_3 }}{1 - \etQ^{(\ell)}_{k \oplus \rvw_4}}\Big)\cdot (N^{(\ell)}_{v} - c_{N}^{(\ell)}) + \log\Big( \frac{\etQ^{(\ell)}_{k \oplus \rvw_1 }}{\etQ^{(\ell)}_{k \oplus \rvw_2}}\Big)\cdot c_{N}^{(\ell)} \bigg)\\
         =&\, \sum_{\ell \in \sL} \bigg( N^{(\ell)}_{v} \cdot \log\Big(1 + \frac{\etQ^{(\ell)}_{k \oplus \rvw_4} - \etQ^{(\ell)}_{k \oplus \rvw_3 }}{1 - \etQ^{(\ell)}_{k \oplus \rvw_4}}\Big) + \log\Big( \frac{\etQ^{(\ell)}_{k \oplus \rvw_1}(1 - \etQ^{(\ell)}_{k \oplus \rvw_4}) }{\etQ^{(\ell)}_{k \oplus \rvw_2}(1 - \etQ^{(\ell)}_{k \oplus \rvw_3})}\Big)\cdot c_{N}^{(\ell)} \bigg) \\
         \lesssim &\, \sum_{\ell \in \sL} N^{(\ell)}_{v} \etQ^{(\ell)}_{k\oplus\rvw_4} + (1 + o(1)) \cdot \sum_{\ell \in \sL}\log\Big( \frac{\etQ^{(\ell)}_{k \oplus \rvw_1} }{\etQ^{(\ell)}_{k \oplus \rvw_2}}\Big) \cdot c_{N}^{(\ell)}.
    \end{align*}

For the first term, note that $\etQ^{(\ell)}_{k \oplus \rvw} \binom{N - 1}{\ell - 1} = \etP^{(\ell)}_{k \oplus \rvw} \cdot q^{(\ell)}_{k \oplus \rvw} \gg 1$ for any $\rvw \in \WC{\ell - 1}{K}$, and $N_v^{(\ell)} = \mismatch_{N} N \binom{N - 2}{\ell - 2} \lesssim N^{\ell - 1}/\rho_{N}$. Since $(\rho_{N})^{-1} = o(1)$ and $\sum_{\ell \in \sL} q^{(\ell)}_{k \oplus \rvw_4} \lesssim \D_{\mathrm{GCH}}(j, k) \cdot q_{N}$ by \Cref{ass:expected_center_separation}, then for any $j\neq k$, the following holds with probability at least $1 - N^{-10}$
\begin{align}
    \sum_{\ell \in \sL} N^{(\ell)}_{v} \etQ^{(\ell)}_{k\oplus\rvw_4} \asymp (\rho_{N})^{-1}\sum_{\ell \in \sL} q^{(\ell)}_{k \oplus \rvw_4} \ll \D_{\mathrm{GCH}}(j, k) \cdot q_{N}.
\end{align}

For the second term, $\etQ^{(\ell)}_{k \oplus \rvw_1}/\etQ^{(\ell)}_{k \oplus \rvw_2} \leq q^{(\ell)}_{\max}/q^{(\ell)}_{\min}$. Since $\log(q^{(\ell)}_{\max}/q^{(\ell)}_{\min})/\log(q^{(\ell)}_{\max}) = o(1)$ by \eqref{eqn:Strong_consistency_prob_condition} and $\sum_{\ell \in \sL}q^{(\ell)}_{\max} \asymp \rho_{N} \asymp \D_{\mathrm{GCH}}(j, k) q_{N}$ by \Cref{ass:expected_center_separation}, we then have
    \begin{align}
        \log\Big( \frac{\etQ^{(\ell)}_{k \oplus \rvw_1} }{\etQ^{(\ell)}_{k \oplus \rvw_2}}\Big) \cdot c_{N}^{(\ell)} \leq \log( q^{(\ell)}_{\max}/q^{(\ell)}_{\min}) \cdot  \frac{4q^{(\ell)}_{\max}}{\log(q^{(\ell)}_{\max})} \ll \D_{\mathrm{GCH}}(j, k) \cdot q_{N}
    \end{align}
    for any $j\neq k$. Therefore $\log U_k \ll \D_{\mathrm{GCH}}(j, k) \cdot q_{N}$, and the desired argument follows. 
\end{proof}

\begin{proof}[Proof of inequality \eqref{eqn:1UpperBound}]
Conditioning on $v\in \gV_k$ and $\rG(v) = 1$, when $v$ is misclassified to $\widehat{\gV}_j$ by MAP, i.e., $\{\widehat{\P}_j(\rD_v = d_v)\cdot \alpha_j > \widehat{\P}_k(\rD_v = d_v)\cdot \alpha_k \}$, then $\{U_j \cdot \alpha_j \P_j(\rD_v = d_v) >  1/U_k \cdot \alpha_k \P_k(\rD_v = d_v)\}$ is true as well, since the following holds by \Cref{lem:boundsofHatProbability}
\begin{align}
    1/U_k \cdot \alpha_k \P_k(\rD_v = d_v) \leq \alpha_k \widehat{\P}_k(\rD_v = d_v) < \alpha_j \widehat{\P}_j(\rD_v = d_v) \leq U_j \cdot \alpha_j \P_j(\rD_v = d_v).
\end{align}
We then expand $\circled{1}$ as follows, where the $\widehat{\rR}_{k|1}$ in \eqref{eqn:ConditionedProbComparation} is written as a union bound
    \begin{align*}
    \circled{1} =&\, \sum \P_{k}(\rD_v = d_v) \cdot \widehat{\rR}_{k|1} \cdot \P(\rG(v) = 1)\\
    =&\, \sum_{d_v} \P(\rG(v) = 1) \cdot \P_k(\rD_v = d_v) \cdot \sum_{j=1, j\neq k}^{K}\P\Big\{  \,\widehat{\P}_j(\rD_v = d_v)\cdot \alpha_j > \widehat{\P}_k(\rD_v = d_v)\cdot \alpha_k \mid \rG(v) = 1 \Big\} \\
    = &\, \sum_{j=1, j\neq k}^{K}\,\, \sum_{d_v} \P(\rG(v) = 1) \cdot \P_k(\rD_v = d_v) \cdot \P\Big\{  \,\widehat{\P}_j(\rD_v = d_v)\cdot \alpha_j > \widehat{\P}_k(\rD_v = d_v)\cdot \alpha_k \mid \rG(v) = 1 \Big\}\\
    \leq &\, \sum_{j=1, j\neq k}^{K}\,\, \sum_{d_v} \P(\rG(v) = 1)\cdot  \P_k(\rD_v = d_v) \cdot \P\Big\{  \,U_j \cdot \alpha_j \P_j(\rD_v = d_v) > 1/U_k \cdot \alpha_k \P_k(\rD_v = d_v) \mid \rG(v) = 1 \Big\}\\
    = &\, \sum_{j=1, j\neq k}^{K}\, \sum_{d_v} \P(\rG(v) = 1) \cdot \P_k(\rD_v = d_v) \cdot \E \indi{ U_j \cdot \alpha_j \P_j(\rD_v = d_v) > 1/U_k \cdot \alpha_k \P_k(\rD_v = d_v) }\\
    \leq &\, \sum_{j=1, j\neq k}^{K}\,\, \sum_{d_v} \P(\rG(v) = 1) \cdot  \P_k(\rD_v = d_v) \cdot \indi{U_j \cdot \alpha_j \P_j(\rD_v = d_v) > 1/U_k \cdot \alpha_k \P_k(\rD_v = d_v) }
    \end{align*}
Suppose $\{  \,U_j \cdot \alpha_j \P_j(\rD_v = d_v) > 1/U_k \cdot \alpha_k \P_k(\rD_v = d_v)\}$, we then have 
\begin{align*}
    \P_k(\rD_v = d_v) \leq \alpha_j U_j U_k /\alpha_k \P_j(\rD_v = d_v).
\end{align*}
Note that $\P(\rG(v) = 1) \leq 1$, $U_j, U_k \gtrsim 1$ and $\alpha_j, \alpha_k \asymp 1$, it then implies 
\begin{align*}
    \P_k(\rD_v = d_v) \lesssim \alpha_j U_j U_k /\alpha_k \cdot \min\{\P_j(\rD_v = d_v), \P_k(\rD_v = d_v)\}
\end{align*}
Then we can bound $\circled{1}$ as follows
    \begin{align*}
        \circled{1} \leq &\, \sum_{j=1, j\neq k}^{K}\,\, \sum_{d_v} \P(\rG(v) = 1) \cdot \frac{\alpha_j}{\alpha_k}\cdot U_j U_k \cdot \min\{\P_j(\rD_v = d_v), \P_k(\rD_v = d_v)\}\\
        \leq &\, \sum_{j=1, j\neq k}^{K} \frac{\alpha_j}{\alpha_k} U_j U_k \sum_{d_v} \min\{\P_j(\rD_v = d_v), \,\, \P_k(\rD_v = d_v)\} \quad (\textnormal{\Cref{lem:probabilityDivergence}})\\
        \leq &\, \sum_{j=1, j\neq k}^{K} \frac{\alpha_j}{\alpha_k} U_j U_k \cdot \exp(- (1 + o(1)) \cdot \D_{\rm{GCH}}(j,k) \cdot q_{N}) \quad (\textnormal{\Cref{lem:boundsofHatProbability}})\\
        \leq &\, (K - 1) \cdot \frac{1 - \alpha_k}{\alpha_k} \cdot \exp\Big(- (1 + o(1)) \cdot q_{N} \cdot \min_{j\in [K], j\neq k}\D_{\rm{GCH}}(j,k) + \log U_k\Big) \\
        \leq &\, (K - 1)\cdot (1 - \alpha_k)/\alpha_k \cdot \exp\Big(- (1 + o(1))\cdot \D_{\rm{GCH}}\cdot q_{N}\Big)
    \end{align*}
where the last inequality holds since $\log U_k \ll \D_{\rm{GCH}}(j,k) \cdot q_{N}$ for any $j\neq k$, and $\D_{\rm{GCH}}= \min_{1\leq j < k \leq K}\D_{\rm{GCH}}(j,k)$.
\end{proof}

\begin{lemma}\label{lem:probabilityDivergence}
For some $j\neq k$ and $j\in [K]$, we have
\begin{align}
        \sum\limits_{d_v} \min \big\{ \P_j(\rD_v = d_v),\, \P_k(\rD_v = d_v)  \big\} \leq \exp(- (1 - o(1)) \cdot \D_{\rm{GCH}}(j,k) \cdot q_{N}).
    \end{align}
\end{lemma}
\begin{proof}[Proof of \Cref{lem:probabilityDivergence}]
Let $\LM$ denote the maximum element of set $\sL$. According to the independence between different layers, for any $t\in [0, 1]$, we have
        \begin{align*}
            &\,\sum_{d_v} \min \Big\{ \P_j(\rD_v = d_v),\, \P_k(\rD_v = d_v)  \Big\}\\
            =&\, \sum_{d^{(2)}_v } \cdots \sum_{ d^{(\LM)}_{v} } \min \Big\{ \prod_{\ell \in \sL} \P_j(\rD^{(\ell)}_v = d^{(\ell)}_v),\, \prod_{\ell \in \sL} \P_k(\rD^{(\ell)}_v = d^{(\ell)}_v)  \Big\}\\
            \leq &\, \sum_{d^{(2)}_v } \cdots \sum_{ d^{(\LM)}_{v} } \Big( \prod_{\ell \in \sL} \P_j(\rD^{(\ell)}_v = d^{(\ell)}_v)\Big)^t \cdot \Big(\prod_{\ell \in \sL}\P_k(\rD^{(\ell)}_v = d^{(\ell)}_v) \Big)^{1 - t}\\
            =&\, \sum_{d^{(2)}_v } \cdots \sum_{ d^{(\LM)}_{v} } \prod_{\ell \in \sL} \Big(\P_j(\rD^{(\ell)}_v = d^{(\ell)}_v)\Big)^t \Big(\P_k(\rD^{(\ell)}_v = d^{(\ell)}_v) \Big)^{1 - t}\\
            = &\, \prod_{\ell \in \sL} \bigg[ \sum_{ d^{(\ell)}_{v} } \Big(\P_j(\rD^{(\ell)}_v = d^{(\ell)}_v)\Big)^t \Big(\P_k(\rD^{(\ell)}_v = d^{(\ell)}_v) \Big)^{1 - t} \bigg], \quad \forall t\in[0, 1].
        \end{align*}
Note that $\rD^{(\ell)}_{v, \rvw}|_{\rY_{v} = k} \sim {\rm{Bin}}(N_{\rvw}, \etQ^{(\ell)}_{k \oplus \rvw})$ with $N_{\rvw}$ denoting the capacity of such $\ell$-hyperedges, then
\begin{align}
    \P_k(\rD^{(\ell)}_v = d^{(\ell)}_v) = \prod_{\rvw \in \WC{\ell - 1}{K}} \P_k[ \rD^{(\ell)}_{v, \rvw} = d^{(\ell)}_{v, \rvw} ].
\end{align}
Note that $d^{(\ell)}_v = \{d^{(\ell)}_{v, \rvw}\}_{\rvw \in \WC{\ell - 1}{K}}$ and $d^{(\ell)}_{v, \rvw}$ will go through all elements in sample space. Consequently for any $t \in [0, 1]$, by exchanging the order of sum and product, we have the
    \begin{align*}
        &\,\sum_{ d^{(\ell)}_{v} } \Big(\P_j(\rD^{(\ell)}_v = d^{(\ell)}_v)\Big)^t \Big(\P_k(\rD^{(\ell)}_v = d^{(\ell)}_v) \Big)^{1 - t}\\
        =&\, \sum_{ d^{(\ell)}_{v} } \prod_{\rvw \in \WC{\ell - 1}{K}} \Big(\P_j[ \rD^{(\ell)}_{v, \rvw} = d^{(\ell)}_{v, \rvw} ] \Big)^{t} \cdot \Big(\P_k[ \rD^{(\ell)}_{v, \rvw} = d^{(\ell)}_{v, \rvw} ] \Big)^{1 - t}\\
        \leq &\, \prod_{\rvw \in \WC{\ell - 1}{K}} \,\, \sum_{d^{(\ell)}_{v, \rvw}}\Big( \P_j[ \rD^{(\ell)}_{v, \rvw} = d^{(\ell)}_{v, \rvw} ] \Big)^{t} \cdot \Big(\P_k[ \rD^{(\ell)}_{v, \rvw} = d^{(\ell)}_{v, \rvw} ] \Big)^{1 - t}.
    \end{align*}
    By plugging in the probability $\P_k[ \rD^{(\ell)}_{v, \rvw} = d^{(\ell)}_{v, \rvw} ] = \binom{N_{\rvw}}{d^{(\ell)}_{v, \rvw} } [\etQ^{(\ell)}_{k \oplus \rvw}]^{d^{(\ell)}_{v, \rvw}} \cdot[ 1 - \etQ^{(\ell)}_{k \oplus \rvw} ]^{N_{\rvw} - d^{(\ell)}_{v, \rvw}}$, we have
    \begin{align*}
        &\,\sum_{d^{(\ell)}_{v, \rvw}} \big( \P_j[ \rD^{(\ell)}_{v, \rvw} = d^{(\ell)}_{v, \rvw} ] \big)^{t} \cdot \big(\P_k[ \rD^{(\ell)}_{v, \rvw} = d^{(\ell)}_{v, \rvw} ] \big)^{1 - t}\\
        =&\, \sum_{d^{(\ell)}_{v, \rvw}} \binom{N_{\rvw}}{d^{(\ell)}_{v, \rvw} } \Big(\frac{ [\etQ^{(\ell)}_{j \oplus \rvw} ]^t }{ [\etQ^{(\ell)}_{k \oplus \rvw} ]^{t-1} } \Big)^{d^{(\ell)}_{v, \rvw}} \cdot \, \Big( \frac{ [1 - \etQ^{(\ell)}_{j \oplus \rvw} ]^t }{ [1 - \etQ^{(\ell)}_{k \oplus \rvw} ]^{t-1} } \Big)^{N_{\rvw} - d^{(\ell)}_{v, \rvw}} \\
        =&\,\Big(  [\etQ^{(\ell)}_{j \oplus \rvw} ]^{t} \cdot [\etQ^{(\ell)}_{k \oplus \rvw} ]^{1 - t} +  [1 - \etQ^{(\ell)}_{j \oplus \rvw}]^t \cdot [ 1 - \etQ^{(\ell)}_{k \oplus \rvw}]^{1-t} \Big)^{N_{\rvw}}, \quad \forall t\in[0, 1],
    \end{align*}
    where the last equality holds by Binomial theorem. Then we put pieces together
    \begin{align*}
        &\, \sum_{d_v} \min \Big\{ \P_j(\rD_v = d_v),\, \P_k(\rD_v = d_v)  \Big\}\\
        \leq &\, \prod_{\ell \in \sL} \prod_{\rvw \in \WC{\ell - 1}{K}} \Big(  [\etQ^{(\ell)}_{j \oplus \rvw} ]^{t} \cdot [\etQ^{(\ell)}_{k \oplus \rvw} ]^{1 - t} +  [1 - \etQ^{(\ell)}_{j \oplus \rvw}]^t \cdot [ 1 - \etQ^{(\ell)}_{k \oplus \rvw}]^{1-t} \Big)^{N_{\rvw}} \\
        = &\, \exp\Bigg\{ \sum_{\ell \in \sL} \sum_{\rvw \in \WC{\ell - 1}{K}}  N_{\rvw} \cdot \log \Big(  [\etQ^{(\ell)}_{j \oplus \rvw} ]^{t} \cdot [\etQ^{(\ell)}_{k \oplus \rvw} ]^{1 - t} +  [1 - \etQ^{(\ell)}_{j \oplus \rvw}]^t \cdot [ 1 - \etQ^{(\ell)}_{k \oplus \rvw}]^{1-t} \Big) \Bigg\}.
    \end{align*}
    Note that for sufficiently small $x$ and $y$, $\log(1 - x) = -x + O(x^2)$ and
    \begin{align*}
         (1 - x)^t(1 - y)^{1-t} =&\,\exp( t\log(1 - x) + (1 - t)\log(1 - y)) = \exp( -tx - (1 - t)y + O(x^2 + y^2))\\
         =&\, 1 - tx - (1-t)y + O(x^2 + y^2).
    \end{align*}
  Note that $\etQ^{(\ell)}_{k \oplus \rvw} = o(1)$ under \eqref{eqn:edge_density}. Let $t^{\star} \in[0, 1]$ be the maximizer of \eqref{eqn:GCH_not_proportional}. By replacing $N_{\rvw}$ with $\overline{N}_{\rvw}$ using \Cref{lem:sizeDeviation}, we then have
    \begin{align*}
        &\,\sum_{d_v} \min \Big\{ \P_j(\rD_v = d_v),\, \P_k(\rD_v = d_v)  \Big\} \\
        \leq &\, \exp\bigg\{ \sum_{\ell \in \sL} \sum_{\rvw \in \WC{\ell - 1}{K}} N_{\rvw} \cdot \log \bigg[ [\etQ^{(\ell)}_{j \oplus \rvw} ]^{t^{\star}} \cdot [\etQ^{(\ell)}_{k \oplus \rvw} ]^{1 - t^{\star}} + 1 - t^{\star}\etQ^{(\ell)}_{j \oplus \rvw}  - (1 - t^{\star})\etQ^{(\ell)}_{k \oplus \rvw}  + O\big( (\etQ^{(\ell)}_{k \oplus \rvw})^2  \big)  \,\bigg] \bigg\}\\
        =&\, \exp\bigg\{ -\sum_{\ell \in \sL} \sum_{\rvw \in \WC{\ell - 1}{K}}  N_{\rvw} \cdot \bigg[ t^{\star}\etQ^{(\ell)}_{j \oplus \rvw}  + (1 - t^{\star})\etQ^{(\ell)}_{k \oplus \rvw} - [\etQ^{(\ell)}_{j \oplus \rvw} ]^{t^{\star}} \cdot [\etQ^{(\ell)}_{k \oplus \rvw} ]^{1 - t^{\star}}  + O\big( (\etQ^{(\ell)}_{k \oplus \rvw})^2  \big) \bigg] \bigg\}\\
        = &\,
        \exp\bigg\{ -\sum_{\ell \in \sL} \sum_{\rvw \in \WC{\ell - 1}{K}}  (1 + o(1))\cdot \overline{N}_{\rvw} \cdot \bigg[ t^{\star}\etQ^{(\ell)}_{j \oplus \rvw}  + (1 - t^{\star})\etQ^{(\ell)}_{k \oplus \rvw} - [\etQ^{(\ell)}_{j \oplus \rvw} ]^{t^{\star}} \cdot [\etQ^{(\ell)}_{k \oplus \rvw} ]^{1 - t^{\star}} \bigg] \bigg\}\\
        =&\, \exp(- (1 + o(1)) \cdot \D_{\rm{GCH}}(j,k) \cdot q_{N}),
    \end{align*}
which completes the proof of the desired argument.
\end{proof}

\begin{proof}[Proof of \Cref{lem:badProbUpperBound}]
For the sake of convenience, denote $c_{N}^{(\ell)} = 4q^{(\ell)}_{\max}/\log(q^{(\ell)}_{\max})$. The event $\{\rG(v) = 0\}$ means that for at least one of $\ell \in \sL$, the cardinality of the edge set $\gE_{\ell}(v, \widehat{\gN}^{(0)})$ (each edge contains $v$ and at least one node from $\widehat{\gN}^{(0)}$ ) is not small, i.e., 
\begin{align*}
    \P_k( \rG(v) = 0) \leq \sum\limits_{\ell \in \sL} \P_k( \rG_{\ell}(v) = 0)  = \sum\limits_{\ell \in \sL} \P\big( |\gE_{\ell}(v, \widehat{\gN}^{(0)})| > c_{N}^{(\ell)}\big).
\end{align*}
Since $\sL$ is a finite set, we only need to focus on extreme event for each $\ell \in \sL$. By \eqref{eqn:weak_consistency_known}, with probability at least $1 - N^{-10}$, the number of misclassified nodes in $\widehat{\rvy}^{(0)}$ is at most $\mismatch_{N} N$ where $\mismatch_{N} \lesssim (\rho_{N})^{-1}$. Then the cardinality of the edge set $\gE_{\ell}(v, \widehat{\gN}^{(0)})$ is at most $\mismatch_{N} N \binom{N - 2}{\ell - 2} \lesssim N^{\ell - 1}/\rho_{N}$, since we can choose the remaining $\ell - 2$ nodes (beside $v$ and the misclassified one) from all $N - 2$ nodes. Since $\rho_{N} \asymp \sum_{\ell \in \sL} q_{\max}^{(\ell)}$, then
\begin{align*}
    \E \Big( |\gE_{\ell}(v, \widehat{\gN}^{(0)})| \Big) \leq \mismatch_{N} N \binom{N - 2}{\ell - 2} \cdot \etQ_{\max}^{(\ell)} \lesssim \frac{q_{\max}^{(\ell)}}{\rho_{N}} \lesssim \frac{1}{q_{\max}^{(\ell)}} \cdot q_{\max}^{(\ell)} \asymp 1.
\end{align*}
We choose $\gamma_{N}^{(\ell)} \coloneqq 1/q_{\max}^{(\ell)}$, and $\zeta_{N}^{(\ell)} \coloneqq 4/\log(q_{\max}^{(\ell)})$, which satisfies $\zeta_{N}^{(\ell)}\geq \gamma_{N}^{(\ell)}$, and
\begin{align*}
    &\, \zeta_{N}^{(\ell)}\log(\zeta_{N}^{(\ell)}/\gamma_{N}^{(\ell)}) + (\gamma_{N}^{(\ell)} - \zeta_{N}^{(\ell)}) \\
    \geq &\, \frac{4}{\log(q_{\max}^{(\ell)})} \bigg( \log\Big( \frac{4q_{\max}^{(\ell)}}{\log(q_{\max}^{(\ell)})} \Big) + \frac{\log(q_{\max}^{(\ell)})}{4q_{\max}^{(\ell)}} - 1 \bigg) \\
    =&\, 4 \bigg( \frac{\log(4q_{\max}^{(\ell)}) - \log \log(q_{\max}^{(\ell)})}{\log(q_{\max}^{(\ell)})}  + \frac{1}{4q_{\max}^{(\ell)}} - \frac{1}{\log(q_{\max}^{(\ell)})}\bigg) \geq 4 + o(1).
\end{align*}
Consequently, according to \Cref{lem:Chernoffvariant}, for any $j\neq k$, we have
\begin{align*}
    \P\big( |\gE_{\ell}(v, \widehat{\gN}^{(0)})| > c_{N}^{(\ell)}\big) \leq &\, e^{-4q_{\max}^{(\ell)}} \leq \exp(- 4\cdot \D_{\rm{GCH}}(j, k) \cdot q_{N}),
\end{align*}
where the last inequality holds since $q_{\max}^{(\ell)} \gtrsim \D_{\rm{GCH}} \cdot q_{N}$ by \eqref{eqn:GCH_not_proportional} and \Cref{ass:expected_center_separation}. The desired argument then follows since $\D_{\rm{GCH}} = \min_{1\leq j < k\leq K} \D_{\rm{GCH}}(j, k)$.
\end{proof}

\section{Technical Lemmas}

\begin{lemma}\label{lem:GCH_max_different_order}
For $a, b > 0$, define the function $g(a, b) \coloneqq \max_{t\in [0, 1]} f(t)$, where
    \begin{align*}
        f(t) = t a + (1 - t) b - a^{t} b^{1-t}.
    \end{align*}
   When $a \asymp b$ but $a/b \neq 1$, $g(a, b) \asymp a$. Otherwise when $a \gg b$, $g(a, b) = a \cdot (1 - o(1))$.
\end{lemma}
\begin{proof}[Proof of \Cref{lem:GCH_max_different_order}]
    By basic calculus, the maximizer $t^{\star}$ is unique, and it satisfies
    \begin{align*}
        \frac{d f}{dt}\Big|_{t = t^{\star}} = a - b - b (a/b)^{t^{\star}} \log(a/b) = 0 \implies t^{\star} = \frac{1}{\log(a/b)} \log\Big( \frac{a/b - 1}{\log(a/b)} \Big).
    \end{align*}
    Consequently,
    \begin{align*}
        g(a, b) &\, = f(t^{\star}) = (a - b) t^{\star} + b\big( 1 - (a/b)^{t^{\star}} \big)\\
                &\, = \frac{a-b}{\log(a/b)} \Big( \log(a/b - 1) - \log\log(a/b)\Big) + b - \frac{a-b}{\log(a/b)}\\
                &\, = a\cdot (1 - b/a) \Big( \frac{\log(a/b - 1)}{\log(a/b)} - \frac{\log\log(a/b) - 1}{\log(a/b)}\Big) + b.
    \end{align*}
    The first part follows obviously. When $a \gg b$, then $b/a = o(1)$, $\log(a/b - 1)/\log(a/b) = 1 - o(1)$, $(\log\log(a/b) - 1)/\log(a/b) = o(1)$, thus $g(a, b) = a (1 - o(1))$.
\end{proof}

\begin{lemma}\label{lem:stirling}
    For integers $n, k \geq 1$, we have
    \begin{align*}
        \log(n!) = &\, n\log(n) - n + \frac{1}{2}\log(2\pi n) + O(n^{-1})\\
        \log \binom{n}{k} =&\, \frac{1}{2}\log \frac{n}{2\pi k (n-k)} + n\log(n) - k\log(k) - (n-k) \log(n-k) + O\Big( \frac{1}{n} + \frac{1}{k} + \frac{1}{n-k} \Big)\\
        \log \binom{n}{k} =&\, k\log\Big( \frac{n}{k} - 1 \Big) -\frac{1}{2}\log(2\pi k) + o(k^{-1})\,, \textnormal{ for } k = \omega(1) \textnormal{ and } \frac{k}{n} = o(1).
    \end{align*}
    Moreover, for any $1\leq k \leq \sqrt{n}$, we have
    \begin{align*}
        \frac{n^{k}}{4 \cdot k!} \leq \binom{n}{k} \leq \frac{n^{k}}{k!}\,,\quad \log \binom{n}{k} \geq k\log \Big(\frac{en}{k} \Big) - \frac{1}{2}\log(k)  - \frac{1}{12k} - \log(4\sqrt{2\pi})\,.
    \end{align*}
\end{lemma}
\begin{proof}[Proof of \Cref{lem:stirling}]
According to Stirling's series \cite{Diaconis1986AnEP}, for any $n \geq 1$, 
\begin{align*}
    \sqrt{2\pi n}\, \left(\frac{n}{e}\right)^{n}e^{\frac{1}{12n+1}} n! \sqrt{2\pi n}\, \left(\frac{n}{e}\right)^{n}e^{\frac{1}{12n}}.
\end{align*}
Then the following asymptotic expansion holds:
\begin{align*}
    \Big| \log(n!) - \Big(\frac{1}{2} \log(2\pi n) + n\log(n) - n \Big) \Big|\leq \frac{1}{12n}.
\end{align*}
Thus for any $k \geq 1$, we have $\log \binom{n}{k} = \log(n!) - \log(k!) - \log((n-k)!)$, as well as
\begin{align*}
    \Big| \log \binom{n}{k} - \Big( \frac{1}{2}\log \frac{n}{2\pi k (n-k)} + n\log(n) - k\log(k) - (n-k) \log(n-k) \Big) \Big| \leq \frac{1}{12} \Big( \frac{1}{n} + \frac{1}{k} + \frac{1}{n-k} \Big).
\end{align*}
When $k = \omega(1)$ but $\frac{k}{n} = o(1)$, we have $\log(1 - \frac{k}{n}) = - \frac{k}{n} + O(N^{-2})$, and 
\begin{align*}
    \log \binom{n}{k} =&\, -\frac{1}{2}\log(2\pi k) - \frac{1}{2}\log\Big(1 - \frac{k}{n} \Big) + n\log(n) - k\log(k) \\
    &\, \quad - (n-k) \log(n) - (n-k)\log\Big(1 - \frac{k}{n} \Big) + o(1)\\
    =&\, k\log(n) - k\log(k) + k + o(k).
\end{align*}
At the same time, we write
    \begin{align*}
        \binom{n}{k} = \frac{n(n-1)\cdots (n - k + 1)}{k!} = \frac{n^{k}}{k!} \cdot \Big(1 - \frac{1}{n} \Big)\Big(1 - \frac{2}{n} \Big) \cdots \Big(1 - \frac{k-1}{n} \Big) \geq \frac{n^{k}}{k!}\Big(1 - \frac{k-1}{n} \Big)^{k-1} 
    \end{align*}
The upper bound is trivial. For the lower bound, let $f(x) = (1 - x/n)^{x}$ denote the function. It is easy to see that $f(x)$ is decreasing with respect to $x$, and 
\begin{align*}
    \Big(1 - \frac{k-1}{n} \Big)^{k-1} \geq \Big(1 - \frac{1}{\sqrt{n}} \Big)^{\sqrt{n}}\,.
\end{align*}
Let $g(t) = (1 - 1/t)^{t}$ and we know that $g(t)$ is increasing when $t \geq 1$, hence $g(t) \geq g(2) = 1/4$ for any $t \geq 2$ and the desired lower bound follows. Then for any $k\leq \sqrt{n}$, we have
\begin{align*}
    \log \binom{n}{k} \geq &\,\log \frac{n^{k}}{4 \cdot k!} \\
    \geq &\, \log \left( \frac{n^{k}}{4 \cdot \sqrt{2\pi k}\ (\frac{k}{e})^{k}e^{\frac {1}{12k}}} \right) = k\log \Big(\frac{en}{k} \Big) - \frac{1}{2}\log(k)  - \frac{1}{12k} - \log(4\sqrt{2\pi})\,.
\end{align*}
\end{proof}

\begin{lemma}[Weyl's inequality, \cite{Weyl1912DasAV}]\label{lem:weyl}
Let $\rmA, \rmE \in \R^{m \times n}$ be two real $m\times n$ matrices, then $|\sigma_i(\rmA + \rmE) - \sigma_i(\rmA)| \leq \|\rmE\|$ for every $1 \leq i \leq \min\{ m, n\}$. Furthermore, if $m = n$ and $\rmA, \rmE \in \R^{n \times n}$ are real symmetric, then $|\lambda_i(\rmA + \rmE) - \lambda_i(\rmA)| \leq \|\rmE\|$ for all $1 \leq i \leq n$.
\end{lemma}

\begin{lemma}[Eckart–Young–Mirsky Theorem \cite{Eckart1936TheAO}]\label{lem:EYM}
Suppose that the matrix $\rmA \in \R^{m \times n}$ ($m\leq n$) adapts the singular value decomposition $\rmA = \rmU \mSigma \rmV^{\sT}$, where $\rmU = [\rvu_1, \ldots, \rvu_{m}] \in \R^{m\times m}$, $\rmV = [\rvv_1, \ldots, \rvv_{n}] \in \R^{n \times n}$ and $\mSigma \in \R^{m\times n}$ contains diagonal elements $\sigma_1 \geq \ldots \geq \sigma_{m}$. Let $\rmA^{(K)} = \sum_{k=1}^{K} \sigma_k \rvu_k \rvv_k ^{\sT}$ be the rank-$K$ approximation of $\rmA$. Then $\|\rmA - \rmA^{(K)}\|_2 = \sigma_{K+1} \leq \|\rmA - \rmB\|_2$ for any matrix $\rmB$ with $\mathrm{rank}(\rmB) = K$.
\end{lemma}

\begin{lemma}[Markov's inequality, {\cite[Proposition $1.2.4$]{Vershynin2018HighDP}}]\label{lem:Markov}
    For any non-negative random variable $\rX$ and $t>0$, we have 
    \begin{align*}
        \P(\rX > t) \leq \E(\rX)/t.
    \end{align*}
\end{lemma}

\begin{lemma}[Hoeffding's inequality, {\cite[Theorem $2.2.6$]{Vershynin2018HighDP}}]\label{lem:Hoeffding}
    Let $\rX_1,\dots, \rX_{N}$ be independent  random variables with $\rX_i \in [a_i, b_i]$, then for any, $t \geq 0$, we have
    \begin{align*}
        \P\Bigg( \bigg|\sum_{i=1}^{N}(\rX_i - \E \rX_i ) \bigg| \geq t \Bigg) \leq 2\exp \Bigg( -\frac{2t^{2} }{ \sum_{i=1}^{N}(b_i - a_i )^{2} } \Bigg)\,.
    \end{align*}
\end{lemma}

\begin{lemma}[Chernoff's inequality, {\cite[Theorem $2.3.1$]{Vershynin2018HighDP}}]\label{lem:Chernoff}
    Let $\rX_i$ be independent Bernoulli random variables with parameters $p_i$. Consider their sum $\rS_{N} = \sum_{i=1}^{N}\rX_i$ and denote its mean by $\mu = \E \rS_{N}$. Then for any $t > \mu$, 
    \begin{align*}
        \P \big( \rS_{N} \geq t \big) \leq e^{-\mu} \left( \frac{e \mu}{t} \right)^{t}\,.
    \end{align*}
\end{lemma}

\begin{lemma}[Bernstein's inequality, {\cite[Theorem $2.8.4$]{Vershynin2018HighDP}}]\label{lem:Bernstein}
    Let $\rX_1,\dots, \rX_{N}$ be independent mean-zero random variables such that $|\rX_i|\leq K$ for all $i$. Let $\sigma^2 = \sum_{i=1}^{N}\E \rX_i^2$. Then for every $t \geq 0$,
    \begin{align*}
        \P \Bigg( \Big|\sum_{i=1}^{N} \rX_i \Big| \geq t \Bigg) \leq 2 \exp \Bigg( - \frac{t^2/2}{\sigma^2 + Kt/3} \Bigg)\,.
    \end{align*}
\end{lemma}

\begin{lemma}[Bennett's inequality, {\cite[Theorem $2.9.2$]{Vershynin2018HighDP} }]\label{lem:Bennett}
    Let $\rX_1,\dots, \rX_{N}$ be independent random variables. Assume that $|\rX_i - \E \rX_i| \leq K$ almost surely for every $i$. Then for any $t>0$, we have
    \begin{align*}
        \P \Bigg( \sum_{i=1}^{N} (\rX_i - \E \rX_i) \geq t \Bigg) \leq \exp \Bigg( - \frac{\sigma^2}{K^2} \cdot h \bigg( \frac{Kt}{\sigma^2} \bigg)\Bigg)\,, \notag 
    \end{align*}
    where $\sigma^2 = \sum_{i=1}^{N}\Var(\rX_i)$ is the variance of the sum and $h(u) := (1 + u)\log(1 + u) - u$.
\end{lemma}

\begin{lemma}[Variants of Chernoff]\label{lem:Chernoffvariant}
Let $\rX_i$ be independent Bernoulli random variable with mean $p_i$. Consider their sum $\rS_{N} \coloneqq \sum_{i=1}^{N}\rX_i$ and denote its mean by $\mu \coloneqq \E \rS_{N}$.
    \begin{enumerate}
        \item Let $\mu = \gamma \cdot q_{N}$ for some $\gamma >0$ and $q_{N} \gg 1$, then for any $\zeta > \gamma$, we have
            \begin{align*}
                \P(\rS_{N} \geq \zeta q_{N} ) \leq \exp\big( - q_{N} \cdot \big[ (\gamma - \zeta) + \zeta \log (\zeta/\gamma) \big] \big)
            \end{align*}
        Specifically, $\P(\rS_{N} \geq t ) \leq N^{-[\zeta \log \frac{\zeta}{\gamma}  +  ( \gamma - \zeta) ]}$ when $q_{N} = \log(N)$.
        \item If $\mu = \gamma_{N} \cdot q_{N}$ for some vanishing sequence $\{\gamma_{N}\}_{N\geq 1}$ with $\gamma_{N} = o(1)$, then for any vanishing sequence $\{\zeta_{N}\}_{ N \geq 1}$ with $\zeta_{N} = \omega(\gamma_{N})$, we have
            \begin{align*}
                \P(\rS_{N} \geq \zeta_{N} \cdot q_{N} ) \leq \exp(-q_{N}\cdot [\zeta_{N} \log \frac{\zeta_{N}}{\gamma_{N}}( 1  +  o(1)\, ) ]).
            \end{align*}
    \end{enumerate}
\end{lemma}
\begin{proof}[Proof of \Cref{lem:Chernoffvariant}]
Chernoff \Cref{lem:Chernoff} gives
\begin{align*}
    \P(\rS_{N} \geq t ) &\,  \leq e^{-\mu} \Big( \frac{e \mu}{t} \Big)^{t} = \exp\Big( -(\mu-t) - t\log \Big(\frac{t}{\mu} \Big) \Big)
\end{align*}
Let $f(t) = t\log(\frac{t}{\mu}) + (\mu - t)$, where $f(\mu) = 0$ and $f^{'}(t) = \log(\frac{t}{\mu}) > 0$ for $t > \mu$, thus $f(t) > 0$ for all $t > \mu$.
\begin{enumerate}
    \item By taking $t = \zeta q_{N} $, we have
        \begin{align*}
            \P(\rS_{N} \geq t ) \leq \exp\big( - (\gamma - \zeta) q_{N} - \zeta \log( \zeta/\gamma) \cdot q_{N} \big).
        \end{align*}
    \item Note that $\gamma_{N} = o(1)$, $\zeta_{N} = o(1)$ but $\frac{\zeta_{N}}{\gamma_{N}} = \omega(1)$, then $\frac{\gamma_{N} - \zeta_{N}}{\zeta_{N} \log(\zeta_{N}/\gamma_{N}) } = o(1)$, hence follows.
\end{enumerate}
\end{proof}